\documentclass[conference,onecolumn,10pt]{IEEEtran}

\usepackage[ruled,vlined,longend]{algorithm2e}
\usepackage{amsfonts, amsmath, amssymb, amsthm}
\usepackage{bbm}
\usepackage{blindtext}
\usepackage{cite}
\usepackage{color}
\usepackage{fancyhdr}
\usepackage{float}
\usepackage{graphicx}
\usepackage{hyperref}
\usepackage{mathdots}
\usepackage{mathtools}
\usepackage{nicefrac}
\usepackage{outlines}
\usepackage{setspace}
\usepackage{titlesec}

\hypersetup{
    urlcolor=blue
}

\mathtoolsset{showonlyrefs}

\graphicspath{{./figures/}}

\newtheorem{lemma}{Lemma}
\newtheorem{theorem}{Theorem}

\theoremstyle{definition}

\theoremstyle{definition}
\newtheorem{remark}{Remark}

\theoremstyle{definition}
\newtheorem{problem}{Problem}

\newcommand{\fatone}{\mathbbm{1}}

\newcommand{\SO}{\mathrm{SO}\left(2\right)}
\newcommand{\SE}{\mathrm{SE}\left(2\right)}

\newcommand{\posei}{\mathbf{x}_{i}}
\newcommand{\posej}{\mathbf{x}_{j}}
\newcommand{\E}{\mathcal{E}}
\newcommand{\fijX}{f_{ij}\left(\mathcal{X}\right)}
\newcommand{\F}{\mathcal{F}}
\newcommand{\FX}{\mathcal{F}\left(\mathcal{X}\right)}
\newcommand{\FXk}{\mathcal{F}\left(\mathcal{X}_{k}\right)}
\newcommand{\FZero}{\mathcal{F}\left(\mathcal{X}_{0}\right)}
\newcommand{\fij}{f_{ij}}
\newcommand{\gij}{\mathbf{g}_{ij}}
\newcommand{\gijX}{\mathbf{g}_{ij}\left(\X\right)}
\newcommand{\hii}{\mathbf{h}_{ii}}
\newcommand{\hij}{\mathbf{h}_{ij}}
\newcommand{\hji}{\mathbf{h}_{ji}}
\newcommand{\hjj}{\mathbf{h}_{jj}}
\newcommand{\fijbar}{\bar{f}_{ij}}
\newcommand{\TBar}{\overline{\mathbf{T}}}
\newcommand{\JBar}{\bar{\mathcal{J}}}
\newcommand{\ebar}{\bar{\mathbf{e}}}
\newcommand{\gbar}{\bar{\mathbf{g}}}
\newcommand{\tbarx}{\bar{\mathbf{t}}_{\mathbf{x}}}
\newcommand{\tbarr}{\bar{\mathbf{t}}_{\mathbf{r}}}
\newcommand{\G}{\mathcal{G}}
\newcommand{\HijRGN}{\mathcal{R}_{ij}}
\newcommand{\Hk}{\mathcal{H}_{k}}
\newcommand{\HRGN}{\mathcal{R}}
\newcommand{\Sk}{\mathcal{S}_{k}}

\newcommand{\M}{\mathcal{M}}
\newcommand{\MN}{\mathcal{M}^{N}}
\newcommand{\TxM}{\mathcal{T}_{\mathbf{x}}\mathcal{M}}
\newcommand{\TxMNormal}{\mathcal{T}_{\mathbf{x}}^{\perp}\mathcal{M}}
\newcommand{\TXMN}{\mathcal{T}_{\mathcal{X}}\mathcal{M}^{N}}
\newcommand{\TXMNNormal}{\mathcal{T}_{\mathcal{X}}^{\perp}\mathcal{M}^{N}}
\newcommand{\TXkM}{\mathcal{T}_{\mathcal{X}_{k}}\mathcal{M}^{N}}
\newcommand{\TXkMN}{\mathcal{T}_{\mathcal{X}_{k}}\mathcal{M}^{N}}

\newcommand{\Ptilde}{\tilde{P}}
\newcommand{\PtildeV}{\tilde{P}_{\mathcal{V}}}

\newcommand{\Hbar}{\bar{\mathcal{H}}}
\newcommand{\Hijbar}{\bar{\mathcal{H}}_{ij}}
\newcommand{\HijbarX}{\bar{\mathcal{H}}_{ij}\left(\X\right)}
\newcommand{\K}{\mathcal{K}}
\newcommand{\U}{\mathcal{U}}
\newcommand{\V}{\mathcal{V}}
\newcommand{\W}{\mathcal{W}}
\newcommand{\X}{\mathcal{X}}
\newcommand{\Y}{\mathcal{Y}}
\newcommand{\XZero}{\mathcal{X}_{0}}
\newcommand{\Xk}{\mathcal{X}_{k}}
\newcommand{\XStar}{\mathcal{X}^{\star}}
\newcommand{\Z}{\mathcal{Z}}
\newcommand{\PX}{\mathcal{P}_{\mathcal{X}}}
\newcommand{\PXk}{\mathcal{P}_{\mathcal{X}_{k}}}

\newcommand{\gradF}{\textnormal{grad}\,\mathcal{F}}
\newcommand{\gradFX}{\textnormal{grad}\,\mathcal{F}\left(\mathcal{X}\right)}
\newcommand{\gradFXk}{\textnormal{grad}\,\mathcal{F}\left(\mathcal{X}_{k}\right)}
\newcommand{\HessF}{\textnormal{Hess}\,\mathcal{F}}
\newcommand{\Hessfij}{\textnormal{Hess}\,f_{ij}}
\newcommand{\HessfijX}{\textnormal{Hess}\,f_{ij}\left(\mathcal{X}\right)}
\newcommand{\HessFX}{\textnormal{Hess}\,\mathcal{F}\left(\mathcal{X}\right)}
\newcommand{\noise}{\eta_{ij}}
\newcommand{\zij}{\tilde{\mathbf{z}}_{ij}}
\newcommand{\eij}{\mathbf{e}_{ij}}
\newcommand{\rij}{{\mathbf{r}}_{ij}}
\newcommand{\Aij}{\mathcal{A}_{ij}}
\newcommand{\Bij}{\mathcal{B}_{ij}}
\newcommand{\Cii}{\mathcal{C}_{ii}}

\newcommand{\Cjj}{\mathcal{C}_{jj}}
\newcommand{\Omegaij}{\Omega_{ij}}
\newcommand{\LieAlgebra}{\mathcal{T}_{\fatone}\mathcal{M}}
\newcommand{\LogOne}{\textnormal{Log}_{\fatone}\!}
\newcommand{\Sublevel}{\left\{\mathcal{X}\mid\mathcal{F}\left(\mathcal{X}\right)\leq\mu\right\}}
\newcommand{\FZeroSublevel}{\left\{\mathcal{X}\mid\mathcal{F}\left(\mathcal{X}\right)\leq\mathcal{F}\left(\mathcal{X}_{0}\right)\right\}}
\newcommand{\Log}{\textnormal{Log}}
\newcommand{\Exp}{\textnormal{Exp}}
\newcommand{\transpose}{\top}
\newcommand{\vect}{\operatorname{vec}}
\newcommand{\MXZero}{M_{\X_{0}}}

\newcommand{\curve}{c}

\DeclareMathOperator*{\argmin}{arg\,min}

\newcommand{\numtrials}{15}

\newcommand{\rhoprime}{10^{-2}}
\newcommand{\deltabar}{10^{6}}
\newcommand{\deltazero}{100}
\newcommand{\epsilong}{10^{-2}}

\newcommand{\quotes}[1]{``#1''}
\newcommand{\vectx}[1]{\operatorname{vec}\left(#1\right)}
\newcommand{\diagx}[1]{\operatorname{diag}\left(#1\right)}
\newcommand{\norm}[1]{{\left\Vert#1\right\Vert}_{2}}
\newcommand{\fnorm}[1]{\left\Vert#1\right\Vert_{F}}
\newcommand{\Logx}[2]{\textnormal{Log}_{#1}\!\left(#2\right)}
\newcommand{\Expx}[2]{\textnormal{Exp}_{#1}\!\left(#2\right)}
\newcommand{\unaryminus}{\scalebox{0.75}[1.0]{\( - \)}}

\definecolor{imp_color}{rgb}{0.0, 0.5, 0.0}
\definecolor{ww_color}{rgb}{0.0, 0.6, 0.0}
\definecolor{hr_color}{rgb}{1.0, 0.0, 0.0}
\definecolor{pg_color}{rgb}{0.5, 0.0, 1.0}
\definecolor{kb_color}{rgb}{1.0, 0.5, 0.0}
\definecolor{mh_color}{rgb}{0.0, 0.0, 1.0}

\newboolean{show_comments}
\setboolean{show_comments}{true}
\ifthenelse{\boolean{show_comments}}{
    \newcommand{\ww}[1]{\textcolor{ww_color}{WW:~#1}} 
    \newcommand{\hr}[1]{\textcolor{hr_color}{HR:~#1}} 
    \newcommand{\pg}[1]{\textcolor{pg_color}{PG:~#1}} 
    \newcommand{\kb}[1]{\textcolor{kb_color}{KB:~#1}} 
    \newcommand{\mh}[1]{\textcolor{mh_color}{MH:~#1}} 
} {
    \newcommand{\ww}[1]{} 
    \newcommand{\hr}[1]{} 
    \newcommand{\pg}[1]{} 
    \newcommand{\kb}[1]{} 
    \newcommand{\mh}[1]{} 
}

\title{Technical Report: Pose Graph Optimization over Planar Unit Dual Quaternions: Improved Accuracy with Provably Convergent Riemannian Optimization}

\makeatletter
\newcommand{\linebreakand}{%
  \end{@IEEEauthorhalign}
  \hfill\mbox{}\par
  \mbox{}\hfill\begin{@IEEEauthorhalign}
}
\makeatother

\author{
    \IEEEauthorblockN{William D. Warke}
    \IEEEauthorblockA{\textit{Mechanical and Aerospace Engineering} \\
    \textit{University of Florida}\\
    Gainesville, FL, USA \\
    \href{mailto:william.warke@ufl.edu}{william.warke@ufl.edu}}
    \and
    \IEEEauthorblockN{J. Humberto Ramos}
    \IEEEauthorblockA{\textit{Mechanical and Aerospace Engineering} \\
    \textit{University of Florida}\\
    Shalimar, FL, USA \\
    \href{mailto:humbertoramos@ufl.edu}{humbertoramos@ufl.edu}}
    \and
    \IEEEauthorblockN{Prashant Ganesh}
    \IEEEauthorblockA{\textit{Principal Research Engineer} \\
    \textit{EpiSci Science, Incorporated}\\
    San Diego, CA, USA \\
    \href{mailto:prashantganesh@episci.com}{prashantganesh@episci.com}}
    \linebreakand
    \IEEEauthorblockN{Kevin M. Brink}
    \IEEEauthorblockA{\textit{Air Force Research Laboratory} \\
    \textit{Eglin Air Force Base}\\
    FL, USA \\
    \href{mailto:kevin.brink@us.af.mil}{kevin.brink@us.af.mil}}
    \and
    \IEEEauthorblockN{Matthew T. Hale}
    \IEEEauthorblockA{\textit{Electrical and Computer Engineering} \\
    \textit{Georgia Institute of Technology}\\
    Atlanta, GA, USA \\
    \href{mailto:mhale30@gatech.edu}{mhale30@gatech.edu}}
}

\begin{document}
\maketitle

\begin{abstract}
It is common in pose graph optimization (PGO) algorithms to assume that noise in the translations and rotations of relative pose measurements is uncorrelated. However, existing work shows that in practice these measurements can be highly correlated, which leads to degradation in the accuracy of PGO solutions that rely on this assumption. Therefore, in this paper we develop a novel algorithm derived from a realistic, correlated model of relative pose uncertainty, and we quantify the resulting improvement in the accuracy of the solutions we obtain relative to state-of-the-art PGO algorithms. Our approach utilizes Riemannian optimization on the planar unit dual quaternion (PUDQ) manifold, and we prove that it converges to first-order stationary points of a Lie-theoretic maximum likelihood objective. Then we show experimentally that, compared to state-of-the-art
PGO algorithms, ours produces estimation errors that are lower by 10\% to 25\% across several orders of magnitude of correlated noise levels and graph sizes. 
\end{abstract}

\thispagestyle{fancy}
\renewcommand{\headrulewidth}{0pt}
\cfoot{\textbf{Distribution Statement A.} Approved for public release: distribution is unlimited.}

\pagebreak
\tableofcontents
\pagebreak

\section{Introduction}\label{sec:intro}

Pose graph optimization (PGO) algorithms aim to optimally reconstruct the trajectory of a mobile agent using a set of uncertain relative measurements that were collected en-route. PGO is a backend component for numerous applications in robotics and computer vision, including simultaneous localization and mapping (SLAM)~\cite{grisetti2010tutorial, sunderhauf2012towards}, bundle adjustment~\cite{bender2013graph}, structure from motion~\cite{havlena2010efficient}, and photogrammetry~\cite{huang2006reassembling}. Additionally, a variety of related practical problems of interest~\cite{thunberg2011distributed, tron2012intrinsic, peters2014sensor, li2019pose} can be transformed into PGO problems, making it a versatile tool for optimization in these fields.

Some well-established PGO frameworks, such as g2o~\cite{grisetti2011g2o}, GTSAM~\cite{dellaert2012factor}, and iSAM~\cite{kaess2008isam}, have addressed the PGO problem using a mix of Euclidean and heuristic optimization techniques. More recently, algorithms based on Riemannian optimization, including SE-Sync~\cite{rosen2019se}, Cartan-Sync~\cite{briales2017cartan}, and CPL-Sync~\cite{fan2019efficient}, have demonstrated that, under certain conditions, the PGO problem admits a semidefinite relaxation whose solution approximates the solution of the original, unrelaxed problem. One condition assumed by the above algorithms (and others) is that uncertainties in position and orientation are modeled by isotropic (uncorrelated) noise.

However, the isotropic noise assumption runs contrary to existing results on uncertainty representations for rigid motion groups, which mathematically encode PGO problems. Specifically, it was shown in 2D~\cite{long2013banana} and in 3D~\cite{barfoot2014associating} that the propagation of uncertainty through compound rigid motions is best captured by a Lie-theoretic model~\cite{sola2018micro}, namely, a Gaussian distribution on the Lie algebra of a rigid motion group. In fact, the authors of~\cite{wheeler2018relative} demonstrated that such a Lie-theoretic model accurately predicted the distribution of a compound rigid motion trajectory where traditional models failed. These Lie-theoretic models are inherently anisotropic, which suggests that a PGO algorithm that incorporates anisotropy may attain improved accuracy. 

Therefore, in this paper, we formulate 2D PGO problems on the manifold of \textit{planar unit dual quaternions} (PUDQs), which 
we use to explicitly incorporate anisotropy in uncertainty models. 
To solve such problems, we use a Riemannian trust region (RTR) algorithm, for which we derive
global convergence guarantees. 
The contributions of this paper are: 
\begin{itemize}
    \item We present what is, to the best of our knowledge, the first provably convergent PGO algorithm that permits arbitrarily large, anisotropic uncertainties.
    \item We prove that the proposed algorithm converges to first-order critical points given \emph{any} initialization.
    \item We show that the resulting pose estimates are always at least~$10$\% more accurate than the state of the art and more than~$25$\% more accurate on high-dimensional problems.
\end{itemize}

The closest related works are~\cite{cheng2016dual,bai2021sparse,li2020improved}. In~\cite{cheng2016dual}, a unit dual quaternion approach to PGO was developed using heuristic optimization techniques without formal guarantees, whereas we employ provably convergent Riemannian-geometric techniques. The authors of~\cite{bai2021sparse} used a Lie-theoretic objective, but did not include convergence guarantees or quantify the accuracy of their solutions. The work in~\cite{li2020improved} uses a similar problem 
formulation to us, though that work was entirely empirical. We differ both by proving convergence and showing improvement in accuracy over a class of Riemannian algorithms that were not studied in~\cite{li2020improved}.

The rest of the paper is organized as follows. Section~\ref{sec:prelim} provides preliminaries, and Section~\ref{sec:prob_form}  provides a formal problem statement. Section~\ref{sec:algorithm} outlines the proposed algorithm, and Section~\ref{sec:conv_analysis} 
proves that it converges. 
Section~\ref{sec:results} contains numerical results, and Section~\ref{sec:conclusion} concludes.

\section{Preliminaries}\label{sec:prelim}
In this section, we include mathematical preliminaries that are necessary for our PUDQ PGO problem formulation. For detailed derivations, see Appendices~\ref{app:construction}-\ref{app:riemannian}.
\subsection{Planar unit dual quaternion construction}\label{sec:pudq_construction}
We construct the PUDQ manifold as a representation of planar rigid motion. Given an orthonormal basis $\{\mathbf{i},\mathbf{j},\mathbf{k}\}$, a planar rigid motion is characterized by a translation, denoted $\mathbf{t}=t_{x}\mathbf{i}+t_{y}\mathbf{j}$, and a rotation about the $\mathbf{k}$ axis by an angle $\theta\in(-\pi,\pi]$. 
The PUDQ parameterization of this motion is given by $\mathbf{x}=\mathbf{x}_{r}+\epsilon\mathbf{x}_{d}$, where $\epsilon$ is a \textit{dual number} satisfying $\epsilon^{2}=0,\epsilon\neq0$. The \textit{real} and \textit{dual} parts of $\mathbf{x}$, denoted $\mathbf{x}_{r}\in\mathbb{S}^{1}$ and $\mathbf{x}_{d}\in\mathbb{R}^{2}$, respectively, are $\mathbf{x}_{r}\triangleq\cos\left(\nicefrac{\theta}{2}\right)+\sin\left(\nicefrac{\theta}{2}\right)\mathbf{k}$ and $\mathbf{x}_{d}\triangleq\frac{1}{2}\boldsymbol{\mathrm{t}}\otimes \mathbf{x}_{r}$, with \quotes{$\otimes$} denoting the Hamilton product~\cite{hamilton1848xi} under the convention $\mathbf{i}^{2}=\mathbf{j}^{2}=\mathbf{k}^{2}=\mathbf{i}\mathbf{j}\mathbf{k}=-1$. Applying the Hamilton product to two PUDQs, denoted $\mathbf{x}$ and $\mathbf{y}$, yields the composition operator \quotes{$\boxplus$}, which can be
expressed as
\begin{equation}\label{eq:QL_QR}
    \mathbf{x}\boxplus\mathbf{y}=\underset{Q_{L}\left(\mathbf{x}\right)}{\underbrace{\Bigg[\begin{smallmatrix}
        x_{0} & -x_{1} & 0 & 0\\
        x_{1} & x_{0} & 0 & 0\\
        x_{2} & x_{3} & x_{0} & -x_{1}\\
        x_{3} & -x_{2} & x_{1} & x_{0}
        \end{smallmatrix}\Bigg]}}\underset{\mathbf{y}}{\underbrace{\Bigg[\begin{smallmatrix}
        y_{0}\\
        y_{1}\\
        y_{2}\\
        y_{3}
        \end{smallmatrix}\Bigg]}}=\underset{Q_{R}\left(\mathbf{y}\right)}{\underbrace{\Bigg[\begin{smallmatrix}
            y_{0} & -y_{1} & 0 & 0\\
            y_{1} & y_{0} & 0 & 0\\
            y_{2} & -y_{3} & y_{0} & y_{1}\\
            y_{3} & y_{2} & -y_{1} & y_{0}
            \end{smallmatrix}\Bigg]}}\underset{\mathbf{x}}{\underbrace{\Bigg[\begin{smallmatrix}
            x_{0}\\
            x_{1}\\
            x_{2}\\
            x_{3}
            \end{smallmatrix}\Bigg]}},
\end{equation}
where $Q_{L}(\cdot)$ and $Q_{R}(\cdot)$ denote the left and right composition maps, respectively. From~\eqref{eq:QL_QR}, we have the identity element $\fatone=\left[1, 0, 0, 0\right]^{\transpose}$ and inverse formula $\mathbf{x}^{-1}=\left[x_{0},-x_{1},-x_{2},-x_{3}\right]^{\transpose}$. The set of  PUDQs forms the smooth manifold $\M\triangleq\mathbb{S}^{1}\rtimes\mathbb{R}^{2}\subset\mathbb{R}^{4}$, which we embed in $\mathbb{R}^{4}$ as
\begin{equation}\label{eq:M_embedding}
    \M\triangleq\left\{\mathbf{x}\in\mathbb{R}^{4}\mid h\left(\mathbf{x}\right)=\mathbf{x}^{\transpose}\Ptilde\mathbf{x}-1=0\right\}\subset\mathbb{R}^{4},
\end{equation}
where $\Ptilde\triangleq\mathrm{diag}(\{1,1,0,0\})$ and $h(\mathbf{x})$ is the \textit{defining function}~\cite{boumal2023introduction} for $\M$. 
PGO algorithms optimize over $N$ poses, so we extend~\eqref{eq:M_embedding} to the $N$-fold product manifold $\MN\triangleq(\mathbb{S}^{1}\rtimes\mathbb{R}^{2})^{N}$. Below, we will use the operator $\vect(\cdot)$, where
\begin{equation}
    \vect((\mathbf{x}_{i})_{i=1}^{N})\triangleq[\mathbf{x}_{1}^{\transpose},\mathbf{x}_{2}^{\transpose},\dots,\mathbf{x}_{N}^{\transpose}]^{\transpose}, 
\end{equation}
with each $\mathbf{x}_{i}\in\M$. Since $(\mathbb{S}^{1}\rtimes\mathbb{R}^{2})^{N}\subset\mathbb{R}^{4\times N}\cong\mathbb{R}^{4N}$, 
we embed $\MN$ in $\mathbb{R}^{4N}$. 
For~$\mathcal{X}, \mathcal{Y} \in \M^N$, 
this embedding lets us write~$\mathcal{X} = \vect((\mathbf{x}_{i})_{i=1}^{N})$ and $\mathcal{Y} = \vect((\mathbf{y}_{i})_{i=1}^{N})$,
where~$\mathbf{x}_i, \mathbf{y}_i \in \M$ for each~$i$. 
This embedding also gives the identity $\fatone^{N}=\vect((\fatone)_{i=1}^{N})$, the inverse formula $\X^{-1}=\vect((\mathbf{x}_{i}^{-1})_{i=1}^{N})$, and the product $\X\boxplus\Y=\vect((\mathbf{x}_{i}\boxplus\mathbf{y}_{i})_{i=1}^{N})$. 
\subsection{Logarithm and exponential maps}\label{sec:log_exp}
The smooth manifold $\M$ with the identity, inverse, and composition operator form a Lie group\cite{sola2018micro} whose Lie algebra is the \textit{tangent space} at the identity element, denoted $\mathcal{T}_{\fatone}\M$. Given $\mathbf{x}\in\M$, the logarithm map at the identity element is $\textnormal{Log}_{\fatone}:\M\rightarrow \mathcal{T}_{\fatone}\M$, given by
\begin{equation}\label{eq:log_1_def}
    \Logx{\fatone}{\mathbf{x}}=\frac{1}{\gamma\left(\mathbf{x}\right)}{\left[x_{1},~x_{2},~x_{3}\right]}^{\transpose},
\end{equation}
with $\gamma\left(\mathbf{x}\right)\triangleq\mathrm{sinc}\left(\phi\left(\mathbf{x}\right)\right)=\sin\left(\phi\left(\mathbf{x}\right)\right)/\phi\left(\mathbf{x}\right)$, where $\phi\left(\mathbf{x}\right)\triangleq\mathrm{wrap}\left(\mathrm{arctan}\left(x_{1},x_{0}\right)\right)$, $\operatorname{arctan}:\mathbb{S}^{1}\rightarrow(-\pi,\pi]$ is the four-quadrant arctangent and
\begin{equation}\label{eq:wrap}
    \mathrm{wrap}\left(\alpha\right)\triangleq\begin{cases}
        \alpha+\pi & \text{if }\alpha\leq-\nicefrac{\pi}{2}\\
        \alpha-\pi & \text{if }\alpha>\nicefrac{\pi}{2}\\
        \alpha & \text{otherwise}.
        \end{cases}
\end{equation}
Here, $\phi:\M\rightarrow \left(-\nicefrac{\pi}{2},\nicefrac{\pi}{2}\right]$ computes the half-angle of rotation about the $\mathbf{k}$-axis encoded by a point on~$\M$. The half-angles $\phi + n\pi$ for all $n\in\mathbb{Z}$ encode the same rotation, so it is valid to wrap $\phi$ to $(-\nicefrac{\pi}{2},\nicefrac{\pi}{2}]$ via~\eqref{eq:wrap}.
 
Given some $\mathbf{x}_{t}=\left[x_{t,1},~x_{t,2},~x_{t,3}\right]^{\transpose}\in\mathcal{T}_{\fatone}\M$, the exponential map at the identity, denoted $\mathrm{Exp}_{\fatone}:T_{\fatone}\M\rightarrow\M$, is given by
    $\Expx{\fatone}{\mathbf{x}_{t}}={\left[\cos\left(x_{t,1}\right), \gamma\left(\mathbf{x}_{t}\right)\mathbf{x}_{t}^{\transpose}\right]}^{\transpose}$,
where $\gamma\left(\mathbf{x}_{t}\right)\triangleq\operatorname{sinc}\left(x_{t,1}\right)$ as above. For any~$\mathbf{x}, \boldsymbol{\mathrm{y}}\in\M$, we also have the point-wise logarithm map
\begin{equation}\label{eq:log_x_def}
    \Logx{\mathbf{x}}{\mathbf{y}}=\mathbf{x}\boxplus[0,\LogOne(\mathbf{x}^{-1}\boxplus\mathbf{y})^{\transpose}]^{\transpose},
\end{equation}
and, for $\mathbf{x}\in\M$, and some $\mathbf{y}_{t}\in \mathcal{T}_{\mathbf{x}}\M$, the point-wise exponential map
\begin{equation}\label{eq:exp_x_def}
    \Expx{\mathbf{x}}{\mathbf{y}_{t}}=\mathbf{x}\boxplus\Expx{\fatone}{\left(\mathbf{x}^{-1}\boxplus\mathbf{y}_{t}\right)_{1:3}},
\end{equation}
where ${\left(\cdot\right)}_{1:3}$ selects the last three entries of a vector. 
For $\X,\Y\in\MN$,~\eqref{eq:log_x_def}-\eqref{eq:exp_x_def} give 
logarithm and exponential maps over the product manifold~$\M^N$, namely
    $\Logx{\X}{\Y}=\vectx{(\Log_{\mathbf{x}_{i}}(\mathbf{y}_{i}))_{i=1}^{N}}$,
and, for any $\Y_{t}=\vectx{(\mathbf{y}_{t,i})_{i=1}^{N}}\in\TXMN$, the mapping
\begin{equation}\label{eq:ExpX_MN}
    \Expx{\X}{\Y_{t}}=\vectx{(\Exp_{\mathbf{x}_{i}}(\mathbf{y}_{t,i}))_{i=1}^{N}},
\end{equation}
with $\Log_{\mathbf{x}_{i}}(\cdot)$ and $\Exp_{\mathbf{x}_{i}}(\cdot)$ given by~\eqref{eq:log_x_def} and~\eqref{eq:exp_x_def}.
\subsection{Pose graph construction}
\begin{figure}
    \vspace{2mm}
    \centering
    \includegraphics[scale=1.0]{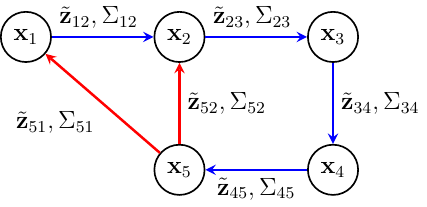}
    \caption{A pose graph with $N=M=5$, labeled with vertex poses $\mathbf{x}_i$, edge measurements $\zij$, and edge covariances $\Sigma_{ij}$.~\textit{Odometry} edges, shown in blue, connect neighboring vertices (i.e., $\left|j-i\right|=1$). \textit{Loop closure} edges, shown in red, connect any non-neighboring vertices (i.e., $\left|j-i\right|>1$).}\label{fig:pose_graph}
\end{figure}

We now address the construction of a pose graph, as exemplified in Figure~\ref{fig:pose_graph}. First, let $\G=(\V,\E)$ be a directed graph with vertex set $\V$ and edge set $\E$ of ordered pairs $(i,j)\in\V\times\V$. Letting $|\mathcal{V}|=N$, we define $\X=\vect((\mathbf{x}_{i})_{i\in\V})\in\MN$ to be the vector of $N$ poses to be estimated, with individual poses denoted $\mathbf{x}_{i}\in\M$. Then, letting $\left|\mathcal{E}\right|=M$, we define $\Z=\vect((\zij)_{\left(i,j\right)\in\mathcal{E}})\in\M^{M}$ to be the vector of $M$ relative pose measurements, where $\zij\in\M$ encodes a measured transformation from $\posei$ to $\posej$, taken in the frame of $\posei$. The noise covariance for $\zij$ is given by the matrix $\Sigma_{ij}$. The corresponding \textit{pose graph} is then constructed by associating the vertex set $\V$ with $\X$, and the edge set $\E$ with $\Z$. 
\section{Problem Formulation}\label{sec:prob_form}
We now derive the problem to be solved. From the perspective of Bayesian inference, PGO algorithms aim to estimate the posterior distribution of poses that best fits a given dataset of relative measurements made along a 
trajectory. Because a prior distribution is not always available, PGO is typically formulated as a \textit{maximum likelihood estimation} (MLE) problem~\cite{grisetti2010tutorial}, and
we use such a formulation here. 

Motivated by~\cite{long2013banana}, we utilize a Lie-theoretic measurement model for $\zij$ in which zero-mean Gaussian noise $\noise$ is mapped from $\mathcal{T}_{\fatone}\mathcal{M}$ to $\mathcal{M}$ via the exponential map, i.e.,
\begin{equation}\label{eq:meas_model}
    \zij=\mathbf{x}_{i}^{-1}\boxplus\mathbf{x}_{j}\boxplus\Expx{\fatone}{\noise},
\end{equation}
with $\noise\in\mathbb{R}^{3}$ and $\noise\sim\mathcal{N}\left(0,\Sigma_{ij}\right)$. As noted in the Introduction,~\eqref{eq:meas_model} gives a realistic model of compound, uncertain transformations. In Appendix~\ref{app:mle}, we show that~\eqref{eq:meas_model} yields the MLE objective $\mathcal{F}:\MN\rightarrow\mathbb{R}$, given by 
\begin{equation}\label{eq:mle_F}
    \mathcal{F}\left(\X\right)=\frac{1}{2}\sum_{\left(i,j\right)\in\mathcal{E}}\fijX,
\end{equation}
where
\begin{equation}\label{eq:fij}
    \fijX\triangleq\left\Vert \mathbf{e}_{ij}\left(\mathbf{x}_{i},\mathbf{x}_{j}\right)\right\Vert _{\Omega_{ij}}^{2}.
\end{equation}
Here, $\Omega_{ij}=\Sigma_{ij}^{-1}$ is the information matrix for edge $(i,j)$, and $\mathbf{e}_{ij}:\M\times\M\rightarrow\LieAlgebra$ is the \textit{tangent} residual given by
\begin{equation}\label{eq:e_ij_def}
    \mathbf{e}_{ij}\left(\mathbf{x}_{i},\mathbf{x}_{j}\right)\triangleq\Logx{\fatone}{\mathbf{r}_{ij}\left(\mathbf{x}_{i},\mathbf{x}_{j}\right)},
\end{equation}
and $\mathbf{r}_{ij}:\mathcal{M}\times\mathcal{M}\rightarrow\mathcal{M}$ is the \textit{geodesic} residual, defined as\footnote{Henceforth, we simply write $\mathbf{e}_{ij}\triangleq\mathbf{e}_{ij}\left(\mathbf{x}_{i},\mathbf{x}_{j}\right)$ and $\mathbf{r}_{ij}\triangleq\mathbf{r}_{ij}\left(\mathbf{x}_{i},\mathbf{x}_{j}\right)$.}\
\begin{equation}\label{eq:r_ij_def}
    \mathbf{r}_{ij}\left(\mathbf{x}_{i},\mathbf{x}_{j}\right)\triangleq\zij^{-1} \boxplus \mathbf{x}_{i}^{-1} \boxplus \mathbf{x}_{j}.
\end{equation}
 In a geometric sense, $\mathbf{r}_{ij}$ encodes the geodesic along $\mathcal{M}$ from a measurement $\zij$ to the estimated relative transformation $\mathbf{x}_{i}^{-1}\boxplus\mathbf{x}_{j}$. The map $\mathbf{e}_{ij}$ then \quotes{unwraps} the geodesic to the Lie algebra. 

We now address \textit{anchoring}, a problem that arises because the objective in~\eqref{eq:mle_F} is invariant to an identical rigid transformation of all poses, i.e., $\mathcal{F}(\X)=\mathcal{F}(\vect((\mathbf{y})_{i\in\V})\boxplus\X)=\mathcal{F}(\X\boxplus\vect((\mathbf{y})_{i\in\V}))$ for any $\mathbf{y}\in\M$. To remedy this, one must \quotes{anchor} at least one vertex by setting $\mathbf{x}_{a}\triangleq\fatone$ for some $a\in\V$, so we assume that this has been done for some node. Given this formulation, we now formally state the problem that we
solve in the remainder of the paper.
\begin{problem}\label{prob:mle}
Given a measurement set $\mathcal{Z}\in\mathcal{M}^{M}$, compute the \textit{maximum likelihood estimate} $\X^{\star} \in \M^N$, where
\begin{equation}\label{eq:X_star}
    \X^{\star}=\underset{\X\in\mathcal{M}^{N}}{\argmin}\ \FX,
\end{equation}
with $\F$ given by~\eqref{eq:mle_F}.
\end{problem}
Problem~\ref{prob:mle} is a nonconvex, nonlinear least squares problem over a Riemannian manifold. In the following section, we employ Riemannian optimization techniques to solve~\eqref{eq:X_star}.

%
\section{Algorithm Description}\label{sec:algorithm}
This section presents the method by which we solve Problem~\ref{prob:mle}, starting with a brief description of the class of algorithms we employ. Trust-region methods~\cite{nocedal2006trust} for optimization in $\mathbb{R}^{n}$ employ a local approximation of the objective function, called the \textit{model}, about each iterate. The model is restricted to a neighborhood of the current iterate, called the \textit{trust region}. At each iteration, a tentative update step is computed, and is accepted to compute the next iterate if the model sufficiently agrees with the objective at the computed point. Riemannian trust region (RTR) methods~\cite[Chapter 7]{absil2008optimization} generalize this idea to Riemannian manifolds, and our proposed algorithm adapts the RTR framework to planar PGO on $\mathcal{M}^{N}$. 
\begin{figure}
    \vspace{2mm}
    \centering
    \includegraphics[scale=0.35]{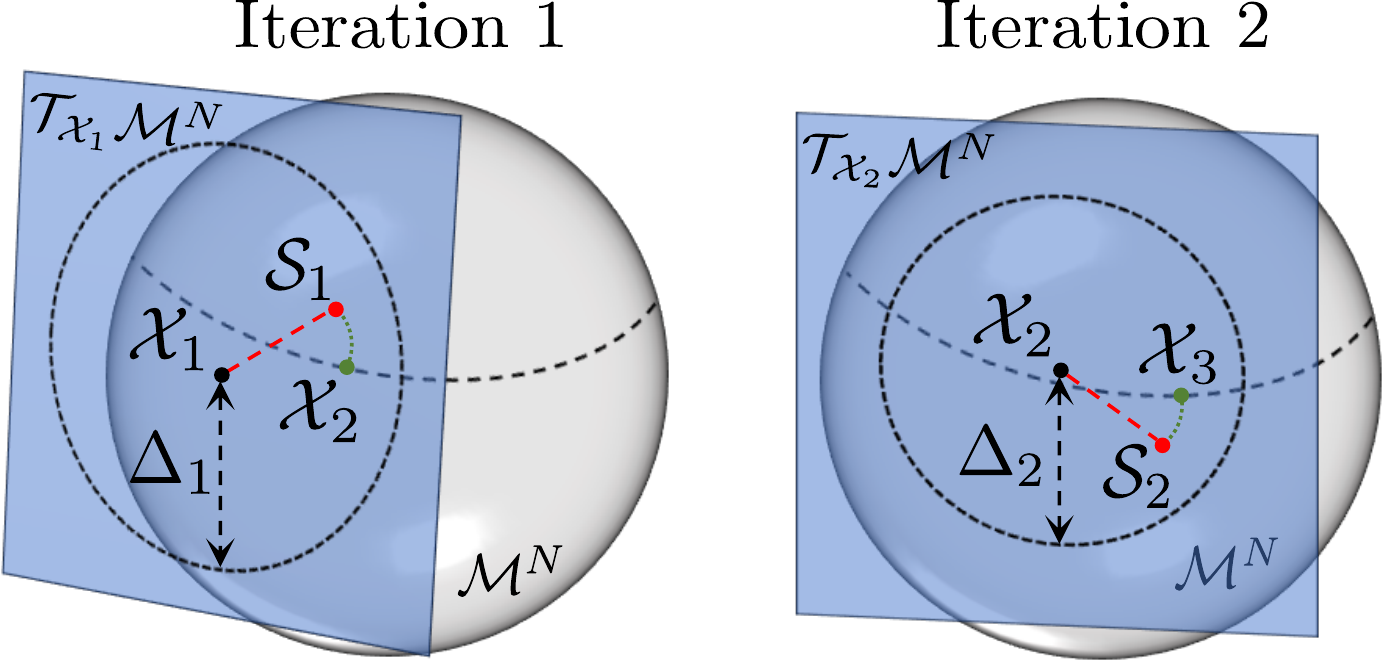}
    \caption{An illustration of two iterations of the RTR algorithm. At each iteration, the algorithm computes a tangent step $\mathcal{S}_{k}\in\mathcal{T}_{\X_{k}}\mathcal{M}^{N}$, shown in red, within a trust region of radius $\Delta_{k}$, which
    is indicated by the dotted circle shown in each tangent space. 
    If the step is accepted (as defined in~\eqref{eq:X_kplusone}), then the next iterate is computed as $\X_{k+1}=\mathrm{Exp}_{\X_{k}}(\mathcal{S}_{k})$, which maps the step from the tangent space back to the manifold itself, as shown in green.}
    \label{fig:rtr_visual}
\end{figure}

An illustration of the proposed RTR algorithm is shown in Figure~\ref{fig:rtr_visual}. At each iteration $k$, instead of approximating the objective $\F$, RTR computes an approximation of $\F$ in the tangent space at $\X_{k}$, called a \textit{pullback}. The pullback is defined as $\hat{\F}_{k}\triangleq\F\circ\textnormal{Exp}_{\Xk}$.\footnote{The pullback can be implemented using any retraction~\cite{absil2007trust,boumal2019global}, and we choose to use the exponential map since it is well-defined on $\mathcal{M}^{N}$ and straightforward to compute.} The approximation takes the form of a second-order model $\hat{m}_{k}:\TXkMN\rightarrow\mathbb{R}$, namely
\begin{equation}\label{eq:m_k}
    \hat{m}_{k}(\mathcal{S})\triangleq\mathcal{F}\left(\X_{k}\right)+ \mathcal{S}^{\transpose}\gradF\left(\X_{k}\right) +\frac{1}{2} \mathcal{S}^{\transpose}\mathcal{H}_{k}\mathcal{S},
\end{equation}
where $\mathcal{S}\in\TXkMN$ is a tangent vector centered at $\X_{k}$, $\gradF:\mathcal{M}^{N}\rightarrow\TXkMN$ is the Riemannian gradient, and $\mathcal{H}_{k}:\mathcal{T}_{\X_{k}}\mathcal{M}^{N}\rightarrow\mathcal{T}_{\X_{k}}\mathcal{M}^{N}$ is a symmetric approximation of the Riemannian Hessian at $\X_{k}$. We include explicit forms for $\gradF$ in Appendix~\ref{app:rgrad} and our choice of $\Hk$ in Appendix~\ref{app:rgn_hess}.

Our procedure corresponds to the RTR update given in~\cite[Chapter 7]{absil2008optimization}. The algorithm is initialized with $\X_{0}\in\mathcal{M}^{N}$ and trust-region radius $\Delta_{0}\in\left(0,\bar{\Delta}\right]$, where $\bar{\Delta}>0$ is the user-specified maximum radius. At iteration~$k$, the tentative step
$\mathcal{S}_{k}$ is computed by solving the inner sub-problem
\begin{equation}\label{eq:inner_subproblem}
\mathcal{S}_{k}=\argmin_{\mathcal{S}\in\TXkMN}\hat{m}_{k}\left(\mathcal{S}\right)\text{~subject~to~}{\left\Vert \mathcal{S}\right\Vert}_{2} \leq\Delta_{k},
\end{equation}
where~$\hat{m}_k$ is from~\eqref{eq:m_k}. To solve~\eqref{eq:inner_subproblem}, we employ the Steihaug-Toint truncated conjugate gradients (tCG) algorithm~\cite{steihaug1983conjugate,toint1981towards}, which offers unique benefits for trust-region sub-problems, including monotonic cost decrease and early termination (thereby approximating~\eqref{eq:inner_subproblem}) in the cases of negative curvature or trust region violation. To measure the agreement between the model and objective functions, we use
\begin{equation}\label{eq:rho_k}
     \rho_{k}=\frac{\hat{\mathcal{F}}_{k}\left(\mathbf{0}\right)-\hat{\mathcal{F}}_{k}\left(\mathcal{S}_{k}\right)}{\hat{m}_{k}\left(\mathbf{0}\right)-\hat{m}_{k}\left(\mathcal{S}_{k}\right)},
\end{equation}
where $\mathbf{0}\in\mathbb{R}^{4N}$ is the zero vector. Based on the level of agreement, the trust-region radius $\Delta_{k}$ is then updated via 
\begin{equation}\label{eq:Delta_kplusone}
     \Delta_{k+1}=\begin{cases}
     \frac{1}{4}\Delta_{k} & \text{if~}\rho_{k}<\frac{1}{4}\\
     \min\left\{2\Delta_{k},\bar{\Delta}\right\} & \text{if~}\rho_{k}>\frac{3}{4}\text{~and~}{\left\Vert \mathcal{S}_{k}\right\Vert}_{2} =\Delta_{k}\\
     \Delta_{k} & \text{otherwise.}
     \end{cases}
\end{equation}
The tentative step $\mathcal{S}_{k}$ is accepted to compute $\X_{k+1}$ only if the model agreement ratio $\rho_{k}$ from~\eqref{eq:rho_k} is greater than a user-defined model agreement threshold $\rho^{\prime}\in\left(0,\nicefrac{1}{4}\right)$, i.e.,
\begin{equation}\label{eq:X_kplusone}
     \X_{k+1}=\begin{cases}
     \Expx{\X_{k}}{\Sk} & \text{if~}\rho_{k}>\rho^{\prime}\\
     \X_{k} & \text{otherwise.}
     \end{cases}
\end{equation}

As summarized in Algorithm~\ref{alg:RTR}, 
the steps from~\eqref{eq:inner_subproblem}-\eqref{eq:X_kplusone} are repeated until the gradient norm crosses below a user-defined threshold $\varepsilon_{g}$, i.e., until ${\left\Vert \operatorname{grad}\mathcal{F}\left(\X_{k}\right)\right\Vert}_{2} \leq\varepsilon_{g}$. 

\begin{algorithm}[htbp]
\SetAlgoLined
\textbf{Input:} Edge measurement set $\mathcal{Z}\in \mathcal{M}^{M}$,\\
Maximum trust-region radius $\bar{\Delta}>0$,\\
Model agreement threshold $\rho^{\prime}\in\left(0,\nicefrac{1}{4}\right]$,\\
Gradient termination threshold $\varepsilon_g>0$.\\
\textbf{Initialize:} $k\leftarrow 0$, $\X_{0}\in \mathcal{M}^{N}$, $\Delta_0 \in \left(0,\bar{\Delta}\right]$
\\
\While{$\norm{\gradFXk} >\varepsilon_{g}$}{
    Compute $\Sk$ from~\eqref{eq:inner_subproblem} using tCG. \\
    Compute $\rho_{k}$ using~\eqref{eq:rho_k}.\\
    Compute $\Delta_{k+1}$ using~\eqref{eq:Delta_kplusone}.\\
    Compute $\X_{k+1}$ using~\eqref{eq:X_kplusone}.\\
    $k\leftarrow k+1$
}
\textbf{return} $\X_{k}$
\caption{RTR for PUDQ PGO}\label{alg:RTR}
\end{algorithm}

\section{Convergence Analysis}\label{sec:conv_analysis}
In this section, we prove that Algorithm~\ref{alg:RTR} is globally convergent. Specifically, given any initialization,
it reaches a first-order critical point to within a user-specified tolerance in finite time. 
The authors of~\cite{boumal2019global} proposed global rates of convergence for the RTR algorithm given a set of assumptions about the problem, so we treat these assumptions as sufficient conditions for convergence. For our proof, we will establish:
\begin{enumerate}
    \item Lower-boundedness of $\mathcal{F}$ on $\mathcal{M}^{N}$.\label{cond:lower_bound}
    \item Sufficient decrease in the model cost at each iteration.\label{cond:model_decrease}
    \item A Lipschitz-type condition for gradients of pullbacks.\label{cond:lipschitz}
    \item Radial linearity and boundedness of $\mathcal{H}_{k}$.\label{cond:rad_lin_bound}
\end{enumerate}
We will make each of these statements mathematically precise in the following analysis. Towards proving Condition~\ref{cond:lower_bound}, we first derive a lemma on continuity of $\F$.
\begin{lemma}\label{lem:continuity}
    The objective $\F$ is continuous on $\MN$.
\end{lemma}
\emph{Proof:} By inspection of~\eqref{eq:log_1_def} and \eqref{eq:mle_F}-\eqref{eq:e_ij_def}, and continuity of \quotes{$\boxplus$} from~\eqref{eq:QL_QR} as a linear map, it suffices to show that $\LogOne$ is continuous on $\M$. While~\eqref{eq:log_1_def} and \eqref{eq:wrap} contain discontinuities independently, we will show that their composition to form~$\textnormal{Log}_{\fatone}$ 
does not. Let $\phi_{1}\triangleq\arctan(r_{1}, r_{0})$ (where $\rij=[r_{0},r_{1},r_{2},r_{3}]^{\transpose}$ denotes the element-wise map), and let $\phi_{2}\triangleq\mathrm{wrap}(\phi_{1})$. Then, we have discontinuities in $\phi_{1}$ at $(r_{0},r_{1})=(-1,0)$, in $\mathrm{wrap}(\phi_{1})$ at $\phi_{1}=\pm\nicefrac{\pi}{2}$, and in $(\gamma(\phi_{2}))^{-1}$  at $\phi_{2}=\pm\pi$. We now observe that $\mathrm{wrap}(-\pi)=\mathrm{wrap}(\pi)=0$, so $\lim_{(r_{0},r_{1})\rightarrow(-1,0)}\mathrm{wrap}(\phi_{1})=0$, thereby nullifying the discontinuities in $\phi_{1}$. Next, $(\mathrm{sinc}(\phi_{2}))^{-1}$ is even and continuous on the domain $[\nicefrac{-\pi}{2},\nicefrac{\pi}{2}]$, so $\lim_{\phi_{2}\rightarrow\nicefrac{-\pi}{2}}(\gamma(\phi_{2}))^{-1}=\lim_{\phi_{2}\rightarrow\nicefrac{\pi}{2}}(\gamma(\phi_{2}))^{-1}=\nicefrac{\pi}{2}$, nullifying the discontinuities in $\phi_{2}$. Finally, because $\lim_{\phi_{2}\rightarrow 0}(\gamma(\phi_{2}))^{-1}=1$ and, by~\eqref{eq:wrap}, $\phi_{2}\in(\nicefrac{-\pi}{2},\nicefrac{-\pi}{2}]$, we conclude that $\LogOne$ is continuous on $\mathcal{M}$, which implies that $\F$ is continuous on $\MN$.\hfill$\blacksquare$

We now show compactness of sublevel sets of $\F$. 
\begin{theorem}\label{thm:compact_sublevels}
    The $\mu$-sublevel sets of $\F$, given by $\Sublevel$, are compact.
\end{theorem}
\emph{Proof:} From~\eqref{eq:ExpX_MN}, for every $\X\in\MN$, $\textnormal{Exp}_{\X}$ is globally defined on $\TXMN$, which implies that $\MN$ is geodesically complete. Therefore, the Hopf-Rinow Theorem~\cite{udriste2013convex} implies that closed and bounded subsets of $\MN$ are compact, so it suffices to show that the sublevel sets are closed and bounded.

From~\eqref{eq:mle_F}-\eqref{eq:fij}, $\FX\geq0$ for all $\X\in\MN$, which implies that the $\mu$-sublevel sets of $\F$ are the preimages of the closed subsets $[0,\mu]$, i.e., $\mu$-sublevel sets are
of the form $\mathcal{F}^{-1}\big([0, \mu]\big)$. 
These sets are closed because $\F$ is continuous by Lemma~\ref{lem:continuity}. 

Turning to boundedness of sublevel sets,~\eqref{eq:M_embedding} implies that $\M$ is unbounded, and therefore $\MN$ is unbounded. Then, by~\cite[Theorem 1]{celledoni2018dissipative}, the $\mu$-sublevel sets are bounded if and only if $\F$ is coercive, i.e., for all $\Y\in\MN$, every sequence $\{\X_{l}\}_{l\in\mathbb{N}}\subset\MN$ such that $\lim_{l\rightarrow\infty}d_{\MN}\left(\X_{l},\Y\right)=\infty$ also satisfies $\lim_{l\rightarrow\infty}F(\X_{l})=\infty$.\footnote{Here, $d_{\MN}(\cdot,\cdot)$ is the geodesic distance on $\MN$ defined in Appendix~\ref{app:geodesic}.} Therefore, it suffices to show that $\F$ is coercive, which we do next.

First, let $\mathcal{X}_{l}=\vectx{(\mathbf{x}_{l,i})_{i\in\mathcal{V}}}$ and $\mathcal{Y}=\vectx{(\mathbf{y}_{i})_{i\in\mathcal{V}}}$, and observe from the definition of $d_{\mathcal{M}^{N}}\left(\mathcal{X}_{l},\mathcal{Y}\right)$ that
\begin{equation}
    \lim_{d_{\mathcal{M}^{N}}\left(\mathcal{X}_{l},\mathcal{Y}\right)\rightarrow\infty}	\max_{i\in\mathcal{V}}\Vert\LogOne(\mathbf{x}_{l,i}^{-1}\boxplus\mathbf{y}_{i})\Vert _{2}^{2} =\infty.
\end{equation}
We now rewrite $\Vert \LogOne(\mathbf{x}_{l,i}^{-1}\boxplus\mathbf{y}_{i})\Vert _{2}^{2}$ as
\begin{equation}
    \Vert \LogOne(\mathbf{x}_{l,i}^{-1}\boxplus\mathbf{y}_{i})\Vert _{2}^{2}=\gamma(\mathbf{x}_{l,i}^{-1}\boxplus\mathbf{y}_{i})^{-2}\mathbf{x}_{l,i}^{\top}M_{LR}^{-}(\mathbf{y}_{i})\mathbf{x}_{l,i},
\end{equation}
where $M_{LR}^{-}(\mathbf{y}_{i})\triangleq Q_{LR}^{-}\left(\mathbf{y}_{i}\right)^{\top}\mathrm{diag}\left(\left\{ 0,I_{3}\right\} \right)Q_{LR}^{-}\left(\mathbf{y}_{i}\right)$, with $Q_{LR}^{-}\left(\mathbf{y}_{i}\right)$ given in Appendix~\ref{app:construction}. Since $\gamma\left(\mathbf{x}\right)\in\left[\nicefrac{-\pi}{2},\nicefrac{\pi}{2}\right]$ for all $\mathbf{x}\in\mathcal{M}$, we have
\begin{equation}\label{eq:log_xinv_y}
    \Vert \LogOne(\mathbf{x}_{l,i}^{-1}\boxplus\mathbf{y}_{i})\Vert _{2}^{2}\leq(\nicefrac{\pi^{2}}{4})\lambda_{\max}(M_{LR}^{-}(\mathbf{y}_{i}))\mathbf{x}_{l,i}^{\top}\mathbf{x}_{l,i},
\end{equation}
where $\lambda_{\max}\left(\cdot\right)$ denotes the maximum eigenvalue of a matrix. Since $\mathbf{y}_{i}$ is constant and $\lambda_{max}(M_{LR}^{-}(\mathbf{y}_{i}))\geq0$,~\eqref{eq:log_xinv_y} implies that $\lim_{\Vert\LogOne(\mathbf{x}_{l,i}^{-1}\boxplus\mathbf{y}_{i})\Vert_{2}^{2}\rightarrow\infty}(\mathbf{x}_{l,i}^{\top}\mathbf{x}_{l,i})=\infty$. The first element of $\mathbf{x}_{l,i}\in\mathcal{M}$ is bounded by $1$, so $\mathbf{x}_{l,i}^{\top}\mathbf{x}_{l,i}-1\leq\Vert\LogOne\left(\mathbf{x}_{l,i}\right)\Vert_{2}^{2}$. Therefore, $\lim_{(\mathbf{x}_{l,i}^{\top}\mathbf{x}_{l,i})\rightarrow\infty}\Vert\LogOne(\mathbf{x}_{l,i})\Vert_{2}^{2}=\infty$. Now, we note that for any $\mathbf{x},\mathbf{y}\in\M$, we can write
\begin{align}
    \left\Vert \LogOne\left(\mathbf{x}\boxplus\mathbf{y}\right)\right\Vert _{2}^{2}&=\gamma\left(\mathbf{x}\boxplus\mathbf{y}\right)^{-1}\mathbf{y}^{\top}M_{L}\left(\mathbf{x}\right)\mathbf{y} \\
	&=\gamma\left(\mathbf{x}\boxplus\mathbf{y}\right)^{-1}\mathbf{x}^{\top}M_{R}\left(\mathbf{y}\right)\mathbf{x},
\end{align}
where $M_{L}\left(\mathbf{x}\right)\triangleq Q_{L}\left(\mathbf{x}\right)^{\top}\mathrm{diag}\left(\left\{ 0,I_{3}\right\} \right)Q_{L}\left(\mathbf{x}\right)$ and $M_{R}\left(\mathbf{y}\right)\triangleq Q_{R}\left(\mathbf{y}\right)^{\top}\mathrm{diag}\left(\left\{ 0,I_{3}\right\} \right)Q_{R}\left(\mathbf{y}\right)$. Because $M_{L}(\cdot),M_{R}(\cdot)\succeq 0$, it holds that, for any $\mathbf{x},\mathbf{y}\in\mathcal{M}$,
\begin{equation}\label{eq:log_xy_lim}
    \lim_{\left\Vert \LogOne\left(\mathbf{x}\boxplus\mathbf{y}\right)\right\Vert _{2}^{2}\rightarrow\infty}\max \left\{\left\Vert \LogOne\left(\mathbf{x}\right)\right\Vert _{2}^{2},\left\Vert \LogOne\left(\mathbf{y}\right)\right\Vert _{2}^{2}\right\} =\infty.
\end{equation}
We now observe that for any two vertices $\mathbf{x}_{i},\mathbf{x}_{j}\in\M$, with $i,j\in\mathcal{V}$ and $i>j$, it follows from connectedness of odometry edges in $\E$ that $\mathbf{x}_{i}=\mathbf{x}_{j}\boxplus\mathbf{c}_{i,j}$, where
\begin{equation}\label{eq:c_ij}
    \mathbf{c}_{i,j}\triangleq\tilde{\mathbf{z}}_{j(j+1)}\boxplus\mathbf{r}_{(j+1)(j+2)}\boxplus\cdots\boxplus\mathbf{\tilde{z}}_{\left(i-1\right)i}\boxplus\mathbf{r}_{\left(i-1\right)i}.
\end{equation}
Equivalently, we have $\mathbf{x}_{j}=\mathbf{x}_{i}\boxplus\mathbf{c}_{i,j}^{-1}$. Per Section~\ref{sec:prob_form}, we have anchored $\mathbf{x}_{a} \triangleq \fatone$ for some $a\in\mathcal{V}$, and since $\LogOne(\mathbf{x}^{-1})=-\LogOne(\mathbf{x})$, it holds that $\Vert\LogOne(\mathbf{x}_{l,m})\Vert_{2}^{2}=\Vert\LogOne(\mathbf{c}_{a,m})\Vert_{2}^{2}$ for any $m\in\mathcal{V}$. Furthermore, because the $\tilde{\mathbf{z}}_{ij}$ terms in~\eqref{eq:c_ij} are constant, applying~\eqref{eq:log_xy_lim} inductively yields, for any $m\in\mathcal{V}$,
\begin{equation}
    \lim_{\left\Vert\LogOne\left(\mathbf{x}_{l,m}\right)\right\Vert _{2}^{2}\rightarrow\infty}\max_{\left(i,j\right)\in\mathcal{E}}\left\Vert \LogOne\left(\mathbf{r}_{i,j}\right)\right\Vert _{2}^{2} = \infty.
\end{equation}
From~\eqref{eq:fij}, $\lambda_{\min}\left(\Omega_{ij}\right)\left\Vert \LogOne\left(\mathbf{r}_{i,j}\right)\right\Vert _{2}^{2}\leq f_{ij}\left(\mathcal{X}_{l}\right)$, where $\lambda_{\min}(\cdot)$ is the minimum eigenvalue, and $\lambda_{\min}(\Omega_{ij})>0$ because~$\Omega_{ij} = \Sigma_{ij}^{-1} \succ 0$. Then $\lim_{\Vert \LogOne(\mathbf{r}_{ij})\Vert _{2}^{2}\rightarrow\infty}f_{ij}(\mathcal{X}_{l})=\infty$, and~\eqref{eq:mle_F} gives
$\lim_{f_{ij}(\mathcal{X}_{l})\rightarrow\infty}\mathcal{F}(\mathcal{X}_{l})=\infty$. 
Then~$\mathcal{F}$ is coercive and the proof is complete. 
\hfill$\blacksquare$\\


Next, we show that the objective $\F$ satisfies Condition~\ref{cond:lower_bound}.
\begin{lemma}\label{lem:F_lower_bound}
    There exists $\mathcal{F}^{\star}\geq0$ such that $\mathcal{F}\left(\X\right)\geq\mathcal{F}^{\star}$ for all $\X\in\mathcal{M}^{N}$.
\end{lemma}
\emph{Proof:} Lemma~\ref{lem:continuity}, Theorem~\ref{thm:compact_sublevels}, and the Weierstrass Theorem~\cite[Prop. A.8]{bertsekas1997nonlinear}
imply the existence of a global minimizer $\XStar\in\MN$, which is the solution to Problem~\ref{prob:mle}. 
Setting $\mathcal{F}^{\star}\triangleq\F\left(\XStar\right)$ completes the proof.~\hfill $\blacksquare$\\

We now show that Algorithm~\ref{alg:RTR} satisfies Condition~\ref{cond:model_decrease}.
\begin{lemma}\label{lem:model_decrease}
    For all $\Xk$ computed by Algorithm~\ref{alg:RTR} such that $\norm{\gradF\left(\Xk\right)}>\varepsilon_{g}$, it holds that the step $\Sk$ satisfies
    \begin{equation}\label{eq:model_decrease}
        \hat{m}_k\left(\mathbf{0}\right)-\hat{m}_k\left(\Sk\right) \geq \frac{1}{2} \min \{\Delta_k, 2\varepsilon_{g}\} \varepsilon_g.
    \end{equation}
\end{lemma}
\emph{Proof:} By design, iterates of the tCG algorithm produce a strict, monotonic decrease of the model cost $\hat{m}_{k}$~\cite{boumal2019global}. For all $k$, the first tCG iterate is the Cauchy step, which satisfies~\eqref{eq:model_decrease} by definition 
and thus completes the proof.\hfill$\blacksquare$\\

The forthcoming analysis in Lemma~\ref{lem:F_decrease}, Theorem~\ref{thm:rgrad_lipschitz}, and Lemma~\ref{lem:conv_lipschitz} addresses Condition~\ref{cond:lipschitz}, namely, Lipschitz continuity of the Riemannian gradient, $\gradF$. First, we use Theorem~\ref{thm:rgrad_lipschitz} to prove its Lipschitz continuity on compact subsets of $\mathcal{M}^{N}$.
\begin{theorem}\label{thm:rgrad_lipschitz}
    The Riemannian gradient, $\gradF$, is $L_{g}$-Lipschitz continuous on any compact subset $\K\subset\mathcal{M}^{N}$.
    That is, there exists $L_{g}>0$ such that for all $\X,\mathcal{Y}\in\K$ we have
    \begin{equation}\label{eq:lipschitz_def}
        \norm{\mathcal{P}_{\X\rightarrow\mathcal{Y}}\gradFX-\gradF\left(\mathcal{Y}\right)} \leq L_{g}d_{\MN}\left(\X,\mathcal{Y}\right),
    \end{equation}
    where $\mathcal{P}_{\X\rightarrow\mathcal{Y}}:\TXMN\rightarrow\mathcal{T}_{\mathcal{Y}}\MN$ is the parallel transport operator defined in Appendix~\ref{app:parallel}.
\end{theorem}
\emph{Proof:} A necessary and sufficient condition for~\eqref{eq:lipschitz_def} is that, for all $\mathcal{X}\in\K$, the Riemannian Hessian, $\HessF$, has operator norm bounded by $L_{g}$, i.e.,
\begin{equation}\label{eq:rhess_op_norm_def}
    \sup_{\mathcal{V}\in\mathcal{T}_{\mathcal{X}}\mathcal{M},\Vert\mathcal{V}\Vert_{2}=1}\left\Vert \text{Hess}~\mathcal{F}(\mathcal{X})[\mathcal{V}]\right\Vert _{2}\leq L_{g}.
\end{equation}
In Appendices~\ref{app:rhess}-\ref{app:lipschitz}, we derive $\HessF$ and derive a constant $L_{g}$ for which~\eqref{eq:rhess_op_norm_def} holds on any compact subset $\K\subset\MN$, completing the proof.\hfill$\blacksquare$

To apply Theorem~\ref{thm:rgrad_lipschitz} to Algorithm~\ref{alg:RTR}, we must first show that the computed iterates remain within the $\FZero$-sublevel set for all $k$, which is accomplished by Lemma~\ref{lem:F_decrease}.

\begin{lemma}\label{lem:F_decrease}
    The objective $\F$ is monotonically decreasing with respect to the iterates of Algorithm~\ref{alg:RTR}.
    In particular, it holds that $\FXk\leq\FZero$ for all $k$.
\end{lemma}
\emph{Proof:} By~\eqref{eq:model_decrease}, we have $\hat{m}_{k}\left(\mathbf{0}\right)-\hat{m}_{k}\left(\Sk\right)>0$ for all $k$. If any $\Sk$ would yield an increase in $\F$, then $\FXk-\F(\mathrm{Exp}_{\Xk}(\Sk))<0$, and~\eqref{eq:rho_k} implies $\rho_{k}<0$. 
By~\eqref{eq:X_kplusone}, such an $\Sk$ is rejected and, therefore the condition $\F(\X_{k+1})=\FXk$ is enforced in such cases. 
Thus, since it cannot occur that~$\mathcal{F}(\mathcal{X}_{k+1}) > \mathcal{F}(\mathcal{X}_k)$, 
we see that~$\F(\X_{k+1})\leq\FXk$ for all $k$. By induction, $\FXk\leq\FZero$ for all $k$, completing the proof.\hfill$\blacksquare$\\

Now, Lemma~\ref{lem:conv_lipschitz} extends Theorem~\ref{thm:rgrad_lipschitz} to any $\Xk$ computed by Algorithm~\ref{alg:RTR}, 
which shows that Condition~\ref{cond:lipschitz} is satisfied.
\begin{lemma}\label{lem:conv_lipschitz}
    For all $\Xk$ computed by Algorithm~\ref{alg:RTR}, there exists $L_{g}\geq0$ such that
\begin{equation}\label{eq:lipschitz_type}
    \left|\mathcal{F}\left(\Expx{\Xk}{\mathcal{S}}\right)-\left(\FXk+\mathcal{S}^{\transpose}\gradF\left(\X_{k}\right) \right)\right|\leq\frac{L_{g}}{2}\norm{\mathcal{S}}^{2}
\end{equation}
for all $\mathcal{S}\in\TXkM$ such that $\norm{\mathcal{S}}\leq\bar{\Delta}$ and for all $k$.
\end{lemma}
\emph{Proof:} Let $\MXZero\triangleq\FZeroSublevel$ and set 
\begin{equation}\label{eq:K_lipschitz}
    \K\triangleq \MXZero\cup\{\Expx{\X}{\mathcal{S}}\mid\X\in \MXZero,\norm{\mathcal{S}}\leq\bar{\Delta}\}.
\end{equation}
Then Theorem~\ref{thm:compact_sublevels} implies that $\MXZero$ is compact, and therefore so is $\K$. Lemma~\ref{lem:F_decrease} implies $\Xk\in \MXZero\subset\K$ for all $k$. By Theorem~\ref{thm:rgrad_lipschitz}, there exists $L_{g}>0$ to which~\eqref{eq:lipschitz_def} applies for all $\Xk\in\K$. From~\cite[Lemma 2.1]{bento2017iteration}, we find that~\eqref{eq:lipschitz_def} implies~\eqref{eq:lipschitz_type}, completing the proof.\hfill $\blacksquare$\\

Lemmas~\ref{lem:Hk_rad_lin} and~\ref{lem:Hk_bounded} address Condition~\ref{cond:rad_lin_bound}, which pertains to properties of $\mathcal{H}_{k}$, the Riemannian Hessian approximation used in~\eqref{eq:inner_subproblem} and spelled out in Appendix~\ref{app:rgn_hess}.

\begin{lemma}\label{lem:Hk_rad_lin}
    The operator $\Hk$ in~\ref{eq:RGN_H_k} is radially linear, i.e., for all $\mathcal{S} \in \TXkM$ and all $\alpha \geq 0$, we have
        $\Hk[\alpha \mathcal{S}]=\alpha \Hk[\mathcal{S}]$.
\end{lemma}
\emph{Proof:} Equation~\eqref{eq:RGN_H_k} is linear by inspection. \hfill$\blacksquare$\\

\begin{lemma}\label{lem:Hk_bounded}
    The operator $\mathcal{H}_{k}$ in~\eqref{eq:RGN_H_k} is bounded
    for all $\Xk$ computed by Algorithm~\ref{alg:RTR}, i.e., there exists $\beta<\infty$ 
    such that 
    \begin{equation}\label{eq:beta_bound}
        \max_{\mathcal{S}}\Big\{\Vert\Hk\mathcal{S}\Vert_{2} \mid \mathcal{S}\in\TXkM,\norm{\mathcal{S}}=1\Big\}\leq\beta.
    \end{equation}
\end{lemma}
\emph{Proof:} First, $\left\Vert \mathcal{S}\right\Vert _{2}=1$ implies $\left\Vert \Hk\mathcal{S}\right\Vert _{2}\leq\left\Vert \Hk\right\Vert _{2}$. Substituting~\eqref{eq:RGN_H_k}, applying the triangle inequality, and using the fact that $\lambda_{max}\left(\mathcal{P}_{\mathcal{X}}\right)=1$ yields
\begin{equation}\label{eq:H_k_bound}
    \left\Vert \Hk\right\Vert _{2}\leq\sum_{\left(i,j\right)\in\mathcal{E}}\left\Vert\PX\HijRGN\PX\right\Vert _{2}\leq\sum_{\left(i,j\right)\in\mathcal{E}}\left\Vert \HijRGN\right\Vert _{2}.
\end{equation}
Since, by definition of~$\|\cdot\|_2$ and~$\|\cdot\|_F$ we have $\Vert\HijRGN\Vert_{2}\leq\Vert\HijRGN\Vert_{F}$, we reach 
\begin{equation}\label{eq:Hij_bound}
    \left\Vert\HijRGN\right\Vert_{2}\leq4\left\Vert\Aij\right\Vert_{F}\left\Vert\Bij\right\Vert_{F}\left\Vert\Omegaij\right\Vert_{F}.
\end{equation}
Now, we set $\K$ as in~\eqref{eq:K_lipschitz} and apply the bounds derived in Appendix~\ref{app:egrad_bounds} for $\left\Vert\Aij\right\Vert_{F}$ and $\left\Vert\Bij\right\Vert_{F}$ on compact subsets of $\MN$. Since every term on 
the right-hand side of~\eqref{eq:Hij_bound} is bounded, we see that
the right-hand side of~\eqref{eq:H_k_bound} is bounded, completing the proof.\hfill$\blacksquare$\\

Our convergence analysis culminates in Theorem~\ref{thm:rtr_conv}.
\begin{theorem}\label{thm:rtr_conv}
    Let $\varepsilon_{g}\leq\nicefrac{\Delta_{0}}{\lambda_{g}}$ be given, 
    where~$\Delta_0$ is from Section~\ref{sec:algorithm}, $\lambda_g\triangleq\nicefrac{1}{4} \min \left\{\nicefrac{1}{\beta}, \nicefrac{1}{2(L_{g}+\beta)}\right\}$, $L_{g}$ is from~\eqref{eq:lipschitz_type}, and $\beta$ is from~\eqref{eq:beta_bound}. 
    Then, for any initialization $\XZero\in\MN$, Algorithm~\ref{alg:RTR} 
    produces an iterate~$\X_k$ that satisfies $\left\Vert\gradF(\X_{k})\right\Vert_{2}\leq\varepsilon_{g}$
    in no more than~$K$ iterations, where
    \begin{equation}\label{eq:iteration_bound}
        K \leq \frac{\FZero-\F\left(\X^{\star}\right)}{\rho^{\prime} \lambda_g} \frac{3}{\varepsilon_{g}^{2}}+\frac{1}{2} \log _2\left(\frac{\Delta_0}{\lambda_g \varepsilon_{g}}\right),
    \end{equation}
    where~$\rho'$ is from~\eqref{eq:X_kplusone} and $\mathcal{X}^{\star}$ is from Lemma~\ref{lem:F_lower_bound}.
\end{theorem}
\emph{Proof:} Lemmas~\ref{lem:F_lower_bound},~\ref{lem:model_decrease}, and~\ref{lem:conv_lipschitz}-\ref{lem:Hk_bounded} 
show the satisfaction of Conditions~\ref{cond:lower_bound}-\ref{cond:rad_lin_bound} in~\cite[Theorem 12]{boumal2019global}, which immediately implies that~\eqref{eq:iteration_bound} holds for Algorithm~\ref{alg:RTR}.\hfill$\blacksquare$\\

Theorem~\ref{thm:rtr_conv} gives provable convergence of Algorithm~\ref{alg:RTR} to approximate first-order critical points of $\F$ \emph{under any initialization} $\XZero$, and we note that the tolerance~$\varepsilon_g$ can be made to take arbitrary values 
by adjusting $\Delta_{0}$. 
\section{Experimental Results}\label{sec:results}
In this section, we validate the accuracy of Algorithm~\ref{alg:RTR} relative to the Riemannian PGO solvers SE-Sync~\cite{rosen2019se} and Cartan-Sync~\cite{briales2017cartan}. Both yield a global minimizer identical to that computed by the class of Riemannian algorithms that use semidefinite relaxations~(e.g., \cite{fan2019efficient,tian2021distributed}), so we omit additional comparisons to those algorithms.

Because an objective comparison necessitated the use of exact ground truth, we opted to adapt three synthetic PGO datasets with diverse vertex and edge counts. The first of these is Grid1000, which we synthesized\footnote{To synthesize the Grid1000 dataset, a ground truth trajectory is computed along a randomized grid resembling the Manhattan datasets created for~\cite{olson2006fast}. Loop closure edges were selected at random, specifically, with $3.0\%$ probability of an edge at Euclidean inter-pose distances of up to $2$ meters.} with $N=1000$ vertices and $M=1250$ edges. The remaining datasets are publicly available, and serve as common benchmarks for PGO evaluations, namely, (i) M3500~\cite{olson2006fast}, with $N=3500$, $M=5598$, 
and (ii) City10000~\cite{kaess2008isam}, with $N=10000$, $M=20687$. To generate PGO trial datasets, we apply calibrated noise to the ground truth dataset for each graph. Each of these datasets, including ground truth, is available at {\small\url{https://github.com/corelabuf/planar_pgo_datasets}}.

\subsection{PGO dataset generation}\label{sec:generation}
\begin{figure}
    \vspace{2mm}
    \centering
    \includegraphics[scale=0.35]{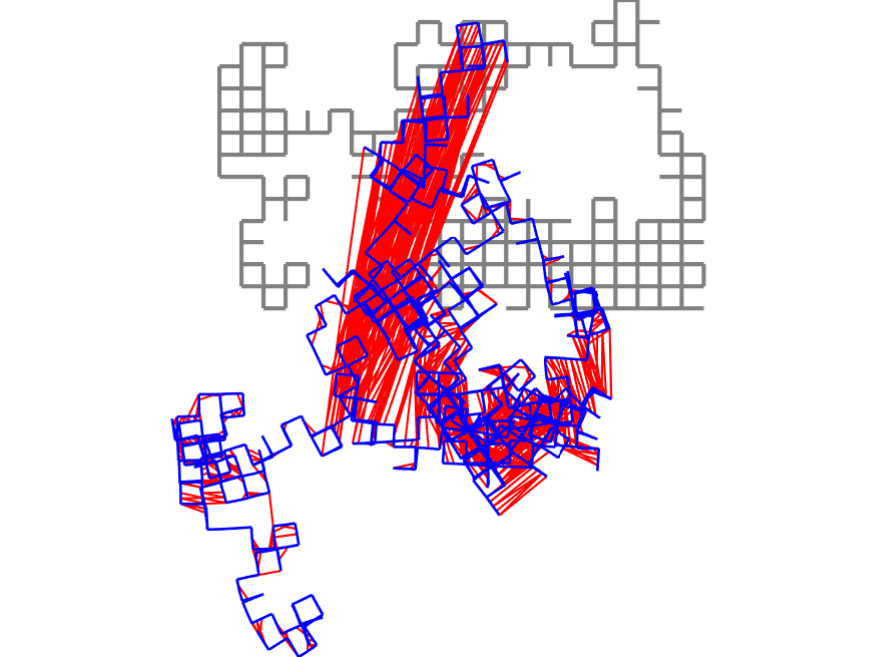}
    \includegraphics[scale=0.35]{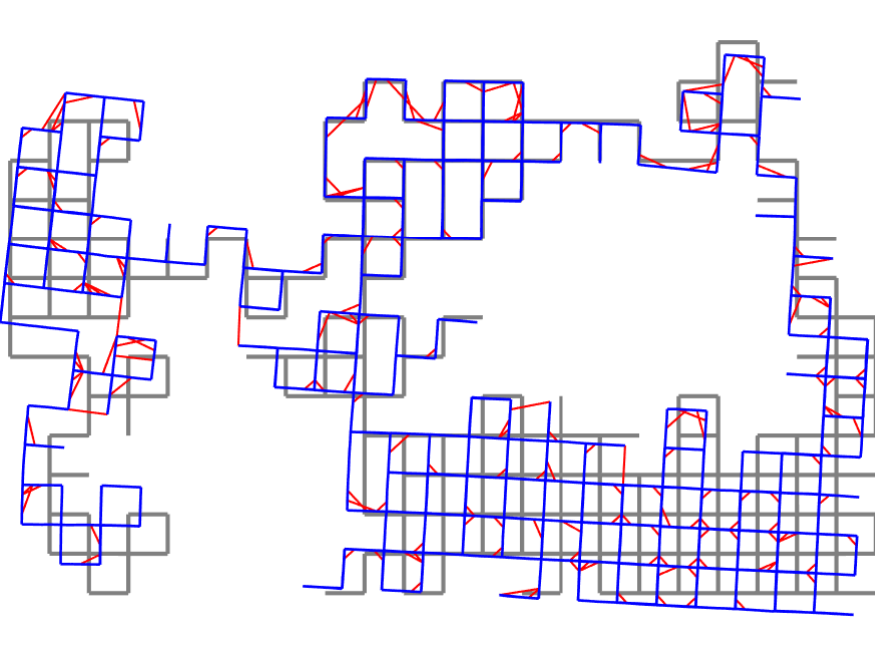}
    \caption{(Left) The M3500 pose graph dataset, corrupted with Lie-theoretic noise. (Right) The estimated
    graph computed by Algorithm~\ref{alg:RTR}. Odometry edges are blue, loop closures are red, and ground truth is shown in gray.}
    \label{fig:random_pg}
\end{figure}

To generate a PGO dataset, the true edge measurements from each dataset are corrupted using the Lie-theoretic noise model from~\eqref{eq:meas_model}. The edge measurement noise covariance, $\Sigma_{ij}$, is computed as $\Sigma_{ij}\sim W_{3}\left(\sigma_{w}\Sigma_{w},10\right)$, where $W_{d}\left(V,n\right)$ is the Wishart distribution with dimension $d$, scale matrix $V$, and $n$ degrees of freedom\footnote{The sample covariance matrix of a multivariate Gaussian random variable is Wishart-distributed\cite{chatfield2018introduction}, making it an apt choice for this application.}. Here, $\sigma_{w}$ is a variance tuning parameter, and $\Sigma_{w}$ is given by $\Sigma_{w}\triangleq J_{3}+\mathrm{diag}\left(\left[u_{1}, u_{2},u_{3}\right]\right)$, where $J_{3}\in\mathbb{R}^{3\times3}$ is a matrix of ones and $u_{i}\sim\mathcal{U}_{(0,1]}$ are uniformly sampled on the interval $(0,1]$. This generates random, positive-definite, anisotropic covariance matrices with $\mathbb{E}[\Sigma_{ij}]=10\sigma_{w}\Sigma_{w}$, which simulates relative pose covariances computed by a Lie-theoretic estimator. Using this approach, we generated 5 trial datasets per ground truth, for a total of 15. Figure~\ref{fig:random_pg} depicts an M3500 variant generated with $\sigma_{w}=5.62\cdot10^{-5}$ alongside the estimate computed by Algorithm~\ref{alg:RTR}.

\subsection{Evaluation methodology}\label{sec:metrics}
Solutions computed by each algorithm were evaluated using the root-mean square relative pose error (RPE) metric. RPE measures total edge deformation with respect to the ground truth, and gives an objective performance metric for SLAM algorithms~\cite{kummerle2009measuring}. We denote $(\mathbf{x}_{i}^{\diamond})_{i=1}^{N}$ to 
be the ground truth poses, and $( \hat{\mathbf{x}}_{i})_{i=1}^{N}$ to be the solution computed by a given algorithm. The Lie-theoretic RPE (RPE-L) is defined as
\begin{equation}\label{eq:R-RPE}
    \text{RPE-L}\triangleq\sqrt{\frac{1}{M}\sum_{\left(i,j\right)\in\mathcal{E}}{\left\Vert \Logx{\fatone}{\hat{\mathbf{z}}_{ij}^{-1}\boxplus\mathbf{z}_{ij}^{\diamond}}\right\Vert}_{2}^{2}},
\end{equation}
where $\hat{\mathbf{z}}_{ij}\triangleq\hat{\mathbf{x}}_{i}^{-1}\boxplus\hat{\mathbf{x}}_{j}$ and $\mathbf{z}_{ij}^{\diamond}\triangleq(\mathbf{x}_{i}^{\diamond})^{-1}\boxplus\mathbf{x}_{j}^{\diamond}$. Now, let $(\hat{\mathbf{t}}_{i},\hat{\theta}_{i})$ and $(\mathbf{t}_{i}^{\diamond},\theta_{i}^{\diamond})$ denote the translations and rotations corresponding to $\hat{\mathbf{x}}_{i}$ and $\mathbf{x}_{i}^{\diamond}$, respectively. The Euclidean RPE (RPE-E) is defined as
\begin{equation}\label{eq:E-RPE}
    \text{RPE-E}\triangleq\sqrt{\frac{1}{M}\sum_{\left(i,j\right)\in\mathcal{E}}\left(\left\Vert \hat{\mathbf{t}}_{ij}-\mathbf{t}_{ij}^{\diamond}\right\Vert ^{2}+ d(\hat{\theta}_{ij},\theta_{ij}^{\diamond})^{2}\right)},
\end{equation}
where $\hat{\mathbf{t}}_{ij}\triangleq R^{\top}(\hat{\theta}_{i})\left(\hat{\mathbf{t}}_{j}-\hat{\mathbf{t}}_{i}\right)$, $\mathbf{t}_{ij}^{\diamond}\triangleq R^{\top}\left(\theta_{i}^{\diamond}\right)\left(\mathbf{t}_{j}^{\diamond}-\mathbf{t}_{i}^{\diamond}\right)$, $d\left(\theta_{1},\theta_{2}\right)$ is the minimal angle between $\theta_{1}$ and $\theta_{2}$, and 
\begin{equation}
    R\left(\theta\right)\triangleq\left[\begin{smallmatrix}
        \cos\left(\theta\right) & -\sin\left(\theta\right)\\
        \sin\left(\theta\right) & \cos\left(\theta\right)
        \end{smallmatrix}\right].
\end{equation}

For evaluation, the variance scaling parameter, $\sigma_{w}$, was varied from $10^{-5}$ to $10^{-2}$, which equated to mean Euclidean covariances with standard deviations ranging from $7.26\cdot 10^{-3}$ to $2.29\cdot 10^{-1}$ meters for translations, and from $4.05\cdot 10^{-1}$ to $12.81$ degrees for rotations. We anchor $\mathbf{x}_{1}\triangleq\fatone$ for all three algorithms. The initial iterate $\XZero$ is computed using the chordal relaxation~\cite{martinec2007robust} method; though not necessary for convergence
of Algorithm~\ref{alg:RTR}, it is the default for both SE-Sync and Cartan-Sync, so we implement it to provide a fair comparison. Algorithm~\ref{alg:RTR} was configured with parameters $\varepsilon_{g}=\epsilong$, $\Delta_{0}=\deltazero$, $\bar{\Delta}=\deltabar$, $\rho^{\prime}=\rhoprime$, and the inner tCG algorithm was implemented with parameters $\kappa=0.05$, $\theta=0.25$, per the notation in~\cite[Section~6.5]{boumal2023introduction}. 

\subsection{Evaluation results}
\begin{figure}
    \vspace{2mm}
    \centering
    \includegraphics[scale=0.8]{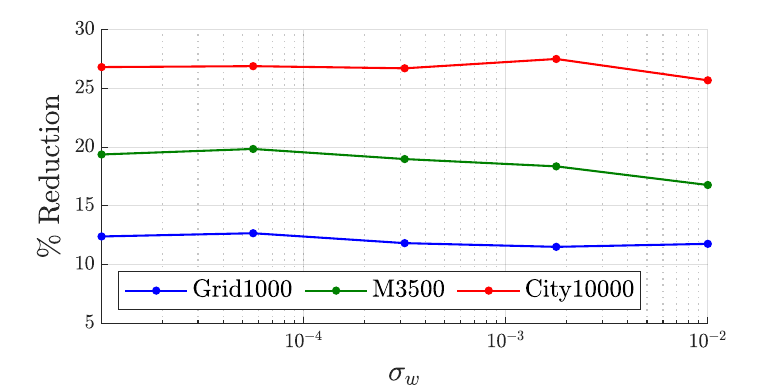}
    \caption{Percent reduction in Lie-theoretic RPE for the solutions computed by Algorithm~\ref{alg:RTR} relative to Cartan-Sync and SE-Sync. Reduction in Euclidean RPE was omitted due to it being indistinguishable from the Riemannian case. We see greater than~$10$\% decrease for the Grid1000 dataset over the entire noise regime, and greater than $15$\% \& $25$\% for the M3500 and City10000 datasets, respectively. In all cases, the improvement in accuracy attained by Algorithm~\ref{alg:RTR} grows with the number of vertices and edges present in a graph.}\label{fig:improvement}
\end{figure}

Algorithm~\ref{alg:RTR} converged to an approximate stationary point in all of the $\numtrials$ pose graphs. 
The RPEs computed for each run according to~\eqref{eq:R-RPE} and~\eqref{eq:E-RPE} are included in Table~\ref{tbl:comparison}, alongside  the percent reduction in RPE attained by Algorithm~\ref{alg:RTR} for each run, which is plotted in Figure~\ref{fig:improvement}. SE-Sync and Cartan-Sync computed identical solutions for each dataset, and exhibited a notable estimation bias across the entire noise spectrum, owing to the assumption of isotropic noise and the resulting approximation error. As shown in Figure~\ref{fig:improvement}, Algorithm~\ref{alg:RTR} demonstrated a consistent reduction in RPE. In fact, the gap in RPE increases with the number of vertices and edges in each graph, highlighting the scalability of our proposed solution.
\begin{table*}[t]
    \vspace{2mm}
  \centering
    \caption{Results of the 2D PGO dataset evaluation. RPE and percent reduction in RPE attained by Algorithm~\ref{alg:RTR} are shown on the right.}\label{tbl:comparison}
  \begin{tabular}{|c|c|c|c|c|c|c|c|}
  \cline{3-8} \cline{4-8} \cline{5-8} \cline{6-8} \cline{7-8} \cline{8-8} 
  \multicolumn{1}{c}{} &  & \multicolumn{2}{c|}{SE-Sync~\cite{rosen2019se}} & \multicolumn{2}{c|}{Cartan-Sync~\cite{briales2017cartan}} & \multicolumn{2}{c|}{Algorithm~\ref{alg:RTR} [ours] (\textcolor{imp_color}{\% Reduction})}\tabularnewline
  \hline 
  Dataset & $\sigma_{w}$ & RPE-L & RPE-E & RPE-L & RPE-E & RPE-L & RPE-E \tabularnewline
  \hline 
  Grid1000 & $1.0\cdot10^{-5}$ & $6.2\cdot10^{-3}$ & $1.2\cdot10^{-2}$ & $6.2\cdot10^{-3}$ & $1.2\cdot10^{-2}$ & $5.4\cdot10^{-3}$ (\textcolor{imp_color}{$\unaryminus12.4$\%}) & $1.1\cdot10^{-2}$ (\textcolor{imp_color}{$\unaryminus12.4$\%}) \tabularnewline
\hline
Grid1000 & $5.6\cdot10^{-5}$ & $1.5\cdot10^{-2}$ & $2.9\cdot10^{-2}$ & $1.5\cdot10^{-2}$ & $2.9\cdot10^{-2}$ & $1.3\cdot10^{-2}$ (\textcolor{imp_color}{$\unaryminus12.7$\%}) & $2.6\cdot10^{-2}$ (\textcolor{imp_color}{$\unaryminus12.7$\%}) \tabularnewline
\hline
Grid1000 & $3.2\cdot10^{-4}$ & $3.5\cdot10^{-2}$ & $7.1\cdot10^{-2}$ & $3.5\cdot10^{-2}$ & $7.1\cdot10^{-2}$ & $3.1\cdot10^{-2}$ (\textcolor{imp_color}{$\unaryminus11.8$\%}) & $6.2\cdot10^{-2}$ (\textcolor{imp_color}{$\unaryminus11.8$\%}) \tabularnewline
\hline
Grid1000 & $1.8\cdot10^{-3}$ & $7.9\cdot10^{-2}$ & $1.6\cdot10^{-1}$ & $7.9\cdot10^{-2}$ & $1.6\cdot10^{-1}$ & $7.0\cdot10^{-2}$ (\textcolor{imp_color}{$\unaryminus11.5$\%}) & $1.4\cdot10^{-1}$ (\textcolor{imp_color}{$\unaryminus11.5$\%}) \tabularnewline
\hline
Grid1000 & $1.0\cdot10^{-2}$ & $1.9\cdot10^{-1}$ & $3.9\cdot10^{-1}$ & $1.9\cdot10^{-1}$ & $3.9\cdot10^{-1}$ & $1.7\cdot10^{-1}$ (\textcolor{imp_color}{$\unaryminus11.8$\%}) & $3.4\cdot10^{-1}$ (\textcolor{imp_color}{$\unaryminus11.7$\%}) \tabularnewline
\hline
M3500 & $1.0\cdot10^{-5}$ & $5.4\cdot10^{-3}$ & $1.1\cdot10^{-2}$ & $5.4\cdot10^{-3}$ & $1.1\cdot10^{-2}$ & $4.4\cdot10^{-3}$ (\textcolor{imp_color}{$\unaryminus19.4$\%}) & $8.7\cdot10^{-3}$ (\textcolor{imp_color}{$\unaryminus19.4$\%}) \tabularnewline
\hline
M3500 & $5.6\cdot10^{-5}$ & $1.3\cdot10^{-2}$ & $2.6\cdot10^{-2}$ & $1.3\cdot10^{-2}$ & $2.6\cdot10^{-2}$ & $1.0\cdot10^{-2}$ (\textcolor{imp_color}{$\unaryminus19.8$\%}) & $2.1\cdot10^{-2}$ (\textcolor{imp_color}{$\unaryminus19.8$\%}) \tabularnewline
\hline
M3500 & $3.2\cdot10^{-4}$ & $3.1\cdot10^{-2}$ & $6.2\cdot10^{-2}$ & $3.1\cdot10^{-2}$ & $6.2\cdot10^{-2}$ & $2.5\cdot10^{-2}$ (\textcolor{imp_color}{$\unaryminus19.0$\%}) & $5.0\cdot10^{-2}$ (\textcolor{imp_color}{$\unaryminus19.0$\%}) \tabularnewline
\hline
M3500 & $1.8\cdot10^{-3}$ & $7.4\cdot10^{-2}$ & $1.5\cdot10^{-1}$ & $7.4\cdot10^{-2}$ & $1.5\cdot10^{-1}$ & $6.0\cdot10^{-2}$ (\textcolor{imp_color}{$\unaryminus18.4$\%}) & $1.2\cdot10^{-1}$ (\textcolor{imp_color}{$\unaryminus18.4$\%}) \tabularnewline
\hline
M3500 & $1.0\cdot10^{-2}$ & $1.7\cdot10^{-1}$ & $3.4\cdot10^{-1}$ & $1.7\cdot10^{-1}$ & $3.4\cdot10^{-1}$ & $1.4\cdot10^{-1}$ (\textcolor{imp_color}{$\unaryminus16.8$\%}) & $2.9\cdot10^{-1}$ (\textcolor{imp_color}{$\unaryminus16.8$\%}) \tabularnewline
\hline
City10k & $1.0\cdot10^{-5}$ & $4.9\cdot10^{-3}$ & $9.7\cdot10^{-3}$ & $4.9\cdot10^{-3}$ & $9.7\cdot10^{-3}$ & $3.6\cdot10^{-3}$ (\textcolor{imp_color}{$\unaryminus26.8$\%}) & $7.1\cdot10^{-3}$ (\textcolor{imp_color}{$\unaryminus26.8$\%}) \tabularnewline
\hline
City10k & $5.6\cdot10^{-5}$ & $1.2\cdot10^{-2}$ & $2.3\cdot10^{-2}$ & $1.2\cdot10^{-2}$ & $2.3\cdot10^{-2}$ & $8.5\cdot10^{-3}$ (\textcolor{imp_color}{$\unaryminus26.9$\%}) & $1.7\cdot10^{-2}$ (\textcolor{imp_color}{$\unaryminus26.9$\%}) \tabularnewline
\hline
City10k & $3.2\cdot10^{-4}$ & $2.8\cdot10^{-2}$ & $5.5\cdot10^{-2}$ & $2.8\cdot10^{-2}$ & $5.5\cdot10^{-2}$ & $2.0\cdot10^{-2}$ (\textcolor{imp_color}{$\unaryminus26.7$\%}) & $4.0\cdot10^{-2}$ (\textcolor{imp_color}{$\unaryminus26.7$\%}) \tabularnewline
\hline
City10k & $1.8\cdot10^{-3}$ & $6.6\cdot10^{-2}$ & $1.3\cdot10^{-1}$ & $6.6\cdot10^{-2}$ & $1.3\cdot10^{-1}$ & $4.8\cdot10^{-2}$ (\textcolor{imp_color}{$\unaryminus27.5$\%}) & $9.5\cdot10^{-2}$ (\textcolor{imp_color}{$\unaryminus27.5$\%}) \tabularnewline
\hline
City10k & $1.0\cdot10^{-2}$ & $1.6\cdot10^{-1}$ & $3.1\cdot10^{-1}$ & $1.6\cdot10^{-1}$ & $3.1\cdot10^{-1}$ & $1.2\cdot10^{-1}$ (\textcolor{imp_color}{$\unaryminus25.7$\%}) & $2.3\cdot10^{-1}$ (\textcolor{imp_color}{$\unaryminus25.7$\%}) \tabularnewline
\hline
  \end{tabular}
  \end{table*}

\section{Conclusion}\label{sec:conclusion}
We presented a novel algorithm for planar PGO derived from a realistic, Lie-theoretic model for uncertainty in sensor measurements. The proposed algorithm was proven to converge in finite-time to approximate first-order stationary points under any initialization, while requiring no additional assumptions about the problem. 
Numerically, the proposed algorithm showed significantly improved accuracy
over the state of the art, and future work will extend the algorithm to the 3D case and distributed/asynchronous implementations.

\pagebreak
\begin{appendices}

\section{Algebraic Construction}\label{app:construction}
Given an orthonormal basis $\{\mathbf{i},\mathbf{j},\mathbf{k}\}$, a rotation in the plane is characterized by a rotation angle $\theta\in(-\pi,\pi]$ about the $\mathbf{k}$ axis. In standard form, we can write the \textit{planar unit quaternion}\footnote{It is noted that a planar unit quaternion is a standard Hamiltonian unit quaternion restricted to a rotation about the $\mathbf{k}$-axis, i.e., $\mathbf{q}\in\mathbb{H},~\mathbf{q}=w+x\mathbf{i}+y\mathbf{j}+z\mathbf{k}$, with $\norm{\mathbf{q}}=1$ and $x=y=0$.} $\mathbf{q}\in\mathbb{S}^{1}$ corresponding to this rotation as
\begin{equation}\label{eq:puq_rotation}
    \mathbf{q}=\cos\left(\nicefrac{\theta}{2}\right)+\mathbf{k}\sin\left(\nicefrac{\theta}{2}\right)=r_{0}+\mathbf{k}r_{1},
\end{equation}
or, in vector form, ${\mathbf{q}=\left[q_{0}, q_{1}\right]}^{\top}$. Let \quotes{$\otimes$} denote the Hamilton product~\cite{hamilton1848xi} under the convention $\mathbf{i}^{2}=\mathbf{j}^{2}=\mathbf{k}^{2}=\mathbf{i}\mathbf{j}\mathbf{k}=-1$. Then, performing the Hamiltonian multiplication of two planar quaternions, denoted $\mathbf{r},\mathbf{s}$, yields
\begin{equation}
    \mathbf{r}\otimes\mathbf{s}=\left(r_{0}+\mathbf{k}r_{1}\right)\left(s_{0}+\mathbf{k}s_{1}\right)=r_{0}s_{0}-r_{1}s_{1}+\mathbf{k}\left(r_{1}s_{0}+r_{0}s_{1}\right)
\end{equation}
\noindent In matrix-vector form, the operation can be written as
\begin{align*}
\mathbf{r}\otimes\mathbf{s}&=\left[\begin{array}{cc}
r_{0} & -r_{1}\\
r_{1} & r_{0}
\end{array}\right]\left[\begin{array}{c}
s_{0}\\
s_{1}
\end{array}\right]=\left[\begin{array}{cc}
s_{0} & -s_{1}\\
s_{1} & s_{0}
\end{array}\right]\left[\begin{array}{c}
r_{0}\\
r_{1}
\end{array}\right].
\end{align*}
A planar rigid motion is characterized by a translation, denoted $\mathbf{t}=t_{x}\mathbf{i}+t_{y}\mathbf{j}$, and a rotation about the $\mathbf{k}$ axis by an angle $\theta\in(-\pi,\pi]$. This can be written in $\mathbb{R}^{3}$ as the Euclidean vector $\mathbf{p}=[\mathbf{t}^{\transpose},\theta]^{\transpose}$. The \textit{planar unit dual quaternion} (PUDQ) parameterization of this motion is given by $\mathbf{x}=\mathbf{x}_{r}+\epsilon\mathbf{x}_{d}$, where $\epsilon$ is a \textit{dual number} satisfying $\epsilon^{2}=0,\epsilon\neq0$. The \textit{real} part of $\mathbf{x}$, denoted $\mathbf{x}_{r}\in\mathbb{S}^{1}$, is a planar unit quaternion of the form 
\begin{equation}\label{eq:xr_def}
    \mathbf{x}_{r}=\cos\left(\nicefrac{\theta}{2}\right)+\sin\left(\nicefrac{\theta}{2}\right)\mathbf{k}=r_{0}+\mathbf{k}r_{1}.
\end{equation}
The \textit{dual} part of $\mathbf{x}$, denoted $\mathbf{x}_{d}\in\mathbb{R}^{2}$, is given by
\begin{equation}\label{eq:xd_def}
    \mathbf{x}_{d}=\frac{1}{2}\mathbf{t}\otimes \mathbf{x}_{r}=\frac{1}{2}\left(t_{x}\mathbf{i}+t_{y}\mathbf{j}\right)\left(r_{0}+\mathbf{k}r_{1}\right)=\frac{1}{2}\left(\left(t_{x}r_{0}+t_{y}r_{1}\right)\mathbf{i}+\left(t_{y}r_{0}-t_{x}r_{1}\right)\mathbf{j}\right).
\end{equation}
In matrix-vector form,~\eqref{eq:xd_def} can be rewritten as
\begin{equation}\label{eq:xd_matvec}
    \mathbf{x}_{d}=\frac{1}{2}\left[\begin{array}{cc}
        t_{x} & t_{y}\\
        t_{y} & -t_{x}
        \end{array}\right]\left[\begin{array}{c}
        r_{0}\\
        r_{1}
        \end{array}\right]=\frac{1}{2}\left[\begin{array}{cc}
        r_{0} & r_{1}\\
        -r_{1} & r_{0}
        \end{array}\right]\left[\begin{array}{c}
        t_{x}\\
        t_{y}
        \end{array}\right].
\end{equation}
In vector form, a PUDQ can be expressed in terms of the bases $\{\mathbf{i}, \mathbf{k}, \epsilon\mathbf{i}, \epsilon\mathbf{j}\}$ as
\begin{equation}
    \mathbf{x}=x_{0}+\mathbf{k}x_{1}+\epsilon\left(\mathbf{i}x_{2}+\mathbf{j}x_{3}\right)={\left[
        x_{0},~x_{1},~x_{2},~x_{3}\right]}^{\top}.
\end{equation}
Equivalently, we can write $\mathbf{x}=[\mathbf{x}_{r}^{\transpose},\mathbf{x}_{d}^{\transpose}]$. Given two PUDQs, $\mathbf{x}={\left[
    x_{0},~x_{1},~x_{2},~x_{3}\right]}^{\top}$ and $\mathbf{y}={\left[
    y_{0},~y_{1},~y_{2},~y_{3}\right]}^{\top}$, we can compute the composition operation “$\boxplus$” by applying Hamiltonian multiplication, which yields
\begin{align}
    \mathbf{x}\boxplus\mathbf{y} & =\left(x_{0}+\mathbf{k}x_{1}+\epsilon\left(\mathbf{i}x_{2}+\mathbf{j}x_{3}\right)\right)\left(y_{0}+\mathbf{k}y_{1}+\epsilon\left(\mathbf{i}y_{2}+\mathbf{j}y_{3}\right)\right)\\
    & =\left(x_{0}y_{0}-x_{1}y_{1}\right)+\mathbf{k}\left(x_{0}y_{1}+x_{1}y_{0}\right)+\epsilon\left(\mathbf{i}\left(x_{0}y_{2}-x_{1}y_{3}+x_{2}y_{0}+x_{3}y_{1}\right)+\mathbf{j}\left(x_{0}y_{3}+x_{1}y_{2}-x_{2}y_{1}+x_{3}y_{0}\right)\right).\label{eq:x_boxplus_y_def}
\end{align}
From~\eqref{eq:x_boxplus_y_def}, we can deduce the identity PUDQ, denoted $\fatone$, to be $\fatone=[1,0,0,0]^{\transpose}$, so that $\fatone\boxplus\mathbf{x}=\mathbf{x}\boxplus\fatone=\mathbf{x}$. Moreover, the inverse of a PUDQ $\mathbf{x}$, denoted $\mathbf{x}^{-1}$, is given by $\mathbf{x}^{-1}=[x_{0},-x_{1},-x_{2},-x_{3}]^{\transpose}$, so that $\mathbf{x}\boxplus\mathbf{x}^{-1}=\mathbf{x}^{-1}\boxplus\mathbf{x}=\fatone$. The operation described by~\eqref{eq:x_boxplus_y_def} is equivalent to the matrix-vector multiplication(s)
\begin{equation}\label{eq:Q_L_R}
    \mathbf{x}\boxplus\mathbf{y}=\underset{Q_{L}\left(\mathbf{x}\right)}{\underbrace{\left[\begin{array}{cccc}
        x_{0} & -x_{1} & 0 & 0\\
        x_{1} & x_{0} & 0 & 0\\
        x_{2} & x_{3} & x_{0} & -x_{1}\\
        x_{3} & -x_{2} & x_{1} & x_{0}
        \end{array}\right]}}\underset{\mathbf{y}}{\underbrace{\left[\begin{array}{c}
        y_{0}\\
        y_{1}\\
        y_{2}\\
        y_{3}
        \end{array}\right]}}=\underset{Q_{R}\left(\mathbf{y}\right)}{\underbrace{\left[\begin{array}{cccc}
        y_{0} & -y_{1} & 0 & 0\\
        y_{1} & y_{0} & 0 & 0\\
        y_{2} & -y_{3} & y_{0} & y_{1}\\
        y_{3} & y_{2} & -y_{1} & y_{0}
        \end{array}\right]}}\underset{\mathbf{x}}{\underbrace{\left[\begin{array}{c}
        x_{0}\\
        x_{1}\\
        x_{2}\\
        x_{3}
        \end{array}\right]}},
\end{equation}
where we have implicitly defined the left and right-handed matrix-valued left and right-hand composition mappings $Q_{L}:\M\rightarrow\mathbb{R}^{4\times4}$ and $Q_{R}:\mathcal{M}\rightarrow\mathbb{R}^{4\times4}$. 
Using $Q_{L}$, we define $Q_{RL}^{-}\left(\mathbf{x}\right):\mathcal{M}\rightarrow\mathbb{R}^{4\times4}$ such that
\begin{equation}\label{eq:Q_R_L_M}
    \mathbf{x}\boxplus\mathbf{y}^{-1}	=\underset{Q_{L}\left(\mathbf{x}\right)}{\underbrace{\left[\begin{array}{cccc}
    x_{0} & -x_{1} & 0 & 0\\
    x_{1} & x_{0} & 0 & 0\\
    x_{2} & x_{3} & x_{0} & -x_{1}\\
    x_{3} & -x_{2} & x_{1} & x_{0}
    \end{array}\right]}}\underset{\mathbf{y}^{-1}}{\underbrace{\left[\begin{array}{c}
    y_{0}\\
    -y_{1}\\
    -y_{2}\\
    -y_{3}
    \end{array}\right]}}=\underset{Q_{RL}^{-}\left(\mathbf{x}\right)}{\underbrace{\left[\begin{array}{cccc}
    x_{0} & x_{1} & 0 & 0\\
    x_{1} & -x_{0} & 0 & 0\\
    x_{2} & -x_{3} & -x_{0} & x_{1}\\
    x_{3} & x_{2} & -x_{1} & -x_{0}
    \end{array}\right]}}\underset{\mathbf{y}}{\underbrace{\left[\begin{array}{c}
    y_{0}\\
    y_{1}\\
    y_{2}\\
    y_{3}
    \end{array}\right]}}
\end{equation}
and $Q_{L}^{--}:\M\rightarrow\mathbb{R}^{4\times4}$ such that
\begin{equation}\label{eq:Q_L_MM}
    \mathbf{x}^{-1}\boxplus\mathbf{y}^{-1}	=\underset{Q_{L}\left(\mathbf{x}^{-1}\right)}{\underbrace{\left[\begin{array}{cccc}
        x_{0} & x_{1} & 0 & 0\\
        -x_{1} & x_{0} & 0 & 0\\
        -x_{2} & -x_{3} & x_{0} & x_{1}\\
        -x_{3} & x_{2} & -x_{1} & x_{0}
        \end{array}\right]}}\underset{\mathbf{y}^{-1}}{\underbrace{\left[\begin{array}{c}
        y_{0}\\
        -y_{1}\\
        -y_{2}\\
        -y_{3}
        \end{array}\right]}}=\underset{Q_{L}^{--}\left(\mathbf{x}\right)}{\underbrace{\left[\begin{array}{cccc}
        x_{0} & -x_{1} & 0 & 0\\
        -x_{1} & -x_{0} & 0 & 0\\
        -x_{2} & x_{3} & -x_{0} & -x_{1}\\
        -x_{3} & -x_{2} & x_{1} & -x_{0}
        \end{array}\right]}}\underset{\mathbf{y}}{\underbrace{\left[\begin{array}{c}
        y_{0}\\
        y_{1}\\
        y_{2}\\
        y_{3}
        \end{array}\right]}}.
\end{equation}
Using $Q_{R}$, we define the mapping $Q_{LR}^{-}:\M\rightarrow\mathbb{R}^{4\times4}$ such that
\begin{equation}\label{eq:Q_L_R_M}
    \mathbf{x}^{-1}\boxplus\mathbf{y}	=\underset{Q_{R}\left(\mathbf{y}\right)}{\underbrace{\left[\begin{array}{cccc}
        y_{0} & -y_{1} & 0 & 0\\
        y_{1} & y_{0} & 0 & 0\\
        y_{2} & -y_{3} & y_{0} & y_{1}\\
        y_{3} & y_{2} & -y_{1} & y_{0}
        \end{array}\right]}}\underset{\mathbf{x}^{-1}}{\underbrace{\left[\begin{array}{c}
        x_{0}\\
        -x_{1}\\
        -x_{2}\\
        -x_{3}
    \end{array}\right]}}=\underset{Q_{LR}^{-}\left(\mathbf{y}\right)}{\underbrace{\left[\begin{array}{cccc}
        y_{0} & y_{1} & 0 & 0\\
        y_{1} & -y_{0} & 0 & 0\\
        y_{2} & y_{3} & -y_{0} & -y_{1}\\
        y_{3} & -y_{2} & y_{1} & -y_{0}
        \end{array}\right]}}\underset{\mathbf{x}}{\underbrace{\left[\begin{array}{c}
        x_{0}\\
        x_{1}\\
        x_{2}\\
        x_{3}
        \end{array}\right]}},
\end{equation}
and $Q_{R}^{--}:\mathcal{M}\rightarrow\mathbb{R}^{4\times4}$ such that
\begin{equation}\label{eq:Q_R_MM}
    \mathbf{x}^{-1}\boxplus\mathbf{y}^{-1}	=\underset{Q_{R}\left(\mathbf{y}^{-1}\right)}{\underbrace{\left[\begin{array}{cccc}
    y_{0} & y_{1} & 0 & 0\\
    -y_{1} & y_{0} & 0 & 0\\
    -y_{2} & y_{3} & y_{0} & -y_{1}\\
    -y_{3} & -y_{2} & y_{1} & y_{0}
    \end{array}\right]}}\underset{\mathbf{x}^{-1}}{\underbrace{\left[\begin{array}{c}
    x_{0}\\
    -x_{1}\\
    -x_{2}\\
    -x_{3}
    \end{array}\right]}}=\underset{Q_{R}^{--}\left(\mathbf{y}\right)}{\underbrace{\left[\begin{array}{cccc}
    y_{0} & -y_{1} & 0 & 0\\
    -y_{1} & -y_{0} & 0 & 0\\
    -y_{2} & -y_{3} & -y_{0} & y_{1}\\
    -y_{3} & y_{2} & -y_{1} & -y_{0}
    \end{array}\right]}}\underset{\mathbf{x}}{\underbrace{\left[\begin{array}{c}
    x_{0}\\
    x_{1}\\
    x_{2}\\
    x_{3}
    \end{array}\right]}}
\end{equation}
The maps $Q_{L}$ and $Q_{R}$ additionally yield the definitions for $Q_{LL}^{-}\left(\mathbf{x}\right)\triangleq Q_{L}\left(\mathbf{x}^{-1}\right)$ and $Q_{RR}^{-}\left(\mathbf{x}\right)\triangleq Q_{R}\left(\mathbf{x}^{-1}\right)$ such that $\mathbf{x}^{-1}\boxplus\mathbf{y}=Q_{LL}^{-}\left(\mathbf{x}\right)\mathbf{y}$ and $\mathbf{x}\boxplus\mathbf{y}^{-1}=Q_{RR}^{-}\left(\mathbf{y}\right)\mathbf{x}$.

\pagebreak
\section{Riemannian Geometry of the Planar Unit Dual Quaternion Manifold}\label{app:riemannian}
In this appendix, we provide derivations relating to the Riemannian geometry of the PUDQ manifold $\M$ and its product manifold extension $\MN$. For a general coverage of these topics, we refer the reader to~\cite{boumal2023introduction}. 

\subsection{Embedded Submanifolds}\label{app:submanifolds}
The set of all PUDQs forms a smooth manifold, denoted $\M$. In this work, we employ an embedding of $\mathcal{M}$ in the ambient Euclidean space $\mathbb{R}^{4}$ with the inner product $\left\langle u,w\right\rangle =u^{\transpose}w$ and induced Euclidean norm $\norm{u}=\sqrt{u^{\transpose}u}$ for all $u,w\in\mathbb{R}^{4}$. This embedding yields the coordinatized definition for $\M$ given by
\begin{equation}\label{eq:M_embedding_app}
    \M\triangleq\left\{\mathbf{x}\in\mathbb{R}^{4}\mid h\left(\mathbf{x}\right)=\mathbf{x}^{\transpose}\Ptilde\mathbf{x}-1=0\right\}\subset\mathbb{R}^{4},
\end{equation}
where $h(\mathbf{x})$ is the \textit{defining function}~\cite{boumal2023introduction} for $\M$ and $\Ptilde\in\mathbb{R}^{4\times4}$ defined as
\begin{equation}\label{eq:Ptilde_def}
    \Ptilde\triangleq\left[\begin{array}{cc}
        I_{2} & \boldsymbol{0}_{2\times2}\\
        \boldsymbol{0}_{2\times2} & \boldsymbol{0}_{2\times2}
        \end{array}\right],
\end{equation}
where $I_{2}\in\mathbb{R}^{2\times2}$ denotes an identity matrix and $\boldsymbol{0}_{2\times2}\in\mathbb{R}^{2\times2}$ denotes a matrix of zeroes. By~\eqref{eq:M_embedding_app}, we have $\M=\mathbb{S}^{1}\rtimes\mathbb{R}^{2}\subset\mathbb{R}^{4}$. We now extend~\eqref{eq:M_embedding_app} to the $N$-fold PUDQ product manifold $\MN\triangleq\M\times\M\times\cdots\times\M=(\mathbb{S}^{1}\rtimes\mathbb{R}^{2})^{N}$. For notational convenience, we define the operator $\vect(\cdot)$ such that 
\begin{equation}\label{eq:vect_app}
    \vect((\mathbf{x}_{i})_{i=1}^{N})\triangleq[\mathbf{x}_{1}^{\transpose},\mathbf{x}_{2}^{\transpose},\dots,\mathbf{x}_{N}^{\transpose}]^{\transpose},
\end{equation}
with each $\mathbf{x}_{i}\in\M$. Since $(\mathbb{S}^{1}\rtimes\mathbb{R}^{2})^{N}\subset\mathbb{R}^{4\times N}\cong\mathbb{R}^{4N}$, we embed $\MN$ in $\mathbb{R}^{4N}$ via the coordinatized definition
\begin{equation}\label{eq:MN_embedding_app}
    \MN\triangleq\left\{\vect(\left(\mathbf{x}_{i}\right)_{i=1}^{N})\mid\mathbf{x}_{i}\in\M\text{~for~all~}i\in\left\{1,\ldots,N\right\}\right\}\subset\mathbb{R}^{4N},
\end{equation}
with $\vect(\cdot)$ given by~\eqref{eq:vect_app}. For~$\mathcal{X}, \mathcal{Y} \in \M^N$, the embedding in~\eqref{eq:MN_embedding_app} lets us write~$\mathcal{X} = \vect((\mathbf{x}_{i})_{i=1}^{N})$ and $\mathcal{Y}=\vect((\mathbf{y}_{i})_{i=1}^{N})$, where~$\mathbf{x}_i, \mathbf{y}_i \in \M$ for each~$i$. Furthermore,~\eqref{eq:MN_embedding_app} admits natural extensions to $\MN$ of the identity $\fatone^{N}=[\fatone^{\transpose},\fatone^{\transpose},\ldots,\fatone^{\transpose}]^{\transpose}$, inverse $\X^{-1}=[\mathbf{x}_{1}^{-\transpose},\mathbf{x}_{2}^{-\transpose},\dots,\mathbf{x}_{N}^{-\transpose}]^{\transpose}$, and, for $\X,\Y\in\MN$, the product $\X\boxplus\Y=\vect((\mathbf{x}_{i}\boxplus\mathbf{y}_{i})_{i=1}^{N})$.

\subsection{Tangent Space and Projection Operators}\label{app:tangent_space}
The tangent space of $\M$ at a point $\mathbf{x}\in\M$, denoted $\TxM$, is the local, Euclidean linearization of $\M$ about $\mathbf{x}$. It is defined as $\TxM\triangleq\ker(\mathrm{D}h(\mathbf{x}))$, where $h(\mathbf{x})$ is the defining function from~\eqref{eq:M_embedding_app}, and $\mathrm{D}h(\mathbf{x})[v]$ is the directional derivative of $h$ along $v\in\mathbb{R}^{4}$ at $\mathbf{x}$. We compute $\mathrm{D}h(\mathbf{x})[v]$ from the definition given in~\cite{boumal2023introduction} as
\begin{align}
    \mathrm{D}h(\mathbf{x})[v]&=\lim_{t\rightarrow0}\frac{h(\mathbf{x}+tv)-h(\mathbf{x})}{t} \\
    &=\lim_{t\rightarrow0}\frac{\left(\mathbf{x}+tv\right)^{\transpose}\Ptilde\left(\mathbf{x}+tv\right)-\mathbf{x}^{\transpose}\Ptilde\mathbf{x}}{t} \\
    &=2\mathbf{x}^{\transpose}\Ptilde v.\label{eq:Dhxv}
\end{align}
Since $\mathcal{T}_{\mathbf{x}}\mathcal{M}\triangleq\ker\left(\mathrm{D}h\left(\mathbf{x}\right)\right)$, it follows from~\eqref{eq:Dhxv} that $\mathbf{x}^{\transpose}\Ptilde v=0$ for all $v\in\mathcal{T}_{\mathbf{x}}\mathcal{M}$. Therefore, $\TxM$ is given by
\begin{equation}\label{eq:TxM_def}
    \mathcal{T}_{\mathbf{x}}\mathcal{M}	=\left\{ v\in\mathbb{R}^{4}\mid\mathbf{x}^{\transpose}\Ptilde v=0\right\}.
\end{equation}
We can then derive the orthogonal projection matrix, denoted $\mathcal{P}_{\mathbf{x}}$, by identifying from~\eqref{eq:TxM_def} that, for any $u\in\mathbb{R}^{4}$, it holds that
\begin{equation}
    \mathrm{proj}_{\mathbf{x}}u=u-\mathrm{proj}_{\Ptilde\mathbf{x}}u,
\end{equation}
where
\begin{equation}
    \mathrm{proj}_{\Ptilde\mathbf{x}}u=(\mathbf{x}^{\transpose}\Ptilde u)\frac{\Ptilde\mathbf{x}}{\Vert\Ptilde\mathbf{x}\Vert_{2}}.
\end{equation}
Since $\Vert\Ptilde\mathbf{x}\Vert_{2}=x_{0}^{2}+x_{1}^{2}=1$ for all $\mathbf{x}=[x_0,x_1,x_2,x_3]^{\transpose}\in\mathcal{M}$, we have $\mathrm{proj}_{\Ptilde\mathbf{x}}u=\Ptilde\mathbf{x}\mathbf{x}^{\transpose}\Ptilde u$. Therefore,
\begin{equation}\label{eq:projx_u}
    \mathrm{proj}_{\mathbf{x}}u	=u-\Ptilde\mathbf{x}\mathbf{x}^{\transpose}\Ptilde u=(I_{4}-\Ptilde\mathbf{x}\mathbf{x}^{\transpose}\Ptilde)u,
\end{equation}
where $I_{m}\in\mathbb{R}^{m\times m}$ denotes the identity matrix. Equation~\eqref{eq:projx_u} yields $\mathrm{proj}_{\mathbf{x}}u=\mathcal{P}_{\mathbf{x}}u$, with $\mathcal{P}_{\mathbf{x}}\in\mathbb{R}^{4\times 4}$ given by the symmetric, idempotent matrix
\begin{equation}\label{eq:Px_def}
    \mathcal{P}_{\mathbf{x}}=I_{4}-\Ptilde\mathbf{x}\mathbf{x}^{\transpose}\Ptilde,
\end{equation}
with $\Ptilde$ given by~\eqref{eq:Ptilde_def}. We also have the normal projection operator, denoted $\mathcal{P}_{\mathbf{x}}^{\perp}\in\mathbb{R}^{4\times4}$, given by
\begin{equation}\label{eq:Px_norm_def}
    \mathcal{P}_{\mathbf{x}}^{\perp}=I_4-\mathcal{P}_{\mathbf{x}}=\Ptilde\mathbf{x}\mathbf{x}^{\transpose}\Ptilde.
\end{equation}
Furthermore, the embedding in~\eqref{eq:MN_embedding_app} gives the orthogonal projector onto $\TXMN$, denoted $\mathcal{P}_{\X}\in\mathbb{R}^{4N\times4N}$, to be
\begin{equation}\label{eq:PXN_def}
    \PX=\mathrm{diag}(\{\mathcal{P}_{\mathbf{x}_{i}}\mid i\in\{1,\ldots,N\}\}),
\end{equation}
with $\mathcal{P}_{\mathbf{x}_{i}}$ given by~\eqref{eq:Px_def}. Finally, the normal projector onto $\mathcal{T}_{\mathcal{X}}^{\perp}\mathcal{M}^{N}$, denoted $\PX^{\perp}\in\mathbb{R}^{4N\times4N}$, is given by
\begin{equation}\label{eq:PXN_norm_def}
    \PX^{\perp}=\mathrm{diag}(\{\mathcal{P}_{\mathbf{x}_{i}}^{\perp}\mid i\in\{1,\ldots,N\}\}),
\end{equation}
with $\mathcal{P}_{\mathbf{x}_{i}}^{\perp}$ given by~\eqref{eq:Px_norm_def}.

\subsection{Riemannian Metrics}\label{app:metrics}
Because we employ the embedding defined in Appendix~\ref{app:submanifolds}, $\mathcal{M}$ inherits the Euclidean metric $g_{\mathbf{x}}\left(u,w\right)=\left\langle u,w\right\rangle _{\mathbf{x}}\triangleq u^{\transpose}w$ and norm ${\left\Vert u\right\Vert}_{\mathbf{x}}\triangleq\norm{u}$ for all $\mathbf{x}\in\mathcal{M}$ and $u,w\in\mathcal{T}_{\mathbf{x}}\mathcal{M}$. Moreover, per~\cite[Section 3.7]{boumal2023introduction}, $\MN$ admits the product metric $g_{\mathcal{X}}\left(\U,\W\right)=\sum_{i=1}^{N}g_{\mathbf{x}_{i}}\left(u_{i},w_{i}\right)=\U^{\transpose}\W$, and norm $\left\Vert \mathcal{U}\right\Vert _{\mathcal{X}}\triangleq\norm{\U}$ for all $\mathcal{X}\in\mathcal{M}^{N}$ and $\U,\W\in\TXMN$.

\subsection{Parallel Transport}\label{app:parallel}
The parallel transport operator maps tangent vectors between tangent spaces. On $\M$, $\mathcal{P}_{\mathbf{x}\rightarrow\mathbf{y}}:\TxM\rightarrow\mathcal{T}_{\mathbf{y}}\M$ denotes the parallel transport from $\TxM$ to $\mathcal{T}_{\mathbf{y}}\M$ for any $\mathbf{x},\mathbf{y}\in\M$. For $u_{\mathbf{x}}\in\mathcal{T}_{\mathbf{x}}\M$, it is given by 
\begin{equation}\label{eq:par_trans_M}
    \mathcal{P}_{\mathbf{x}\rightarrow\mathbf{y}}(u_{\mathbf{x}})=\mathbf{y}\boxplus(\mathbf{x}^{-1}\boxplus u_{\mathbf{x}}).
\end{equation}
Extending this definition to $\MN$ yields $\mathcal{P}_{\X\rightarrow\mathcal{Y}}:\TXMN\rightarrow\mathcal{T}_{\mathcal{Y}}\MN$ to be 
\begin{equation}\label{eq:par_trans_MN}
    \mathcal{P}_{\X\rightarrow\mathcal{Y}}(\mathcal{U}_{\mathcal{X}})=\vectx{(\mathbf{y}_{i}\boxplus(\mathbf{x}_{i}^{-1}\boxplus u_{i}))_{i=1}^{N}}
\end{equation}
for $\mathcal{U}_{\mathcal{X}}=\vect((u_{i})_{i=1}^{N})\in\mathcal{T}_{\mathcal{Y}}\MN$, $\X=\vect((\mathbf{x}_{i})_{i=1}^{N})\in\MN$, and $\Y=\vect((\mathbf{y}_{i})_{i=1}^{N})\in\MN$.

\subsection{Logarithm and Exponential}\label{app:log_exp}
Here, we derive the logarithm and exponential maps for $\M$ and $\MN$. The smooth manifold $\M$ with the identity, inverse, and composition operator form a Lie group\cite{sola2018micro} whose Lie algebra is the tangent space~\eqref{eq:TxM_def} at the identity element, denoted $\mathcal{T}_{\fatone}\M$. The geometry of screw motions encoded by elements of the PUDQ group is a consequence of Chasles' Theorem~\cite{chen1991screw}, which states that any rigid transformation can be modeled as a rotation and translation about a singular axis, termed the \textit{screw axis}. In~\cite{daniilidis1999hand}, the logarithm map at the identity for the unit dual quaternion (UDQ) group $\mathbb{D}\mathbb{H}$ was derived for rigid transformations in 3D in terms of four screw parameters: the rotation angle $\theta$, pitch $d$, direction vector $\boldsymbol{l}$, and moment $\boldsymbol{m}$. Given an orthonormal basis $\{\mathbf{i},\mathbf{j},\mathbf{k}\}$, we define a translation $\mathbf{t}\triangleq t_{x}\mathbf{i}+t_{y}\mathbf{j}+t_{z}\mathbf{k}$ and direction vector $\boldsymbol{l}\triangleq l_{x}\mathbf{i}+l_{y}\mathbf{j}+l_{z}\mathbf{k}$. Then, the pitch is $d$ given by $d=\mathbf{t}^{\transpose}\boldsymbol{l}$, and the moment $\boldsymbol{m}$ is given by
\begin{equation}\label{eq:udq_moment}
    \boldsymbol{m}=\frac{1}{2}\left(\mathbf{t}\times\boldsymbol{l}+\cot\left(\frac{\theta}{2}\right)\boldsymbol{l}\times(\mathbf{t}\times\boldsymbol{l})\right),
\end{equation}
where \quotes{$\times$} denotes the standard cross product in 3D. Following the methodology in~\cite{li2020improved}, we can treat the PUDQ group $\M$ as the degenerate, planar case of the UDQ group $\mathbb{D}\mathbb{H}$, in which case $\theta$ remains unchanged, $\mathbf{t}=t_{x}\mathbf{i}+t_{y}\mathbf{j}$, and $\boldsymbol{l}=\mathbf{k}$. Furthermore, for planar rigid motions, $\mathbf{t}$ and $\mathbf{k}$ are orthogonal vectors, so $d=\mathbf{t}^{\top}\boldsymbol{l}=\mathbf{t}^{\top}\mathbf{k}=0$. Moreover, applying these planar definitions to~\eqref{eq:udq_moment} and simplifying yields the planar moment $\boldsymbol{m}$ to be
\begin{align}
    \boldsymbol{m}&=\frac{1}{2}\left(\mathbf{t}\times\mathbf{k}+\cot\left(\frac{\theta}{2}\right)\mathbf{k}\times(\mathbf{t}\times\mathbf{k})\right) \\
    &=\frac{1}{2}\left(\left(t_{y}+\cot\left(\frac{\theta}{2}\right)t_{x}\right)\mathbf{i}+\left(\cot\left(\frac{\theta}{2}\right)t_{y}-t_{x}\right)\mathbf{j}\right).\label{eq:pudq_moment}
\end{align}
Finally, substituting~\eqref{eq:pudq_moment} and the preceding planar definitions into the UDQ logarithm map derived in~\cite{daniilidis1999hand} and simplifying yields the PUDQ logarithm map at the identity for $\mathbf{x}\in\M$ to be
\begin{equation}\label{eq:log_1_symbolic}
    \Logx{\fatone}{\mathbf{x}}=\frac{1}{2}\left(\theta+\varepsilon d\right)\left(\boldsymbol{l}+\varepsilon\boldsymbol{m}\right)=\frac{\theta}{2}\left(\mathbf{k}+\varepsilon\boldsymbol{m}\right).
\end{equation}
We can express~\eqref{eq:log_1_symbolic} in vector form according to the basis $\{\mathbf{k},\varepsilon\mathbf{i},\varepsilon\mathbf{j}\}$ as $\Logx{\fatone}{\mathbf{x}}=[\frac{1}{2},\frac{1}{2}\boldsymbol{m}^{\transpose}]^{\transpose}$. Then, substituting~\eqref{eq:pudq_moment}, letting $\phi\triangleq\nicefrac{\theta}{2}$, and applying the definition of $\cot(\cdot)$ yields the vector expression
\begin{equation}
    \Logx{\fatone}{\mathbf{x}}=\left[\phi,\frac{\phi}{2}\left(t_{y}+\frac{\cos\left(\phi\right)}{\sin\left(\phi\right)}t_{x}\right),\frac{\phi}{2}\left(\frac{\cos\left(\phi\right)}{\sin\left(\phi\right)}t_{y}-t_{x}\right)\right]^{\transpose},
\end{equation}
which simplifies to 
\begin{equation}\label{eq:log_unsimp}
    \Logx{\fatone}{\mathbf{x}}=\left[\phi,\frac{1}{2}\frac{\phi}{\sin\left(\phi\right)}\left(\left[\begin{array}{cc}
        \cos\left(\phi\right) & \sin\left(\phi\right)\\
        -\sin\left(\phi\right) & \cos\left(\phi\right)
        \end{array}\right]\left[\begin{array}{c}
        t_{x}\\
        t_{y}
        \end{array}\right]\right)^{\transpose}\right].
\end{equation}
Now, we write $\mathbf{x}=[\mathbf{x}_{r}^{\transpose},\mathbf{x}_{d}^{\transpose}]^{\transpose}=[x_{0},x_{1},x_{2},x_{3}]^{\transpose}$, and note that, from~\eqref{eq:xd_matvec}, we have
\begin{equation}\label{eq:xd_log}
    \mathbf{x}_{d}=\frac{1}{2}\left[\begin{array}{cc}
        \cos\left(\phi\right) & \sin\left(\phi\right)\\
        -\sin\left(\phi\right) & \cos\left(\phi\right)
        \end{array}\right]\left[\begin{array}{c}
        t_{x}\\
        t_{y}
        \end{array}\right].
\end{equation}
Finally, substituting~\eqref{eq:xd_log} into~\eqref{eq:log_unsimp} and simplifying with $x_{1}=\sin(\phi)$ yields
\begin{equation}
    \Logx{\fatone}{\mathbf{x}}=\frac{\phi}{\sin\left(\phi\right)}\left[x_{1},x_{2},x_{3}\right].
\end{equation}
Therefore, given $\mathbf{x}\in\M$, the logarithm map at the identity, denoted $\textnormal{Log}_{\fatone}:\M\rightarrow \mathcal{T}_{\fatone}\M$, is given by
\begin{equation}\label{eq:log_1_app}
    \Logx{\fatone}{\mathbf{x}}=\frac{1}{\gamma\left(\mathbf{x}\right)}{\left[x_{1},~x_{2},~x_{3}\right]}^{\transpose},
\end{equation}
where 
\begin{equation}\label{eq:gamma_app}
    \gamma\left(\mathbf{x}\right)\triangleq\mathrm{sinc}\left(\phi\left(\mathbf{x}\right)\right)=\frac{\sin\left(\phi\left(\mathbf{x}\right)\right)}{\phi\left(\mathbf{x}\right)},
\end{equation}
with
\begin{equation}\label{eq:phi_def_app}
    \phi\left(\mathbf{x}\right)\triangleq\mathrm{wrap}\left(\arctan\left(x_{1},x_{0}\right)\right),
\end{equation}
where $\arctan:\mathbb{S}^{1}\rightarrow(-\pi,\pi]$ is the four-quadrant arctangent and
\begin{equation}\label{eq:wrap_app}
    \mathrm{wrap}\left(\alpha\right)\triangleq\begin{cases}
        \alpha+\pi & \text{if }\alpha\leq-\nicefrac{\pi}{2}\\
        \alpha-\pi & \text{if }\alpha>\nicefrac{\pi}{2}\\
        \alpha & \text{otherwise}.
        \end{cases}
\end{equation}
Here, $\phi:\M\rightarrow \left(-\nicefrac{\pi}{2},\nicefrac{\pi}{2}\right]$ computes the half-angle of rotation about the $\mathbf{k}$-axis encoded by a point on~$\M$. The half-angles $\phi + n\pi$ for all $n\in\mathbb{Z}$ encode the same rotation, so it is valid to wrap $\phi$ to $(-\nicefrac{\pi}{2},\nicefrac{\pi}{2}]$ via~\eqref{eq:wrap_app}. Moreover, the exponential map at the identity, denoted $\mathrm{Exp}_{\fatone}:T_{\fatone}\M\rightarrow\M$, is the inverse of~\eqref{eq:log_1_app}. Given $\mathbf{x}_{t}=\left[x_{t,1},~x_{t,2},~x_{t,3}\right]^{\transpose}\in\mathcal{T}_{\fatone}\M$, it is given by
\begin{equation}\label{eq:exp_1_app}
    \Expx{\fatone}{\mathbf{x}_{t}}={\left[\cos\left(x_{t,1}\right), \gamma\left(\mathbf{x}_{t}\right)\mathbf{x}_{t}^{\transpose}\right]}^{\transpose},
\end{equation}
where $\gamma\left(\mathbf{x}_{t}\right)\triangleq\operatorname{sinc}\left(x_{t,1}\right)$ from~\eqref{eq:gamma_app}. For context,~\eqref{eq:log_1_app} and~\eqref{eq:exp_1_app} constitute the Lie-theoretic logarithm and exponential maps on $\M$, when treated as a Lie group. By equipping $\M$ with the Riemannian metric derived in Appendix~\ref{app:metrics}, we can treat $\M$ as a Riemannian manifold, in which case~\eqref{eq:log_1_app} and~\eqref{eq:exp_1_app} define the logarithm and exponential maps evaluated at the identity. Furthermore, we can apply the parallel transport operator on $\M$ from~\eqref{eq:par_trans_M} to extend~\eqref{eq:log_1_app} and~\eqref{eq:exp_1_app} to arbitrary points on $\M$ as in~\cite{busam2017camera}. This yields, for any~$\mathbf{x},\mathbf{y}\in\M$, the pointwise logarithm map
\begin{equation}\label{eq:log_x_app}
    \Logx{\mathbf{x}}{\mathbf{y}}=\mathbf{x}\boxplus[0,\LogOne(\mathbf{x}^{-1}\boxplus\mathbf{y})^{\transpose}]^{\transpose},
\end{equation}
and, for $\mathbf{x}\in\M,\mathbf{y}_{t}\in \mathcal{T}_{\mathbf{x}}\M$, the pointwise exponential map
\begin{equation}\label{eq:exp_x_app}
    \Expx{\mathbf{x}}{\mathbf{y}_{t}}=\mathbf{x}\boxplus\Expx{\fatone}{\left(\mathbf{x}^{-1}\boxplus\mathbf{y}_{t}\right)_{1:3}},
\end{equation}
where ${\left(\cdot\right)}_{1:3}$ selects the last three elements of a vector. For the product manifold $\MN$,~\eqref{eq:log_1_app}-\eqref{eq:exp_x_app} yield, for $\X,\Y\in\MN$,
\begin{equation}\label{eq:LogX_MN_app}
    \Logx{\X}{\Y}=\vectx{(\Log_{\mathbf{x}_{i}}(\mathbf{y}_{i}))_{i=1}^{N}},
\end{equation}
and, for $\Y_{t}=\vectx{(\mathbf{y}_{t,i})_{i=1}^{N}}\in\TXMN$,
\begin{equation}\label{eq:ExpX_MN_app}
    \Expx{\X}{\Y_{t}}=\vectx{(\Exp_{\mathbf{x}_{i}}(\mathbf{y}_{t,i}))_{i=1}^{N}},
\end{equation}
with $\Log_{\mathbf{x}_{i}}(\cdot)$ and $\Exp_{\mathbf{x}_{i}}(\cdot)$ given by~\eqref{eq:log_x_app} and~\eqref{eq:exp_x_app}.

\subsection{Geodesic Distance}\label{app:geodesic}
The geodesic distance metric extends the Riemannian metric to measure the lengths of minimal curves between points on manifolds. The geodesic distance on $\M$ is given by 
\begin{equation}
    d_{\M}(\mathbf{x},\mathbf{y})=\left\Vert\LogOne(\mathbf{x}^{-1}\boxplus\mathbf{y})\right\Vert_{2}
\end{equation}
for $\mathbf{x},\mathbf{y}\in\M$. For the product manifold $\MN$, it is given by 
\begin{equation}
    d_{\MN}(\X,\Y)=\sqrt{\sum_{i=1}^{N}\left\Vert\LogOne(\mathbf{x}_{i}^{-1}\boxplus\mathbf{y}_{i})\right\Vert_{2}^{2}}
\end{equation}
for $\X=\vect((\mathbf{x}_{i})_{i=1}^{N})\in\MN$, and $\Y=\vect((\mathbf{y}_{i})_{i=1}^{N})\in\MN$.

\subsection{Weingarten Map}
The Weingarten map describes the extrinsic curvature of a manifold. Here, we derive the Weingarten maps for $\M$ and $\MN$, which will be used in our derivation of the Riemannian Hessian in Appendix~\ref{app:rhess}. From~\cite{absil2013extrinsic}, the Weingarten map at $\mathbf{x}\in\M$, denoted $\mathfrak{A}_{\mathbf{x}}:\TxM\times\mathcal{T}_{\mathbf{x}}^{\perp}\M\rightarrow\TxM$, is given by, for $u\in\mathcal{T}_{\mathbf{x}}\mathcal{M}$, $w\in\mathcal{T}_{\mathbf{x}}^{\perp}\mathcal{M}$,\footnote{It is noted that $u\in\mathcal{T}_{\mathbf{x}}\mathcal{M}$ implies $\mathcal{P}_{\mathbf{x}}u=u$ and $w\in\mathcal{T}_{\mathbf{x}}^{\perp}\mathcal{M}$ implies $\mathcal{P}_{\mathbf{x}}^{\perp}w=w$.}
\begin{equation}\label{eq:wein_def}
    \mathfrak{A}_{\mathbf{x}}\left(u,w\right)=\mathcal{P}_{\mathbf{x}}\mathrm{D}_{u}\mathcal{P}_{\mathbf{x}}w,
\end{equation}
with $\mathcal{P}_{\mathbf{x}}$ given by~\eqref{eq:Px_def}. In~\eqref{eq:wein_def}, $D_{u}$ denotes the directional derivative along $u$ at $\mathbf{x}$, which is defined for any function $f$ on $\M$ into a vector space, and for any $u\in\TxM$ as
\begin{equation}\label{eq:DuF}
    \mathrm{D}_{u}f\left(\mathbf{x}\right)=\underset{t\rightarrow0}{\lim}\ \frac{d}{dt}f\left(\curve\left(t\right)\right),
\end{equation}
where $\curve$ is any curve on $\mathcal{M}$ with $\curve\left(0\right)=\mathbf{x}$ and $\curve^{\prime}\left(0\right)=u$. Applying~\eqref{eq:Px_def} to~\eqref{eq:DuF} and letting $\mathbf{x}=\curve\left(t\right)$ yields
\begin{equation}\label{eq:wein_DuP_prelim}
    \mathrm{D}_{u}\mathcal{P}_{\mathbf{x}}=\underset{t\rightarrow0}{\lim}\ \frac{d}{dt}\mathcal{P}_{\curve\left(t\right)}=-\Ptilde\left(\curve^{\prime}\left(0\right){\left(\curve\left(0\right)\right)}^{\transpose}+\curve\left(0\right){\left(\curve^{\prime}\left(0\right)\right)}^{\transpose}\right)\Ptilde,
\end{equation}
which simplifies to
\begin{equation}\label{eq:wein_DuP}
    \mathrm{D}_{u}\mathcal{P}_{\mathbf{x}}=-\left(\Ptilde u\mathbf{x}^{\transpose}\Ptilde+\Ptilde\mathbf{x}u^{\transpose}\Ptilde\right),
\end{equation}
with $\Ptilde$ given by~\eqref{eq:Ptilde_def}. Substituting~\eqref{eq:wein_DuP} into~\eqref{eq:wein_def} yields
\begin{equation}\label{eq:wein_incomplete_1}
    \mathfrak{A}_{\mathbf{x}}\left(u,w\right)=\mathcal{P}_{\mathbf{x}}\mathrm{D}_{u}\mathcal{P}_{\mathbf{x}}w=-\mathcal{P}_{\mathbf{x}}\left(\Ptilde u\mathbf{x}^{\transpose}\Ptilde+\Ptilde\mathbf{x}v^{\transpose}\Ptilde\right)w=-\mathcal{P}_{\mathbf{x}}\Ptilde u\mathbf{x}^{\transpose}\Ptilde w-\mathcal{P}_{\mathbf{x}}\Ptilde\mathbf{x}u^{\transpose}\Ptilde w.
\end{equation}
The following two lemmas allow us to further simplify~\eqref{eq:wein_incomplete_1}.
\begin{lemma}\label{lemma:P_tilde_w}
    For all $w\in T_{\mathbf{x}}^{\perp}\mathcal{M}$ and for all $\mathbf{x}\in\mathcal{M}$, it holds that $\Ptilde w=w$.
\end{lemma}
\emph{Proof:} Since $w\in T_{\mathbf{x}}^{\perp}\mathcal{M}$, it holds that $\mathcal{P}_{\mathbf{x}}^{\perp}w=w$. Therefore, since $\Ptilde$ is idempotent (i.e., $\Ptilde\Ptilde=\Ptilde$) we have
\begin{equation}
    \Ptilde w=\Ptilde\mathcal{P}_{\mathbf{x}}^{\perp}w=\Ptilde\left(\Ptilde\mathbf{x}\mathbf{x}^{\transpose}\Ptilde\right)w=\Ptilde\mathbf{x}\mathbf{x}^{\transpose}\Ptilde w=\mathcal{P}_{\mathbf{x}}^{\perp}w=w,
\end{equation}
completing the proof.\hfill $\blacksquare$
\begin{lemma}\label{lemma:P_tilde_P_x}
    For all $\mathbf{x}\in\mathcal{M}$, $\Ptilde\mathcal{P}_{\mathbf{x}}=\mathcal{P}_{\mathbf{x}}\Ptilde$.
\end{lemma}
\emph{Proof:}
Since $\Ptilde$ is idempotent, we have
\begin{equation}
    \Ptilde\mathcal{P}_{\mathbf{x}}	=\Ptilde\left(I-\Ptilde\mathbf{x}\mathbf{x}^{\transpose}\Ptilde\right)=\Ptilde-\Ptilde\Ptilde\mathbf{x}\mathbf{x}^{\transpose}\Ptilde=\Ptilde-\Ptilde\mathbf{x}\mathbf{x}^{\transpose}\Ptilde\Ptilde=\left(I-\Ptilde\mathbf{x}\mathbf{x}^{\transpose}\Ptilde\right)\Ptilde=\mathcal{P}_{\mathbf{x}}\Ptilde,
\end{equation} 
completing the proof.\hfill $\blacksquare$ \\

Applying Lemmas~\ref{lemma:P_tilde_w} and~\eqref{lemma:P_tilde_P_x} to~\eqref{eq:wein_incomplete_1} yields
\begin{equation}\label{eq:wein_preorth}
    \mathfrak{A}_{\mathbf{x}}\left(u,w\right)=-\Ptilde\mathcal{P}_{\mathbf{x}}u\mathbf{x}^{\transpose}w-\Ptilde\mathcal{P}_{\mathbf{x}}\mathbf{x}u^{\transpose}w.
\end{equation}
Finally, since $u\in T_{\mathbf{x}}\mathcal{M}$ and $w\in T_{\mathbf{x}}^{\perp}\mathcal{M}$, it follows that $u$ and $w$ are orthogonal and therefore $u^{\transpose}w=0$. Applying this to~\eqref{eq:wein_preorth} yields
\begin{equation}\label{eq:wein_pudq}
    \mathfrak{A}_{\mathbf{x}}\left(u,w\right)=-\mathcal{P}_{\mathbf{x}}\Ptilde u\mathbf{x}^{\transpose}w,
\end{equation}
which gives the Weingarten map for $\M$. We now extend~\eqref{eq:wein_pudq} to derive the Weingarten map for $\mathcal{M}^{N}$, denoted $\mathfrak{A}_{\mathcal{X}}:\mathcal{T}_{\mathcal{X}}\mathcal{M}^{N}\times\mathcal{T}_{\mathcal{X}}^{\perp}\mathcal{M}^{N}\rightarrow\mathcal{T}_{\mathcal{X}}\mathcal{M}^{N}$. First, given $\X\in\MN$, $\mathcal{U}\in\mathcal{T}_{\mathcal{X}}\mathcal{M}^{N}$, and $\mathcal{W}\in\mathcal{T}_{\mathcal{X}}^{\perp}\mathcal{M}^{N}$, we define $C\left(t\right)={\left[\curve_{1}^{\transpose}\left(t\right),\curve_{2}^{\transpose}\left(t\right),\ldots,\curve_{N}^{\transpose}\left(t\right)\right]}^{\transpose}$ such that $C\left(0\right)=\mathcal{X}$ and $C^{\prime}\left(0\right)=\mathcal{U}$. Rewriting~\eqref{eq:wein_DuP} in terms of $\MN$ yields the Weingarten map at $\X\in\MN$ to be
\begin{equation}\label{eq:wein_prod_def}
    \mathfrak{A}_{\X}\left(\U,\W\right)=\mathcal{P}_{\X}\mathrm{D}_{\U}\mathcal{P}_{\X}\W,
\end{equation}
for any $\mathcal{U}=\vect((u_{i})_{i=1}^{N})\in \mathcal{T}_{\mathcal{X}}\mathcal{M}^{N}$, $\mathcal{W}=\vect((w_{i})_{i=1}^{N})\in \mathcal{T}_{\mathcal{X}}^{\perp}\mathcal{M}^{N}$, with $\PX$ given by~\eqref{eq:PXN_def}. From the definition in~\eqref{eq:DuF}, we now derive $\mathrm{D}_{\mathcal{\mathcal{U}}}\mathcal{P}_{\mathcal{X}}$ to be
\begin{equation}\label{eq:DUP_prelim}
    \mathrm{D}_{\U}\mathcal{P}_{\X}=\underset{t\rightarrow0}{\lim}\ \frac{d}{dt}\mathcal{\mathcal{P}}_{C\left(t\right)} =\frac{d}{dt}\lim_{t\rightarrow0}\left(\mathrm{diag}\left(\left\{ \mathcal{P}_{c_{i}\left(t\right)}\mid i\in\{1,\ldots,N\}\right\} \right)\right)=\mathrm{diag}\Big(\Big\{ \frac{d}{dt}\lim_{t\rightarrow0}\mathcal{P}_{c_{i}\left(t\right)}\mid i\in\{1,\ldots,N\}\Big\} \Big).
\end{equation}
Now, using~\eqref{eq:wein_DuP_prelim}, we see that~\eqref{eq:DUP_prelim} simplifies to
\begin{equation}\label{eq:DUP}
    \mathrm{D}_{\U}\mathcal{P}_{\X}=\mathrm{diag}\left(\left\{ \mathrm{D}_{u_{i}}\mathcal{P}_{\mathbf{x}_{i}}\mid i\in\{1,\ldots,N\}\right\} \right).
\end{equation}
Substituting~\eqref{eq:DUP} into~\eqref{eq:wein_prod_def} yields
\begin{equation}
    \mathfrak{A}_{\mathcal{X}}\left(\mathcal{U},\mathcal{W}\right)=\mathcal{P}_{\mathcal{X}}\mathrm{diag}\left(\left\{ \mathrm{D}_{u_{i}}\mathcal{P}_{\mathbf{x}_{i}}\mid i\in\{1,\ldots,N\}\right\} \right)\mathcal{W},
\end{equation}
which simplifies to 
\begin{equation}\label{eq:wein_prod_prelim}
    \mathfrak{A}_{\mathcal{X}}\left(\mathcal{U},\mathcal{W}\right)=\vectx{\left(\mathcal{P}_{\mathbf{x}_{i}}\mathrm{D}_{u_{i}}\mathcal{P}_{\mathbf{x}_{i}}w_{i}\right)_{i=1}^{N}}.
\end{equation}
Finally, noting that~\eqref{eq:wein_def} gives $\mathcal{P}_{\mathbf{x}_{i}}\mathrm{D}_{u_{i}}\mathcal{P}_{\mathbf{x_{i}}}w_{i}=\mathfrak{A}_{\mathbf{x_{i}}}\left(u_{i},w_{i}\right)$, we observe that
\begin{equation}\label{eq:wein_prod}
    \mathfrak{A}_{\mathcal{X}}\left(\mathcal{U},\mathcal{W}\right)=\vectx{\left(\mathfrak{A}_{\mathbf{x_{i}}}\left(u_{i},w_{i}\right)\right)_{i=1}^{N}},
\end{equation}
which gives the Weingarten map for $\MN$.

\section{Maximum Likelihood Objective Derivation}\label{app:mle}
Here, we derive the MLE objective $\mathcal{F}$ for PGO on the PUDQ product manifold. First, let $\G=(\V,\E)$ be a (directed) pose graph with vertex set $\V$ and edge set $\E$ consisting of ordered pairs $(i,j)\in\V\times\V$. Let $\X=\vect((\mathbf{x}_{i})_{i\in\V})\in\MN$ denote $N$ poses to be estimated. The $M$ relative pose measurements are denoted $\Z=\vect((\zij)_{(i,j)\in\E})\in\M^{M}$, where each $\zij\in\M$ encodes a measured transformation from $\posei$ to $\posej$, taken in the frame of $\posei$. We utilize a Lie-theoretic measurement model for $\zij$ in which zero-mean Gaussian noise $\noise$ is mapped from $\mathcal{T}_{\fatone}\mathcal{M}$ to $\mathcal{M}$ via the exponential map, i.e.,
\begin{equation}\label{eq:meas_model_app}
    \mathbf{z}_{ij}=\mathbf{x}_{i}^{-1}\boxplus\mathbf{x}_{j}\boxplus\Expx{\fatone}{\noise},\ \eta_{ij}\in\mathbb{R}^{3},\ \eta_{ij}\sim\mathcal{N}\left(0,\Sigma_{ij}\right).
\end{equation}
Rearranging terms and noting that $\LogOne\left(\mathbf{x}^{-1}\right)=-\LogOne\left(\mathbf{x}\right)$ and $\left(\mathbf{x}\boxplus \mathbf{y}\right)^{-1}=\mathbf{y}^{-1}\boxplus \mathbf{x}^{-1}$, we see that~\eqref{eq:meas_model_app} gives the likelihood function $\mathcal{L}(\mathcal{X}\mid\mathcal{Z})=\mathbb{P}(Z=\Z\mid\X)$ (where $Z$ denotes the random variable corresponding to observation $\Z$), with
\begin{equation}\label{eq:likelihood_1}
    \mathcal{L}\left(\mathcal{X}\mid\mathcal{Z}\right)=\prod_{\left(i,j\right)\in\mathcal{E}}\frac{1}{\sqrt{\left(2\pi\right)^{3}\det(\Sigma_{ij})}}\exp\left(-\frac{1}{2}\mathrm{Log}_{\fatone}\left(\zij^{-1}\boxplus\mathbf{x}_{i}^{-1}\boxplus\mathbf{x}_{j}\right)^{\top}\Sigma_{ij}^{-1}\LogOne\left(\zij^{-1}\boxplus\mathbf{x}_{i}^{-1}\boxplus\mathbf{x}_{j}\right)\right),
\end{equation}
whose maximizer over $\mathcal{X}\in\mathcal{M}^{N}$ is the \textit{maximum likelihood estimate}, denoted $\mathcal{X}^\star$. Equivalently, $\mathcal{X}^\star$ is the minimizer of the negative likelihood $-(\mathcal{L}(\X\mid\Z))$. Now, taking the natural logarithm of $-\log\left(\mathcal{L}\left(\mathcal{X}\mid\mathcal{Z}\right)\right)$ and simplifying yields
\begin{equation}\label{eq:likelihood_2}
    -\log\left(\mathcal{L}\left(\mathcal{X}\mid\mathcal{Z}\right)\right)=\sum_{\left(i,j\right)\in\mathcal{E}}\log\left(\frac{1}{\sqrt{\left(2\pi\right)^{3}\det\left(\Sigma_{ij}\right)}}\right)+\sum_{\left(i,j\right)\in\mathcal{E}}\frac{1}{2}\left(\mathrm{Log}_{\fatone}\left(\zij^{-1}\boxplus\mathbf{x}_{i}^{-1}\boxplus\mathbf{x}_{j}\right)^{\top}\Sigma_{ij}^{-1}\LogOne\left(\zij^{-1}\boxplus\mathbf{x}_{i}^{-1}\boxplus\mathbf{x}_{j}\right)\right).
\end{equation}
We now observe from~\eqref{eq:likelihood_2} that 
\begin{equation}
    \argmin_{\X} \left(-\log\left(\mathcal{L}\left(\mathcal{X}\mid\mathcal{Z}\right)\right)\right) = \argmin_{\X}\FX,
\end{equation}
where the maximum likelihood objective, denoted $\mathcal{F}\left(\mathcal{X}\right)$, is given by
\begin{equation}\label{eq:mle_F_app}
    \mathcal{F}\left(\mathcal{X}\right)=\sum_{\left(i,j\right)\in\mathcal{E}}f_{ij}\left(\X\right),
\end{equation}
where 
\begin{equation}\label{eq:f_ij_app}
    f_{ij}\left(\X\right)=\frac{1}{2}\mathbf{e}_{ij}\left(\mathbf{x}_{i},\mathbf{x}_{j}\right)^{\top}\Omega_{ij}\mathbf{e}_{ij}\left(\mathbf{x}_{i},\mathbf{x}_{j}\right)=\frac{1}{2}\left\Vert \mathbf{e}_{ij}\left(\mathbf{x}_{i},\mathbf{x}_{j}\right)\right\Vert _{\Omega_{ij}}^{2}.
\end{equation}
In~\eqref{eq:f_ij_app}, $\Omega_{ij}=\Sigma_{ij}^{-1}$ is the information matrix for edge $(i,j)$, and $\mathbf{e}_{ij}:\M\times\M\rightarrow\LieAlgebra$ is the \textit{tangent} residual given by
\begin{equation}\label{eq:e_ij_app}
    \mathbf{e}_{ij}\left(\mathbf{x}_{i},\mathbf{x}_{j}\right)\triangleq\mathrm{Log}_{1}\left(\mathbf{r}_{ij}\left(\mathbf{x}_{i},\mathbf{x}_{j}\right)\right)=\mathrm{Log}_{\fatone}\left(\zij^{-1}\boxplus\mathbf{x}_{i}^{-1}\boxplus\mathbf{x}_{j}\right),
\end{equation}
where we have implicitly defined the \textit{geodesic} residual $\mathbf{r}_{ij}:\mathcal{M}\times\mathcal{M}\rightarrow\mathcal{M}$ as $\mathbf{r}_{ij}\left(\mathbf{x}_{i},\mathbf{x}_{j}\right)\triangleq\zij^{-1}\boxplus\mathbf{x}_{i}^{-1}\boxplus\mathbf{x}_{j}$.\footnote{Henceforth, we omit the dependency on $\left(\mathbf{x}_{i},\mathbf{x}_{j}\right)$ from our notation, i.e., $\mathbf{e}_{ij}\triangleq\mathbf{e}_{ij}\left(\mathbf{x}_{i},\mathbf{x}_{j}\right)$, $\mathbf{r}_{ij}\triangleq\mathbf{r}_{ij}\left(\mathbf{x}_{i},\mathbf{x}_{j}\right)$.}
\section{Transformations of Pose Parameterizations and Uncertainties}\label{app:covariance}
In this appendix, we derive transformations of poses and pose uncertainties between three parameterizations of planar rigid motion, namely, Euclidean space, denoted $\mathbb{R}^{3}$, and the planar unit dual quaternion group, denoted $\M$, and the planar special Euclidean group, denoted $\SE$.

\subsection{Pose Transformations}
Here, we derive transformations of poses between the three aforementioned parameterizations. We first define a planar pose represented in an orthonormal basis $(\boldsymbol{\mathrm{i}},\boldsymbol{\mathrm{j}},\boldsymbol{\mathrm{k}})$ and characterized by a translation $\boldsymbol{\mathrm{t}}=t_{x}\boldsymbol{\mathrm{i}}+t_{y}\boldsymbol{\mathrm{j}}$ and a rotation angle $\theta\in\left(-\pi,\pi\right]$ about the $\mathbf{k}$ axis. In Euclidean space, such a pose is given by the vector $\mathbf{p}=\left[\mathbf{t}^{\top},\theta\right]^{\top}\in\mathbb{R}^{3}$, with no additional structure applied. An alternative planar pose parameterization is that of the planar unit dual quaternion group (as detailed in Section~\ref{sec:pudq_construction}), which we denote $\M$. Letting $\phi\triangleq\nicefrac{\theta}{2}$, Euclidean poses are mapped to $\M$ via $\psi_{p}:\mathbb{R}^{3}\rightarrow\M$, defined as
\begin{equation}\label{eq:psi_p}
    \psi_{p}\left(\mathbf{p}\right)\triangleq\left[\cos\left(\phi\right),\sin\left(\phi\right),\frac{1}{2}(R_{\phi}\mathbf{t})^{\transpose}\right]^{\transpose},
\end{equation}
where, letting $\mathbf{x}={[x_{0},x_{1},x_{2},x_{3}]}^{\transpose}$, $R_{\phi}$ is given by
\begin{equation}\label{eq:R_phi}
    R_{\phi}\triangleq\left[\begin{array}{cc}
        \cos\left(\phi\right) & \sin\left(\phi\right)\\
        -\sin\left(\phi\right) & \cos\left(\phi\right)
        \end{array}\right]=\left[\begin{array}{cc}
            x_{0} & x_{1}\\
            -x_{1} & x_{0}
            \end{array}\right].
\end{equation}
The inverse map, $\psi_{p}^{-1}:\M\rightarrow\mathbb{R}^{3}$, is defined as
\begin{equation}\label{eq:psi_p_inv}
    \psi_{p}^{-1}\left(\mathbf{x}\right)=2\left[(R_{\phi}^{\transpose}\Ptilde\mathbf{x})^{\transpose},\phi(\mathbf{x})\right]^{\transpose},
\end{equation}
with $\phi(\mathbf{x})$ given by~\eqref{eq:phi_def_app}. Another common pose parameterization is the planar special Euclidean group, denoted $\SE$, which is defined as $\SE\triangleq\SO\rtimes\mathbb{R}^{2}$, where \quotes{$\rtimes$} denotes the semidirect product, and $\SO$ is the special orthogonal group, i.e., the set of all rotation matrices, which is given by
\begin{equation}
    \SO\triangleq\left\{R\in\mathbb{R}^{2\times 2}\mid R^{\transpose}R=RR^{\transpose}=I_{2},\operatorname{det}\left(R\right)=1\right\},
\end{equation}
where $I_{2}\in\mathbb{R}^{2\times 2}$ is an identity matrix. $\SE$ is traditionally coordinatized using the homogeneous transformation matrix (HTM) representation, which gives the definition
\begin{equation}\label{eq:SE2_def}
    \SE\triangleq\left\{ \left[\begin{array}{cc}
        \boldsymbol{R} & \boldsymbol{t}\\
        \boldsymbol{0} & 1
        \end{array}\right]\in\mathbb{R}^{3\times3}\mid\boldsymbol{R}\in \SO,\boldsymbol{t}\in\mathbb{R}^{2}\right\}.
\end{equation}
In this work, we equate $\SE$ with its HTM representation in~\eqref{eq:SE2_def}. Given a planar Euclidean pose, $\mathbf{p}\in\mathbb{R}^{3}$, the mapping from Euclidean space to $\SE$, which we denote $\psi_{s}:\mathbb{R}^{3}\rightarrow\SE$, is then given by 
\begin{equation}\label{eq:psi_s}
    \psi_{s}\left(\mathbf{p}\right)=\left[\begin{array}{ccc}
        \cos\left(\theta\right) & -\sin\left(\theta\right) & t_{x}\\
        \sin\left(\theta\right) & \cos\left(\theta\right) & t_{y}\\
        0 & 0 & 1
        \end{array}\right].
\end{equation}
Moreover, given a planar special Euclidean pose $T\in\SE$, the mapping from $\SE$ back to Euclidean space is given by the inverse mapping $\psi_{s}^{-1}:\SE\rightarrow\mathbb{R}^{3}$, which is given by
\begin{equation}\label{eq:psi_s_inv}
    \psi_{s}^{-1}\left(T\right)=\left[T_{13},T_{23},\arctan\left(T_{21},T_{11}\right)\right]^{\top},
\end{equation}
where $T_{ij}$ denotes the entry of the matrix $T$ at row $i$, column $j$, and $\arctan(\cdot)$ denotes the four-quadrant arctangent function. Furthermore, poses can be mapped between $\SE$ and $\M$ using compositions of~\eqref{eq:psi_s}-\eqref{eq:psi_s_inv},~\eqref{eq:psi_p}, and~\eqref{eq:psi_p_inv}, i.e., $\psi_{s}\circ\psi_{p}^{-1}:\M\rightarrow\SE$ and $\psi_{p}\circ\psi_{s}^{-1}:\SE\rightarrow\M$.

\subsection{Pose Covariance Transformations}
We now derive transformations between uncertainties of poses corresponding to random variables in~$\mathbb{R}^{3}$, $\M$, and $\SE$. These transformations presume that pose uncertainties in $\M$ and $\SE$ are modeled as Gaussian distributions in the Lie algebras of their respective groups. First, let $\mathbf{x}_{e}\triangleq[t_{x},t_{y},\theta]^{\top}\in\mathbb{R}^{3}$ be a planar Euclidean pose. Then, given $v_{p}\in\mathcal{T}_{\fatone}\mathcal{M}$, where $v_{p}\triangleq\mathrm{Log}_{\fatone}\left(\psi_{p}\left(\mathbf{x}_{e}\right)\right)$, and noting that $\theta\triangleq2\phi$, it holds from~\eqref{eq:log_1_app} that
\begin{equation}
    v_{p}=\left[\frac{\theta}{2},\frac{\left(x_{0}t_{x}+x_{1}t_{y}\right)}{2\mathrm{sinc}\left(\nicefrac{\theta}{2}\right)},\frac{\left(x_{0}t_{y}-x_{1}t_{x}\right)}{2\mathrm{sinc}\left(\nicefrac{\theta}{2}\right)}\right]^{\top},
\end{equation}
which simplifies to
\begin{equation}\label{eq:v_p_map}
   v_{p}=\frac{1}{2}B_{p}M_{p}\left(\theta\right)\mathbf{x}_{e},
\end{equation}
where
\begin{equation}
    B_{p}\triangleq\left[\begin{array}{ccc}
        0 & 0 & 1\\
        1 & 0 & 0\\
        0 & 1 & 0
        \end{array}\right],\text{~and~}M_{p}\triangleq\left[\begin{array}{ccc}
        \omega\left(\theta\right) & \nicefrac{\theta}{2} & 0\\
        -\nicefrac{\theta}{2} & \omega\left(\theta\right) & 0\\
        0 & 0 & 1
        \end{array}\right],\text{~with~}\omega(\theta)\triangleq\frac{\cos\left(\nicefrac{\theta}{2}\right)}{\mathrm{sinc}\left(\nicefrac{\theta}{2}\right)}.
\end{equation}
Here,~\eqref{eq:v_p_map} gives an invertible map from $\mathbb{R}^{3}$ to $\mathcal{T}_{\fatone}\mathcal{M}$. Now, let $\mathbf{x}_{e}\sim\mathcal{N}\left(0,\Sigma_{e}\right)$ be a random vector. Letting $\Sigma_{p}\triangleq\mathrm{Cov}\left[v_{p}\right]$ and 
applying~\eqref{eq:v_p_map} yields
\begin{equation}\label{eq:eucl_cov_map}
    \Sigma_{p}=\frac{1}{4}B_{p}M_{p}\left(\theta\right)\Sigma_{e}M_{p}^{\top}\left(\theta\right)B_{p}^{\top}.
\end{equation}
Additionally, letting $\Omega_{p}\triangleq\Sigma_{p}^{-1}$, $\Omega_{e}\triangleq\Sigma_{e}^{-1}$, and noting that $B_{p}^{-1}=B_{p}^{\transpose}$, we have 
\begin{equation}\label{eq:eucl_info_map}
    \Omega_{p}=4B_{p}M_{p}^{-\top}\left(\theta\right)\Omega_{e}M_{p}^{-1}B_{p}^{\top}\left(\theta\right).
\end{equation}
Equations~\eqref{eq:eucl_cov_map} and~\eqref{eq:eucl_info_map} give invertible maps, and thus transform the covariance and information matrices of Gaussian random variables between $\mathbb{R}^{3}$ and $\mathcal{T}_{\fatone}\mathcal{M}$. However, this requires \textit{a priori} knowledge of $\theta$, which is not always available. Moreover, given a vector in the Lie algebra of $\SE$, denoted $v_{s}\in se(2)$, with $v_{s}=\psi_{s}\left(\mathbf{x}_{e}\right)$ (where $\psi_{s}:\mathbb{R}^{3}\rightarrow \SE$ is derived in~\cite{long2013banana}), it holds that
\begin{equation}\label{eq:v_s_map}
\mathbf{x}_{e}=M_{s}\left(\theta\right)v_{s}^{\vee},
\end{equation}
where
\begin{equation}\label{eq:M_S}
    M_{s}\left(\theta\right)\triangleq\left[\begin{array}{ccc}
        \mathrm{sinc}\left(\theta\right) & \frac{\cos\theta-1}{\theta} & 0\\
        \frac{1-\cos\theta}{\theta} & \mathrm{sinc}\left(\theta\right) & 0\\
        0 & 0 & 1
        \end{array}\right],
\end{equation}
and the operator $\vee:se(2)\rightarrow\mathbb{R}^{3}$ maps from the Lie algebra to its Euclidean representation. 
Combining~\eqref{eq:v_p_map} and~\eqref{eq:v_s_map} yields the mapping from $se(2)$ to $\mathcal{T}_{\fatone}\mathcal{M}$ to be
\begin{equation}\label{eq:vs_to_vp_exp}
    v_{p}=\frac{1}{2}B_{p}M_{p}\left(\theta\right)M_{s}\left(\theta\right)v_{s}^{\vee},
\end{equation}
and since $M_{p}\left(\theta\right)M_{s}\left(\theta\right)=I_{3}$,~\eqref{eq:vs_to_vp_exp} reduces to
\begin{equation}\label{eq:vs_to_vp}
    v_{p}=\frac{1}{2}B_{p}v_{s}^{\vee},
\end{equation}
which gives an invertible vector map from $\LieAlgebra$ to $se(2)$ that is independent of $\theta$. Now, consider $v_{s}^{\vee}\sim\mathcal{N}\left(0,\Sigma_{s}\right)$. From~\eqref{eq:vs_to_vp}, we have
\begin{equation}\label{eq:se2_cov_map}
    \Sigma_{p}=\mathrm{Cov}\left[\frac{1}{2}B_{p}v_{s}^{\vee}\right]=\frac{1}{4}B_{p}\Sigma_{s}B_{p}^{\transpose}.
\end{equation}
Letting $\Omega_{s}\triangleq\Sigma_{s}^{-1}$, we also have
\begin{equation}\label{eq:se2_info_map}
\Omega_{p}	=\left(\frac{1}{4}B_{p}\Sigma_{s}B_{p}^{\transpose}\right)^{-1}=4B_{p}\Omega_{s}B_{p}^{\transpose}.
\end{equation}
The maps in~\eqref{eq:se2_cov_map} and~\eqref{eq:se2_info_map} are also invertible, and thus transform the covariance and information matrices of Gaussian random variables between $\mathcal{T}_{\fatone}\mathcal{M}$ and $se(2)$.
\section{Riemannian Gradient Derivation}\label{app:rgrad}
In this appendix, we derive the Riemannian gradient for the maximum likelihood objective $\F$ given by~\eqref{eq:mle_F_app}. Because $\MN$ is an Riemannian submanifold of a Euclidean space~\cite{boumal2023introduction}, the Riemannian gradient at $\X\in\MN$, denoted $\gradF\left(\mathcal{X}\right)$, is computed by projecting the Euclidean gradient at $\X$, denoted $\partial\bar{\mathcal{F}}\left(\mathcal{X}\right)$, onto $\TXMN$, i.e., 
\begin{equation}\label{eq:rgrad_F}
    \gradFX=\mathcal{P}_{\mathcal{X}}\partial\bar{\mathcal{F}}\left(\mathcal{X}\right), 
\end{equation}
with $\mathcal{P}_{\mathcal{X}}$ given by equation~\eqref{eq:Px_def}. Thus, the remainder of this appendix serves to derive the Euclidean gradient of $\F$.

\subsection{Euclidean Gradient}\label{app:egrad}
The Euclidean gradient of $\mathcal{F}$, denoted $\partial\bar{\mathcal{F}}$, is derived by omitting the manifold constraint from equation~\eqref{eq:mle_F_app} and computing the gradient  of $\F$ in $\mathbb{R}^{4N}$ with respect to $\mathcal{X}$. Differentiating~\eqref{eq:mle_F_app} in this manner and simplifying yields
\begin{equation}\label{eq:eucl_grad_def}
    \partial\bar{\mathcal{F}}\left(\mathcal{X}\right)=\frac{\partial\mathcal{F}\left(\mathcal{X}\right)}{\partial\mathcal{X}}=\frac{\partial}{\partial\mathcal{X}}\sum_{\left(i,j\right)\in\mathcal{E}}f_{ij}\left(\X\right)=\sum_{\left(i,j\right)\in\mathcal{E}}\frac{\partial f_{ij}\left(\X\right)}{\partial\mathcal{X}}.
\end{equation}
Since 
\begin{equation}\label{eq:fij_diff_1}
    \frac{\partial f_{ij}\left(\X\right)}{\partial\mathbf{x}_{l}}=\begin{cases}
    \frac{\partial}{\partial\mathbf{x}_{i}}f_{ij}\left(\X\right) & l=i,\\
    \frac{\partial}{\partial\mathbf{x}_{j}}f_{ij}\left(\X\right) & l=j,\\
    0 & \mathrm{otherwise},
    \end{cases}
\end{equation}
it suffices to compute the partial derivatives of $f_{ij}$ with respect to $\mathbf{x}_{i}$ and $\mathbf{x}_{j}$. Omitting the arguments $\left(\mathbf{x}_{i},\mathbf{x}_{j}\right)$ from $\mathbf{e}_{ij}$ and applying the chain rule to~\eqref{eq:fij_diff_1}, we have
\begin{equation}\label{eq:fij_diff_i}
    \frac{\partial f_{ij}\left(\X\right)}{\partial\mathbf{x}_{i}}=\frac{\partial}{\partial\mathbf{x}_{i}}\left(\mathbf{e}_{ij}^{\top}\Omega_{ij}\mathbf{e}_{ij}\right)={\left(\frac{\partial\mathbf{e}_{ij}}{\partial\mathbf{x}_{i}}\right)}^{\top}\Omega_{ij}\mathbf{e}_{ij}.
\end{equation}
Similarly,
\begin{equation}\label{eq:fij_diff_j}
    \frac{\partial f_{ij}\left(\X\right)}{\partial\mathbf{x}_{j}}={\left(\frac{\partial\mathbf{e}_{ij}}{\partial\mathbf{x}_{j}}\right)}^{\top}\Omega_{ij}\mathbf{e}_{ij}.
\end{equation}
Now, we denote $\mathcal{A}_{ij}\triangleq\frac{\partial}{\partial\mathbf{x}_{i}}\mathbf{e}_{ij}$ and $\mathcal{B}_{ij}\triangleq\frac{\partial}{\partial\mathbf{x}_{j}}\mathbf{e}_{ij}$ to be the Jacobians of $\eij$, which we derive in Appendix~\ref{app:egrad_jacobians}. Applying these definitions to~\eqref{eq:fij_diff_i} and~\eqref{eq:fij_diff_j} yields
\begin{equation}
    \frac{\partial f_{ij}\left(\X\right)}{\partial\mathbf{x}_{i}}=\mathcal{A}_{ij}^{\top}\Omega_{ij}\mathbf{e}_{ij}\text{~~and~~}\frac{\partial f_{ij}\left(\X\right)}{\partial\mathbf{x}_{j}}=\mathcal{B}_{ij}^{\top}\Omega_{ij}\mathbf{e}_{ij}.
\end{equation}
For each $f_{ij}$, with $\left(i,j\right)\in\mathcal{E}$, we have the block column vector 
\begin{equation}\label{eq:egrad_gij}
    \gijX\triangleq\frac{\partial \fijX}{\partial\mathcal{X}}={\left[g_{ij,1}^{\top},g_{ij,2}^{\top},\ldots,g_{ij,N}^{\top}\right]}^{\top},
\end{equation}
where
\begin{equation}\label{eq:g_ijl}
    g_{ij,l}=\begin{cases}
        \mathcal{A}_{ij}^{\top}\Omega_{ij}\mathbf{e}_{ij} & l=i,\\
        \mathcal{B}_{ij}^{\top}\Omega_{ij}\mathbf{e}_{ij} & l=j,\\
        \boldsymbol{0}_{4\times1} & \mathrm{otherwise},
        \end{cases}
\end{equation}
with each $g_{ij,l}\in\mathbb{R}^{4}$. Therefore, the Euclidean gradient of $\mathcal{F}$ is given by
\begin{equation}\label{eq:eucl_grad}
    \partial\bar{\mathcal{F}}\left(\mathcal{X}\right)=\sum_{\left(i,j\right)\in\mathcal{E}}\gijX,
\end{equation}
with $\mathbf{g}_{ij}$ given by~\eqref{eq:egrad_gij}.
\section{Riemannian Hessian Derivation}\label{app:rhess}
In this appendix, we derive the Riemannian Hessian for the maximum likelihood objective $\F$ given by~\eqref{eq:mle_F_app}. Additionally, in Appendix~\ref{app:rgn_hess}, we derive Riemannian Gauss-Newton Hessian approximation utilized in Section~\ref{sec:algorithm}. Towards deriving the Riemannian Hessian, we note that the embedding in Appendix~\ref{app:submanifolds} gives $\mathcal{M}^{N}$ as a Riemannian submanifold of the ambient Euclidean space $\mathbb{R}^{4N}$, and thus $\mathcal{M}^{N}$ takes on an extrinsic definition within the confines of this work. Leveraging this fact, we utilize the derivation proposed in~\cite{absil2013extrinsic}, in which, given $\mathcal{X}\in\mathcal{M}^{N}$ and $\mathcal{U}\in\mathcal{T}_{\mathcal{X}}\mathcal{M}^{N}$, the Riemannian Hessian is derived to be
\begin{equation}\label{eq:riem_hess_def}
    \HessF\left(\mathcal{X}\right)\left[\mathcal{U}\right]=\mathcal{P}_{\mathcal{X}}\partial^{2}\bar{\mathcal{F}}(\mathcal{X})\mathcal{U}+\mathfrak{A}_{\mathcal{X}}\left(\mathcal{U},\mathcal{P}_{\mathcal{X}}^{\perp}\partial\bar{\mathcal{F}}\left(\mathcal{X}\right)\right),
\end{equation}
where $\partial^{2}\bar{\mathcal{F}}$ is the Euclidean Hessian of $\mathcal{F}$ which we derive in Appendix~\ref{app:ehess}, $\partial\bar{\mathcal{F}}$ is the Euclidean gradient of $\mathcal{F}$ given by~\eqref{eq:eucl_grad}, $\mathcal{P}_{\mathcal{X}}$ is the orthogonal projector onto $\mathcal{T}_{\mathcal{X}}\mathcal{M}$ given by~\eqref{eq:PXN_def}, $\mathcal{P}_{\mathcal{X}}^{\perp}$ is the orthogonal projector onto $\mathcal{T}_{\mathcal{X}}^{\perp}\mathcal{M}$ given by~\eqref{eq:PXN_norm_def}, and $\mathfrak{A}_{\mathcal{X}}$ is the Weingarten map for $\MN$ given by~\eqref{eq:wein_prod}. To simplify~\eqref{eq:riem_hess_def}, we first separate the equation in terms of individual edges $\left(i,j\right)\in\mathcal{E}$. Substituting~\eqref{eq:mle_F_app} into~\eqref{eq:riem_hess_def} and simplifying yields
\begin{equation}\label{eq:rhess_expanded}
    \HessF\left(\mathcal{X}\right)\left[\mathcal{U}\right]=\sum_{\left(i,j\right)\in\mathcal{E}}\mathcal{P}_{\mathcal{X}}\HijbarX\mathcal{U}+\mathfrak{A}_{\mathcal{X}}\bigg(\mathcal{U},\sum_{\left(i,j\right)\in\mathcal{E}}\mathcal{P}_{\mathcal{X}}^{\perp}\gijX\bigg),
\end{equation}
where $\gij$ denotes the Euclidean gradient of $\fij$ given by~\eqref{eq:egrad_gij}, and $\Hijbar$ denotes the Euclidean Hessian of $\fij$, which we derive in Appendix~\ref{app:ehess}. To further simplify~\eqref{eq:rhess_expanded}, we prove in the following lemma that the Weingarten map on $\MN$ is linear in its second argument.
\begin{lemma}\label{lem:wein_prod_linear}
    Given $\X\in\MN$, $\mathcal{U}\in\mathcal{T}_{\mathcal{X}}\mathcal{M}^{N}$, and $\mathcal{W}\in\mathcal{T}_{\mathcal{X}}^{\perp}\mathcal{M}^{N}$, the Weingarten map $\mathfrak{A}_{\mathcal{X}}\left(\mathcal{U},\mathcal{W}\right)$ on $\mathcal{M}^{N}$ is linear in $\mathcal{W}$.
\end{lemma}
\emph{Proof:} First, we observe from~\eqref{eq:wein_pudq} that for any $\mathbf{x}\in\M$, $u\in\TxM$, $w,y\in\TxMNormal$, and $\alpha,\beta\in\mathbb{R}$, it holds that
\begin{equation}\label{eq:wein_pudq_linear}
    \mathfrak{A}_{\mathbf{x}}\left(u,\alpha w+\beta y\right)=-\mathcal{P}_{\mathbf{x}}\tilde{P}u\mathbf{x}^{\top}\left(\alpha w+\beta y\right)=-\alpha\mathcal{P}_{\mathbf{x}}\tilde{P}v\mathbf{x}^{\top}w-\beta\mathcal{P}_{\mathbf{x}}\tilde{P}v\mathbf{x}^{\top}y=\alpha\mathfrak{A}_{\mathbf{x}}\left(u,w\right)+\beta\mathfrak{A}_{\mathbf{x}}\left(u,y\right),
\end{equation}
which implies linearity of $\mathfrak{A}_{\mathbf{x}}\left(v,w\right)$ in $w$ on $\mathcal{M}$. It then suffices to show that linearity of $\mathfrak{A}_{\mathcal{X}}\left(\mathcal{U},\mathcal{W}\right)$ in $\mathcal{W}$ on $\mathcal{M}^{N}$ then follows from linearity of $\mathfrak{A}_{\mathbf{x}}\left(v,w\right)$ in $w$ on $\mathcal{M}$. Applying~\eqref{eq:wein_pudq_linear} to~\eqref{eq:wein_prod} yields, for any $\X\in\MN$, $\U\in\TXMN$, $\W,\Y\in\TXMNNormal$, and $\alpha,\beta\in\mathbb{R}$,
\begin{equation}\label{eq:wein_prod_lin_1}
    \mathfrak{A}_{\mathcal{X}}\left(\mathcal{U},\alpha\mathcal{W}+\beta\mathcal{Y}\right)=\vectx{\left(\mathfrak{A}_{\mathbf{x}_{i}}\left(u_{i},\alpha w_{i}+\beta y_{i}\right)_{i=1}^{N}\right)}=\vectx{\left(\alpha\mathfrak{A}_{\mathbf{x}_{i}}\left(u_{i}, w_{i}\right)+\beta\mathfrak{A}_{\mathbf{x}_{i}}\left(u_{i}, y_{i}\right)\right)_{i=1}^{N}}.
\end{equation}
Furthermore,~\eqref{eq:wein_prod_lin_1} implies that
\begin{equation}\label{eq:wein_prod_linear}
    \mathfrak{A}_{\mathcal{X}}\left(\mathcal{U},\alpha\mathcal{W}+\beta\mathcal{Y}\right)=\alpha\mathfrak{A}_{\mathcal{X}}\left(\mathcal{U},\mathcal{W}\right)+\beta\mathfrak{A}_{\mathcal{X}}\left(\mathcal{U},\mathcal{Y}\right),
\end{equation}
which gives linearity of $\mathfrak{A}_{\mathcal{X}}$ in $\W$, completing the proof.~\hfill $\blacksquare$ \\

Now, applying Lemma~\ref{lem:wein_prod_linear} to equation~\eqref{eq:rhess_expanded} yields
\begin{align}
    \HessFX\left[\mathcal{U}\right]&=\sum_{\left(i,j\right)\in\mathcal{E}}\mathcal{P}_{\mathcal{X}}\HijbarX\mathcal{U}+\sum_{\left(i,j\right)\in\mathcal{E}}\mathfrak{A}_{\mathcal{X}}\left(\mathcal{U},\mathcal{P}_{\mathcal{X}}^{\perp}\gijX\right) \\
    &=\sum_{\left(i,j\right)\in\mathcal{E}}\left(\mathcal{P}_{\mathcal{X}}\HijbarX\mathcal{U}+\mathfrak{A}_{\mathcal{X}}\left(\mathcal{U},\mathcal{P}_{\mathcal{X}}^{\perp}\gijX\right)\right)\label{eq:rhess_expanded_2}
\end{align}
Moreover, applying~\eqref{eq:riem_hess_def} to $f_{ij}$ from~\eqref{eq:f_ij_app} yields the Riemannian Hessian of $f_{ij}$ to be
\begin{equation}\label{eq:hess_fij_def}
    \Hessfij\left(\X\right)\left[\mathcal{U}\right]=\mathcal{P}_{\mathcal{X}}\HijbarX\mathcal{U}+\mathfrak{A}_{\mathcal{X}}\left(\mathcal{U},\mathcal{P}_{\mathcal{X}}^{\perp}\gijX\right),
\end{equation}
and substituting~\eqref{eq:hess_fij_def} into~\eqref{eq:rhess_expanded_2} yields
\begin{equation}\label{eq:rhess_f_sum}
    \HessFX\left[\mathcal{U}\right]=\sum_{\left(i,j\right)\in\mathcal{E}}\Hessfij\left(\X\right)\left[\mathcal{U}\right].
\end{equation}
Therefore, it suffices to derive $\Hessfij$ in order to derive $\HessF$. Towards accomplishing this, we first expand the Weingarten map term in~\eqref{eq:hess_fij_def} as
\begin{equation}\label{eq:pudq_wein_exp_1}
    \mathfrak{A}_{\mathcal{X}}\left(\mathcal{U},\mathcal{P}_{\mathcal{X}}^{\perp}\gijX\right)=\vectx{\left(\mathfrak{A}_{\mathbf{x}_{l}}\left(u_{l},\mathcal{P}_{\mathbf{x}_{l}}^{\perp}g_{ij,l}\right)\right)_{l\in\V}}=\vectx{\big(-\mathcal{P}_{\mathbf{x}_{l}}\tilde{P}u_{l}\mathbf{x}_{l}^{\top}\mathcal{P}_{\mathbf{x}_{l}}^{\perp}g_{ij,l}\big)_{l\in\V}},
\end{equation}
with $g_{ij,l}$ given in~\eqref{eq:g_ijl}. Because the $\mathbf{x}_{l}^{\top}\mathcal{P}_{\mathbf{x}_{l}}^{\perp}\mathbf{g}_{l}$ terms in~\eqref{eq:pudq_wein_exp_1} are scalars, it holds that
\begin{equation}\label{eq:pudq_wein_exp_2}
    \mathfrak{A}_{\mathcal{X}}\left(\mathcal{U},\mathcal{P}_{\mathcal{X}}^{\perp}\gijX\right)=\vectx{\left(-\mathcal{P}_{\mathbf{x}_{l}}\tilde{P}\mathbf{x}_{l}^{\top}\mathcal{P}_{\mathbf{x}_{l}}^{\perp}g_{ij,l}u_{l}\right)}=-\PX\PtildeV\diagx{\left\{\mathbf{x}_{l}^{\top}\mathcal{P}_{\mathbf{x}_{l}}^{\perp}g_{ij,l}I_{4}\right\}_{l\in\V}}\U,
\end{equation}
where $\tilde{P}_{\V}\triangleq\mathrm{diag}(\{\tilde{P}\}_{l\in\V}) \in\mathbb{R}^{4N\times4N}$. Furthermore,~\eqref{eq:pudq_wein_exp_2} is equivalent to the expression
\begin{equation}\label{eq:wein_prod_fij}
    \mathfrak{A}_{\mathcal{X}}\left(\mathcal{U},\mathcal{P}_{\mathcal{X}}^{\perp}\gijX\right)=-\mathcal{P}_{\mathcal{X}}\PtildeV\left(\left(J_{4}\otimes I_{N}\right)\circ\mathcal{X}^{\top}\mathcal{P}_{\X}^{\perp}\partial\fijbar\left(\X\right)\right)\U,
\end{equation}
where $J_{4}\in\mathbb{R}^{4\times4}$ is a matrix of ones, $\otimes$ is the Kronecker product, and $\circ$ is the Hadamard product. Substituting~\eqref{eq:wein_prod_fij} into~\eqref{eq:hess_fij_def} and simplifying yields the operator form of $\Hessfij$ to be
\begin{equation}\label{eq:hess_fij_operator}
    \Hessfij\left(\X\right)\left[\mathcal{U}\right]=\mathcal{P}_{\mathcal{X}}\left(\HijbarX-\PtildeV\left(\left(J_{4}\otimes I_{N}\right)\circ\mathcal{X}^{\top}\mathcal{P}_{\mathcal{X}}^{\perp}\gijX\right)\right)\mathcal{U},
\end{equation}
and since~\eqref{eq:hess_fij_operator} gives a matrix-vector multiplication, we deduce the matrix form of $\Hessfij$ to be
\begin{equation}\label{eq:rhess_fij}
    \Hessfij\left(\X\right)=\mathcal{P}_{\mathcal{X}}\left(\HijbarX-\PtildeV\left(\left(\boldsymbol{1}_{4}\otimes I_{N}\right)\circ\mathcal{X}^{\top}\mathcal{P}_{\mathcal{X}}^{\perp}\gijX\right)\right).
\end{equation}
Finally, substituting~\eqref{eq:rhess_fij} into~\eqref{eq:rhess_f_sum} and simplifying yields the Riemannian Hessian of $\mathcal{F}$ in matrix form to be
\begin{equation}\label{eq:hess_F}
    \HessFX=\sum_{\left(i,j\right)\in\mathcal{E}}\Hessfij\left(\X\right),
\end{equation}
with $\Hessfij\left(\X\right)$ given by~\eqref{eq:rhess_fij}.
\subsection{Euclidean Hessian}\label{app:ehess}
The Euclidean Hessian of $\F$, denoted $\partial^{2}\bar{\mathcal{F}}$, is computed by differentiating the Euclidean gradient of $\F$ from~\eqref{eq:eucl_grad} with respect to $\X$, i.e.,
\begin{equation}\label{eq:ehess_F}
    \partial^{2}\bar{\mathcal{F}}=\frac{\partial}{\partial\mathcal{X}}\left(\frac{\partial\mathcal{F}(\mathcal{X})}{\partial\mathcal{X}}\right)=\sum_{(i,j)\in\mathcal{E}}\frac{\partial}{\partial\mathcal{X}}\left(\frac{\partial f_{ij}(\mathcal{X})}{\partial\mathcal{X}}\right)=\sum_{(i,j)\in\mathcal{E}}\HijbarX,
\end{equation}
with $\HijbarX\triangleq\frac{\partial}{\partial\mathcal{X}}\gijX$ denoting the Euclidean Hessian of $\fij$, and with $\mathbf{g}_{ij}$ given by~\eqref{eq:egrad_gij}. From~\eqref{eq:egrad_gij}, we observe that
\begin{equation}\label{eq:partial_gij}
    \frac{\partial}{\partial\mathbf{x}_{m}}\left(\frac{\partial f_{ij}\left(\mathcal{X}\right)}{\partial\mathbf{x}_{l}}\right)=\begin{cases}
        \frac{\partial}{\partial\mathbf{x}_{i}}\left(\mathcal{A}_{ij}^{\top}\Omega_{ij}\mathbf{e}_{ij}\right) & m=l=i\\
        \frac{\partial}{\partial\mathbf{x}_{j}}\left(\mathcal{A}_{ij}^{\top}\Omega_{ij}\mathbf{e}_{ij}\right) & m=i,l=j\\
        \frac{\partial}{\partial\mathbf{x}_{i}}\left(\mathcal{B}_{ij}^{\top}\Omega_{ij}\mathbf{e}_{ij}\right) & m=j,l=i\\
        \frac{\partial}{\partial\mathbf{x}_{j}}\left(\mathcal{B}_{ij}^{\top}\Omega_{ij}\mathbf{e}_{ij}\right) & m=l=j\\
        \mathbf{0}_{4\times1} & \text{ otherwise. }
        \end{cases}
\end{equation}
We now denote the four nonzero blocks in~\eqref{eq:partial_gij} to be 
\begin{equation}\label{eq:partial_gij_blocks}
    \mathbf{h}_{ii}\triangleq\frac{\partial}{\partial x_{i}}\left(\mathcal{A}_{ij}^{\top}\Omega_{ij}\mathbf{e}_{ij}\right),~\mathbf{h}_{ij}\triangleq\frac{\partial}{\partial x_{j}}\left(\mathcal{A}_{ij}^{\top}\Omega_{ij}\mathbf{e}_{ij}\right),~\mathbf{h}_{ji}\triangleq\frac{\partial}{\partial x_{i}}\left(\mathcal{B}_{ij}^{\top}\Omega_{ij}\mathbf{e}_{ij}\right),~\mathbf{h}_{jj}\triangleq\frac{\partial}{\partial x_{j}}\left(\mathcal{B}_{ij}^{\top}\Omega_{ij}\mathbf{e}_{ij}\right). 
\end{equation}
Applying the product rule to compute the expressions in~\eqref{eq:partial_gij_blocks} yields
\begin{align}
    \mathbf{h}_{ii}&=\frac{\partial}{\partial\mathbf{x}_{i}}\left(\mathcal{A}_{ij}^{\top}\right)\Omega_{ij}\mathbf{e}_{ij}+\mathcal{A}_{ij}^{\top}\Omega_{ij}\mathcal{A}_{ij}, \\
    \mathbf{h}_{ij}&=\frac{\partial}{\partial\mathbf{x}_{j}}\left(\mathcal{A}_{ij}^{\top}\right)\Omega_{ij}\mathbf{e}_{ij}+\mathcal{A}_{ij}^{\top}\Omega_{ij}\mathcal{B}_{ij}, \\
    \mathbf{h}_{ji}&=\frac{\partial}{\partial\mathbf{x}_{i}}\left(\mathcal{B}_{ij}^{\top}\right)\Omega_{ij}\mathbf{e}_{ij}+\mathcal{B}_{ij}^{\top}\Omega_{ij}\mathcal{A}_{ij}, \\
    \mathbf{h}_{jj}&=\frac{\partial}{\partial\mathbf{x}_{j}}\left(\mathcal{B}_{ij}^{\top}\right)\Omega_{ij}\mathbf{e}_{ij}+\mathcal{B}_{ij}^{\top}\Omega_{ij}\mathcal{B}_{ij}.
\end{align}
Letting $\mathcal{C}_{ii}\triangleq\frac{\partial}{\partial\mathbf{x}_{i}}\left(\mathcal{A}_{ij}\right)^{\top}\Omega_{ij}\mathbf{e}_{ij}$, $\mathcal{C}_{ij}\triangleq\frac{\partial}{\partial\mathbf{x}_{j}}\left(\mathcal{A}_{ij}\right)^{\top}\Omega_{ij}\mathbf{e}_{ij}$, $\mathcal{C}_{ji}\triangleq\frac{\partial}{\partial\mathbf{x}_{i}}\left(\mathcal{B}_{ij}\right)^{\top}\Omega_{ij}\mathbf{e}_{ij}$, and $\mathcal{C}_{jj}\triangleq\frac{\partial}{\partial\mathbf{x}_{j}}\left(\mathcal{B}_{ij}\right)^{\top}\Omega_{ij}\mathbf{e}_{ij}$ gives\footnote{Expressions for $\mathcal{C}_{ii},\mathcal{C}_{ij},\mathcal{C}_{ji}$, and $\mathcal{C}_{jj}$ are derived in Appendix~\ref{app:ehess_combined_bounds}.}
\begin{equation}\label{eq:h_ii-h_jj}
    \mathbf{h}_{ii}=\mathcal{C}_{ii}+\mathcal{A}_{ij}^{\top}\Omega_{ij}\mathcal{A}_{ij},~\mathbf{h}_{ij}=\mathcal{C}_{ij}+\mathcal{A}_{ij}^{\top}\Omega_{ij}\mathcal{B}_{ij},~\mathbf{h}_{ji}=\mathcal{C}_{ji}+\mathcal{B}_{ij}^{\top}\Omega_{ij}\mathcal{A}_{ij},\text{~and~}\mathbf{h}_{jj}=\mathcal{C}_{jj}+\mathcal{B}_{ij}^{\top}\Omega_{ij}\mathcal{B}_{ij}.
\end{equation}
From~\eqref{eq:partial_gij}, it holds that the matrix $\Hijbar$ has only 4 nonzero blocks, which we now define in terms of~\eqref{eq:h_ii-h_jj}. Denoting block indices {\small$\mathtt{i}\triangleq 4i+1:4i+4$} and {\small$\mathtt{j}\triangleq 4j+1:4j+4$}, they are given by 
\begin{equation}\label{eq:ehess_blocks}
    \Hbar_{ij\left[\mathtt{i},\mathtt{i}\right]}=\mathbf{h}_{ii},~\Hbar_{ij\left[\mathtt{i},\mathtt{j}\right]}=\mathbf{h}_{ij},~\Hbar_{ij\left[\mathtt{j},\mathtt{i}\right]}=\mathbf{h}_{ji},\text{~and~}\Hbar_{ij\left[\mathtt{j},\mathtt{j}\right]}=\mathbf{h}_{jj},
\end{equation}
with $\mathbf{h}_{ii}$, $\mathbf{h}_{ij}$, $\mathbf{h}_{ji}$ and $\mathbf{h}_{jj}$ given by~\eqref{eq:h_ii-h_jj}. It then follows that the Euclidean Hessian of $\F$ is computed by applying~\eqref{eq:ehess_blocks} to~\eqref{eq:ehess_F}.
\begin{remark}
We note that, as expected, $\mathbf{h}_{ii}=\mathbf{h}_{ii}^{\top}$, $\mathbf{h}_{ji}=\mathbf{h}_{ij}^{\top}$, and $\mathbf{h}_{jj}=\mathbf{h}_{jj}^{\top}$, so $\bar{\mathcal{H}}_{ij}$ is symmetric, and therefore the Euclidean Hessian $\partial^{2}\bar{\mathcal{F}}$ in~\eqref{eq:ehess_F} is symmetric.
\end{remark}
\subsection{Riemannian Gauss-Newton Hessian}\label{app:rgn_hess}
In~\eqref{eq:m_k}, $\mathcal{H}_{k}:\mathcal{T}_{\X_{k}}\mathcal{M}^{N}\rightarrow\mathcal{T}_{\X_{k}}\mathcal{M}^{N}$ is the Riemannian Gauss-Newton (RGN) approximation of the Riemannian Hessian at $\X_{k}$, which we now derive. First, because $\Omega_{ij}$ is an information matrix, we have $\Omega_{ij}\succ0$, and can write $\Omega_{ij}=\Omega_{ij}^{\nicefrac{1}{2}}\Omega_{ij}^{\nicefrac{1}{2}}$, with $\Omega_{ij}^{\nicefrac{1}{2}}=(\Omega_{ij}^{\nicefrac{1}{2}})^{\top}$. Applying this to the definition of $f_{ij}(\X)$ given in~\eqref{eq:f_ij_app}, we can then write 
\begin{equation}
    f_{ij}\left(\mathcal{X}\right)	=\left\Vert F_{ij}\left(\mathcal{X}\right)\right\Vert _{2}^{2}=\left\langle F_{ij}\left(\mathcal{X}\right),F_{ij}\left(\mathcal{X}\right)\right\rangle, 
\end{equation}
where $F_{ij}\left(\mathcal{X}\right)\triangleq\Omega_{ij}^{\nicefrac{1}{2}}\mathbf{e}_{ij}$. From~\cite[Section~8.4]{absil2008optimization}, the RGN approximation of $\Hessfij$, denoted $\tilde{\mathcal{H}}_{ij}$, is given by
\begin{equation}
    \HessfijX\left[\xi,\eta\right]\approx\tilde{\mathcal{H}}_{ij}\left(\X\right)\left[\xi,\eta\right]\triangleq\left\langle \mathrm{D}F_{ij}\left(\mathcal{X}\right)\left[\xi\right],\mathrm{D}F_{ij}\left(\mathcal{X}\right)\left[\eta\right]\right\rangle,
\end{equation}
for $\xi,\eta\in\mathcal{T}_{\mathcal{X}}\mathcal{M}^{N}$. Applying the inner product definition from Appendix~\ref{app:riemannian} yields 
\begin{equation}
    \tilde{\mathcal{H}}_{ij}\left(\X\right)\left[\xi,\eta\right]=\xi^{\top}\left(\mathrm{D}F_{ij}\left(\mathcal{X}\right)\right)^{\top}\mathrm{D}F_{ij}\left(\mathcal{X}\right)\eta,
\end{equation}
from which we deduce that 
\begin{equation}\label{eq:rgnhess_df}
    \tilde{\mathcal{H}}_{ij}\left(\X\right)=\left(\mathrm{D}F_{ij}\left(\mathcal{X}\right)\right)^{\top}\mathrm{D}F_{ij}\left(\mathcal{X}\right).
\end{equation}
Furthermore, it holds from~\cite[Section~8.4]{absil2008optimization} that
\begin{equation}\label{eq:gradfij_rgn_adj}
    \textnormal{grad}\,f_{ij}\left(\mathcal{X}\right)=\left(\mathrm{D}F_{ij}\left(\mathcal{X}\right)\right)^{*}\left[F_{ij}\left(\mathcal{X}\right)\right],
\end{equation}
where $\left(\cdot\right)^{*}$ is the \textit{adjoint} operator, which we now define. Given two Euclidean spaces, denoted $\mathcal{O}$ and $\mathcal{Q}$, and an operator $T:\mathcal{O}\rightarrow\mathcal{Q}$, the adjoint of $T$ is the operator $T^{*}:\mathcal{Q}\rightarrow\mathcal{O}$ satisfying $\left\langle T\left[\mathcal{U}\right],\mathcal{W}\right\rangle =\left\langle \mathcal{U},T^{*}\left[\mathcal{W}\right]\right\rangle$  for all $\mathcal{U}\in\mathcal{O}$ and all $\mathcal{W}\in\mathcal{Q}$~\cite[Appendix~A]{absil2008optimization}. Applying the inner product definition yields $\left\langle T\left[\mathcal{U}\right],\mathcal{W}\right\rangle =\mathcal{U}^{\top}T^{\top}\mathcal{W}$, from which it follows that $T^{*}=T^{\top}$. Applying this to~\eqref{eq:gradfij_rgn_adj} yields
\begin{equation}\label{eq:gradfij_rgn}
    \textnormal{grad}\,f_{ij}\left(\mathcal{X}\right)=\left(\mathrm{D}F_{ij}\left(\mathcal{X}\right)\right)^{\top}F_{ij}\left(\mathcal{X}\right),
\end{equation}
and equating~\eqref{eq:gradfij_rgn} with~\eqref{eq:rgrad_F} yields 
\begin{equation}
    \left(\mathrm{D}F_{ij}\left(\X\right)\right)^{\top}=\mathcal{P}_{\mathcal{X}}\tilde{\mathbf{g}}_{ij}\left(\X\right), 
\end{equation}
with $\mathcal{P}_{\mathcal{X}}$ from~\eqref{eq:PXN_def} and with $\tilde{\mathbf{g}}_{ij}\left(\X\right)\triangleq[\tilde{g}_{ij,1}^{\transpose},\tilde{g}_{ij,2}^{\transpose},\ldots,\tilde{g}_{ij,N}^{\transpose}]^{\transpose}$, where
\begin{equation}\label{eq:gtilde_ijk}
    \tilde{g}_{ij,k}=\begin{cases}
        \Aij^{\transpose}\Omega_{ij}^{\nicefrac{1}{2}} & i=k,\\
        \mathcal{B}_{ij}^{\transpose}\Omega_{ij}^{\nicefrac{1}{2}} & j=k,\\
        \boldsymbol{0}_{4\times3} & \text{otherwise}.
        \end{cases}
\end{equation}
Substituting this into~\eqref{eq:rgnhess_df} and noting that $\PX=\PX^{\transpose}$ yields
\begin{equation}\label{eq:H_ij}
    \tilde{\mathcal{H}}_{ij}\left(\X\right)=\mathcal{P}_{\mathcal{X}}\mathcal{R}_{ij}\left(\X\right)\mathcal{P}_{\mathcal{X}},
\end{equation}
where $\mathcal{R}_{ij}\triangleq\tilde{\mathbf{g}}_{ij}\tilde{\mathbf{g}}_{ij}^{\top}$ (with argument $\left(\X\right)$ omitted from $\tilde{\mathbf{g}}_{ij}\left(\X\right)$). The matrix $\mathcal{R}_{ij}\in\mathbb{R}^{4N\times4N}$ in~\eqref{eq:H_ij} has only four nonzero blocks, which we now define. Denoting block indices {\small$\mathtt{i}\triangleq 4i+1:4i+4$} and {\small$\mathtt{j}\triangleq 4j+1:4j+4$}, they are given by 
\begin{equation}\label{eq:rgnhess_blocks}
    \HRGN_{ij\left[\mathtt{i},\mathtt{i}\right]}=\Aij^{\transpose}\Omega_{ij}\Aij,~\HRGN_{ij\left[\mathtt{i},\mathtt{j}\right]}=\Aij^{\transpose}\Omega_{ij}\mathcal{B}_{ij},~\HRGN_{ij\left[\mathtt{j},\mathtt{i}\right]}=\mathcal{B}_{ij}^{\transpose}\Omega_{ij}\Aij,\text{~and~}\HRGN_{ij\left[\mathtt{j},\mathtt{j}\right]}=\mathcal{B}_{ij}^{\transpose}\Omega_{ij}\mathcal{B}_{ij}.
\end{equation}
Moreover, applying~\eqref{eq:H_ij} to the definition of $\F$ given in~\eqref{eq:mle_F_app} yields the RGN Hessian approximation for $\HessF$ at $\Xk\in\MN$, denoted $\mathcal{H}_{k}$, to be
\begin{equation}\label{eq:RGN_H_k}
    \mathcal{H}_{k}=\sum_{\left(i,j\right)\in\mathcal{E}}\tilde{\mathcal{H}}_{ij}\left(\Xk\right)=\sum_{\left(i,j\right)\in\mathcal{E}}\PXk\mathcal{R}_{ij}\left(\Xk\right)\PXk.
\end{equation}
\begin{remark}
    As evidenced by comparing the Riemannian Hessian blocks in~\eqref{eq:h_ii-h_jj} to the RGN Hessian blocks in~\eqref{eq:rgnhess_blocks}, $\tilde{\mathcal{H}}_{ij}$ closely approximates $\HessF$ when the $\mathcal{C}_{ii}$-$\mathcal{C}_{jj}$ terms are negligible.
\end{remark}

\section{Lipschitz Continuity of the Riemannian Gradient}\label{app:lipschitz}
In this appendix, we derive a Lipschitz constant for the Riemannian gradient of the maximum likelihood objective given by~\eqref{eq:rgrad_F}, and our derivation serves as a proof for Theorem~\ref{thm:rgrad_lipschitz}. From~\cite{ferreira2019gradient} (see also~\cite{boumal2023introduction, han2023riemannian}), if $\mathcal{F}:\K\rightarrow\mathbb{R}$ is twice continuously differentiable on $\K$, then its Riemannian gradient is Lipchitz continuous on $\K$ with constant $L_{g}>0$ if and only if $\operatorname{Hess}\mathcal{F}\left(\mathcal{X}\right)$ has operator norm bounded by $L_{g}$ for all $\mathcal{X}\in\K$, that is, if for all $\mathcal{X}\in\K$, we have
\begin{equation}\label{eq:rhess_op_norm_lipschitz}
    \|\HessFX\|_{\mathcal{X}}=\sup\left\{ \|\HessFX[\mathcal{U}]\|_{\mathcal{X}}\mid\mathcal{U}\in\mathcal{T}_{\mathcal{X}}\mathcal{M},\|\mathcal{U}\|_{\mathcal{X}}=1\right\} \leq L_{g},
\end{equation}
where $\left\Vert \cdot\right\Vert _{\mathcal{X}}$ is the norm induced by the Riemannian metric at $\mathcal{X}$ on $\mathcal{M}$. Here, $\K\subset\MN$ is any compact subset of $\MN$, and the results we derive in this appendix hold for all $\X\in\K$. Using the inherited Riemannian metric and induced norm included in Appendix~\ref{app:metrics}, we first rewrite the operator norm from equation~\eqref{eq:rhess_op_norm_def} as
\begin{equation}\label{eq:rhess_op_norm_prod}
    \left\Vert \HessFX\right\Vert _{2}	=\sup\left\{ \left\Vert \HessFX\left[\mathcal{U}\right]\right\Vert _{2}\mid\mathcal{U}\in\mathcal{T}_{\mathcal{X}}\mathcal{M},\left\Vert \mathcal{U}\right\Vert _{2}=1\right\}.
\end{equation}
Next, we rewrite~\eqref{eq:rhess_op_norm_prod} in terms of $\Hessfij$ according to equation~\eqref{eq:rhess_f_sum}, which gives
\begin{equation}\label{eq:rhess_op_1}
    \left\Vert \HessFX\right\Vert _{2}=\sup\bigg\{ \Big\Vert \sum_{\left(i,j\right)\in\mathcal{E}}\operatorname{Hess}f_{ij}(\mathcal{X})\left[\mathcal{U}\right]\Big\Vert _{2}\mid\mathcal{U}\in\mathcal{T}_{\mathcal{X}}\mathcal{M},\|\mathcal{U}\|_{2}=1\bigg\}.
\end{equation}
Applying the triangle inequality to~\eqref{eq:rhess_op_1} yields
\begin{equation}
    \left\Vert \HessFX\right\Vert _{2}\leq\sup\bigg\{ \sum_{\left(i,j\right)\in\mathcal{E}}\big\Vert \operatorname{Hess}f_{ij}(\mathcal{X})\left[\mathcal{U}\right]\big\Vert _{2}\mid\mathcal{U}\in\mathcal{T}_{\mathcal{X}}\mathcal{M},\|\mathcal{U}\|_{2}=1\bigg\},
\end{equation}
and since $\sup\{x+y\}\leq\sup\{x\}+\sup\{y\}$, we observe that
\begin{equation}\label{eq:rhess_fij_opnorm_bound}
    \left\Vert \HessFX\right\Vert _{2}\leq\sum_{\left(i,j\right)\in\mathcal{E}}\sup\left\{ \left\Vert \operatorname{Hess}f_{ij}(\mathcal{X})\left[\mathcal{U}\right]\right\Vert _{2}\mid\mathcal{U}\in\mathcal{T}_{\mathcal{X}}\mathcal{M},\|\mathcal{U}\|_{2}=1\right\}.
\end{equation} 
We will now bound~\eqref{eq:rhess_fij_opnorm_bound} by bounding the $\operatorname{Hess}f_{ij}$ operator norms individually. First, it follows from~\eqref{eq:rhess_fij} and~\eqref{eq:ehess_blocks} that $\Hessfij$ has four nonzero blocks. Letting $\mathbf{H}_{ij}\triangleq\Hessfij(\X)$, and denoting block indices {\small$\mathtt{i}\triangleq 4i+1:4i+4$} and {\small$\mathtt{j}\triangleq 4j+1:4j+4$}, they are given by 
\begin{align}
    \mathbf{H}_{ij\left[\mathtt{i},\mathtt{i}\right]}&=\mathcal{P}_{\mathbf{x}_{i}}\left(\mathbf{h}_{ii}-\tilde{P}\mathbf{x}_{i}^{\top}\mathcal{P}_{\mathbf{x}_{i}}^{\perp}g_{ij,i}\right), \\
    \mathbf{H}_{ij\left[\mathtt{i},\mathtt{i}\right]}&=\mathcal{P}_{\mathbf{x}_{i}}\mathbf{h}_{ij},\\
    \mathbf{H}_{ij\left[\mathtt{i},\mathtt{i}\right]}&=\mathcal{P}_{\mathbf{x}_{j}}\mathbf{h}_{ji},\\
    \mathbf{H}_{ij\left[\mathtt{i},\mathtt{i}\right]}&=\mathcal{P}_{\mathbf{x}_{j}}\left(\mathbf{h}_{jj}-\tilde{P}\mathbf{x}_{j}^{\top}\mathcal{P}_{\mathbf{x}_{j}}^{\perp}g_{ij,j}\right),
\end{align}
with $g_{ij,i}$, $g_{ij,j}$ from~\eqref{eq:g_ijl} and $\mathbf{h}_{ii}$-$\mathbf{h}_{jj}$ from~\eqref{eq:h_ii-h_jj}. Then, given $\X\in\MN$ and $\U=\vect((u_{l})_{l\in\V})\in\TXMN$, we have
\begin{equation}\label{eq:hess_fij_hijl}
    \Hessfij\left(\X\right)\left[\U\right]=\vectx{\left(h_{ij,l}\right)_{i\in\V}},
\end{equation} 
with
\begin{equation}\label{eq:h_ijl}
    h_{ij,l}\triangleq\begin{cases}
        \mathcal{P}_{\mathbf{x}_{i}}\left(\mathbf{h}_{ii}u_{i}+\mathbf{h}_{ij}u_{j}-\tilde{P}\mathbf{x}_{i}^{\top}\mathcal{P}_{\mathbf{x}_{i}}^{\perp}g_{ij,i}u_{i}\right) & l=i,\\
        \mathcal{P}_{\mathbf{x}_{j}}\left(\mathbf{h}_{ji}u_{i}+\mathbf{h}_{jj}u_{j}-\tilde{P}\mathbf{x}_{j}^{\top}\mathcal{P}_{\mathbf{x}_{j}}^{\perp}g_{ij,j}u_{j}\right) & l=j,\\
        \mathbf{0}_{4\times1} & \text{otherwise}.
        \end{cases}
\end{equation}
Using~\eqref{eq:hess_fij_hijl} and~\eqref{eq:h_ijl}, we observe that
\begin{equation}\label{eq:rhess_fij_expansion}
    \left\Vert \HessfijX\left[\mathcal{U}\right]\right\Vert_{2}=\sqrt{\left\Vert H_{i}\left[\mathcal{U}\right]\right\Vert _{2}^{2}+\left\Vert H_{j}\left[\mathcal{U}\right]\right\Vert _{2}^{2}},
\end{equation}
where
\begin{align}
    H_{i}\left[\mathcal{U}\right]&=\mathcal{P}_{\mathbf{x}_{i}}\left(\mathbf{h}_{ii}u_{i}+\mathbf{h}_{ij}u_{j}-\mathbf{d}_{ii}u_{i}\right),\label{eq:H_i}\\
    H_{j}\left[\mathcal{U}\right]&=\mathcal{P}_{\mathbf{x}_{j}}\left(\mathbf{h}_{ji}u_{i}+\mathbf{h}_{jj}u_{j}-\mathbf{d}_{jj}u_{j}\right),\label{eq:H_j}
\end{align}
with
\begin{align}
    \mathbf{d}_{ii}&\triangleq\tilde{P}\mathbf{x}_{i}^{\top}\mathcal{P}_{\mathbf{x}_{i}}^{\perp}g_{ij,i},\label{eq:d_ii}\\
    \mathbf{d}_{jj}&\triangleq\tilde{P}\mathbf{x}_{j}^{\top}\mathcal{P}_{\mathbf{x}_{j}}^{\perp}g_{ij,j}.\label{eq:d_jj}
\end{align}
Substituting~\eqref{eq:rhess_fij_expansion} into~\eqref{eq:rhess_fij_opnorm_bound} yields
\begin{equation}\label{eq:rhess_sup_bound}
    \left\Vert \HessFX\right\Vert_{2}\leq \sum_{\left(i,j\right)\in\mathcal{E}}\sup\left\{ \sqrt{\left\Vert H_{i}\left[\mathcal{U}\right]\right\Vert _{2}^{2}+\left\Vert H_{j}\left[\mathcal{U}\right]\right\Vert _{2}^{2}} \mid \mathcal{U}\in\mathcal{T}_{\mathcal{X}}\mathcal{M},\|\mathcal{U}\|_{2}=1\right\},
\end{equation}
which implies that boundedness of $\left\Vert H_{i}\left[\mathcal{U}\right]\right\Vert _{2}$ and $\left\Vert H_{j}\left[\mathcal{U}\right]\right\Vert _{2}$ 
for all $\X\in\K\subset\MN$ implies boundedness of $\norm{\HessFX}$ for all $\X\in\K\subset\MN$, which we will now show. First, we observe that symmetricity and idempotence of $\mathcal{P}_{\mathbf{x}_{i}}$ implies
\begin{equation}\label{eq:H_i_bound_lambda_P_x}
    \left\Vert H_{i}\left[\mathcal{U}\right]\right\Vert _{2}^{2}=\left\|\mathbf{h}_{i i} u_i+\mathbf{h}_{i j} u_j-\mathbf{d}_{i i} u_i\right\|_{\mathcal{P}_{\mathbf{x}_i}}^2\leq\lambda_{\max}\left(\mathcal{P}_{\mathbf{x}_{i}}\right)\left\Vert \mathbf{h}_{ii}u_{i}+\mathbf{h}_{ij}u_{j}-\mathbf{d}_{ii}u_{i}\right\Vert _{2}^{2}
\end{equation}
where $\lambda_{\max}\left(\mathcal{P}_{\mathbf{x}_{i}}\right)$ denotes the maximum eigenvalue of $\mathcal{P}_{\mathbf{x}_{i}}$, which we compute in the following lemma.

\begin{lemma}\label{eq:P_ortho_eigenvalues}
    For all $\mathbf{x}\in\mathcal{M}$, $\lambda_{\max}\left(\mathcal{P}_{\mathbf{x}}\right)=1$.
\end{lemma}
\emph{Proof:} Letting $\mathbf{x}=[\cos(\phi), \sin(\phi), x_{2}, x_{3}]^{\transpose}$, we observe from~\eqref{eq:Px_def} that
\begin{equation}
    \mathcal{P}_{\mathbf{x}}=\left[\begin{array}{cc}
        \begin{array}{cc}
        \sin\left(\phi\right)^{2} & -\sin\left(\phi\right)\cos\left(\phi\right)\\
        -\sin\left(\phi\right)\cos\left(\phi\right) & \cos\left(\phi\right)^{2}
        \end{array} & \boldsymbol{0}_{2\times2}\\
        \boldsymbol{0}_{2\times2} & I_{2}
        \end{array}\right].
\end{equation}
The characteristic polynomial of $\mathcal{P}_{\mathbf{x}}$, denoted $f\left(\lambda\right)$ is then given by
\begin{equation}
    f\left(\lambda\right)=\left|\lambda I-\mathcal{P}_{\mathbf{x}}\right|=\left(\lambda-\sin^{2}\left(\phi\right)\right)\left(\lambda-\cos^{2}\left(\phi\right)\right)\left(\lambda-1\right)^{2}-\sin^{2}\left(\phi\right)\cos^{2}\left(\phi\right)\left(\lambda-1\right)^{2},
\end{equation}
which simplifies to $f\left(\lambda\right)=\lambda\left(\lambda-1\right)^{3}$. Therefore, the eigenvalues of $\mathcal{P}_{\mathbf{x}}$ are $\left\{ 0,1,1,1\right\}$ for all $\mathbf{x}\in\M$ and $\lambda_{\max}\left(\mathcal{P}_{\mathbf{x}}\right)=1$, completing the proof.~\hfill $\blacksquare$\\

Applying Lemma~\eqref{eq:P_ortho_eigenvalues} and the triangle inequality to~\eqref{eq:H_i_bound_lambda_P_x} yields
\begin{equation}
    \left\Vert H_{i}\left[\mathcal{U}\right]\right\Vert _{2}^{2}\leq\left\Vert \mathbf{h}_{ii}u_{i}+\mathbf{h}_{ij}u_{j}-\mathbf{d}_{ii}u_{i}\right\Vert_{2}^{2}\leq\left(\left\Vert \mathbf{h}_{ii}u_{i}\right\Vert _{2}+\left\Vert \mathbf{h}_{ij}u_{j}\right\Vert _{2}+\left\Vert \mathbf{d}_{ii}u_{i}\right\Vert _{2}\right)^{2}
\end{equation}
and further simplifying gives
\begin{equation}\label{eq:Hi_normbound_pre_1}
    \left\Vert H_{i}\left[\mathcal{U}\right]\right\Vert _{2}^{2}\leq\left(\left\Vert \mathbf{h}_{ii}\right\Vert _{2}\left\Vert u_{i}\right\Vert _{2}+\left\Vert \mathbf{h}_{ij}\right\Vert _{2}\left\Vert u_{j}\right\Vert _{2}+\left\Vert \mathbf{d}_{ii}u_{i}\right\Vert _{2}\right)^{2}
\end{equation}
First, we observe that 
\begin{equation}\label{eq:V_norm_bound}
    \|\mathcal{U}\|_{2}^{2}	=\sum_{l\in\V}\left\Vert u_{l}\right\Vert_{2}^{2}=1,
\end{equation}
which implies that $\left\Vert u_{l}\right\Vert^{2}\leq1$ for all $l\in\mathcal{V}$. Applying this and the fact that $\Vert\cdot\Vert_{2}\leq\Vert\cdot\Vert_{F}$ to~\eqref{eq:Hi_normbound_pre_1} yields
\begin{equation}\label{eq:Hi_normbound_pre_2}
    \left\Vert H_{i}\left[\mathcal{U}\right]\right\Vert _{2}^{2}\leq\left(\left\Vert \mathbf{h}_{ii}\right\Vert _{2}+\left\Vert \mathbf{h}_{ij}\right\Vert _{2}+\left\Vert \mathbf{d}_{ii}u_{i}\right\Vert _{2}\right)^{2}\leq\left(\left\Vert \mathbf{h}_{ii}\right\Vert _{F}+\left\Vert \mathbf{h}_{ij}\right\Vert _{F}+\left\Vert \mathbf{d}_{ii}u_{i}\right\Vert _{2}\right)^{2}.
\end{equation}
To further bound~\eqref{eq:Hi_normbound_pre_2}, we will derive a bound for $\norm{\mathbf{d}_{ii}u_{i}}$, with $\mathbf{d}_{ii}$ given by~\eqref{eq:d_ii}. Letting $\mathbf{x}=[\cos(\phi), \sin(\phi), x_{2}, x_{3}]^{\transpose}$, we observe from~\eqref{eq:Px_norm_def} that
\begin{equation}
    \mathcal{P}_{\mathbf{x}}^{\perp}=\left[\begin{array}{cc}
        \begin{array}{cc}
        \cos^{2}(\phi) & \sin(\phi)\cos(\phi)\\
        \sin(\phi)\cos(\phi) & \sin^{2}(\phi)
        \end{array} & \boldsymbol{0}_{2\times2}\\
        \boldsymbol{0}_{2\times2} & \boldsymbol{0}_{2\times2}
        \end{array}\right]
\end{equation}
Therefore,
\begin{equation}
    \mathbf{x}^{\top}\mathcal{P}_{\mathbf{x}}^{\perp}=\left[\cos(\phi)\left(\cos^{2}(\phi)+\sin^{2}(\phi)\right),\sin(\phi)\left(\cos^{2}(\phi)+\sin^{2}(\phi)\right),0,0\right],
\end{equation}
and simplifying with $\cos^{2}(\phi)+\sin^{2}(\phi)=1$ yields
\begin{equation}
    \mathbf{x}^{\top}\mathcal{P}_{\mathbf{x}}^{\perp}=\left[\cos(\phi),\sin(\phi),0,0\right],
\end{equation}
which holds for all $\mathbf{x}\in\M$. Now, letting $\mathbf{x}=[\cos(\phi_{i}), \sin(\phi_{i}), x_{i,2}, x_{i,3}]^{\transpose}$ and $g_{ij,i}=\left[g_{i,0},g_{i,1},g_{i,2},g_{i,3}\right]^{\top}$, we observe that
\begin{equation}
    \mathbf{d}_{ii}=\tilde{P}\mathbf{x}_{i}^{\top}g_{ij,i}=\operatorname{diag}\left(\left\{
        g_{i,0}\cos\left(\phi_{i}\right)+g_{i,1}\sin\left(\phi_{i}\right), g_{i,0}\cos\left(\phi_{i}\right)+g_{i,1}\sin\left(\phi_{i}\right), 0, 0 \right\}\right).
\end{equation}
Then, letting $u_{i}=\left[u_{i,0}, u_{i,1}, u_{i,2}, u_{i,3}\right]^{\top}$, it holds that
\begin{equation}
    \mathbf{d}_{ii}u_{i}=\left[\left(g_{i,0}\cos\left(\phi_{i}\right)+g_{i,1}\sin\left(\phi_{i}\right)\right)u_{i,0}, \left(g_{i,0}\cos\left(\phi_{i}\right)+g_{i,1}\sin\left(\phi_{i}\right)\right)u_{i,1}, 0, 0\right]^{\top},
\end{equation}
which implies that
\begin{equation}\label{eq:dii_ui_bound_1}
    \left\Vert \mathbf{d}_{ii}u_{i}\right\Vert _{2}=\sqrt{u_{i}^{\top}\mathbf{d}_{ii}^{\top}\mathbf{d}_{ii}u_{i}}=\sqrt{\left(g_{i,0}\cos\left(\phi_{i}\right)+g_{i,1}\sin\left(\phi_{i}\right)\right)^{2}\left(u_{i,0}^{2}+u_{i,1}^{2}\right)}.
\end{equation}
Now, because~\eqref{eq:V_norm_bound} implies that $u_{i,0}^{2}+u_{i,1}^{2}\leq1$ for all $i$,~\eqref{eq:dii_ui_bound_1} simplifies to
\begin{equation}\label{eq:dii_ui_bound}
    \left\Vert \mathbf{d}_{ii}u_{i}\right\Vert _{2}\leq\left|g_{i,0}\cos\left(\phi_{i}\right)+g_{i,1}\sin\left(\phi_{i}\right)\right|\leq\left|g_{i,0}\right|+\left|g_{i,1}\right|,
\end{equation}
and applying~\eqref{eq:dii_ui_bound} to~\eqref{eq:Hi_normbound_pre_2} yields
\begin{equation}\label{eq:H_i_bound}
    \left\Vert H_{i}\left[\mathcal{U}\right]\right\Vert_{2}^{2}	\leq\left(\left\Vert \mathbf{h}_{ii}\right\Vert _{F}+\left\Vert \mathbf{h}_{ij}\right\Vert _{F}+\left|g_{i,0}\right|+\left|g_{i,1}\right|\right)^{2}.
\end{equation}
Furthermore, following the derivation of~\eqref{eq:H_i_bound} for $H_{j}[\U]$ and letting $g_{ij,j}=\left[g_{j,0},g_{j,1},g_{j,2},g_{j,3}\right]^{\top}$ yields
\begin{equation}\label{eq:H_j_bound}
    \left\Vert H_{j}\left[\mathcal{U}\right]\right\Vert_{2}^{2}	\leq\left(\left\Vert \mathbf{h}_{ji}\right\Vert _{F}+\left\Vert \mathbf{h}_{jj}\right\Vert _{F}+\left|g_{j,0}\right|+\left|g_{j,1}\right|\right)^{2}.
\end{equation}
In Appendix~\ref{app:egrad_bounds}, we derive bounds for the Euclidean gradient terms appearing in~\eqref{eq:H_i_bound} and~\eqref{eq:H_j_bound}, namely, $\left|g_{i,0}\right|$, $\left|g_{i,1}\right|$, $\left|g_{j,0}\right|$, and $\left|g_{j,1}\right|$ that hold for all $\X\in\K\subset\MN$, with $\K$ compact. Specifically, from~\eqref{eq:gij_bounds} we have
\begin{equation}\label{eq:gi_gj_bounds}
    \left|g_{i,0}\right|+\left|g_{i,1}\right|\leq2\gbar\text{~and~}\left|g_{j,0}\right|+\left|g_{j,1}\right|\leq2\gbar,
\end{equation}
with constant $\gbar$ given by~\eqref{eq:gbar_def}. Furthermore, in Appendix~\ref{app:ehess_bounds}, we derive bounds for the Euclidean Hessian terms appearing in~\eqref{eq:H_i_bound} and~\eqref{eq:H_j_bound}, namely, $\fnorm{\hii}$, $\fnorm{\hij}$, $\fnorm{\hji}$, and $\fnorm{\hjj}$, that also hold for all $\X\in\K\subset\MN$, with $\K$ compact. Specifically,~\eqref{eq:hii_hjj_bound} and~\eqref{eq:hij_hji_bound} give
\begin{align}
    \fnorm{\hii},\fnorm{\hjj}&\leq\overline{\mathbf{h}}_{ii}\fnorm{\Omegaij}\label{eq:hii_hjj_bound_restate}\\
    \fnorm{\hij},\fnorm{\hji}&\leq\overline{\mathbf{h}}_{ij}\fnorm{\Omegaij}\label{eq:hij_hji_bound_restate},
\end{align}
with $\overline{\mathbf{h}}_{ii}$  and $\overline{\mathbf{h}}_{ij}$ defined in~\eqref{eq:hijbar_def} and~\eqref{eq:hiibar_def}. Applying~\eqref{eq:gi_gj_bounds} and~\eqref{eq:hii_hjj_bound_restate}-~\eqref{eq:hij_hji_bound_restate} to~\eqref{eq:H_i_bound} and~\eqref{eq:H_j_bound} and substituting the result into~\eqref{eq:rhess_fij_expansion} yields
\begin{equation}
    \sup\left\{ \left\Vert \operatorname{Hess}f_{ij}(\mathcal{X})\left[\mathcal{U}\right]\right\Vert _{2}\mid\mathcal{U}\in\mathcal{T}_{\mathcal{X}}\mathcal{M},\|\mathcal{U}\|_{2}=1\right\} \leq\sqrt{2}\left(\left(\overline{\mathbf{h}}_{ii}+\overline{\mathbf{h}}_{ij}\right)\left\Vert \Omega_{ij}\right\Vert _{F}+2\mathbf{g}\right),
\end{equation}
and applying this to~\eqref{eq:rhess_sup_bound} and simplifying yields $\norm{\HessFX}\leq L_{g}$, with
\begin{equation}\label{eq:L_g}
    L_{g}\triangleq\sqrt{2}\left(\overline{\mathbf{h}}_{ii}+\overline{\mathbf{h}}_{ij}\right)\overline{\boldsymbol{\Omega}}+2M\mathbf{g},
\end{equation}
where $\M=|\E|$ and $\overline{\boldsymbol{\Omega}}\triangleq\sum_{\left(i,j\right)\in\mathcal{E}}\left\Vert \Omega_{ij}\right\Vert _{F}$. Equation~\eqref{eq:L_g} gives a Lipschitz constant $L_{g}$ satisfying~\eqref{eq:rhess_op_norm_lipschitz} which holds for all $\X\in\K\subset\MN$, with $\K$ compact. Therefore, the Riemannian gradient from~\eqref{eq:rgrad_F} is Lipschitz continuous on any compact subset of $\MN$, completing our derivation.
\section{Derivation of Euclidean Gradient Jacobians}\label{app:egrad_jacobians}
As derived in Appendix~\ref{app:egrad}, the Jacobians of the tangent residual $\mathbf{e}_{ij}$ in~\eqref{eq:e_ij_app} with respect to $\mathbf{x}_{i},\mathbf{x}_{j}$ are necessary to derive the Euclidean gradient of $\mathcal{F}$ from~\eqref{eq:eucl_grad}. In vector form, we denote $\mathbf{x}_{i}={\left[x_{i,0},~x_{i,1},~x_{i,2},~x_{i,3}\right]}^{\top}$, $\mathbf{x}_{j}={\left[x_{j,0},~x_{j,1},~x_{j,2},~x_{j,3}\right]}^{\top}$, and $\mathbf{e}_{ij}={\left[e_{0},~e_{1},~e_{2}\right]}^{\top}$. In this appendix, we derive the Jacobian matrices $\mathcal{A}_{ij},\mathcal{B}_{ij}\in\mathbb{R}^{3\times4}$, defined as
\begin{equation}
    \mathcal{A}_{ij}\triangleq\left[\begin{array}{ccc}
        \mathcal{A}_{11} & \cdots & \mathcal{A}_{14}\\
        \vdots & \ddots & \vdots\\
        \mathcal{A}_{31} & \cdots & \mathcal{A}_{34}
        \end{array}\right]=\left[\begin{array}{ccc}
        \frac{\partial e_{0}}{\partial x_{i,0}} & \cdots & \frac{\partial e_{0}}{\partial x_{i,3}}\\
        \vdots & \ddots & \vdots\\
        \frac{\partial e_{2}}{\partial x_{i,0}} & \cdots & \frac{\partial e_{2}}{\partial x_{i,3}}
        \end{array}\right],\ \ \ \mathcal{B}_{ij}\triangleq\left[\begin{array}{ccc}
        \mathcal{B}_{11} & \cdots & \mathcal{B}_{14}\\
        \vdots & \ddots & \vdots\\
        \mathcal{B}_{31} & \cdots & \mathcal{B}_{34}
        \end{array}\right]=\left[\begin{array}{ccc}
        \frac{\partial e_{0}}{\partial x_{j,0}} & \cdots & \frac{\partial e_{0}}{\partial x_{j,3}}\\
        \vdots & \ddots & \vdots\\
        \frac{\partial e_{2}}{\partial x_{j,0}} & \cdots & \frac{\partial e_{2}}{\partial x_{j,3}}
        \end{array}\right].
\end{equation}
We first rewrite $\mathbf{e}_{ij}$ in a manner that is conducive to differentiation with respect to $\mathbf{x}_i$ and $\mathbf{x}_j$. Using~\eqref{eq:Q_L_R}-\eqref{eq:Q_L_MM}, the residual $\rij=\zij^{-1}\boxplus\mathbf{x}_{i}^{-1}\boxplus\mathbf{x}_{j}$ can be rewritten as two equivalent expressions, which are given by
\begin{equation}\label{eq:r_ij_expansion}
    \rij=\zij^{-1}\boxplus\mathbf{x}_{i}^{-1}\boxplus\mathbf{x}_{j}=Q_{R}\left(\mathbf{x}_{j}\right)Q_{L}^{--}\left(\zij\right)\mathbf{x}_{i}=Q_{LL}^{-}\left(\zij\right)Q_{LL}^{-}\left(\mathbf{x}_{i}\right)\mathbf{x}_{j}.
\end{equation}
We now define $Q_{i}\triangleq Q_{R}\left(\mathbf{x}_{j}\right)Q_{L}^{--}\left(\zij\right)$ and $Q_{j}\triangleq Q_{LL}^{-}\left(\zij\right)Q_{LL}^{-}\left(\mathbf{x}_{i}\right)$, such that $\rij=Q_{i}\mathbf{x}_{i}=Q_{j}\mathbf{x}_{j}$, and write these matrices in the form
\begin{equation}
    Q_{i}=\left[\begin{array}{cccc}
        \mu_{i} & \omega_{i} & 0 & 0\\
        \eta_{i} & \kappa_{i} & 0 & 0\\
        \alpha_{1} & \beta_{1} & \xi_{1} & \zeta_{1}\\
        \alpha_{2} & \beta_{2} & -\zeta_{1} & \xi_{1}
    \end{array}\right],~~~Q_{j}=\left[\begin{array}{cccc}
        \mu_{j} & \omega_{j} & 0 & 0\\
        \eta_{j} & \kappa_{j} & 0 & 0\\
        \alpha_{3} & \beta_{3} & \kappa_{j} & -\eta_{j}\\
        \beta_{3} & -\alpha_{3} & \eta_{j} & \kappa_{j}
    \end{array}\right],
\label{eq:Q_i_j}
\end{equation}
where the element-wise definitions for $Q_i$ are given by
\begin{align}
    \mu_{i}&\triangleq z_{0}x_{j,0}+z_{1}x_{j,1}, \label{eq:mu_i} \\
    \omega_{i}&\triangleq -z_{1}x_{j,0}+z_{0}x_{j,1}, \label{eq:omega_i} \\
    \eta_{i}&\triangleq -z_{1}x_{j,0}+z_{0}x_{j,1}, \label{eq:eta_i} \\
    \kappa_{i}&\triangleq -z_{0}x_{j,0}-z_{1}x_{j,1}, \label{eq:kappa_i} \\
    \alpha_{1}&\triangleq -z_{2}x_{j,0}-z_{3}x_{j,1}+z_{0}x_{j,2}+z_{1}x_{j,3}, \label{eq:alpha_1} \\
    \beta_{1}&\triangleq z_{3}x_{j,0}-z_{2}x_{j,1}-z_{1}x_{j,2}+z_{0}x_{j,3}, \label{eq:beta_1} \\
    \xi_{1}&\triangleq -z_{0}x_{j,0}+z_{1}x_{j,1}, \label{eq:xi_1} \\
    \zeta_{1}&\triangleq -z_{1}x_{j,0}-z_{0}x_{j,1}, \label{eq:zeta_1} \\
    \alpha_{2}&\triangleq -z_{3}x_{j,0}+z_{2}x_{j,1}-z_{1}x_{j,2}+z_{0}x_{j,3} \label{eq:alpha_3} \\
    \beta_{2}&\triangleq -z_{2}x_{j,0}-z_{3}x_{j,1}-z_{0}x_{j,2}-z_{1}x_{j,3}, \label{eq:beta_3} 
\end{align}
and for $Q_{j}$,
\begin{align}
    \mu_{j}&\triangleq z_{0}x_{i,0}-z_{1}x_{i,1}, \label{eq:mu_j} \\
    \omega_{j}&\triangleq z_{1}x_{i,0}+z_{0}x_{i,1}, \label{eq:omega_j} \\
    \eta_{j}&\triangleq -z_{1}x_{i,0}-z_{0}x_{i,1}, \label{eq:eta_j} \\
    \kappa_{j}&\triangleq z_{0}x_{i,0}-z_{1}x_{i,1}, \label{eq:kappa_j} \\
    \alpha_{3}&\triangleq -z_{2}x_{i,0}+z_{3}x_{i,1}-z_{0}x_{i2}-z_{1}x_{i,3}, \label{eq:alpha_2} \\
    \beta_{3}&\triangleq -z_{3}x_{i,0}-z_{2}x_{i,1}+z_{1}x_{i,2}-z_{0}x_{i,3}. \label{eq:beta_2} 
\end{align}
Letting $\rij=\left[r_0,r_1,r_2,r_3\right]^{\top}$, we can substitute~\eqref{eq:r_ij_expansion}-\eqref{eq:Q_i_j} to expand each term of $\rij$ as
\begin{align}
    r_{0}&=\mu_{i}x_{i,0}+\omega_{i}x_{i,1}=\mu_{j}x_{j,0}+\omega_{j}x_{j,1}\label{eq:r_0_i_j} \\
    r_{1}&=\eta_{i}x_{i,0}+\kappa_{i}x_{i,1}=\eta_{j}x_{j,0}+\kappa_{j}x_{j,1},\label{eq:r_1_i_j} \\
    r_{2}&=\alpha_{1}x_{i,0}+\beta_{1}x_{i,1}+\xi_{1}x_{i,2}+\zeta_{1}x_{i,3}=\alpha_{3}x_{j,0}+\beta_{3}x_{j,1}+\kappa_{j}x_{j,2}-\eta_{j}x_{j,3},\label{eq:r_2_i_j} \\
    r_{3}&=\alpha_{2}x_{i,0}+\beta_{2}x_{i,1}-\zeta_{1}x_{i,2}+\xi_{1}x_{i,3}=\beta_{3}x_{j,0}-\alpha_{3}x_{j,1}+\eta_{j}x_{j,2}+\kappa_{j}x_{j,3},\label{eq:r_3_i_j}
\end{align}
which simplifies the calculation of $\frac{\partial r}{\partial x_{l,m}}$ for any entry $r_{l}$ of $\rij$ and any entry $x_{l,m}$ of $\mathbf{x}_{i},\mathbf{x}_{j}$. From~\eqref{eq:e_ij_def}, letting $\gamma\triangleq\gamma\left(\phi\left(\rij\right)\right)$ yields the element-wise definitions of $\mathbf{e}_{ij}$ to begin
\begin{equation}\label{eq:eij_element_def}
    e_{0}=\frac{r_{1}}{\gamma},\ e_{1}=\frac{r_{2}}{\gamma},\ e_{2}=\frac{r_{3}}{\gamma}.
\end{equation}
Before differentiating $\mathbf{e}_{ij}$, we precompute a general form for partial derivatives of $\gamma$ with respect to any entry $x_{l,m}$ of $\mathbf{x}_{i},\mathbf{x}_{j}$. Letting $\phi\triangleq\phi\left(\rij\right)$ and applying the chain rule to~\eqref{eq:gamma_app} yields
\begin{equation}\label{eq:dgamma_dx_chain}
    \frac{\partial\gamma}{\partial x_{l,m}}=\frac{\partial\gamma}{\partial\phi}\frac{\partial\phi}{\partial x_{l,m}}.
\end{equation}
The term $\frac{\partial\gamma}{\partial\phi}$ is computed by applying the quotient rule to differentiate~\eqref{eq:gamma_app}, yielding
\begin{equation}\label{eq:dgamma_dphi}
    \frac{\partial\gamma}{\partial\phi}	=\frac{\partial}{\partial\phi}\left(\frac{\sin\left(\phi\right)}{\phi}\right)=\frac{\phi\cos\left(\phi\right)-\sin\left(\phi\right)}{\phi^{2}}=\frac{\phi r_{0}-r_{1}}{\phi^{2}}.
\end{equation}
Given the definition of $\phi$ from~\eqref{eq:phi_def_app}, applying the chain rule yields
\begin{equation}\label{eq:dphi_dx_chain}
    \frac{\partial\phi}{\partial x_{l,m}}=\frac{\partial\phi}{\partial r_{0}}\frac{\partial r_{0}}{\partial x_{l,m}}+\frac{\partial\phi}{\partial r_{1}}\frac{\partial r_{1}}{\partial x_{l,m}}.
\end{equation}
We now observe that in~\eqref{eq:phi_def_app}, $\nicefrac{\partial\mathrm{wrap}\left(u\right)}{\partial u}=1$ for all $u\in\left(-\nicefrac{\pi}{2},\nicefrac{\pi}{2}\right)$, and $\phi$ is continuously differentiable on $\left(-\nicefrac{\pi}{2},\nicefrac{\pi}{2}\right)$, with
\begin{equation}
    \frac{\partial}{\partial r_{l}}\left(\arctan\left(r_{1},r_{0}\right)\right)=\frac{\partial}{\partial r_{l}}\left(\arctan\left(\frac{r_{1}}{r_{0}}\right)\right),
\end{equation}
where $\arctan(\nicefrac{u}{v})$ is the two-quadrant arctangent, so we have
\begin{equation}\label{eq:dphi_drij}
    \frac{\partial\phi}{\partial r_{0}}	=-\frac{r_{1}}{r_{0}^{2}+r_{1}^{2}},\quad\frac{\partial\phi}{\partial r_{1}}=\frac{r_{0}}{r_{0}^{2}+r_{1}^{2}}.
\end{equation}
Substituting~\eqref{eq:dphi_drij} into~\eqref{eq:dphi_dx_chain} then gives
\begin{equation}\label{eq:dphi_dx}
    \frac{\partial\phi}{\partial x_{l,m}}	=\left(\frac{1}{r_{0}^{2}+r_{1}^{2}}\right)\left(\frac{\partial r_{1}}{\partial x_{l,m}}r_{0}-\frac{\partial r_{0}}{\partial x_{l,m}}r_{1}\right).
\end{equation}
Substituting~\eqref{eq:dgamma_dphi} and~\eqref{eq:dphi_dx} into~\eqref{eq:dgamma_dx_chain} yields the general form for $\frac{\partial\gamma}{\partial x_{l,m}}$ to be
\begin{equation}\label{eq:dgamma_dx_general}
    \frac{\partial\gamma}{\partial x_{l,m}}=\left(\frac{\phi r_{0}-r_{1}}{\phi^{2}}\right)\left(\frac{1}{r_{0}^{2}+r_{1}^{2}}\right)\left(\frac{\partial r_{1}}{\partial x_{l,m}}r_{0}-\frac{\partial r_{0}}{\partial x_{l,m}}r_{1}\right).
\end{equation}
Using~\eqref{eq:dgamma_dx_general}, it is straightforward to further compute general forms for partial derivatives of $\mathbf{e}_{ij}$ with respect to $\mathbf{x}_{i}, \mathbf{x}_{j}$. For example, applying the quotient rule to differentiate $e_{0}$ from~\eqref{eq:eij_element_def} with respect to any entry $x_{l,m}$ of $\mathbf{x}_{i},\mathbf{x}_{j}$ yields
\begin{equation}
    \frac{\partial e_{0}}{\partial x_{l,m}}=\frac{\partial}{\partial x_{l,m}}\left(\frac{r_{1}}{\gamma}\right)=\frac{\frac{\partial r_{1}}{\partial x_{l,m}}\gamma-r_{1}\frac{\partial\gamma}{\partial x_{l,m}}}{\gamma^{2}},
\end{equation}
and substituting~\eqref{eq:dgamma_dx_general} and simplifying yields
\begin{equation}
    \frac{\partial e_{0}}{\partial x_{l,m}}	=\frac{\frac{\partial r_{1}}{\partial x_{l,m}}}{\gamma}+\frac{r_{1}}{r_{0}^{2}+r_{1}^{2}}\left(\frac{\partial r_{1}}{\partial x_{l,m}}r_{0}-\frac{\partial r_{0}}{\partial x_{l,m}}r_{1}\right)\left(\frac{r_{1}-\phi r_{0}}{\gamma^{2}\phi^{2}}\right),
\end{equation}
which can be further simplified by the fact that $\gamma^{2}\phi^{2}=\sin^{2}\left(\phi\right)=r_{1}^{2}$. Applying this simplification gives the expression
\begin{equation}\label{eq:de0_dx_nof1}
    \frac{\partial e_{0}}{\partial x_{l,m}}	=\frac{\frac{\partial r_{1}}{\partial x_{l,m}}}{\gamma}+\frac{r_{1}}{r_{0}^{2}+r_{1}^{2}}\left(\frac{\partial r_{1}}{\partial x_{l,m}}r_{0}-\frac{\partial r_{0}}{\partial x_{l,m}}r_{1}\right)\left(\frac{r_{1}-\phi r_{0}}{r_{1}^{2}}\right).
\end{equation}
To simplify~\eqref{eq:de0_dx_nof1}, we define the function $f_{1}:\mathbb{R}\rightarrow\mathbb{R}$ as
\begin{equation}\label{eq:f_1_def}
    f_{1}\left(\phi\right)	\triangleq\frac{r_{1}-\phi r_{0}}{r_{1}^{2}}=\frac{\sin\left(\phi\right)-\phi\cos\left(\phi\right)}{\sin^{2}\left(\phi\right)}=\csc^{2}\left(\phi\right)\left(\sin\left(\phi\right)-\phi\cos\left(\phi\right)\right).
\end{equation}
Letting $r_{0}=\cos\left(\phi\right)$ and $r_{1}=\sin\left(\phi\right)$ yields the equivalence
\begin{equation}\label{eq:f_1_equivalence}
    \frac{r_{1}-\phi r_{0}}{r_{1}^{2}}=\frac{\sin\left(\phi\right)-\phi\cos\left(\phi\right)}{\sin^{2}\left(\phi\right)}=f_{1}\left(\phi\right).
\end{equation}
Letting $f_{1}\triangleq f_{1}\left(\phi\right)$ and substituting~\eqref{eq:f_1_equivalence} into~\eqref{eq:de0_dx_nof1} yields the general form for $\frac{\partial e_{0}}{\partial x_{l,m}}$ to be
\begin{equation}\label{eq:de0_dx_general}
    \frac{\partial e_{0}}{\partial x_{l,m}}=\frac{\frac{\partial r_{1}}{\partial x_{l,m}}}{\gamma}+\frac{r_{1}}{r_{0}^{2}+r_{1}^{2}}\left(\frac{\partial r_{1}}{\partial x_{l,m}}r_{0}-\frac{\partial r_{0}}{\partial x_{l,m}}r_{1}\right)f_{1}.
\end{equation}
From~\eqref{eq:r_0_i_j}, it is straightforward to compute the derivatives
\begin{equation}\label{eq:dr0_dxi}
    \frac{\partial r_{0}}{\partial x_{i,0}}	=\mu_{i}\ \frac{\partial r_{0}}{x_{i,1}}=\omega_{i},\ \frac{\partial r_{0}}{\partial x_{i,2}}=\frac{\partial r_{0}}{\partial x_{i,3}}=0,
\end{equation}
and
\begin{equation}\label{eq:dr0_dxj}
    \frac{\partial r_{0}}{\partial x_{j,0}}	=\mu_{j},\ \frac{\partial r_{0}}{\partial x_{j,1}}=\omega_{j},\ \frac{\partial r_{0}}{\partial x_{j,2}}=\frac{\partial r_{0}}{\partial x_{j,3}}=0,
\end{equation}
Similarly, differentiating~\eqref{eq:r_1_i_j} gives
\begin{equation}\label{eq:dr1_dxi}
    \frac{\partial r_{1}}{\partial x_{i,0}}	=\eta_{i},\ \frac{\partial r_{1}}{\partial x_{i,1}}=\kappa_{i},\ \frac{\partial r_{1}}{\partial x_{i,2}}=\frac{\partial r_{1}}{\partial x_{i,3}}=0,
\end{equation}
and
\begin{equation}\label{eq:dr1_dxj}
    \frac{\partial r_{1}}{\partial x_{j,0}}=\eta_{j},\ \frac{\partial r_{1}}{\partial x_{j,1}}=\kappa_{j},\ \frac{\partial r_{1}}{\partial x_{j,2}}=\frac{\partial r_{1}}{\partial x_{j,3}}=0.
\end{equation}
Substituting~\eqref{eq:dr0_dxi}-\eqref{eq:dr1_dxj} into the general form given by~\eqref{eq:de0_dx_general} yields $\mathcal{A}_{11}$-$\mathcal{A}_{14}$ and $\mathcal{B}_{11}$-$\mathcal{B}_{14}$ to be
\begin{align}
    \mathcal{A}_{11}=\frac{\partial e_{0}}{\partial x_{i,0}}&=\frac{\eta_{i}}{\gamma}+\frac{r_{1}}{r_{0}^{2}+r_{1}^{2}}\left(\eta_{i}r_{0}-\mu_{i}r_{1}\right)f_{1}, \label{eq:A_11}\\
    \mathcal{A}_{12}=\frac{\partial e_{0}}{\partial x_{i,1}}&=\frac{\kappa_{i}}{\gamma}+\frac{r_{1}}{r_{0}^{2}+r_{1}^{2}}\left(\kappa_{i}r_{0}-\omega_{i}r_{1}\right)f_{1},\label{eq:A_12}\\
    \mathcal{A}_{13}=\frac{\partial e_{0}}{\partial x_{i,2}}&=0,\label{eq:A_13}\\
    \mathcal{A}_{14}=\frac{\partial e_{0}}{\partial x_{i,3}}&=0,\label{eq:A_14}\\
    \mathcal{B}_{11}=\frac{\partial e_{0}}{\partial x_{j,0}}&=\frac{\eta_{j}}{\gamma}+\frac{r_{1}}{r_{0}^{2}+r_{1}^{2}}\left(\eta_{j}r_{0}-\mu_{j}r_{1}\right)f_{1},\label{eq:B_11}\\
    \mathcal{B}_{12}=\frac{\partial e_{0}}{\partial x_{j,1}}&=\frac{\kappa_{j}}{\gamma}+\frac{r_{1}}{r_{0}^{2}+r_{1}^{2}}\left(\kappa_{j}r_{0}-\omega_{j}r_{1}\right)f_{1},\label{eq:B_12}\\
    \mathcal{B}_{13}=\frac{\partial e_{0}}{\partial x_{j,2}}&=0,\label{eq:B_13}\\
    \mathcal{B}_{14}=\frac{\partial e_{0}}{\partial x_{j,3}}&=0.\label{eq:B_14}
\end{align}
Because $\frac{\partial e_{1}}{\partial x_{l,m}}$ has the same structure as $\frac{\partial e_{0}}{\partial x_{l,m}}$, its general form is computed to be
\begin{equation}
    \frac{\partial e_{1}}{\partial x_{l,m}}	=\frac{\partial}{\partial x_{l,m}}\left(\frac{r_{2}}{\gamma}\right)=\frac{\frac{\partial r_{2}}{\partial x_{l,m}}}{\gamma}+\frac{r_{2}}{r_{0}^{2}+r_{1}^{2}}\left(\frac{\partial r_{1}}{\partial x_{l,m}}r_{0}-\frac{\partial r_{0}}{\partial x_{l,m}}r_{1}\right)f_{1}.
\label{eq:de1_dx_general}
\end{equation}
From~\eqref{eq:r_2_i_j}, we have the derivatives
\begin{equation}\label{eq:dr2_dxi}
    \frac{\partial r_{2}}{x_{i,0}}	=\alpha_{1},\ \frac{\partial r_{2}}{x_{i,1}}=\beta_{1},\ \frac{\partial r_{2}}{x_{i,2}}=\xi_{1},\ \frac{\partial r_{2}}{x_{i,3}}=\zeta_{1},
\end{equation}
and
\begin{equation}\label{eq:dr2_dxj}
    \frac{\partial r_{2}}{x_{j0}}	=\alpha_{3},\ \frac{\partial r_{2}}{x_{j1}}=\beta_{3},\ \frac{\partial r_{2}}{x_{j2}}=\kappa_{j},\ \frac{\partial r_{2}}{x_{j3}}=-\eta_{j}.
\end{equation}
The terms $\mathcal{A}_{21}-\mathcal{A}_{24}$ and $\mathcal{B}_{21}-\mathcal{B}_{24}$ are then computed by substituting~\eqref{eq:dr0_dxi}-\eqref{eq:dr1_dxj} and~\eqref{eq:dr2_dxi}-\eqref{eq:dr2_dxj} into~\eqref{eq:de1_dx_general}, yielding
\begin{align}
    \mathcal{A}_{21}=\frac{\partial e_{1}}{\partial x_{i,0}}&=\frac{\alpha_{1}}{\gamma}+\frac{r_{2}}{r_{0}^{2}+r_{1}^{2}}\ensuremath{\left(\eta_{i}r_{0}-\mu_{i}r_{1}\right)f_{1}}\label{eq:A_21}\\
    \mathcal{A}_{22}=\frac{\partial e_{1}}{\partial x_{i,1}}&=\frac{\beta_{1}}{\gamma}+\frac{r_{2}}{r_{0}^{2}+r_{1}^{2}}\ensuremath{\left(\kappa_{i}r_{0}-\omega_{i}r_{1}\right)}f_{1}\label{eq:A_22}\\
    \mathcal{A}_{23}=\frac{\partial e_{1}}{\partial x_{i,2}}&=\frac{\xi_{1}}{\gamma}\label{eq:A_23}\\
    \mathcal{A}_{24}=\frac{\partial e_{1}}{\partial x_{i,3}}&=\frac{\zeta_{1}}{\gamma}\label{eq:A_24}\\
    \mathcal{B}_{21}=\frac{\partial e_{1}}{\partial x_{j,0}}&=\frac{\alpha_{3}}{\gamma}+\frac{r_{2}}{r_{0}^{2}+r_{1}^{2}}\ensuremath{\left(\eta_{j}r_{0}-\mu_{j}r_{1}\right)}f_{1}\label{eq:B_21}\\
    \mathcal{B}_{22}=\frac{\partial e_{1}}{\partial x_{j,1}}&=\frac{\beta_{3}}{\gamma}+\frac{r_{2}}{r_{0}^{2}+r_{1}^{2}}\ensuremath{\left(\kappa_{j}r_{0}-\omega_{j}r_{1}\right)f_{1}}\label{eq:B_22}\\
    \mathcal{B}_{23}=\frac{\partial e_{1}}{\partial x_{j,2}}&=\frac{\kappa_{j}}{\gamma}\label{eq:B_23}\\
    \mathcal{B}_{24}=\frac{\partial e_{1}}{\partial x_{j,3}}&=-\frac{\eta_{j}}{\gamma}\label{eq:B_24}\\
\end{align}
The final derivative, $\frac{\partial e_{2}}{\partial x_{l,m}}$, also has the same structure as $\frac{\partial e_{0}}{\partial x_{l,m}}$, so its general form is given by
\begin{equation}
    \frac{\partial e_{2}}{\partial x_{l,m}}	=\frac{\partial}{\partial x_{l,m}}\left(\frac{r_{3}}{\gamma}\right)=\frac{\frac{\partial r_{3}}{\partial x_{l,m}}}{\gamma}+\frac{r_{3}}{r_{0}^{2}+r_{1}^{2}}\left(\frac{\partial r_{1}}{\partial x_{l,m}}r_{0}-\frac{\partial r_{0}}{\partial x_{l,m}}r_{1}\right)f_{1}.
\label{eq:de2_dx_general}
\end{equation}
From equations~\eqref{eq:r_3_i_j}, we have the derivatives
\begin{equation}
    \frac{\partial r_{3}}{x_{i,0}}	=\alpha_{2},\ \frac{\partial r_{3}}{x_{i,1}}=\beta_{2},\ \frac{\partial r_{3}}{x_{i,2}}=-\zeta_{1},\ \frac{\partial r_{3}}{x_{i,3}}=\xi_{1},
\label{eq:dr3_dxi}
\end{equation}
and
\begin{equation}
    \frac{\partial r_{3}}{x_{j,0}}	=\beta_{3},\ \frac{\partial r_{3}}{x_{j,1}}=-\alpha_{3},\ \frac{\partial r_{3}}{x_{j,2}}=\eta_{j},\ \frac{\partial r_{3}}{x_{j,3}}=\kappa_{j}.
\label{eq:dr3_dxj}
\end{equation}
Finally, the terms $\mathcal{A}_{31}-\mathcal{A}_{34}$ and $\mathcal{B}_{31}-\mathcal{B}_{34}$ are computed by substituting~\eqref{eq:dr0_dxi}-\eqref{eq:dr1_dxj} and~\eqref{eq:dr3_dxi}-\eqref{eq:dr3_dxj} into~\eqref{eq:de2_dx_general}, yielding
\begin{align}
    \mathcal{A}_{31}=\frac{\partial e_{2}}{\partial x_{i,0}}&=\frac{\alpha_{2}}{\gamma}+\frac{r_{3}}{r_{0}^{2}+r_{1}^{2}}\ensuremath{\left(\eta_{i}r_{0}-\mu_{i}r_{1}\right)}f_{1},\label{eq:A_31}\\
    \mathcal{A}_{32}=\frac{\partial e_{2}}{\partial x_{i,1}}&=\frac{\beta_{2}}{\gamma}+\frac{r_{3}}{r_{0}^{2}+r_{1}^{2}}\ensuremath{\left(\kappa_{i}r_{0}-\omega_{i}r_{1}\right)f_{1}},\label{eq:A_32}\\
    \mathcal{A}_{33}=\frac{\partial e_{2}}{\partial x_{i,2}}&=-\frac{\zeta_{1}}{\gamma},\label{eq:A_33}\\
    \mathcal{A}_{34}=\frac{\partial e_{2}}{\partial x_{i,3}}&=\frac{\xi_{1}}{\gamma},\label{eq:A_34}\\
    \mathcal{B}_{31}=\frac{\partial e_{2}}{\partial x_{j,0}}&=\frac{\beta_{3}}{\gamma}+\frac{r_{3}}{r_{0}^{2}+r_{1}^{2}}\ensuremath{\left(\eta_{j}r_{0}-\mu_{j}r_{1}\right)}f_{1},\label{eq:B_31}\\
    \mathcal{B}_{32}=\frac{\partial e_{2}}{\partial x_{j,1}}&=-\frac{\alpha_{3}}{\gamma}+\frac{r_{3}}{r_{0}^{2}+r_{1}^{2}}\ensuremath{\left(\kappa_{j}r_{0}-\omega_{j}r_{1}\right)f_{1},}\label{eq:B_32}\\
    \mathcal{B}_{33}=\frac{\partial e_{2}}{\partial x_{j,2}}&=\frac{\eta_{j}}{\gamma},\label{eq:B_33}\\
    \mathcal{B}_{34}=\frac{\partial e_{2}}{\partial x_{j,3}}&=\frac{\kappa_{j}}{\gamma}.\label{eq:B_34}\\
\end{align}
which concludes the derivation of Jacobians $\mathcal{A}_{ij}$ and $\mathcal{B}_{ij}$.
\section{Derivation of Euclidean Hessian Tensors}\label{app:ehess_tensors}
Here we compute the quantities $\frac{\partial}{\partial\mathbf{x}_{i}}\mathcal{A}_{ij}$, $\frac{\partial}{\partial\mathbf{x}_{j}}\mathcal{A}_{ij}$, $\frac{\partial}{\partial\mathbf{x}_{i}}\mathcal{B}_{ij}$, and $\frac{\partial}{\partial\mathbf{x}_{j}}\mathcal{B}_{ij}$. Because we are differentiating a matrix in $\mathbb{R}^{3\times4}$ with respect to a vector in $\mathbb{R}^{4}$, each of these quantities represents a tensor in $\mathbb{R}^{3\times4\times4}$, in which the third dimension encodes the index of a respective entry in $\mathbf{x}_{i}$ or $\mathbf{x}_{j}$. We note that since further derivatives will not be taken, we are directly computing the implementation form of each of the expressions in this section.
\subsection{Partial Derivatives of $\mathcal{A}_{ij}$}\label{app:dAij_dx}
We begin by deriving a general form for differentiating $\mathcal{A}_{11}$, which is given by~\eqref{eq:A_11}, with respect to any entry $x_{l,m}$ of $\mathbf{x}_{i},\mathbf{x}_{j}$. We first separate the derivative as
\begin{equation}
    \frac{\partial\mathcal{A}_{11}}{\partial x_{l,m}}	=\frac{\partial}{\partial x_{l,m}}\left(\frac{\eta_{i}}{\gamma}+\frac{r_{1}}{r_{0}^{2}+r_{1}^{2}}\left(\eta_{i}r_{0}-\mu_{i}r_{1}\right)f_{1}\right)=\frac{\partial}{\partial x_{l,m}}\left(\frac{\eta_{i}}{\gamma}\right)+\frac{\partial}{\partial x_{l,m}}\left(\frac{r_{1}}{r_{0}^{2}+r_{1}^{2}}\left(\eta_{i}r_{0}-\mu_{i}r_{1}\right)f_{1}\right).
    \label{eq:dA11_sep}
\end{equation}
We first examine the left-hand derivative in equation~\eqref{eq:dA11_sep}. Applying the quotient rule yields
\begin{equation}
    \frac{\partial}{\partial x_{l,m}}\left(\frac{\eta_{i}}{\gamma}\right)=\frac{1}{\gamma^{2}}\left(\frac{\partial\eta_{i}}{\partial x_{l,m}}\gamma-\eta_{i}\frac{\partial\gamma}{\partial x_{l,m}}\right).
\label{eq:dA11_left_stmt}
\end{equation}
We now substitute~\eqref{eq:dgamma_dx_general} into~\eqref{eq:dA11_left_stmt} and simplify to obtain
\begin{equation}
    \frac{\partial}{\partial x_{l,m}}\left(\frac{\eta_{i}}{\gamma}\right)=\frac{\frac{\partial\eta_{i}}{\partial x_{l,m}}}{\gamma}+\frac{\eta_{i}}{r_{0}^{2}+r_{1}^{2}}\left(\frac{\partial r_{1}}{\partial x_{l,m}}r_{0}-\frac{\partial r_{0}}{\partial x_{l,m}}r_{1}\right)f_{1}.
\label{eq:dA11_left_unsimp}
\end{equation}
Since we are solving for the implemention form directly, we can subtitute $r_{0}^{2}+r_{1}^{2}=1$ into~\eqref{eq:dA11_left_stmt} to obtain
\begin{equation}
    \frac{\partial}{\partial x_{l,m}}\left(\frac{\eta_{i}}{\gamma}\right)	=\frac{\frac{\partial\eta_{i}}{\partial x_{l,m}}}{\gamma}+\eta_{i}\left(\frac{\partial r_{1}}{\partial x_{l,m}}r_{0}-\frac{\partial r_{0}}{\partial x_{l,m}}r_{1}\right)f_{1}.
\label{eq:dA11_left}
\end{equation}
We now address the right-hand derivative from equation~\eqref{eq:dA11_sep}. Applying the product rule twice yields 
\begin{align}
    \frac{\partial}{\partial x_{l,m}}\left(\frac{r_{1}}{r_{0}^{2}+r_{1}^{2}}\left(\eta_{i}r_{0}-\mu_{i}r_{1}\right)f_{1}\right)=&~\frac{\partial}{\partial x_{l,m}}\left(\frac{r_{1}}{r_{0}^{2}+r_{1}^{2}}\right)\left(\eta_{i}r_{0}-\mu_{i}r_{1}\right)f_{1} \\
    &+\frac{r_{1}}{r_{0}^{2}+r_{1}^{2}}\frac{\partial}{\partial x_{l,m}}\left(\eta_{i}r_{0}-\mu_{i}r_{1}\right)f_{1}\\
    &+\frac{r_{1}}{r_{0}^{2}+r_{1}^{2}}\left(\eta_{i}r_{0}-\mu_{i}r_{1}\right)\frac{\partial f_{1}}{\partial x_{l,m}}.\label{eq:dA11_right_stmt}
\end{align}
The expression given by~\eqref{eq:dA11_right_1} have three derivative terms, which we will now compute invidually. For the first term from the top, applying the quotient rule and simplifying yields
\begin{equation}
    \frac{\partial}{\partial x_{l,m}}\left(\frac{r_{1}}{r_{0}^{2}+r_{1}^{2}}\right)=\frac{\frac{\partial r_{1}}{\partial x_{l,m}}}{r_{0}^{2}+r_{1}^{2}}-2\frac{r_{1}}{\left(r_{0}^{2}+r_{1}^{2}\right)^{2}}\left(\frac{\partial r_{0}}{\partial x_{l,m}}r_{0}+\frac{\partial r_{1}}{\partial x_{l,m}}r_{1}\right).
\end{equation}
Applying the constraint equation $r_{0}^{2}+r_{1}^{2}=1$ then yields
\begin{equation}\label{eq:dA11_right_1}
    \frac{\partial}{\partial x_{l,m}}\left(\frac{r_{1}}{r_{0}^{2}+r_{1}^{2}}\right)=\frac{\partial r_{1}}{\partial x_{l,m}}-2r_{1}\left(\frac{\partial r_{0}}{\partial x_{l,m}}r_{0}+\frac{\partial r_{1}}{\partial x_{l,m}}r_{1}\right).
\end{equation}
For the second term from the top of~\eqref{eq:dA11_right_1}, we simply distribute and apply the product rule, which gives
\begin{equation}\label{eq:dA11_right_2}
    \frac{\partial}{\partial x_{l,m}}\left(\eta_{i}r_{0}-\mu_{i}r_{1}\right)=\eta_{i}\frac{\partial r_{0}}{\partial x_{l,m}}-\mu_{i}\frac{\partial r_{1}}{\partial x_{l,m}}+\frac{\partial\eta_{i}}{\partial x_{l,m}}r_{0}-\frac{\partial\mu_{i}}{\partial x_{l,m}}r_{1}.
\end{equation}
To compute the third term, we apply the chain rule to write
\begin{equation}\label{eq:df_1_stmt}
    \frac{\partial f_{1}}{\partial x_{l,m}}=\frac{\partial f_{1}}{\partial\phi}\frac{\partial\phi}{\partial x_{l,m}},
\end{equation}
where $\nicefrac{\partial\phi}{\partial x_{l,m}}$ is given by~\eqref{eq:dphi_dx}. For $\nicefrac{\partial f_{1}}{\partial\phi}$, with $f_{1}$ given by~\eqref{eq:f_1_def}, a combination of quotient, chain, and product rules and trigonometric simplifications is applied to write
\begin{align}
    \frac{\partial f_{1}}{\partial\phi}&=\frac{\partial}{\partial\phi}\left(\frac{\sin\left(\phi\right)-\phi\cos\left(\phi\right)}{\sin^{2}\left(\phi\right)}\right)\\
	&=\left(\frac{1}{\sin^{4}\left(\phi\right)}\right)\left(\frac{\partial}{\partial\phi}\left(\sin\left(\phi\right)-\phi\cos\left(\phi\right)\right)\sin^{2}\left(\phi\right)-\left(\sin\left(\phi\right)-\phi\cos\left(\phi\right)\right)\frac{\partial}{\partial\phi}\sin^{2}\left(\phi\right)\right)\\
	&=\left(\frac{1}{\sin^{4}\left(\phi\right)}\right)\left(\left(\phi\sin\left(\phi\right)\right)\sin^{2}\left(\phi\right)-\left(\sin\left(\phi\right)-\phi\cos\left(\phi\right)\right)\left(2\sin\left(\phi\right)\cos\left(\phi\right)\right)\right)\\
	&=\left(\frac{1}{\sin\left(\phi\right)}\right)\left(\phi-2\frac{\cos\left(\phi\right)}{\sin\left(\phi\right)}+2\phi\frac{\cos^{2}\left(\phi\right)}{\sin^{2}\left(\phi\right)}\right)\\
	&=\csc\left(\phi\right)\left(\phi-2\cot\left(\phi\right)+2\phi\cot^{2}\left(\phi\right)\right).
\end{align}
We now define the function $f_{2}:\mathbb{R}\rightarrow\mathbb{R}$ as
\begin{equation}\label{eq:f_2_def}
    f_{2}\left(\phi\right)	\triangleq\csc\left(\phi\right)\left(\phi-2\cot\left(\phi\right)+2\phi\cot^{2}\left(\phi\right)\right),
\end{equation}
so that $\nicefrac{\partial f_{1}}{\partial\phi}=f_{2}$. Substituting equations~\eqref{eq:f_2_def} and~\eqref{eq:dphi_dx} into equation~\eqref{eq:df_1_stmt} now gives
\begin{equation}\label{eq:dA11_right_3}
    \frac{\partial f_{1}}{\partial x_{l,m}}=\left(\frac{1}{r_{0}^{2}+r_{1}^{2}}\right)\left(\frac{\partial r_{1}}{\partial x_{l,m}}r_{0}-\frac{\partial r_{0}}{\partial x_{l,m}}r_{1}\right)f_{2},
\end{equation}
Substituting~\eqref{eq:dA11_right_1},~\eqref{eq:dA11_right_2}, and~\eqref{eq:dA11_right_3} into equation~\eqref{eq:dA11_right_stmt}, and letting $r_{0}^{2}+r_{1}^{2}=1$ yields 
\begin{align}
    \frac{\partial}{\partial x_{l,m}}\left(\frac{r_{1}}{r_{0}^{2}+r_{1}^{2}}\left(\eta_{i}r_{0}-\mu_{i}r_{1}\right)f_{1}\right)=&\left(\frac{\partial r_{1}}{\partial x_{l,m}}-2r_{1}\left(\frac{\partial r_{0}}{\partial x_{l,m}}r_{0}+\frac{\partial r_{1}}{\partial x_{l,m}}r_{1}\right)\right)\left(\eta_{i}r_{0}-\mu_{i}r_{1}\right)f_{1}  \\
    &+r_{1}\left(\frac{\partial r_{0}}{\partial x_{l,m}}\eta_{i}-\frac{\partial r_{1}}{\partial x_{l,m}}\mu_{i}+\frac{\partial\eta_{i}}{\partial x_{l,m}}r_{0}-\frac{\partial\mu_{i}}{\partial x_{l,m}}r_{1}\right)f_{1}  \\
    &+r_{1}\left(\eta_{i}r_{0}-\mu_{i}r_{1}\right)\left(\frac{\partial r_{1}}{\partial x_{l,m}}r_{0}-\frac{\partial r_{0}}{\partial x_{l,m}}r_{1}\right)f_{2}\label{eq:dA11_right}
\end{align}
Finally, substituting~\eqref{eq:dA11_left} and~\eqref{eq:dA11_right} back into equation~\eqref{eq:dA11_sep} and simplifying yields the general form for derivatives of $\mathcal{A}_{11}$ as
\begin{align}
    \frac{\partial\mathcal{A}_{11}}{\partial x_{l,m}}=&~\frac{\frac{\partial\eta_{i}}{\partial x_{l,m}}}{\gamma}+\left(\eta_{i}\left(\frac{\partial r_{1}}{\partial x_{l,m}}r_{0}-\frac{\partial r_{0}}{\partial x_{l,m}}r_{1}\right)+\left(\frac{\partial r_{1}}{\partial x_{l,m}}-2r_{1}\left(\frac{\partial r_{0}}{\partial x_{l,m}}r_{0}+\frac{\partial r_{1}}{\partial x_{l,m}}r_{1}\right)\right)\left(\eta_{i}r_{0}-\mu_{i}r_{1}\right)\right)f_{1}  \\
	&+r_{1}\left(\frac{\partial r_{0}}{\partial x_{l,m}}\eta_{i}-\frac{\partial r_{1}}{\partial x_{l,m}}\mu_{i}+\frac{\partial\eta_{i}}{\partial x_{l,m}}r_{0}-\frac{\partial\mu_{i}}{\partial x_{l,m}}r_{1}\right)f_{1}+r_{1}\left(\eta_{i}r_{0}-\mu_{i}r_{1}\right)\left(\frac{\partial r_{1}}{\partial x_{l,m}}r_{0}-\frac{\partial r_{0}}{\partial x_{l,m}}r_{1}\right)f_{2}. \label{eq:dA11_dx_general}
\end{align}
Now, from~\eqref{eq:eta_i}, it is straightforward to compute
\begin{equation}\label{eq:detai_dxi}
    \frac{\partial\eta_{i}}{\partial x_{i,0}}=\frac{\partial\eta_{i}}{\partial x_{i,1}}=\frac{\partial\eta_{i}}{\partial x_{i,2}}=\frac{\partial\eta_{i}}{\partial x_{i,3}}=0,
\end{equation}
and
\begin{equation}\label{eq:detai_dxj}
    \frac{\partial\eta_{i}}{\partial x_{j,0}}=-z_{1},\ \frac{\partial\eta_{i}}{\partial x_{j,1}}=z_{0},\ \frac{\partial\eta_{i}}{\partial x_{j,2}}=\frac{\partial\eta_{i}}{\partial x_{j,3}}=0.
\end{equation}
From~\eqref{eq:mu_i}, we have
\begin{equation}\label{eq:dmui_dxi}
    \frac{\partial\mu_{i}}{\partial x_{i,0}}=\frac{\partial\mu_{i}}{\partial x_{i,1}}=\frac{\partial\mu_{i}}{\partial x_{i,2}}=\frac{\partial\mu_{i}}{\partial x_{i,3}}=0,
\end{equation}
and
\begin{equation}\label{eq:dmui_dxj}
    \frac{\partial\mu_{i}}{\partial x_{j,0}}=z_{0},\ \frac{\partial\mu_{i}}{\partial x_{j,1}}=z_{1},\ \frac{\partial\mu_{i}}{\partial x_{j,2}}=\frac{\partial\mu_{i}}{\partial x_{j,3}}=0.
\end{equation}
Substituting equations~\eqref{eq:dr0_dxi}-\eqref{eq:dr1_dxj} and~\eqref{eq:detai_dxi}-\eqref{eq:dmui_dxj} into~\eqref{eq:dA11_dx_general} yields the following expressions for $\nicefrac{\partial\mathcal{A}_{11}}{\partial x_{l,m}}$.
\begin{align}
    \frac{\partial\mathcal{A}_{11}}{\partial x_{i,0}}=&~2\left(\eta_{i}-r_{1}\left(\mu_{i}r_{0}+\eta_{i}r_{1}\right)\right)\left(\eta_{i}r_{0}-\mu_{i}r_{1}\right)f_{1}+r_{1}\left(\eta_{i}r_{0}-\mu_{i}r_{1}\right)^{2}f_{2}, \label{eq:dA11_dxi0}\\
    \frac{\partial\mathcal{A}_{11}}{\partial x_{i,1}}=&~\left(\eta_{i}\left(\kappa_{i}r_{0}-\omega_{i}r_{1}\right)+\left(\kappa_{i}-2r_{1}\left(\omega_{i}r_{0}+\kappa_{i}r_{1}\right)\right)\left(\eta_{i}r_{0}-\mu_{i}r_{1}\right)+r_{1}\right)f_{1}  \\
    &+r_{1}\left(\eta_{i}r_{0}-\mu_{i}r_{1}\right)\left(\kappa_{i}r_{0}-\omega_{i}r_{1}\right)f_{2}, \label{eq:dA11_dxi1}\\
    \frac{\partial\mathcal{A}_{11}}{\partial x_{i,2}}=&~\frac{\partial\mathcal{A}_{11}}{\partial x_{i,3}}=0, \label{eq:dA11_dxi2_xi3}\\
    \frac{\partial\mathcal{A}_{11}}{\partial x_{j,0}}=&-\frac{z_{1}}{\gamma}+\left(\eta_{i}\left(\eta_{j}r_{0}-\mu_{j}r_{1}\right)+\left(\eta_{j}-2r_{1}\left(\mu_{j}r_{0}+\eta_{j}r_{1}\right)\right)\left(\eta_{i}r_{0}-\mu_{i}r_{1}\right)\right)f_{1}  \\
    &+r_{1}\left(\mu_{j}\eta_{i}-\eta_{j}\mu_{i}-z_{1}r_{0}-z_{0}r_{1}\right)f_{1}  \\
    &+r_{1}\left(\eta_{i}r_{0}-\mu_{i}r_{1}\right)\left(\eta_{j}r_{0}-\mu_{j}r_{1}\right)f_{2}, \label{eq:dA11_dxj0}\\
    \frac{\partial\mathcal{A}_{11}}{\partial x_{j,1}}=&~\frac{z_{0}}{\gamma}+\left(\eta_{i}\left(\kappa_{j}r_{0}-\omega_{j}r_{1}\right)+\left(\kappa_{j}-2r_{1}\left(\omega_{j}r_{0}+\kappa_{j}r_{1}\right)\right)\left(\eta_{i}r_{0}-\mu_{i}r_{1}\right)\right)f_{1}  \\
    &+r_{1}\left(\omega_{j}\eta_{i}-\kappa_{j}\mu_{i}+z_{0}r_{0}-z_{1}r_{1}\right)f_{1}  \\
    &+r_{1}\left(\eta_{i}r_{0}-\mu_{i}r_{1}\right)\left(\kappa_{j}r_{0}-\omega_{j}r_{1}\right)f_{2},\label{eq:dA11_dxj1}\\
    \frac{\partial\mathcal{A}_{11}}{\partial x_{j,2}}=&~\frac{\partial\mathcal{A}_{11}}{\partial x_{j,3}}=0,\label{eq:dA11_dxj2_xj3}
\end{align}
where we have additionally used the fact that
\begin{equation}
    \omega_{i}\eta_{i}-\kappa_{i}\mu_{i}=\cos\left(\phi_{j}-\phi_{z}\right)^{2}+\sin\left(\phi_{j}-\phi_{z}\right)^{2}=1
    \label{eq:omega_eta_minus_kappa_mu}
\end{equation}
to simplify~\eqref{eq:dA11_dxi1}. Furthermore, since $\mathcal{A}_{12}$, which is given by~\eqref{eq:A_12}, has identical structure to $\mathcal{A}_{11}$, the general form for its partial derivatives is computed as
\begin{align}
    \frac{\partial\mathcal{A}_{12}}{\partial x_{l,m}}=&~\frac{\frac{\partial\kappa_{i}}{\partial x_{l,m}}}{\gamma}+\left(\kappa_{i}\left(\frac{\partial r_{1}}{\partial x_{l,m}}r_{0}-\frac{\partial r_{0}}{\partial x_{l,m}}r_{1}\right)+\left(\frac{\partial r_{1}}{\partial x_{l,m}}-2r_{1}\left(\frac{\partial r_{0}}{\partial x_{l,m}}r_{0}+\frac{\partial r_{1}}{\partial x_{l,m}}r_{1}\right)\right)\left(\kappa_{i}r_{0}-\omega_{i}r_{1}\right)\right)f_{1}  \\
	&+r_{1}\left(\frac{\partial r_{0}}{\partial x_{l,m}}\kappa_{i}-\frac{\partial r_{1}}{\partial x_{l,m}}\omega_{i}+\frac{\partial\kappa_{i}}{\partial x_{l,m}}r_{0}-\frac{\partial\omega_{i}}{\partial x_{l,m}}r_{1}\right)f_{1}  \\
	&+r_{1}\left(\kappa_{i}r_{0}-\omega_{i}r_{1}\right)\left(\frac{\partial r_{1}}{\partial x_{l,m}}r_{0}-\frac{\partial r_{0}}{\partial x_{l,m}}r_{1}\right)f_{2}. \label{eq:dA12_dx_general}
\end{align}
From~\eqref{eq:kappa_i}, it is straightforward to compute
\begin{equation}
    \frac{\partial\kappa_{i}}{\partial x_{i,0}}	=\frac{\partial\kappa_{i}}{\partial x_{i,1}}=\frac{\partial\kappa_{i}}{\partial x_{i,2}}=\frac{\partial\kappa_{i}}{\partial x_{i,3}}=0,
\label{eq:dkappai_dxi}
\end{equation}
and
\begin{equation}
    \frac{\partial\kappa_{i}}{\partial x_{j,0}}	=-z_{0},\ \frac{\partial\kappa_{i}}{\partial x_{j,1}}=-z_{1},\ \frac{\partial\kappa_{i}}{\partial x_{j,2}}=\frac{\partial\kappa_{i}}{\partial x_{j,3}}=0,
\label{eq:dkappai_dxj}
\end{equation}
and from~\eqref{eq:omega_i}
\begin{equation}
    \frac{\partial\omega_{i}}{\partial x_{i,0}}	=\frac{\partial\omega_{i}}{\partial x_{i,1}}=\frac{\partial\omega_{i}}{\partial x_{i,2}}=\frac{\partial\omega_{i}}{\partial x_{i,3}}=0,
\label{eq:domegai_dxi}
\end{equation}
and
\begin{equation}
    \frac{\partial\omega_{i}}{\partial x_{j,0}}	=-z_{1},\ \frac{\partial\omega_{i}}{\partial x_{j,1}}=z_{0},\ \frac{\partial\omega_{i}}{\partial x_{j,2}}=\frac{\partial\omega_{i}}{\partial x_{j,3}}=0.
\label{eq:domegai_dxj}
\end{equation}
Substituting equations~\eqref{eq:dr0_dxi}-\eqref{eq:dr1_dxj} and~\eqref{eq:dkappai_dxi}-\eqref{eq:domegai_dxj} into~\eqref{eq:dA12_dx_general} yields the following expressions for $\nicefrac{\partial\mathcal{A}_{12}}{\partial x_{l,m}}$.
\begin{align}
    \frac{\partial\mathcal{A}_{12}}{\partial x_{i,0}}=&\left(\kappa_{i}\left(\eta_{i}r_{0}-\mu_{i}r_{1}\right)+\left(\eta_{i}-2r_{1}\left(\mu_{i}r_{0}+\eta_{i}r_{1}\right)\right)\left(\kappa_{i}r_{0}-\omega_{i}r_{1}\right)-r_{1}\right)f_{1}  \\
	&+r_{1}\left(\kappa_{i}r_{0}-\omega_{i}r_{1}\right)\left(\eta_{i}r_{0}-\mu_{i}r_{1}\right)f_{2}, \label{eq:dA12_dxi0} \\
    \frac{\partial\mathcal{A}_{12}}{\partial x_{i,1}}=&~2\left(\kappa_{i}-r_{1}\left(\omega_{i}r_{0}+\kappa_{i}r_{1}\right)\right)\left(\kappa_{i}r_{0}-\omega_{i}r_{1}\right)f_{1}+r_{1}\left(\kappa_{i}r_{0}-\omega_{i}r_{1}\right)^{2}f_{2}, \label{eq:dA12_dxi1} \\
    \frac{\partial\mathcal{A}_{12}}{\partial x_{i,2}}=&~\frac{\partial\mathcal{A}_{12}}{\partial x_{i,3}}=0, \label{eq:dA12_dxi2_xi3} \\
    \frac{\partial\mathcal{A}_{12}}{\partial x_{j,0}}=&-\frac{z_{0}}{\gamma}+\left(\kappa_{i}\left(\eta_{j}r_{0}-\mu_{j}r_{1}\right)+\left(\eta_{j}-2r_{1}\left(\mu_{j}r_{0}+\eta_{j}r_{1}\right)\right)\left(\kappa_{i}r_{0}-\omega_{i}r_{1}\right)\right)f_{1}  \\
	&+r_{1}\left(\mu_{j}\kappa_{i}-\eta_{j}\omega_{i}-z_{0}r_{0}+z_{1}r_{1}\right)f_{1}  \\
	&+r_{1}\left(\kappa_{i}r_{0}-\omega_{i}r_{1}\right)\left(\eta_{j}r_{0}-\mu_{j}r_{1}\right)f_{2}, \label{eq:dA12_dxj0} \\
    \frac{\partial\mathcal{A}_{12}}{\partial x_{j,1}}=&-\frac{z_{1}}{\gamma}+\left(\kappa_{i}\left(\kappa_{j}r_{0}-\omega_{j}r_{1}\right)+\left(\kappa_{j}-2r_{1}\left(\omega_{j}r_{0}+\kappa_{j}r_{1}\right)\right)\left(\kappa_{i}r_{0}-\omega_{i}r_{1}\right)\right)f_{1}  \\
	&+r_{1}\left(\omega_{j}\kappa_{i}-\kappa_{j}\omega_{i}+-z_{1}r_{0}-z_{0}r_{1}\right)f_{1}  \\
	&+r_{1}\left(\kappa_{i}r_{0}-\omega_{i}r_{1}\right)\left(\kappa_{j}r_{0}-\omega_{j}r_{1}\right)f_{2}, \label{eq:dA12_dxj1} \\
    \frac{\partial\mathcal{A}_{12}}{\partial x_{j,2}}=&~\frac{\partial\mathcal{A}_{12}}{\partial x_{j,3}}=0, \label{eq:dA12_dxj2_xj3}
\end{align}
where we have used~\eqref{eq:omega_eta_minus_kappa_mu} to simplify~\eqref{eq:dA12_dxi0}. Because $\mathcal{A}_{13}=\mathcal{A}_{14}=0$, we have
\begin{equation}\label{eq:dA13_dx_all}
    \frac{\partial\mathcal{A}_{13}}{\partial x_{i,0}}=\frac{\partial\mathcal{A}_{13}}{\partial x_{i,1}}=\frac{\partial\mathcal{A}_{13}}{\partial x_{i,2}}=\frac{\partial\mathcal{A}_{13}}{\partial x_{i,3}}=\frac{\partial\mathcal{A}_{13}}{\partial x_{j,0}}=\frac{\partial\mathcal{A}_{13}}{\partial x_{j,1}}=\frac{\partial\mathcal{A}_{13}}{\partial x_{j,2}}=\frac{\partial\mathcal{A}_{13}}{\partial x_{j,3}}	=0,
\end{equation}
and
\begin{equation}\label{eq:dA14_dx_all}
    \frac{\partial\mathcal{A}_{14}}{\partial x_{i,0}}=\frac{\partial\mathcal{A}_{14}}{\partial x_{i,1}}=\frac{\partial\mathcal{A}_{14}}{\partial x_{i,2}}=\frac{\partial\mathcal{A}_{14}}{\partial x_{i,3}}	=\frac{\partial\mathcal{A}_{14}}{\partial x_{j,0}}=\frac{\partial\mathcal{A}_{14}}{\partial x_{j,1}}=\frac{\partial\mathcal{A}_{14}}{\partial x_{j,2}}=\frac{\partial\mathcal{A}_{14}}{\partial x_{j,3}}=0.
\end{equation}
Because $\mathcal{A}_{21}$ from~\eqref{eq:A_21} again follows the same general structure as $\mathcal{A}_{11}$, the general form for its derivatives is given by
\begin{align}
    \frac{\partial\mathcal{A}_{21}}{\partial x_{l,m}}=&~\frac{\frac{\partial\alpha_{1}}{\partial x_{l,m}}}{\gamma}+\left(\alpha_{1}\left(\frac{\partial r_{1}}{\partial x_{l,m}}r_{0}-\frac{\partial r_{0}}{\partial x_{l,m}}r_{1}\right)+\left(\frac{\partial r_{2}}{\partial x_{l,m}}-2r_{2}\left(\frac{\partial r_{0}}{\partial x_{l,m}}r_{0}+\frac{\partial r_{1}}{\partial x_{l,m}}r_{1}\right)\right)\left(\eta_{i}r_{0}-\mu_{i}r_{1}\right)\right)f_{1}  \\
	&+r_{2}\left(\frac{\partial r_{0}}{\partial x_{l,m}}\eta_{i}-\frac{\partial r_{1}}{\partial x_{l,m}}\mu_{i}+\frac{\partial\eta_{i}}{\partial x_{l,m}}r_{0}-\frac{\partial\mu_{i}}{\partial x_{l,m}}r_{1}\right)f_{1}  \\
	&+r_{2}\left(\eta_{i}r_{0}-\mu_{i}r_{1}\right)\left(\frac{\partial r_{1}}{\partial x_{l,m}}r_{0}-\frac{\partial r_{0}}{\partial x_{l,m}}r_{1}\right)f_{2}. \label{eq:dA21_dx_general}
\end{align}
From~\eqref{eq:alpha_1}, it is straightforward to compute
\begin{equation}
    \frac{\partial\alpha_{1}}{\partial x_{i,0}}	=\frac{\partial\alpha_{1}}{\partial x_{i,1}}=\frac{\partial\alpha_{1}}{\partial x_{i,2}}=\frac{\partial\alpha_{1}}{\partial x_{i,3}}=0,
\label{eq:dalpha_1_dxi}
\end{equation}
and
\begin{equation}
    \frac{\partial\alpha_{1}}{\partial x_{j,0}}	=-z_{2},\ \frac{\partial\alpha_{1}}{\partial x_{j,1}}=-z_{3},\ \frac{\partial\alpha_{1}}{\partial x_{j,2}}=z_{0},\ \frac{\partial\alpha_{1}}{\partial x_{j,3}}=z_{1}.
    \label{eq:dalpha_1_dxj}
\end{equation}
Substituting equations~\eqref{eq:dr0_dxi}-\eqref{eq:dr1_dxj},~\eqref{eq:dr2_dxi}-\eqref{eq:dr2_dxj}, and~\eqref{eq:dalpha_1_dxi}-\eqref{eq:dalpha_1_dxj} into~\eqref{eq:dA21_dx_general} yields the following expressions for $\nicefrac{\partial\mathcal{A}_{21}}{\partial x_{l,m}}$.
\begin{align}
    \frac{\partial\mathcal{A}_{21}}{\partial x_{i,0}}=&~2\left(\alpha_{1}-r_{2}\left(\mu_{i}r_{0}+\eta_{i}r_{1}\right)\right)\left(\eta_{i}r_{0}-\mu_{i}r_{1}\right)f_{1}+r_{2}\left(\eta_{i}r_{0}-\mu_{i}r_{1}\right)^{2}f_{2}, \label{eq:dA21_dxi0} \\
    \frac{\partial\mathcal{A}_{21}}{\partial x_{i,1}}=&	\left(\alpha_{1}\left(\kappa_{i}r_{0}-\omega_{i}r_{1}\right)+\left(\beta_{1}-2r_{2}\left(\omega_{i}r_{0}+\kappa_{i}r_{1}\right)\right)\left(\eta_{i}r_{0}-\mu_{i}r_{1}\right)+r_{2}\right)f_{1}  \\
    &+r_{2}\left(\eta_{i}r_{0}-\mu_{i}r_{1}\right)\left(\kappa_{i}r_{0}-\omega_{i}r_{1}\right)f_{2}, \label{eq:dA21_dxi1} \\
    \frac{\partial\mathcal{A}_{21}}{\partial x_{i,2}}=&~\xi_{1}\left(\eta_{i}r_{0}-\mu_{i}r_{1}\right)f_{1}, \label{eq:dA21_dxi2} \\
    \frac{\partial\mathcal{A}_{21}}{\partial x_{i,3}}=&~\zeta_{1}\left(\eta_{i}r_{0}-\mu_{i}r_{1}\right)f_{1}, \label{eq:dA21_dxi3} \\
    \frac{\partial\mathcal{A}_{21}}{\partial x_{j,0}}=&	-\frac{z_{2}}{\gamma}+\left(\alpha_{1}\left(\eta_{j}r_{0}-\mu_{j}r_{1}\right)+\left(\alpha_{3}-2r_{2}\left(\mu_{j}r_{0}+\eta_{j}r_{1}\right)\right)\left(\eta_{i}r_{0}-\mu_{i}r_{1}\right)\right)f_{1}  \\
    &+r_{2}\left(\mu_{j}\eta_{i}-\eta_{j}\mu_{i}-z_{1}r_{0}-z_{0}r_{1}\right)f_{1}  \\
    &+r_{2}\left(\eta_{i}r_{0}-\mu_{i}r_{1}\right)\left(\eta_{j}r_{0}-\mu_{j}r_{1}\right)f_{2}, \label{eq:dA21_dxj0} \\
    \frac{\partial\mathcal{A}_{21}}{\partial x_{j,1}}=&	-\frac{z_{3}}{\gamma}+\left(\alpha_{1}\left(\kappa_{j}r_{0}-\omega_{j}r_{1}\right)+\left(\beta_{3}-2r_{2}\left(\omega_{j}r_{0}+\kappa_{j}r_{1}\right)\right)\left(\eta_{i}r_{0}-\mu_{i}r_{1}\right)\right)f_{1}  \\
    &+r_{2}\left(\omega_{j}\eta_{i}-\kappa_{j}\mu_{i}+z_{0}r_{0}-z_{1}r_{1}\right)f_{1}  \\
    &+r_{2}\left(\eta_{i}r_{0}-\mu_{i}r_{1}\right)\left(\kappa_{j}r_{0}-\omega_{j}r_{1}\right)f_{2}, \label{eq:dA21_dxj1} \\
    \frac{\partial\mathcal{A}_{21}}{\partial x_{j,2}}=&~\frac{z_{0}}{\gamma}+\kappa_{j}\left(\eta_{i}r_{0}-\mu_{i}r_{1}\right)f_{1}, \label{eq:dA21_dxj2} \\
    \frac{\partial\mathcal{A}_{21}}{\partial x_{j,3}}=&~\frac{z_{1}}{\gamma}-\eta_{j}\left(\eta_{i}r_{0}-\mu_{i}r_{1}\right)f_{1}, \label{eq:dA21_dxj3} 
\end{align}
where we have used~\eqref{eq:omega_eta_minus_kappa_mu} to simplify~\eqref{eq:dA21_dxi1}. Because $\mathcal{A}_{22}$ from~\eqref{eq:A_22} again follows the same general structure as $\mathcal{A}_{11}$, the general form for its derivatives is given by
\begin{align}
    \frac{\partial\mathcal{A}_{22}}{\partial x_{l,m}}=&~\frac{\frac{\partial\beta_{1}}{\partial x_{l,m}}}{\gamma}+\left(\beta_{1}\left(\frac{\partial r_{1}}{\partial x_{l,m}}r_{0}-\frac{\partial r_{0}}{\partial x_{l,m}}r_{1}\right)+\left(\frac{\partial r_{2}}{\partial x_{l,m}}-2r_{2}\left(\frac{\partial r_{0}}{\partial x_{l,m}}r_{0}+\frac{\partial r_{1}}{\partial x_{l,m}}r_{1}\right)\right)\left(\kappa_{i}r_{0}-\omega_{i}r_{1}\right)\right)f_{1}  \\
	&+r_{2}\left(\kappa_{i}\frac{\partial r_{0}}{\partial x_{l,m}}-\omega_{i}\frac{\partial r_{1}}{\partial x_{l,m}}+\frac{\partial\kappa_{i}}{\partial x_{l,m}}r_{0}-\frac{\partial\omega_{i}}{\partial x_{l,m}}r_{1}\right)f_{1}  \\
	&+r_{2}\left(\kappa_{i}r_{0}-\omega_{i}r_{1}\right)\left(\frac{\partial r_{1}}{\partial x_{l,m}}r_{0}-\frac{\partial r_{0}}{\partial x_{l,m}}r_{1}\right)f_{2}. \label{dA22_dx_general}
\end{align}
From~\eqref{eq:beta_1}, it is straightforward to compute
\begin{equation}
    \frac{\partial\beta_{1}}{\partial x_{i,0}}	=\frac{\partial\beta_{1}}{\partial x_{i,1}}=\frac{\partial\beta_{1}}{\partial x_{i,2}}=\frac{\partial\beta_{1}}{\partial x_{i,3}}=0,
\label{eq:dbeta1_dxi}
\end{equation}
and
\begin{equation}
    \frac{\partial\beta_{1}}{\partial x_{j,0}}	=z_{3},\ \frac{\partial\beta_{1}}{\partial x_{j,1}}=-z_{2},\ \frac{\partial\beta_{1}}{\partial x_{j,2}}=-z_{1},\ \frac{\partial\beta_{1}}{\partial x_{j,3}}=z_{0}.
\label{eq:dbeta1_dxj}
\end{equation}
Substituting equations~\eqref{eq:dr0_dxi}-\eqref{eq:dr1_dxj},~\eqref{eq:dr2_dxi}-\eqref{eq:dr2_dxj}, and~\eqref{eq:dbeta1_dxi}-\eqref{eq:dbeta1_dxj} into~\eqref{eq:dA21_dx_general} yields the following expressions for $\nicefrac{\partial\mathcal{A}_{22}}{\partial x_{l,m}}$.
\begin{align}
    \frac{\partial\mathcal{A}_{22}}{\partial x_{i,0}}=&\left(\beta_{1}\left(\eta_{i}r_{0}-\mu_{i}r_{1}\right)+\left(\alpha_{1}-2r_{2}\left(\mu_{i}r_{0}+\eta_{i}r_{1}\right)\right)\left(\kappa_{i}r_{0}-\omega_{i}r_{1}\right)-r_{2}\right)f_{1}  \\
	&+r_{2}\left(\kappa_{i}r_{0}-\omega_{i}r_{1}\right)\left(\eta_{i}r_{0}-\mu_{i}r_{1}\right)f_{2}, \label{eq:dA22_dxi0} \\
    \frac{\partial\mathcal{A}_{22}}{\partial x_{i,1}}=&~2\left(\beta_{1}-r_{2}\left(\omega_{i}r_{0}+\kappa_{i}r_{1}\right)\right)\left(\kappa_{i}r_{0}-\omega_{i}r_{1}\right)f_{1}+r_{2}\left(\kappa_{i}r_{0}-\omega_{i}r_{1}\right)^{2}f_{2}, \label{eq:dA22_dxi1} \\
    \frac{\partial\mathcal{A}_{22}}{\partial x_{i,2}}=&~\xi_{1}\left(\kappa_{i}r_{0}-\omega_{i}r_{1}\right)f_{1}, \label{eq:dA22_dxi2} \\
    \frac{\partial\mathcal{A}_{22}}{\partial x_{i,3}}=&~\zeta_{1}\left(\kappa_{i}r_{0}-\omega_{i}r_{1}\right)f_{1}, \label{eq:dA22_dxi3} \\
    \frac{\partial\mathcal{A}_{22}}{\partial x_{j,0}}=&~\frac{z_{3}}{\gamma}+\left(\beta_{1}\left(\eta_{j}r_{0}-\mu_{j}r_{1}\right)+\left(\alpha_{3}-2r_{2}\left(\mu_{j}r_{0}+\eta_{j}r_{1}\right)\right)\left(\kappa_{i}r_{0}-\omega_{i}r_{1}\right)\right)f_{1}  \\
	&+r_{2}\left(\kappa_{i}\mu_{j}-\omega_{i}\eta_{j}-z_{0}r_{0}+z_{1}r_{1}\right)f_{1}  \\
	&+r_{2}\left(\kappa_{i}r_{0}-\omega_{i}r_{1}\right)\left(\eta_{j}r_{0}-\mu_{j}r_{1}\right)f_{2}, \label{eq:dA22_dxj0} \\
    \frac{\partial\mathcal{A}_{22}}{\partial x_{j,1}}=&-\frac{z_{2}}{\gamma}+\left(\beta_{1}\left(\kappa_{j}r_{0}-\omega_{j}r_{1}\right)+\left(\beta_{3}-2r_{2}\left(\omega_{j}r_{0}+\kappa_{j}r_{1}\right)\right)\left(\kappa_{i}r_{0}-\omega_{i}r_{1}\right)\right)f_{1}  \\
	&+r_{2}\left(\kappa_{i}\omega_{j}-\omega_{i}\kappa_{j}-z_{1}r_{0}-z_{0}r_{1}\right)f_{1}  \\
	&+r_{2}\left(\kappa_{i}r_{0}-\omega_{i}r_{1}\right)\left(\kappa_{j}r_{0}-\omega_{j}r_{1}\right)f_{2}, \label{eq:dA22_dxj1} \\
    \frac{\partial\mathcal{A}_{22}}{\partial x_{j,2}}=&-\frac{z_{1}}{\gamma}+\kappa_{j}\left(\kappa_{i}r_{0}-\omega_{i}r_{1}\right)f_{1}, \label{eq:dA22_dxj2} \\
    \frac{\partial\mathcal{A}_{22}}{\partial x_{j,3}}=&~\frac{z_{0}}{\gamma}-\eta_{j}\left(\kappa_{i}r_{0}-\omega_{i}r_{1}\right)f_{1}, \label{eq:dA22_dxj3} 
\end{align}
where we have again used~\eqref{eq:omega_eta_minus_kappa_mu} to simplify~\eqref{eq:dA22_dxi0}. To compute derivatives of $\mathcal{A}_{23}$, which is given by~\eqref{eq:A_23}, we follow the derivation of~\eqref{eq:dA11_left} to derive the general form
\begin{equation}\label{eq:dA23_dx_general}
    \frac{\partial\mathcal{A}_{23}}{\partial x_{l,m}}	=\frac{\frac{\partial\xi_{1}}{\partial x_{l,m}}}{\gamma}+\xi_{1}\left(\frac{\partial r_{1}}{\partial x_{l,m}}r_{0}-\frac{\partial r_{0}}{\partial x_{l,m}}r_{1}\right)f_{1}.
\end{equation}
From~\eqref{eq:xi_1}, we have
\begin{equation}\label{eq:dxi1_dxi}
    \frac{\partial\xi_{1}}{\partial x_{i,0}}=\frac{\partial\xi_{1}}{\partial x_{i,1}}=\frac{\partial\xi_{1}}{\partial x_{i,2}}=\frac{\partial\xi_{1}}{\partial x_{i,3}}=0,
\end{equation}
and
\begin{equation}\label{eq:dxi1_dxj}
    \frac{\partial\xi_{1}}{\partial x_{j,0}}=-z_{0},\ \frac{\partial\xi_{1}}{\partial x_{j,1}}=z_{1},\ \frac{\partial\xi_{1}}{\partial x_{j,2}}=\frac{\partial\xi_{1}}{\partial x_{j,3}}=0.
\end{equation}
Substituting equations~\eqref{eq:dr0_dxi}-\eqref{eq:dr1_dxj} and~\eqref{eq:dxi1_dxi}-\eqref{eq:dxi1_dxj} into~\eqref{eq:dA23_dx_general} yields the following expressions for $\nicefrac{\partial\mathcal{A}_{23}}{\partial x_{l,m}}$.
\begin{align}
    \frac{\partial\mathcal{A}_{23}}{\partial x_{i,0}}&=\xi_{1}\left(\eta_{i}r_{0}-\mu_{i}r_{1}\right)f_{1}, \label{eq:dA23_dxi0} \\
    \frac{\partial\mathcal{A}_{23}}{\partial x_{i,1}}&=\xi_{1}\left(\kappa_{i}r_{0}-\omega_{i}r_{1}\right)f_{1}, \label{eq:dA23_dxi1} \\
    \frac{\partial\mathcal{A}_{23}}{\partial x_{i,2}}&=\frac{\partial\mathcal{A}_{23}}{\partial x_{i,3}}=0, \label{eq:dA23_dxi2_xi3} \\
    \frac{\partial\mathcal{A}_{23}}{\partial x_{j,0}}&=-\frac{z_{0}}{\gamma}+\xi_{1}\left(\eta_{j}r_{0}-\mu_{j}r_{1}\right)f_{1}, \label{eq:dA23_dxj0} \\
    \frac{\partial\mathcal{A}_{23}}{\partial x_{j,1}}&=\frac{z_{1}}{\gamma}+\xi_{1}\left(\kappa_{j}r_{0}-\omega_{j}r_{1}\right)f_{1}, \label{eq:dA23_dxj1} \\
    \frac{\partial\mathcal{A}_{23}}{\partial x_{j,2}}&=\frac{\partial\mathcal{A}_{23}}{\partial x_{j,3}}=0, \label{eq:dA23_dxj2_xj3}
\end{align}
Similarly, the general form for derivatives of $\mathcal{A}_{24}$ from~\eqref{eq:A_24} is given by
\begin{equation}\label{eq:dA24_dx_general}
    \frac{\partial\mathcal{A}_{24}}{\partial x_{l,m}}=\frac{\frac{\partial\zeta_{1}}{\partial x_{l,m}}}{\gamma}+\zeta_{1}\left(\frac{\partial r_{1}}{\partial x_{l,m}}r_{0}-\frac{\partial r_{0}}{\partial x_{l,m}}r_{1}\right)f_{1}.
\end{equation}
From~\eqref{eq:zeta_1}, we have
\begin{equation}\label{eq:dzeta1_dxi}
    \frac{\partial\zeta_{1}}{\partial x_{i,0}}=\frac{\partial\zeta_{1}}{\partial x_{i,1}}=\frac{\partial\zeta_{1}}{\partial x_{i,2}}=\frac{\partial\zeta_{1}}{\partial x_{i,3}}=0,
\end{equation}
and
\begin{equation}\label{eq:dzeta1_dxj}
    \frac{\partial\zeta_{1}}{\partial x_{j,0}}=-z_{1},\ \frac{\partial\zeta_{1}}{\partial x_{j,1}}=-z_{0},\ \frac{\partial\zeta_{1}}{\partial x_{j,2}}=\frac{\partial\zeta_{1}}{\partial x_{j,3}}=0.
\end{equation}
Substituting equations~\eqref{eq:dr0_dxi}-\eqref{eq:dr1_dxj} and~\eqref{eq:dzeta1_dxi}-\eqref{eq:dzeta1_dxj} into~\eqref{eq:dA24_dx_general} yields the following expressions for $\nicefrac{\partial\mathcal{A}_{24}}{\partial x_{l,m}}$.
\begin{align}
    \frac{\partial\mathcal{A}_{24}}{\partial x_{i,0}}&=\zeta_{1}\left(\eta_{i}r_{0}-\mu_{i}r_{1}\right)f_{1}, \label{eq:dA24_dxi0} \\
    \frac{\partial\mathcal{A}_{24}}{\partial x_{i,1}}&=\zeta_{1}\left(\kappa_{i}r_{0}-\omega_{i}r_{1}\right)f_{1}, \label{eq:dA24_dxi1} \\
    \frac{\partial\mathcal{A}_{24}}{\partial x_{i,2}}&=\frac{\partial\mathcal{A}_{24}}{\partial x_{i,3}}=0, \label{eq:dA24_dxi2_xi3} \\
    \frac{\partial\mathcal{A}_{24}}{\partial x_{j,0}}&=-\frac{z_{1}}{\gamma}+\zeta_{1}\left(\eta_{j}r_{0}-\mu_{j}r_{1}\right)f_{1}, \label{eq:dA24_dxj0} \\
    \frac{\partial\mathcal{A}_{24}}{\partial x_{j,1}}&=-\frac{z_{0}}{\gamma}+\zeta_{1}\left(\kappa_{j}r_{0}-\omega_{j}r_{1}\right)f_{1}, \label{eq:dA24_dxj1} \\
    \frac{\partial\mathcal{A}_{24}}{\partial x_{j,2}}&=\frac{\partial\mathcal{A}_{24}}{\partial x_{j,3}}=0, \label{eq:dA24_dxj2_xj3}
\end{align}
Again following a similar derivation to~\eqref{eq:dA11_dx_general}, the general form for derivatives of $\mathcal{A}_{31}$ from~\eqref{eq:A_31} is derived to be
\begin{align}
    \frac{\partial\mathcal{A}_{31}}{\partial x_{l,m}}=&~\frac{\frac{\partial\alpha_{2}}{\partial x_{l,m}}}{\gamma}+\left(\alpha_{2}\left(\frac{\partial r_{1}}{\partial x_{l,m}}r_{0}-\frac{\partial r_{0}}{\partial x_{l,m}}r_{1}\right)+\left(\frac{\partial r_{3}}{\partial x_{l,m}}-2r_{3}\left(\frac{\partial r_{0}}{\partial x_{l,m}}r_{0}+\frac{\partial r_{1}}{\partial x_{l,m}}r_{1}\right)\right)\left(\eta_{i}r_{0}-\mu_{i}r_{1}\right)\right)f_{1}\\
	&+r_{3}\left(\eta_{i}\frac{\partial r_{0}}{\partial x_{l,m}}-\mu_{i}\frac{\partial r_{1}}{\partial x_{l,m}}+\frac{\partial\eta_{i}}{\partial x_{l,m}}r_{0}-\frac{\partial\mu_{i}}{\partial x_{l,m}}r_{1}\right)f_{1}  \\
	&+r_{3}\left(\eta_{i}r_{0}-\mu_{i}r_{1}\right)\left(\frac{\partial r_{1}}{\partial x_{l,m}}r_{0}-\frac{\partial r_{0}}{\partial x_{l,m}}r_{1}\right)f_{2}. \label{eq:dA31_dx_general}
\end{align}
From~\eqref{eq:alpha_3}, it is straightforward to compute
\begin{equation}\label{eq:dalpha3_dxi}
    \frac{\partial\alpha_{2}}{\partial x_{i,0}}	=\frac{\partial\alpha_{2}}{\partial x_{i,1}}=\frac{\partial\alpha_{2}}{\partial x_{i,2}}=\frac{\partial\alpha_{2}}{\partial x_{i,3}}=0,
\end{equation}
and
\begin{equation}\label{eq:dalpha3_dxj}
    \frac{\partial\alpha_{2}}{\partial x_{j,0}}	=-z_{3},\ \frac{\partial\alpha_{2}}{\partial x_{j,1}}=z_{2},\ \frac{\partial\alpha_{2}}{\partial x_{j,2}}=-z_{1},\ \frac{\partial\alpha_{2}}{\partial x_{j,3}}=z_{0}.
\end{equation}
Substituting equations~\eqref{eq:dr0_dxi}-\eqref{eq:dr1_dxj},~\eqref{eq:dr3_dxi}-\eqref{eq:dr3_dxj}, and~\eqref{eq:dalpha3_dxi}-\eqref{eq:dalpha3_dxj} into~\eqref{eq:dA31_dx_general} yields the following expressions for $\nicefrac{\partial\mathcal{A}_{22}}{\partial x_{l,m}}$.
\begin{align}
    \frac{\partial\mathcal{A}_{31}}{\partial x_{i,0}}=&~2\left(\alpha_{2}-r_{3}\left(\mu_{i}r_{0}+\eta_{i}r_{1}\right)\right)\left(\eta_{i}r_{0}-\mu_{i}r_{1}\right)f_{1}+r_{3}\left(\eta_{i}r_{0}-\mu_{i}r_{1}\right)^{2}f_{2}, \label{eq:dA31_dxi0} \\
    \frac{\partial\mathcal{A}_{31}}{\partial x_{i,1}}=&\left(\alpha_{2}\left(\kappa_{i}r_{0}-\omega_{i}r_{1}\right)+\left(\beta_{2}-2r_{3}\left(\omega_{i}r_{0}+\kappa_{i}r_{1}\right)\right)\left(\eta_{i}r_{0}-\mu_{i}r_{1}\right)+r_{3}\right)f_{1}  \\
	&+r_{3}\left(\eta_{i}r_{0}-\mu_{i}r_{1}\right)\left(\kappa_{i}r_{0}-\omega_{i}r_{1}\right)f_{2},\label{eq:dA31_dxi1} \\
    \frac{\partial\mathcal{A}_{31}}{\partial x_{i,2}}=&-\zeta_{1}\left(\eta_{i}r_{0}-\mu_{i}r_{1}\right)f_{1},\label{eq:dA31_dxi2} \\
    \frac{\partial\mathcal{A}_{31}}{\partial x_{i,3}}=&~\xi_{1}\left(\eta_{i}r_{0}-\mu_{i}r_{1}\right)f_{1},\label{eq:dA31_dxi3} \\
    \frac{\partial\mathcal{A}_{31}}{\partial x_{j,0}}=&-\frac{z_{3}}{\gamma}+\left(\alpha_{2}\left(\eta_{j}r_{0}-\mu_{j}r_{1}\right)+\left(\beta_{3}-2r_{3}\left(\mu_{j}r_{0}+\eta_{j}r_{1}\right)\right)\left(\eta_{i}r_{0}-\mu_{i}r_{1}\right)\right)f_{1}  \\
	&+r_{3}\left(\eta_{i}\mu_{j}-\mu_{i}\eta_{j}-z_{1}r_{0}-z_{0}r_{1}\right)f_{1}  \\
	&+r_{3}\left(\eta_{i}r_{0}-\mu_{i}r_{1}\right)\left(\eta_{j}r_{0}-\mu_{j}r_{1}\right)f_{2},\label{eq:dA31_dxj0} \\
    \frac{\partial\mathcal{A}_{31}}{\partial x_{j,1}}=&~\frac{z_{2}}{\gamma}+\left(\alpha_{2}\left(\kappa_{j}r_{0}-\omega_{j}r_{1}\right)-\left(\alpha_{3}+2r_{3}\left(\omega_{j}r_{0}+\kappa_{j}r_{1}\right)\right)\left(\eta_{i}r_{0}-\mu_{i}r_{1}\right)\right)f_{1}  \\
	&+r_{3}\left(\eta_{i}\omega_{j}-\mu_{i}\kappa_{j}+z_{0}r_{0}-z_{1}r_{1}\right)f_{1}  \\
	&+r_{3}\left(\eta_{i}r_{0}-\mu_{i}r_{1}\right)\left(\kappa_{j}r_{0}-\omega_{j}r_{1}\right)f_{2},\label{eq:dA31_dxj1} \\
    \frac{\partial\mathcal{A}_{31}}{\partial x_{j,2}}=&-\frac{z_{1}}{\gamma}+\eta_{j}\left(\eta_{i}r_{0}-\mu_{i}r_{1}\right)f_{1},\label{eq:dA31_dxj2} \\
    \frac{\partial\mathcal{A}_{31}}{\partial x_{j,3}}=&~\frac{z_{0}}{\gamma}+\kappa_{j}\left(\eta_{i}r_{0}-\mu_{i}r_{1}\right)f_{1},\label{eq:dA31_dxj3} 
\end{align}
where we have again used~\eqref{eq:omega_eta_minus_kappa_mu} to simplify~\eqref{eq:dA31_dxi1}. Again following a similar derivation to~\eqref{eq:dA11_dx_general}, the general form for derivatives of $\mathcal{A}_{32}$ from~\eqref{eq:A_32} is derived to be
\begin{align}
    \frac{\partial\mathcal{A}_{32}}{\partial x_{l,m}}=&~\frac{\frac{\partial\beta_{2}}{\partial x_{l,m}}}{\gamma}+\left(\beta_{2}\left(\frac{\partial r_{1}}{\partial x_{l,m}}r_{0}-\frac{\partial r_{0}}{\partial x_{l,m}}r_{1}\right)+\left(\frac{\partial r_{3}}{\partial x_{l,m}}-2r_{3}\left(\frac{\partial r_{0}}{\partial x_{l,m}}r_{0}+\frac{\partial r_{1}}{\partial x_{l,m}}r_{1}\right)\right)\left(\kappa_{i}r_{0}-\omega_{i}r_{1}\right)\right)f_{1}  \\
	&+r_{3}\left(\frac{\partial r_{0}}{\partial x_{l,m}}\kappa_{i}-\frac{\partial r_{1}}{\partial x_{l,m}}\omega_{i}+\frac{\partial\kappa_{i}}{\partial x_{l,m}}r_{0}-\frac{\partial\omega_{i}}{\partial x_{l,m}}r_{1}\right)f_{1}  \\
	&+r_{3}\left(\kappa_{i}r_{0}-\omega_{i}r_{1}\right)\left(\frac{\partial r_{1}}{\partial x_{l,m}}r_{0}-\frac{\partial r_{0}}{\partial x_{l,m}}r_{1}\right)f_{2} \label{eq:dA32_dx_general}
\end{align}
From~\eqref{eq:beta_3}, it is straightforward to compute
\begin{equation}
    \frac{\partial\beta_{2}}{\partial x_{i,0}}	=\frac{\partial\beta_{2}}{\partial x_{i,1}}=\frac{\partial\beta_{2}}{\partial x_{i,2}}=\frac{\partial\beta_{2}}{\partial x_{i,3}}=0,
    \label{eq:dbeta3_dxi}
\end{equation}
and
\begin{equation}
    \frac{\partial\beta_{2}}{\partial x_{j,0}}	=-z_{2},\ \frac{\partial\beta_{2}}{\partial x_{j,1}}=-z_{3},\ \frac{\partial\beta_{2}}{\partial x_{j,2}}=-z_{0},\ \frac{\partial\beta_{2}}{\partial x_{j,3}}=-z_{1}.
    \label{eq:dbeta3_dxj}
\end{equation}
Substituting equations~\eqref{eq:dr0_dxi}-\eqref{eq:dr1_dxj},~\eqref{eq:dr3_dxi}-\eqref{eq:dr3_dxj}, and~\eqref{eq:dbeta3_dxi}-\eqref{eq:dbeta3_dxj} into~\eqref{eq:dA32_dx_general} yields the following expressions for $\nicefrac{\partial\mathcal{A}_{32}}{\partial x_{l,m}}$.
\begin{align}
    \frac{\partial\mathcal{A}_{32}}{\partial x_{i,0}}=&\left(\beta_{2}\left(\eta_{i}r_{0}-\mu_{i}r_{1}\right)+\left(\alpha_{2}-2r_{3}\left(\mu_{i}r_{0}+\eta_{i}r_{1}\right)\right)\left(\kappa_{i}r_{0}-\omega_{i}r_{1}\right)-r_{3}\right)f_{1}  \\
	&+r_{3}\left(\kappa_{i}r_{0}-\omega_{i}r_{1}\right)\left(\eta_{i}r_{0}-\mu_{i}r_{1}\right)f_{2}, \label{eq:dA32_dxi0} \\
    \frac{\partial\mathcal{A}_{32}}{\partial x_{i,1}}=&2\left(\beta_{2}-r_{3}\left(\omega_{i}r_{0}+\kappa_{i}r_{1}\right)\right)\left(\kappa_{i}r_{0}-\omega_{i}r_{1}\right)f_{1}+r_{3}\left(\kappa_{i}r_{0}-\omega_{i}r_{1}\right)^{2}f_{2},\label{eq:dA32_dxi1} \\
    \frac{\partial\mathcal{A}_{32}}{\partial x_{i,2}}=&-\zeta_{1}\left(\kappa_{i}r_{0}-\omega_{i}r_{1}\right)f_{1},\label{eq:dA32_dxi2} \\
    \frac{\partial\mathcal{A}_{32}}{\partial x_{i,3}}=&\xi_{1}\left(\kappa_{i}r_{0}-\omega_{i}r_{1}\right)f_{1},\label{eq:dA32_dxi3} \\
    \frac{\partial\mathcal{A}_{32}}{\partial x_{j,0}}=&-\frac{z_{2}}{\gamma}+\left(\beta_{2}\left(\eta_{j}r_{0}-\mu_{j}r_{1}\right)+\left(\beta_{3}-2r_{3}\left(\mu_{j}r_{0}+\eta_{j}r_{1}\right)\right)\left(\kappa_{i}r_{0}-\omega_{i}r_{1}\right)\right)f_{1}  \\
	&+r_{3}\left(\kappa_{i}\mu_{j}-\omega_{i}\eta_{j}-z_{0}r_{0}+z_{1}r_{1}\right)f_{1}  \\
	&+r_{3}\left(\kappa_{i}r_{0}-\omega_{i}r_{1}\right)\left(\eta_{j}r_{0}-\mu_{j}r_{1}\right)f_{2},\label{eq:dA32_dxj0} \\
    \frac{\partial\mathcal{A}_{32}}{\partial x_{j,1}}=&-\frac{z_{3}}{\gamma}+\left(\beta_{2}\left(\kappa_{j}r_{0}-\omega_{j}r_{1}\right)-\left(\alpha_{3}+2r_{3}\left(\omega_{j}r_{0}+\kappa_{j}r_{1}\right)\right)\left(\kappa_{i}r_{0}-\omega_{i}r_{1}\right)\right)f_{1}  \\
	&+r_{3}\left(\kappa_{i}\omega_{j}-\omega_{i}\kappa_{j}-z_{1}r_{0}-z_{0}r_{1}\right)f_{1}  \\
	&+r_{3}\left(\kappa_{i}r_{0}-\omega_{i}r_{1}\right)\left(\kappa_{j}r_{0}-\omega_{j}r_{1}\right)f_{2},\label{eq:dA32_dxj1} \\
    \frac{\partial\mathcal{A}_{32}}{\partial x_{j,2}}=&-\frac{z_{0}}{\gamma}+\eta_{j}\left(\kappa_{i}r_{0}-\omega_{i}r_{1}\right)f_{1},\label{eq:dA32_dxj2} \\
    \frac{\partial\mathcal{A}_{32}}{\partial x_{j,3}}=&-\frac{z_{1}}{\gamma}+\kappa_{j}\left(\kappa_{i}r_{0}-\omega_{i}r_{1}\right)f_{1},\label{eq:dA32_dxj3}
\end{align}
where we have again used~\eqref{eq:omega_eta_minus_kappa_mu} to simplify~\eqref{eq:dA32_dxi0}. To compute derivatives of $\mathcal{A}_{33}$ from~\eqref{eq:A_33}, we again follow the derivation of~\eqref{eq:dA11_left} to derive the general form
\begin{equation}
    \frac{\partial\mathcal{A}_{33}}{\partial x_{l,m}}	=-\frac{\frac{\partial\zeta_{1}}{\partial x_{l,m}}}{\gamma}-\zeta_{1}\left(\frac{\partial r_{1}}{\partial x_{l,m}}r_{0}-\frac{\partial r_{0}}{\partial x_{l,m}}r_{1}\right)f_{1}
\label{eq:dA33_dx_general}
\end{equation}
Substituting equations~\eqref{eq:dr0_dxi}-\eqref{eq:dr1_dxj} and~\eqref{eq:dzeta1_dxi}-\eqref{eq:dzeta1_dxj} into~\eqref{eq:dA33_dx_general} yields the following expressions for $\nicefrac{\partial\mathcal{A}_{33}}{\partial x_{l,m}}$.
\begin{align}
    \frac{\partial\mathcal{A}_{33}}{\partial x_{i,0}}&=-\zeta_{1}\left(\eta_{i}r_{0}-\mu_{i}r_{1}\right)f_{1}, \label{eq:dA33_dxi0} \\
    \frac{\partial\mathcal{A}_{33}}{\partial x_{i,1}}&=-\zeta_{1}\left(\kappa_{i}r_{0}-\omega_{i}r_{1}\right)f_{1}, \label{eq:dA33_dxi1} \\
    \frac{\partial\mathcal{A}_{33}}{\partial x_{i,2}}&=\frac{\partial\mathcal{A}_{33}}{\partial x_{i,3}}=0, \label{eq:dA33_dxi2_xi3} \\
    \frac{\partial\mathcal{A}_{33}}{\partial x_{j,0}}&=\frac{z_{1}}{\gamma}-\zeta_{1}\left(\eta_{j}r_{0}-\mu_{j}r_{1}\right)f_{1}, \label{eq:dA33_dxj0} \\
    \frac{\partial\mathcal{A}_{33}}{\partial x_{j,1}}&=\frac{z_{0}}{\gamma}-\zeta_{1}\left(\kappa_{j}r_{0}-\omega_{j}r_{1}\right)f_{1}, \label{eq:dA33_dxj1} \\
    \frac{\partial\mathcal{A}_{33}}{\partial x_{j,2}}&=\frac{\partial\mathcal{A}_{33}}{\partial x_{j,3}}=0. \label{eq:dA33_dxj2_xj3}
\end{align}
Similarly, the general form for derivatives of $\mathcal{A}_{34}$ from~\eqref{eq:A_34} is given by
\begin{equation}
    \frac{\partial\mathcal{A}_{34}}{\partial x_{l,m}}	=\frac{\frac{\partial\xi_{1}}{\partial x_{l,m}}}{\gamma}+\xi_{1}\left(\frac{\partial r_{1}}{\partial x_{l,m}}r_{0}-\frac{\partial r_{0}}{\partial x_{l,m}}r_{1}\right)f_{1}
\label{eq:dA34_dx_general}
\end{equation}
Substituting equations~\eqref{eq:dr0_dxi}-\eqref{eq:dr1_dxj} and~\eqref{eq:dxi1_dxi}-\eqref{eq:dxi1_dxj} into~\eqref{eq:dA34_dx_general} yields the following expressions for $\nicefrac{\partial\mathcal{A}_{34}}{\partial x_{l,m}}$.
\begin{align}
    \frac{\partial\mathcal{A}_{34}}{\partial x_{i,0}}&=\xi_{1}\left(\eta_{i}r_{0}-\mu_{i}r_{1}\right)f_{1}, \label{eq:dA34_dxi0} \\
    \frac{\partial\mathcal{A}_{34}}{\partial x_{i,1}}&=\xi_{1}\left(\kappa_{i}r_{0}-\omega_{i}r_{1}\right)f_{1}, \label{eq:dA34_dxi1} \\
    \frac{\partial\mathcal{A}_{34}}{\partial x_{i,2}}&=\frac{\partial\mathcal{A}_{34}}{\partial x_{i,3}}=0, \label{eq:dA34_dxi2_xi3} \\
    \frac{\partial\mathcal{A}_{34}}{\partial x_{j,0}}&=-\frac{z_{0}}{\gamma}+\xi_{1}\left(\eta_{j}r_{0}-\mu_{j}r_{1}\right)f_{1}, \label{eq:dA34_dxj0} \\
    \frac{\partial\mathcal{A}_{34}}{\partial x_{j,1}}&=\frac{z_{1}}{\gamma}+\xi_{1}\left(\kappa_{j}r_{0}-\omega_{j}r_{1}\right)f_{1}, \label{eq:dA34_dxj1} \\
    \frac{\partial\mathcal{A}_{34}}{\partial x_{j,2}}&=\frac{\partial\mathcal{A}_{34}}{\partial x_{j,3}}=0. \label{eq:dA34_dxj2_xj3}
\end{align}
\subsection{Partial Derivatives of $\mathcal{B}_{ij}$}\label{app:dBij_dx}
Partial derivatives of $\mathcal{B}_{ij}$ with respect to $x\in\mathbf{x}_{i},\mathbf{x}_{j}$ are computed in a similar manner. For example, following the derivation from equations~\eqref{eq:dA11_sep}-\eqref{eq:dA11_dx_general} with respect to the structure of $\mathcal{B}_{11}$ from~\eqref{eq:B_11}, the general form for its derivatives is given by
\begin{align}
    \frac{\partial\mathcal{B}_{11}}{\partial x_{l,m}}=&~\frac{\frac{\partial\eta_{j}}{\partial x_{l,m}}}{\gamma}+\left(\eta_{j}\left(\frac{\partial r_{1}}{\partial x_{l,m}}r_{0}-\frac{\partial r_{0}}{\partial x_{l,m}}r_{1}\right)+\left(\frac{\partial r_{1}}{\partial x_{l,m}}-2r_{1}\left(\frac{\partial r_{0}}{\partial x_{l,m}}r_{0}+\frac{\partial r_{1}}{\partial x_{l,m}}r_{1}\right)\right)\left(\eta_{j}r_{0}-\mu_{j}r_{1}\right)\right)f_{1}  \\
	&+r_{1}\left(\frac{\partial r_{0}}{\partial x_{l,m}}\eta_{j}-\frac{\partial r_{1}}{\partial x_{l,m}}\mu_{j}+\frac{\partial\eta_{j}}{\partial x_{l,m}}r_{0}-\frac{\partial\mu_{j}}{\partial x_{l,m}}r_{1}\right)f_{1}+r_{1}\left(\eta_{j}r_{0}-\mu_{j}r_{1}\right)\left(\frac{\partial r_{1}}{\partial x_{l,m}}r_{0}-\frac{\partial r_{0}}{\partial x_{l,m}}r_{1}\right)f_{2}. \label{eq:dB11_dx_general}
\end{align}
From~\eqref{eq:eta_j}, it is straightforward to compute
\begin{equation}
    \frac{\partial\eta_{j}}{\partial x_{i,0}}=-z_{1},\ \frac{\partial\eta_{j}}{\partial x_{i,1}}=-z_{0},\ \frac{\partial\eta_{j}}{\partial x_{i,2}}=\frac{\partial\eta_{j}}{\partial x_{i,3}}=0,
    \label{eq:detaj_dxi}
\end{equation}
and
\begin{equation}
    \frac{\partial\eta_{j}}{\partial x_{j,0}}=\frac{\partial\eta_{j}}{\partial x_{j,1}}=\frac{\partial\eta_{j}}{\partial x_{j,2}}=\frac{\partial\eta_{j}}{\partial x_{j,3}}=0,
    \label{eq:detaj_dxj}
\end{equation}
and from~\eqref{eq:mu_j}, we have
\begin{equation}
    \frac{\partial\mu_{j}}{\partial x_{i,0}}=z_{0},\ \frac{\partial\mu_{j}}{\partial x_{i,1}}=-z_{1},\ \frac{\partial\mu_{j}}{\partial x_{i,2}}=\frac{\partial\mu_{j}}{\partial x_{i,3}}=0,
    \label{eq:dmuj_dxi}
\end{equation}
and
\begin{equation}
    \frac{\partial\mu_{j}}{\partial x_{j,0}}=\frac{\partial\mu_{j}}{\partial x_{j,1}}=\frac{\partial\mu_{j}}{\partial x_{j,2}}=\frac{\partial\mu_{j}}{\partial x_{j,3}}=0,
    \label{eq:dmuj_dxj}
\end{equation}
Substituting equations~\eqref{eq:dr0_dxi}-\eqref{eq:dr1_dxj} and~\eqref{eq:detaj_dxi}-\eqref{eq:dmuj_dxj} into~\eqref{eq:dB11_dx_general} yields the following expressions for $\nicefrac{\partial\mathcal{B}_{11}}{\partial x_{l,m}}$.
\begin{align}
    \frac{\partial\mathcal{B}_{11}}{\partial x_{i,0}}=&-\frac{z_{1}}{\gamma}+\left(\eta_{j}\left(\eta_{i}r_{0}-\mu_{i}r_{1}\right)+\left(\eta_{i}-2r_{1}\left(\mu_{i}r_{0}+\eta_{i}r_{1}\right)\right)\left(\eta_{j}r_{0}-\mu_{j}r_{1}\right)\right)f_{1}  \\
	&+r_{1}\left(\mu_{i}\eta_{j}-\eta_{i}\mu_{j}-z_{1}r_{0}-z_{0}r_{1}\right)f_{1}+r_{1}\left(\eta_{j}r_{0}-\mu_{j}r_{1}\right)\left(\eta_{i}r_{0}-\mu_{i}r_{1}\right)f_{2}, \label{eq:dB11_dxi0} \\
    \frac{\partial\mathcal{B}_{11}}{\partial x_{i,1}}=&-\frac{z_{0}}{\gamma}+\left(\eta_{j}\left(\kappa_{i}r_{0}-\omega_{i}r_{1}\right)+\left(\kappa_{i}-2r_{1}\left(\omega_{i}r_{0}+\kappa_{i}r_{1}\right)\right)\left(\eta_{j}r_{0}-\mu_{j}r_{1}\right)\right)f_{1}  \\
	&+r_{1}\left(\omega_{i}\eta_{j}-\kappa_{i}\mu_{j}-z_{0}r_{0}+z_{1}r_{1}\right)f_{1}+r_{1}\left(\eta_{j}r_{0}-\mu_{j}r_{1}\right)\left(\kappa_{i}r_{0}-\omega_{i}r_{1}\right)f_{2}, \label{eq:dB11_dxi1} \\
    \frac{\partial\mathcal{B}_{11}}{\partial x_{i,2}}=&~\frac{\partial\mathcal{B}_{11}}{\partial x_{i,3}}=0, \label{eq:dB11_dxi2_xi3} \\
    \frac{\partial\mathcal{B}_{11}}{\partial x_{j,0}}=&~2\left(\eta_{j}-r_{1}\left(\mu_{j}r_{0}+\eta_{j}r_{1}\right)\right)\left(\eta_{j}r_{0}-\mu_{j}r_{1}\right)f_{1}+r_{1}\left(\eta_{j}r_{0}-\mu_{j}r_{1}\right)^{2}f_{2}, \label{eq:dB11_dxj0} \\
    \frac{\partial\mathcal{B}_{11}}{\partial x_{j,1}}=&\left(\eta_{j}\left(\kappa_{j}r_{0}-\omega_{j}r_{1}\right)+\left(\kappa_{j}-2r_{1}\left(\omega_{j}r_{0}+\kappa_{j}r_{1}\right)\right)\left(\eta_{j}r_{0}-\mu_{j}r_{1}\right)-r_{1}\right)f_{1}  \\
	&+r_{1}\left(\eta_{j}r_{0}-\mu_{j}r_{1}\right)\left(\kappa_{j}r_{0}-\omega_{j}r_{1}\right)f_{2}, \label{eq:dB11_dxj1} \\
    \frac{\partial\mathcal{B}_{11}}{\partial x_{j,2}}=&~\frac{\partial\mathcal{B}_{11}}{\partial x_{j,3}}=0, \label{eq:dB11_dxj2_xj3}
\end{align}
where we have used the fact that
\begin{equation}\label{eq:muj_kappaj_min_etaj_muj}
    \mu_{j}\kappa_{j}-\eta_{j}\omega_{j}=\sin^{2}\left(\phi_{i}+\phi_{z}\right)+\cos^{2}\left(\phi_{i}+\phi_{z}\right)=1
\end{equation}
to simplify~\eqref{eq:dB11_dxj1}. Furthermore, since $\mathcal{B}_{12}$ from~\eqref{eq:B_12} has identical structure to $\mathcal{B}_{11}$, the general form for its partial derivatives is computed as
\begin{align}
    \frac{\partial\mathcal{B}_{12}}{\partial x_{l,m}}=&\frac{\frac{\partial\kappa_{j}}{\partial x_{l,m}}}{\gamma}+\left(\kappa_{j}\left(\frac{\partial r_{1}}{\partial x_{l,m}}r_{0}-\frac{\partial r_{0}}{\partial x_{l,m}}r_{1}\right)+\left(\frac{\partial r_{1}}{\partial x_{l,m}}-2r_{1}\left(\frac{\partial r_{0}}{\partial x_{l,m}}r_{0}+\frac{\partial r_{1}}{\partial x_{l,m}}r_{1}\right)\right)\left(\kappa_{j}r_{0}-\omega_{j}r_{1}\right)\right)f_{1}  \\
	&+r_{1}\left(\frac{\partial r_{0}}{\partial x_{l,m}}\kappa_{j}-\frac{\partial r_{1}}{\partial x_{l,m}}\omega_{j}+\frac{\partial\kappa_{j}}{\partial x_{l,m}}r_{0}-\frac{\partial\omega_{j}}{\partial x_{l,m}}r_{1}\right)f_{1}+r_{1}\left(\kappa_{j}r_{0}-\omega_{j}r_{1}\right)\left(\frac{\partial r_{1}}{\partial x_{l,m}}r_{0}-\frac{\partial r_{0}}{\partial x_{l,m}}r_{1}\right)f_{2}. \label{eq:dB12_dx_general}
\end{align}
From~\eqref{eq:kappa_j}, it is straightforward to compute
\begin{equation}
    \frac{\partial\kappa_{j}}{\partial x_{i,0}}=z_{0},\ \frac{\partial\kappa_{j}}{\partial x_{i,1}}=-z_{1},\ \frac{\partial\kappa_{j}}{\partial x_{i,2}}=\frac{\partial\kappa_{j}}{\partial x_{i,3}}=0,
\label{eq:dkappaj_dxi}
\end{equation}
and
\begin{equation}
    \frac{\partial\kappa_{j}}{\partial x_{j,0}}=\frac{\partial\kappa_{j}}{\partial x_{j,1}}=\frac{\partial\kappa_{j}}{\partial x_{j,2}}=\frac{\partial\kappa_{j}}{\partial x_{j,3}}=0,
    \label{eq:dkappaj_dxj}
\end{equation}
and from~\eqref{eq:omega_j}, we have
\begin{equation}
    \frac{\partial\omega_{j}}{\partial x_{i,0}}=z_{1},\ \frac{\partial\omega_{j}}{\partial x_{i,1}}=z_{0},\ \frac{\partial\omega_{j}}{\partial x_{i,2}}=\frac{\partial\omega_{j}}{\partial x_{i,3}}=0,
    \label{eq:domegaj_dxi}
\end{equation}
and
\begin{equation}
    \frac{\partial\omega_{j}}{\partial x_{j,0}}=\frac{\partial\omega_{j}}{\partial x_{j,1}}=\frac{\partial\omega_{j}}{\partial x_{j,2}}=\frac{\partial\omega_{j}}{\partial x_{j,3}}=0,
    \label{eq:domegaj_dxj}
\end{equation}
Substituting equations~\eqref{eq:dr0_dxi}-\eqref{eq:dr1_dxj} and~\eqref{eq:dkappaj_dxi}-\eqref{eq:domegaj_dxj} into~\eqref{eq:dB12_dx_general} yields the following expressions for $\nicefrac{\partial\mathcal{B}_{12}}{\partial x_{l,m}}$.
\begin{align}
    \frac{\partial\mathcal{B}_{12}}{\partial x_{i,0}}=&~\frac{z_{0}}{\gamma}+\left(\kappa_{j}\left(\eta_{i}r_{0}-\mu_{i}r_{1}\right)+\left(\eta_{i}-2r_{1}\left(\mu_{i}r_{0}+\eta_{i}r_{1}\right)\right)\left(\kappa_{j}r_{0}-\omega_{j}r_{1}\right)\right)f_{1}  \\
	&+r_{1}\left(\mu_{i}\kappa_{j}-\eta_{i}\omega_{j}+z_{0}r_{0}-z_{1}r_{1}\right)f_{1}+r_{1}\left(\kappa_{j}r_{0}-\omega_{j}r_{1}\right)\left(\eta_{i}r_{0}-\mu_{i}r_{1}\right)f_{2}, \label{eq:dB12_dxi0} \\
    \frac{\partial\mathcal{B}_{12}}{\partial x_{i,1}}=&-\frac{z_{1}}{\gamma}+\left(\kappa_{j}\left(\kappa_{i}r_{0}-\omega_{i}r_{1}\right)+\left(\kappa_{i}-2r_{1}\left(\omega_{i}r_{0}+\kappa_{i}r_{1}\right)\right)\left(\kappa_{j}r_{0}-\omega_{j}r_{1}\right)\right)f_{1}  \\
	&+r_{1}\left(\omega_{i}\kappa_{j}-\kappa_{i}\omega_{j}-z_{1}r_{0}-z_{0}r_{1}\right)f_{1}+r_{1}\left(\kappa_{j}r_{0}-\omega_{j}r_{1}\right)\left(\kappa_{i}r_{0}-\omega_{i}r_{1}\right)f_{2}, \label{eq:dB12_dxi1} \\
    \frac{\partial\mathcal{B}_{12}}{\partial x_{i,2}}=&~\frac{\partial\mathcal{B}_{12}}{\partial x_{i,3}}=0,\label{eq:dB12_dxi2_xi3} \\
    \frac{\partial\mathcal{B}_{12}}{\partial x_{j,0}}=&\left(\kappa_{j}\left(\eta_{j}r_{0}-\mu_{j}r_{1}\right)+\left(\eta_{j}-2r_{1}\left(\mu_{j}r_{0}+\eta_{j}r_{1}\right)\right)\left(\kappa_{j}r_{0}-\omega_{j}r_{1}\right)+r_{1}\right)f_{1}  \\
	&+r_{1}\left(\kappa_{j}r_{0}-\omega_{j}r_{1}\right)\left(\eta_{j}r_{0}-\mu_{j}r_{1}\right)f_{2}, \label{eq:dB12_dxj0} \\
    \frac{\partial\mathcal{B}_{12}}{\partial x_{j,1}}=&~2\left(\kappa_{j}-r_{1}\left(\omega_{j}r_{0}+\kappa_{j}r_{1}\right)\right)\left(\kappa_{j}r_{0}-\omega_{j}r_{1}\right)f_{1}+r_{1}\left(\kappa_{j}r_{0}-\omega_{j}r_{1}\right)^{2}f_{2},\label{eq:dB12_dxj1} \\
    \frac{\partial\mathcal{B}_{12}}{\partial x_{j,2}}=&~\frac{\partial\mathcal{B}_{12}}{\partial x_{j,3}}=0,\label{eq:dB12_dxj2_xj3}
\end{align}
where we have again used~\eqref{eq:muj_kappaj_min_etaj_muj} to simplify~\eqref{eq:dB12_dxj0}. Because $\mathcal{B}_{13}=\mathcal{B}_{14}=0$, we have
\begin{equation}\label{eq:dB13_dx_all}
    \frac{\partial\mathcal{B}_{13}}{\partial x_{i,0}}=\frac{\partial\mathcal{B}_{13}}{\partial x_{i,1}}=\frac{\partial\mathcal{B}_{13}}{\partial x_{i,2}}=\frac{\partial\mathcal{B}_{13}}{\partial x_{i,3}}=\frac{\partial\mathcal{B}_{13}}{\partial x_{j,0}}=\frac{\partial\mathcal{B}_{13}}{\partial x_{j,1}}=\frac{\partial\mathcal{B}_{13}}{\partial x_{j,2}}=\frac{\partial\mathcal{B}_{13}}{\partial x_{j,3}}	=0,
\end{equation}
and
\begin{equation}\label{eq:dB14_dx_all}
    \frac{\partial\mathcal{B}_{14}}{\partial x_{i,0}}=\frac{\partial\mathcal{B}_{14}}{\partial x_{i,1}}=\frac{\partial\mathcal{B}_{14}}{\partial x_{i,2}}=\frac{\partial\mathcal{B}_{14}}{\partial x_{i,3}}	=\frac{\partial\mathcal{B}_{14}}{\partial x_{j,0}}=\frac{\partial\mathcal{B}_{14}}{\partial x_{j,1}}=\frac{\partial\mathcal{B}_{14}}{\partial x_{j,2}}=\frac{\partial\mathcal{B}_{14}}{\partial x_{j,3}}=0.
\end{equation}
Because $\mathcal{B}_{21}$ from~\eqref{eq:B_21} again follows the same general structure as $\mathcal{B}_{11}$, the general form for its derivatives is given by
\begin{align}
    \frac{\partial\mathcal{B}_{21}}{\partial x_{l,m}}=&\frac{\frac{\partial\alpha_{3}}{\partial x_{l,m}}}{\gamma}+\left(\alpha_{3}\left(\frac{\partial r_{1}}{\partial x_{l,m}}r_{0}-\frac{\partial r_{0}}{\partial x_{l,m}}r_{1}\right)+\left(\frac{\partial r_{2}}{\partial x_{l,m}}-2r_{2}\left(\frac{\partial r_{0}}{\partial x_{l,m}}r_{0}+\frac{\partial r_{1}}{\partial x_{l,m}}r_{1}\right)\right)\left(\eta_{j}r_{0}-\mu_{j}r_{1}\right)\right)f_{1}  \\
	&+r_{2}\left(\frac{\partial r_{0}}{\partial x_{l,m}}\eta_{j}-\frac{\partial r_{1}}{\partial x_{l,m}}\mu_{j}+\frac{\partial\eta_{j}}{\partial x_{l,m}}r_{0}-\frac{\partial\mu_{j}}{\partial x_{l,m}}r_{1}\right)f_{1}+r_{2}\left(\eta_{j}r_{0}-\mu_{j}r_{1}\right)\left(\frac{\partial r_{1}}{\partial x_{l,m}}r_{0}-\frac{\partial r_{0}}{\partial x_{l,m}}r_{1}\right)f_{2}. \label{eq:dB21_dx_general}
\end{align}
From~\eqref{eq:alpha_2}, it is straightforward to compute
\begin{equation}
    \frac{\partial\alpha_{3}}{\partial x_{i,0}}	=-z_{2},\ \frac{\partial\alpha_{3}}{\partial x_{i,1}}=z_{3},\frac{\partial\alpha_{3}}{\partial x_{i,2}}=-z_{0},\ \frac{\partial\alpha_{3}}{\partial x_{i,3}}=-z_{1},
    \label{eq:dalpha2_dxi}
\end{equation}
and
\begin{equation}
    \frac{\partial\alpha_{3}}{\partial x_{j,0}}	=\frac{\partial\alpha_{3}}{\partial x_{j,1}}=\frac{\partial\alpha_{3}}{\partial x_{j,2}}=\frac{\partial\alpha_{3}}{\partial x_{j,3}}=0.
    \label{eq:dalpha2_dxj}
\end{equation}
Substituting equations~\eqref{eq:dr0_dxi}-\eqref{eq:dr1_dxj},~\eqref{eq:dr2_dxi}-\eqref{eq:dr2_dxj}, and~\eqref{eq:dalpha2_dxi}-\eqref{eq:dalpha2_dxj} into~\eqref{eq:dB21_dx_general} yields the following expressions for $\nicefrac{\partial\mathcal{B}_{21}}{\partial x_{l,m}}$.
\begin{align}
    \frac{\partial\mathcal{B}_{21}}{\partial x_{i,0}}=&-\frac{z_{2}}{\gamma}+\left(\alpha_{3}\left(\eta_{i}r_{0}-\mu_{i}r_{1}\right)+\left(\alpha_{1}-2r_{2}\left(\mu_{i}r_{0}+\eta_{i}r_{1}\right)\right)\left(\eta_{j}r_{0}-\mu_{j}r_{1}\right)\right)f_{1}  \\
	&+r_{2}\left(\mu_{i}\eta_{j}-\eta_{i}\mu_{j}-z_{1}r_{0}-z_{0}r_{1}\right)f_{1}+r_{2}\left(\eta_{j}r_{0}-\mu_{j}r_{1}\right)\left(\eta_{i}r_{0}-\mu_{i}r_{1}\right)f_{2}, \label{eq:dB21_dxi0} \\
    \frac{\partial\mathcal{B}_{21}}{\partial x_{i,1}}=&~\frac{z_{3}}{\gamma}+\left(\alpha_{3}\left(\kappa_{i}r_{0}-\omega_{i}r_{1}\right)+\left(\beta_{1}-2r_{2}\left(\omega_{i}r_{0}+\kappa_{i}r_{1}\right)\right)\left(\eta_{j}r_{0}-\mu_{j}r_{1}\right)\right)f_{1}  \\
	&+r_{2}\left(\omega_{i}\eta_{j}-\kappa_{i}\mu_{j}-z_{0}r_{0}+z_{1}r_{1}\right)f_{1}+r_{2}\left(\eta_{j}r_{0}-\mu_{j}r_{1}\right)\left(\kappa_{i}r_{0}-\omega_{i}r_{1}\right)f_{2}, \label{eq:dB21_dxi1} \\
    \frac{\partial\mathcal{B}_{21}}{\partial x_{i,2}}=&-\frac{z_{0}}{\gamma}+\xi_{1}\left(\eta_{j}r_{0}-\mu_{j}r_{1}\right)f_{1}, \label{eq:dB21_dxi2} \\
    \frac{\partial\mathcal{B}_{21}}{\partial x_{i,3}}=&-\frac{z_{1}}{\gamma}+\zeta_{1}\left(\eta_{j}r_{0}-\mu_{j}r_{1}\right)f_{1}, \label{eq:dB21_dxi3} \\
    \frac{\partial\mathcal{B}_{21}}{\partial x_{j,0}}=&~2\left(\alpha_{3}-r_{2}\left(\mu_{j}r_{0}+\eta_{j}r_{1}\right)\right)\left(\eta_{j}r_{0}-\mu_{j}r_{1}\right)f_{1}+r_{2}\left(\eta_{j}r_{0}-\mu_{j}r_{1}\right)^{2}f_{2}, \label{eq:dB21_dxj0} \\
    \frac{\partial\mathcal{B}_{21}}{\partial x_{j,1}}=&\left(\alpha_{3}\left(\kappa_{j}r_{0}-\omega_{j}r_{1}\right)+\left(\beta_{3}-2r_{2}\left(\omega_{j}r_{0}+\kappa_{j}r_{1}\right)\right)\left(\eta_{j}r_{0}-\mu_{j}r_{1}\right)-r_{2}\right)f_{1}  \\
	&+r_{2}\left(\eta_{j}r_{0}-\mu_{j}r_{1}\right)\left(\kappa_{j}r_{0}-\omega_{j}r_{1}\right)f_{2}, \label{eq:dB21_dxj1} \\
    \frac{\partial\mathcal{B}_{21}}{\partial x_{j,2}}=&~\kappa_{j}\left(\eta_{j}r_{0}-\mu_{j}r_{1}\right)f_{1}, \label{eq:dB21_dxj2} \\
    \frac{\partial\mathcal{B}_{21}}{\partial x_{j,3}}=&-\eta_{j}\left(\eta_{j}r_{0}-\mu_{j}r_{1}\right)f_{1}, \label{eq:dB21_dxj3}
\end{align}
where we have again used~\eqref{eq:muj_kappaj_min_etaj_muj} to simplify~\eqref{eq:dB21_dxj1}.
Because $\mathcal{B}_{22}$ from~\eqref{eq:B_22} again follows the same general structure as $\mathcal{B}_{11}$, the general form for its derivatives is given by
\begin{align}
    \frac{\partial\mathcal{B}_{22}}{\partial x_{l,m}}=&~\frac{\frac{\partial\beta_{3}}{\partial x_{l,m}}}{\gamma}+\left(\beta_{3}\left(\frac{\partial r_{1}}{\partial x_{l,m}}r_{0}-\frac{\partial r_{0}}{\partial x_{l,m}}r_{1}\right)+\left(\frac{\partial r_{2}}{\partial x_{l,m}}-2r_{2}\left(\frac{\partial r_{0}}{\partial x_{l,m}}r_{0}+\frac{\partial r_{1}}{\partial x_{l,m}}r_{1}\right)\right)\left(\kappa_{j}r_{0}-\omega_{j}r_{1}\right)\right)f_{1}  \\
	&+r_{2}\left(\frac{\partial r_{0}}{\partial x_{l,m}}\kappa_{j}-\frac{\partial r_{1}}{\partial x_{l,m}}\omega_{j}+\frac{\partial\kappa_{j}}{\partial x_{l,m}}r_{0}-\frac{\partial\omega_{j}}{\partial x_{l,m}}r_{1}\right)f_{1}+r_{2}\left(\kappa_{j}r_{0}-\omega_{j}r_{1}\right)\left(\frac{\partial r_{1}}{\partial x_{l,m}}r_{0}-\frac{\partial r_{0}}{\partial x_{l,m}}r_{1}\right)f_{2}. \label{eq:dB22_dx_general}
\end{align}
From~\eqref{eq:beta_2}, it is straightforward to compute
\begin{equation}
    \frac{\partial\beta_{3}}{\partial x_{i,0}}=-z_{3},\ \frac{\partial\beta_{3}}{\partial x_{i,1}}=-z_{2},\frac{\partial\beta_{3}}{\partial x_{i,2}}=z_{1},\ \frac{\partial\beta_{3}}{\partial x_{i,3}}=-z_{0},
    \label{eq:dbeta2_dxi}
\end{equation}
and
\begin{equation}
    \frac{\partial\beta_{3}}{\partial x_{j,0}}=\frac{\partial\beta_{3}}{\partial x_{j,1}}=\frac{\partial\beta_{3}}{\partial x_{j,2}}=\frac{\partial\beta_{3}}{\partial x_{j,3}}=0.
    \label{eq:dbeta2_dxj}
\end{equation}
Substituting equations~\eqref{eq:dr0_dxi}-\eqref{eq:dr1_dxj},~\eqref{eq:dr2_dxi}-\eqref{eq:dr2_dxj}, and~\eqref{eq:dbeta2_dxi}-\eqref{eq:dbeta2_dxj} into~\eqref{eq:dB22_dx_general} yields the following expressions for $\nicefrac{\partial\mathcal{B}_{22}}{\partial x_{l,m}}$.
\begin{align}
    \frac{\partial\mathcal{B}_{22}}{\partial x_{i,0}}=&-\frac{z_{3}}{\gamma}+\left(\beta_{3}\left(\eta_{i}r_{0}-\mu_{i}r_{1}\right)+\left(\alpha_{1}-2r_{2}\left(\mu_{i}r_{0}+\eta_{i}r_{1}\right)\right)\left(\kappa_{j}r_{0}-\omega_{j}r_{1}\right)\right)f_{1}  \\
	&+r_{2}\left(\mu_{i}\kappa_{j}-\eta_{i}\omega_{j}+z_{0}r_{0}-z_{1}r_{1}\right)f_{1}+r_{2}\left(\kappa_{j}r_{0}-\omega_{j}r_{1}\right)\left(\eta_{i}r_{0}-\mu_{i}r_{1}\right)f_{2}, \label{eq:dB22_dxi0} \\
    \frac{\partial\mathcal{B}_{22}}{\partial x_{i,1}}=&-\frac{z_{2}}{\gamma}+\left(\beta_{3}\left(\kappa_{i}r_{0}-\omega_{i}r_{1}\right)+\left(\beta_{1}-2r_{2}\left(\omega_{i}r_{0}+\kappa_{i}r_{1}\right)\right)\left(\kappa_{j}r_{0}-\omega_{j}r_{1}\right)\right)f_{1}  \\
	&+r_{2}\left(\omega_{i}\kappa_{j}-\kappa_{i}\omega_{j}-z_{1}r_{0}-z_{0}r_{1}\right)f_{1}+r_{2}\left(\kappa_{j}r_{0}-\omega_{j}r_{1}\right)\left(\kappa_{i}r_{0}-\omega_{i}r_{1}\right)f_{2}, \label{eq:dB22_dxi1} \\
    \frac{\partial\mathcal{B}_{22}}{\partial x_{i,2}}=&\frac{z_{1}}{\gamma}+\xi_{1}\left(\kappa_{j}r_{0}-\omega_{j}r_{1}\right)f_{1}, \label{eq:dB22_dxi2} \\
    \frac{\partial\mathcal{B}_{22}}{\partial x_{i,3}}=&-\frac{z_{0}}{\gamma}+\zeta_{1}\left(\kappa_{j}r_{0}-\omega_{j}r_{1}\right)f_{1}, \label{eq:dB22_dxi3} \\
    \frac{\partial\mathcal{B}_{22}}{\partial x_{j,0}}=&\left(\beta_{3}\left(\eta_{j}r_{0}-\mu_{j}r_{1}\right)+\left(\alpha_{3}-2r_{2}\left(\mu_{j}r_{0}+\eta_{j}r_{1}\right)\right)\left(\kappa_{j}r_{0}-\omega_{j}r_{1}\right)+r_{2}\right)f_{1}  \\
	&+r_{2}\left(\kappa_{j}r_{0}-\omega_{j}r_{1}\right)\left(\eta_{j}r_{0}-\mu_{j}r_{1}\right)f_{2}, \label{eq:dB22_dxj0} \\
    \frac{\partial\mathcal{B}_{22}}{\partial x_{j,1}}=&2\left(\beta_{3}-r_{2}\left(\omega_{j}r_{0}+\kappa_{j}r_{1}\right)\right)\left(\kappa_{j}r_{0}-\omega_{j}r_{1}\right)f_{1}+r_{2}\left(\kappa_{j}r_{0}-\omega_{j}r_{1}\right)^{2}f_{2},\label{eq:dB22_dxj1} \\
    \frac{\partial\mathcal{B}_{22}}{\partial x_{j,2}}=&\kappa_{j}\left(\kappa_{j}r_{0}-\omega_{j}r_{1}\right)f_{1},\label{eq:dB22_dxj2} \\
    \frac{\partial\mathcal{B}_{22}}{\partial x_{j,3}}=&-\eta_{j}\left(\kappa_{j}r_{0}-\omega_{j}r_{1}\right)f_{1},\label{eq:dB22_dxj3}
\end{align}
where we have again used~\eqref{eq:muj_kappaj_min_etaj_muj} to simplify~\eqref{eq:dB22_dxj0}. To compute derivatives of $\mathcal{B}_{23}$ from~\eqref{eq:B_23}, we  follow the derivation of~\eqref{eq:dA11_left} to derive the general form
\begin{equation}
    \frac{\partial\mathcal{B}_{23}}{\partial x_{l,m}}	=\frac{\frac{\partial\kappa_{j}}{\partial x_{l,m}}}{\gamma}+\kappa_{j}\left(\frac{\partial r_{1}}{\partial x_{l,m}}r_{0}-\frac{\partial r_{0}}{\partial x_{l,m}}r_{1}\right)f_{1}.
    \label{eq:dB23_dx_general}
\end{equation}
Substituting equations~\eqref{eq:dr0_dxi}-\eqref{eq:dr1_dxj} and~\eqref{eq:dkappai_dxi}-\eqref{eq:dkappai_dxj} into~\eqref{eq:dB23_dx_general} yields the following expressions for $\nicefrac{\partial\mathcal{B}_{23}}{\partial x_{l,m}}$.
\begin{align}
    \frac{\partial\mathcal{B}_{23}}{\partial x_{i,0}}&=\frac{z_{0}}{\gamma}+\kappa_{j}\left(\eta_{i}r_{0}-\mu_{i}r_{1}\right)f_{1}, \label{eq:dB23_dxi0} \\
    \frac{\partial\mathcal{B}_{23}}{\partial x_{i,1}}&=-\frac{z_{1}}{\gamma}+\kappa_{j}\left(\kappa_{i}r_{0}-\omega_{i}r_{1}\right)f_{1}, \label{eq:dB23_dxi1} \\
    \frac{\partial\mathcal{B}_{23}}{\partial x_{i,2}}&=\frac{\partial\mathcal{B}_{23}}{\partial x_{i,3}}=0, \label{eq:dB23_dxi2_xi3} \\
    \frac{\partial\mathcal{B}_{23}}{\partial x_{j,0}}&=\kappa_{j}\left(\eta_{j}r_{0}-\mu_{j}r_{1}\right)f_{1}, \label{eq:dB23_dxj0} \\
    \frac{\partial\mathcal{B}_{23}}{\partial x_{j,1}}&=\kappa_{j}\left(\kappa_{j}r_{0}-\omega_{j}r_{1}\right)f_{1}, \label{eq:dB23_dxj1} \\
    \frac{\partial\mathcal{B}_{23}}{\partial x_{j,2}}&=\frac{\partial\mathcal{B}_{23}}{\partial x_{j,3}}=0. \label{eq:dB23_dxj2_xj3}
\end{align}
Derivatives of $\mathcal{B}_{24}$ from~\eqref{eq:B_24} again follow the derivation of~\eqref{eq:dA11_left}, allowing us to derive the general form
\begin{equation}
    \frac{\partial\mathcal{B}_{24}}{\partial x_{l,m}}	=-\frac{\frac{\partial\eta_{j}}{\partial x_{l,m}}}{\gamma}-\eta_{j}\left(\frac{\partial r_{1}}{\partial x_{l,m}}r_{0}-\frac{\partial r_{0}}{\partial x_{l,m}}r_{1}\right)f_{1}.
    \label{eq:dB24_dx_general}
\end{equation}
Substituting equations~\eqref{eq:dr0_dxi}-\eqref{eq:dr1_dxj} and~\eqref{eq:detaj_dxi}-\eqref{eq:detaj_dxj} into~\eqref{eq:dB23_dx_general} yields the following expressions for $\nicefrac{\partial\mathcal{B}_{24}}{\partial x_{l,m}}$.
\begin{align}
    \frac{\partial\mathcal{B}_{24}}{\partial x_{i,0}}&=\frac{z_{1}}{\gamma}-\eta_{j}\left(\eta_{i}r_{0}-\mu_{i}r_{1}\right)f_{1}, \label{eq:dB24_dxi0} \\
    \frac{\partial\mathcal{B}_{24}}{\partial x_{i,1}}&=\frac{z_{0}}{\gamma}-\eta_{j}\left(\kappa_{i}r_{0}-\omega_{i}r_{1}\right)f_{1}, \label{eq:dB24_dxi1} \\
    \frac{\partial\mathcal{B}_{24}}{\partial x_{i,2}}&=\frac{\partial\mathcal{B}_{24}}{\partial x_{i,3}}=0, \label{eq:dB24_dxi2_xi3} \\
    \frac{\partial\mathcal{B}_{24}}{\partial x_{j,0}}&=-\eta_{j}\left(\eta_{j}r_{0}-\mu_{j}r_{1}\right)f_{1}, \label{eq:dB24_dxj0} \\
    \frac{\partial\mathcal{B}_{24}}{\partial x_{j,1}}&=-\eta_{j}\left(\kappa_{j}r_{0}-\omega_{j}r_{1}\right)f_{1}, \label{eq:dB24_dxj1} \\
    \frac{\partial\mathcal{B}_{24}}{\partial x_{j,2}}&=\frac{\partial\mathcal{B}_{24}}{\partial x_{j,3}}=0. \label{eq:dB24_dxj2_xj3}
\end{align}
Since $\mathcal{B}_{31}$ from~\eqref{eq:B_31} matches the structure of $\mathcal{A}_{11}$, its general form is given by
\begin{align}
    \frac{\partial\mathcal{B}_{31}}{\partial x_{l,m}}=&\frac{\frac{\partial\beta_{3}}{\partial x_{l,m}}}{\gamma}+\left(\beta_{3}\left(\frac{\partial r_{1}}{\partial x_{l,m}}r_{0}-\frac{\partial r_{0}}{\partial x_{l,m}}r_{1}\right)+\left(\frac{\partial r_{3}}{\partial x_{l,m}}-2r_{3}\left(\frac{\partial r_{0}}{\partial x_{l,m}}r_{0}+\frac{\partial r_{1}}{\partial x_{l,m}}r_{1}\right)\right)\left(\eta_{j}r_{0}-\mu_{j}r_{1}\right)\right)f_{1}  \\
	&+r_{3}\left(\frac{\partial r_{0}}{\partial x_{l,m}}\eta_{j}-\frac{\partial r_{1}}{\partial x_{l,m}}\mu_{j}+\frac{\partial\eta_{j}}{\partial x_{l,m}}r_{0}-\frac{\partial\mu_{j}}{\partial x_{l,m}}r_{1}\right)f_{1}+r_{3}\left(\eta_{j}r_{0}-\mu_{j}r_{1}\right)\left(\frac{\partial r_{1}}{\partial x_{l,m}}r_{0}-\frac{\partial r_{0}}{\partial x_{l,m}}r_{1}\right)f_{2}. \label{eq:dB31_dx_general}
\end{align}
Substituting equations~\eqref{eq:dr0_dxi}-\eqref{eq:dr1_dxj},~\eqref{eq:dr3_dxi}-\eqref{eq:dr3_dxj},~\eqref{eq:detaj_dxi}-\eqref{eq:dmuj_dxj}, and~\eqref{eq:dbeta2_dxi}-\eqref{eq:dbeta2_dxj} into~\eqref{eq:dB31_dx_general} yields the following expressions for $\nicefrac{\partial\mathcal{B}_{31}}{\partial x_{l,m}}$.
\begin{align}
    \frac{\partial\mathcal{B}_{31}}{\partial x_{i,0}}&=-\frac{z_{3}}{\gamma}+\left(\beta_{3}\left(\eta_{i}r_{0}-\mu_{i}r_{1}\right)+\left(\alpha_{2}-2r_{3}\left(\mu_{i}r_{0}+\eta_{i}r_{1}\right)\right)\left(\eta_{j}r_{0}-\mu_{j}r_{1}\right)\right)f_{1}  \\
	&+r_{3}\left(\mu_{i}\eta_{j}-\eta_{i}\mu_{j}-z_{1}r_{0}-z_{0}r_{1}\right)f_{1}+r_{3}\left(\eta_{j}r_{0}-\mu_{j}r_{1}\right)\left(\eta_{i}r_{0}-\mu_{i}r_{1}\right)f_{2}, \label{eq:dB31_dxi0} \\
    \frac{\partial\mathcal{B}_{31}}{\partial x_{i,1}}&=-\frac{z_{2}}{\gamma}+\left(\beta_{3}\left(\kappa_{i}r_{0}-\omega_{j}r_{1}\right)+\left(\beta_{2}-2r_{3}\left(\omega_{j}r_{0}+\kappa_{i}r_{1}\right)\right)\left(\eta_{j}r_{0}-\mu_{j}r_{1}\right)\right)f_{1}  \\
	&+r_{3}\left(\omega_{j}\eta_{j}-\kappa_{i}\mu_{j}-z_{0}r_{0}+z_{1}r_{1}\right)f_{1}+r_{3}\left(\eta_{j}r_{0}-\mu_{j}r_{1}\right)\left(\kappa_{i}r_{0}-\omega_{j}r_{1}\right)f_{2}, \label{eq:dB31_dxi1} \\
    \frac{\partial\mathcal{B}_{31}}{\partial x_{i,2}}&=\frac{z_{1}}{\gamma}-\zeta_{1}\left(\eta_{j}r_{0}-\mu_{j}r_{1}\right)f_{1}, \label{eq:dB31_dxi2} \\
    \frac{\partial\mathcal{B}_{31}}{\partial x_{i,3}}&=-\frac{z_{0}}{\gamma}+\xi_{1}\left(\eta_{j}r_{0}-\mu_{j}r_{1}\right)f_{1}, \label{eq:dB31_dxi3} \\
    \frac{\partial\mathcal{B}_{31}}{\partial x_{j,0}}&=2\left(\beta_{3}-r_{3}\left(\mu_{j}r_{0}+\eta_{j}r_{1}\right)\right)\left(\eta_{j}r_{0}-\mu_{j}r_{1}\right)f_{1}+r_{3}\left(\eta_{j}r_{0}-\mu_{j}r_{1}\right)^{2}f_{2}, \label{eq:dB31_dxj0} \\
    \frac{\partial\mathcal{B}_{31}}{\partial x_{j,1}}&=\left(\beta_{3}\left(\kappa_{j}r_{0}-\omega_{j}r_{1}\right)-\left(\alpha_{3}+2r_{3}\left(\omega_{j}r_{0}+\kappa_{j}r_{1}\right)\right)\left(\eta_{j}r_{0}-\mu_{j}r_{1}\right)-r_{3}\right)f_{1}  \\
	&+r_{3}\left(\eta_{j}r_{0}-\mu_{j}r_{1}\right)\left(\kappa_{j}r_{0}-\omega_{j}r_{1}\right)f_{2}, \label{eq:dB31_dxj1} \\
    \frac{\partial\mathcal{B}_{31}}{\partial x_{j,2}}&=\eta_{j}\left(\eta_{j}r_{0}-\mu_{j}r_{1}\right)f_{1}, \label{eq:dB31_dxj2} \\
    \frac{\partial\mathcal{B}_{31}}{\partial x_{j,3}}&=\kappa_{j}\left(\eta_{j}r_{0}-\mu_{j}r_{1}\right)f_{1}, \label{eq:dB31_dxj3}
\end{align}
where we have again used~\eqref{eq:muj_kappaj_min_etaj_muj} to simplify~\eqref{eq:dB31_dxj1}. $\mathcal{B}_{32}$ from~\eqref{eq:B_32} also matches the structure of $\mathcal{A}_{11}$, so its general form is given by
\begin{align}
    \frac{\partial\mathcal{B}_{32}}{\partial x_{l,m}}=&-\frac{\frac{\partial\alpha_{3}}{\partial x_{l,m}}}{\gamma}+\left(-\alpha_{3}\left(\frac{\partial r_{1}}{\partial x_{l,m}}r_{0}-\frac{\partial r_{0}}{\partial x_{l,m}}r_{1}\right)+\left(\frac{\partial r_{3}}{\partial x_{l,m}}-2r_{3}\left(\frac{\partial r_{0}}{\partial x_{l,m}}r_{0}+\frac{\partial r_{1}}{\partial x_{l,m}}r_{1}\right)\right)\left(\kappa_{j}r_{0}-\omega_{j}r_{1}\right)\right)f_{1}  \\
	&+r_{3}\left(\frac{\partial r_{0}}{\partial x_{l,m}}\kappa_{j}-\frac{\partial r_{1}}{\partial x_{l,m}}\omega_{j}+\frac{\partial\kappa_{j}}{\partial x_{l,m}}r_{0}-\frac{\partial\omega_{j}}{\partial x_{l,m}}r_{1}\right)f_{1}+r_{3}\left(\kappa_{j}r_{0}-\omega_{j}r_{1}\right)\left(\frac{\partial r_{1}}{\partial x_{l,m}}r_{0}-\frac{\partial r_{0}}{\partial x_{l,m}}r_{1}\right)f_{2}.\label{eq:dB32_dx_general}
\end{align}
Substituting equations~\eqref{eq:dr0_dxi}-\eqref{eq:dr1_dxj},~\eqref{eq:dr3_dxi}-\eqref{eq:dr3_dxj},~\eqref{eq:dkappaj_dxi}-\eqref{eq:domegaj_dxj}, and~\eqref{eq:dalpha2_dxi}-\eqref{eq:dalpha2_dxj} into~\eqref{eq:dB32_dx_general} yields the following expressions for $\nicefrac{\partial\mathcal{B}_{32}}{\partial x_{l,m}}$.
\begin{align}
    \frac{\partial\mathcal{B}_{32}}{\partial x_{i,0}}=&\frac{z_{2}}{\gamma}+\left(-\alpha_{3}\left(\eta_{i}r_{0}-\mu_{i}r_{1}\right)+\left(\alpha_{2}-2r_{3}\left(\mu_{i}r_{0}+\eta_{i}r_{1}\right)\right)\left(\kappa_{j}r_{0}-\omega_{j}r_{1}\right)\right)f_{1}  \\
	&+r_{3}\left(\mu_{i}\kappa_{j}-\eta_{i}\omega_{j}+z_{0}r_{0}-z_{1}r_{1}\right)f_{1}+r_{3}\left(\kappa_{j}r_{0}-\omega_{j}r_{1}\right)\left(\eta_{i}r_{0}-\mu_{i}r_{1}\right)f_{2}, \label{eq:dB32_dxi0} \\
    \frac{\partial\mathcal{B}_{32}}{\partial x_{i,1}}=&-\frac{z_{3}}{\gamma}+\left(-\alpha_{3}\left(\kappa_{i}r_{0}-\omega_{j}r_{1}\right)+\left(\beta_{2}-2r_{3}\left(\omega_{j}r_{0}+\kappa_{i}r_{1}\right)\right)\left(\kappa_{j}r_{0}-\omega_{j}r_{1}\right)\right)f_{1}  \\
	&+r_{3}\left(\omega_{j}\kappa_{j}-\kappa_{i}\omega_{j}-z_{1}r_{0}-z_{0}r_{1}\right)f_{1}+r_{3}\left(\kappa_{j}r_{0}-\omega_{j}r_{1}\right)\left(\kappa_{i}r_{0}-\omega_{j}r_{1}\right)f_{2}, \label{eq:dB32_dxi1} \\
    \frac{\partial\mathcal{B}_{32}}{\partial x_{i,2}}=&\frac{z_{0}}{\gamma}-\zeta_{1}\left(\kappa_{j}r_{0}-\omega_{j}r_{1}\right)f_{1}, \label{eq:dB32_dxi2} \\
    \frac{\partial\mathcal{B}_{32}}{\partial x_{i,3}}=&\frac{z_{1}}{\gamma}+\xi_{1}\left(\kappa_{j}r_{0}-\omega_{j}r_{1}\right)f_{1}, \label{eq:dB32_dxi3} \\
    \frac{\partial\mathcal{B}_{32}}{\partial x_{j,0}}=&\left(-\alpha_{3}\left(\eta_{j}r_{0}-\mu_{j}r_{1}\right)+\left(\beta_{3}-2r_{3}\left(\mu_{j}r_{0}+\eta_{j}r_{1}\right)\right)\left(\kappa_{j}r_{0}-\omega_{j}r_{1}\right)+r_{3}\right)f_{1}  \\
	&+r_{3}\left(\kappa_{j}r_{0}-\omega_{j}r_{1}\right)\left(\eta_{j}r_{0}-\mu_{j}r_{1}\right)f_{2}, \label{eq:dB32_dxj0} \\
    \frac{\partial\mathcal{B}_{32}}{\partial x_{j,1}}=&-2\left(\alpha_{3}+r_{3}\left(\omega_{j}r_{0}+\kappa_{j}r_{1}\right)\right)\left(\kappa_{j}r_{0}-\omega_{j}r_{1}\right)f_{1}+r_{3}\left(\kappa_{j}r_{0}-\omega_{j}r_{1}\right)^{2}f_{2}, \label{eq:dB32_dxj1} \\
    \frac{\partial\mathcal{B}_{32}}{\partial x_{j,2}}=&\eta_{j}\left(\kappa_{j}r_{0}-\omega_{j}r_{1}\right)f_{1}, \label{eq:dB32_dxj2} \\
    \frac{\partial\mathcal{B}_{32}}{\partial x_{j,3}}=&\kappa_{j}\left(\kappa_{j}r_{0}-\omega_{j}r_{1}\right)f_{1}, \label{eq:dB32_dxj3}
\end{align}
where we have again used~\eqref{eq:muj_kappaj_min_etaj_muj} to simplify~\eqref{eq:dB32_dxj0}.
Derivatives of $\mathcal{B}_{33}$ from~\eqref{eq:B_33} follow the derivation of~\eqref{eq:dA11_left}, allowing us to derive the general form
\begin{equation}
    \frac{\partial\mathcal{B}_{33}}{\partial x_{l,m}}	=\frac{\frac{\partial\eta_{j}}{\partial x_{l,m}}}{\gamma}+\eta_{j}\left(\frac{\partial r_{1}}{\partial x_{l,m}}r_{0}-\frac{\partial r_{0}}{\partial x_{l,m}}r_{1}\right)f_{1}.
    \label{eq:dB33_dx_general}
\end{equation}
Substituting equations~\eqref{eq:dr0_dxi}-\eqref{eq:dr1_dxj} and~\eqref{eq:detaj_dxi}-\eqref{eq:detaj_dxj} into~\eqref{eq:dB33_dx_general} yields the following expressions for $\nicefrac{\partial\mathcal{B}_{33}}{\partial x_{l,m}}$.
\begin{align}
    \frac{\partial\mathcal{B}_{33}}{\partial x_{i,0}}&=-\frac{z_{1}}{\gamma}+\eta_{j}\left(\eta_{i}r_{0}-\mu_{i}r_{1}\right)f_{1}, \label{eq:dB33_dxi0} \\
    \frac{\partial\mathcal{B}_{33}}{\partial x_{i,1}}&=-\frac{z_{0}}{\gamma}+\eta_{j}\left(\kappa_{i}r_{0}-\omega_{i}r_{1}\right)f_{1}, \label{eq:dB33_dxi1} \\
    \frac{\partial\mathcal{B}_{33}}{\partial x_{i,2}}&=\frac{\partial\mathcal{B}_{33}}{\partial x_{i,3}}=0, \label{eq:dB33_dxi2_xi3} \\
    \frac{\partial\mathcal{B}_{33}}{\partial x_{j,0}}&=\eta_{j}\left(\eta_{j}r_{0}-\mu_{j}r_{1}\right)f_{1}, \label{eq:dB33_dxj0} \\
    \frac{\partial\mathcal{B}_{33}}{\partial x_{j,1}}&=\eta_{j}\left(\kappa_{j}r_{0}-\omega_{j}r_{1}\right)f_{1}, \label{eq:dB33_dxj1} \\
    \frac{\partial\mathcal{B}_{33}}{\partial x_{j,2}}&=\frac{\partial\mathcal{B}_{33}}{\partial x_{j,3}}=0. \label{eq:dB33_dxj2_xj3} \\
\end{align}
Derivatives of $\mathcal{B}_{34}$ from~\eqref{eq:B_34} again follow the derivation of~\eqref{eq:dA11_left}, yielding the general form
\begin{equation}
    \frac{\partial\mathcal{B}_{34}}{\partial x_{l,m}}=\frac{\frac{\partial\kappa_{j}}{\partial x_{l,m}}}{\gamma}+\kappa_{j}\left(\frac{\partial r_{1}}{\partial x_{l,m}}r_{0}-\frac{\partial r_{0}}{\partial x_{l,m}}r_{1}\right)f_{1}.
    \label{eq:dB34_dx_general}
\end{equation}
Substituting equations~\eqref{eq:dr0_dxi}-\eqref{eq:dr1_dxj} and~\eqref{eq:dkappaj_dxi}-\eqref{eq:dkappaj_dxj} into~\eqref{eq:dB34_dx_general} yields the following expressions for $\nicefrac{\partial\mathcal{B}_{34}}{\partial x_{l,m}}$.
\begin{align}
    \frac{\partial\mathcal{B}_{34}}{\partial x_{i,0}}&=\frac{z_{0}}{\gamma}+\kappa_{j}\left(\eta_{i}r_{0}-\mu_{i}r_{1}\right)f_{1}, \label{eq:dB34_dxi0} \\
    \frac{\partial\mathcal{B}_{34}}{\partial x_{i,1}}&=-\frac{z_{1}}{\gamma}+\kappa_{j}\left(\kappa_{i}r_{0}-\omega_{i}r_{1}\right)f_{1}, \label{eq:dB34_dxi1} \\
    \frac{\partial\mathcal{B}_{34}}{\partial x_{i,2}}&=\frac{\partial\mathcal{B}_{34}}{\partial x_{i,3}}=0, \label{eq:dB34_dxi2_xi3} \\
    \frac{\partial\mathcal{B}_{34}}{\partial x_{j,0}}&=\kappa_{j}\left(\eta_{j}r_{0}-\mu_{j}r_{1}\right)f_{1}, \label{eq:dB34_dxj0} \\
    \frac{\partial\mathcal{B}_{34}}{\partial x_{j,1}}&=\kappa_{j}\left(\kappa_{j}r_{0}-\omega_{j}r_{1}\right)f_{1}, \label{eq:dB34_dxj1} \\
    \frac{\partial\mathcal{B}_{34}}{\partial x_{j,2}}&=\frac{\partial\mathcal{B}_{34}}{\partial x_{j,3}}=0, \label{eq:dB34_dxj2_xj3} \\
\end{align}
concluding the derivation of Hessian tensors $\frac{\partial}{\partial\mathbf{x}_{i}}\mathcal{A}_{ij}$, $\frac{\partial}{\partial\mathbf{x}_{j}}\mathcal{A}_{ij}$, $\frac{\partial}{\partial\mathbf{x}_{i}}\mathcal{B}_{ij}$, and $\frac{\partial}{\partial\mathbf{x}_{j}}\mathcal{B}_{ij}$.
\section{Derivation of Euclidean Gradient Bounds}\label{app:egrad_bounds}
In this appendix, we derive bounds for components of the Euclidean gradient that are necessary for the proof of Lipschitz continuity of the Riemannian gradient in Appendix~\ref{app:lipschitz}. Specifically, we show that $\Vert \mathbf{e}_{ij}\Vert _{2}$, $\Vert \mathcal{A}_{ij}\Vert _{F}$, and $\Vert \mathcal{B}_{ij}\Vert _{F}$ are bounded for all $(i,j)\in\E$, given that $\X\in\K$, where $\K$ is a compact subset of $\MN$. We begin by providing preliminary derivations that will serve as a reference for the subsequent analysis in this appendix as well as in Appendix~\ref{app:ehess_bounds}.

\subsection{Preliminaries}\label{app:egrad_bounds_prelim}
For reference, we first include definitions for the Frobenius norm and matrix 2-norm. Given a matrix $A\in\mathbb{R}^{n\times m}$ with entries $a_{ij}$, the Frobenius norm of $A$, denoted $\left\Vert A\right\Vert _{F}$, is computed as 
\begin{equation}\label{eq:frob_norm_def}
    \left\Vert A\right\Vert _{F}=\sqrt{\sum_{i=1}^{n}\sum_{j=1}^{m}\left|a_{ij}\right|^{2}}.
\end{equation}
The matrix 2-norm of $A$, denoted $\Vert A\Vert_{2}$, is given by
\begin{equation}\label{eq:matrix_two_norm}
    \Vert A\Vert_{2}=\sqrt{\lambda_{\max}(A^{\top}A)},
\end{equation}
where $\lambda_{\max}(\cdot)$ denotes the maximum eigenvalue of a matrix. We now define the notion of the Euclidean norm, denoted $\Vert\cdot\Vert_{2}$, on $\M$ and $\mathcal{M}^{N}$. Following from the embedding of $\M$ in $\mathbb{R}^{4}$ given by~\eqref{eq:M_embedding_app}, we have, for $\mathbf{x}\in\mathcal{M}$, $\left\Vert \mathbf{x}\right\Vert _{2}=\sqrt{\mathbf{x}^{\top}\mathbf{x}}$. Moreover, from the embedding of $\MN$ in $\mathbb{R}^{4N}$ given by~\eqref{eq:MN_embedding_app}, we have, for $\mathcal{X}\in\mathcal{M}^{N}$, $\left\Vert \mathcal{X}\right\Vert _{2}=\sqrt{\mathcal{X}^{\top}\mathcal{X}}$. Using these definitions, we now derive a lemma on the boundedness of the translational components of poses and manifold residuals associated with pose graphs whose poses are limited to compact subsets of $\MN$.
\begin{lemma}\label{lem:compact_trans_bounds}
    Let $\mathcal{G}=\left(\mathcal{V},\mathcal{E}\right)$ be a pose graph, with associated poses $\X=\vect((\mathbf{x}_{i})_{i\in\mathcal{V}})$ and relative edge measurements $\Z=\vect({(\zij)}_{(i,j)\in\E})$, with $|\V|=N$ and $|\E|=M$. Now, let $\mathbf{x}_{i}=[\mathbf{x}_{i,r}^{\top},\mathbf{x}_{i,d}^{\top}]^{\top}\in\M$ and $\mathbf{r}_{ij}=[\mathbf{r}_{ij,r}^{\top},\mathbf{r}_{ij,d}^{\top}]^{\top}\in\M$, for all $i\in\V$, and for all $(i,j)\in\E$, denote the poses and manifold residuals associated with $\mathcal{G}$, respectively, represented in vector form with explicit rotational and translational (dual) components. Then, given any compact subset $\K\subset\MN$, it holds for all $\X\in\K$ that for all $i\in\mathcal{V}$, and for all $(i,j)\in\mathcal{E}$, that $\left\Vert \mathbf{x}_{i,d}\right\Vert _{2}\leq\bar{\mathbf{t}}_{\mathbf{x}}$ and $\left\Vert \mathbf{r}_{ij,d}\right\Vert _{2}\leq\bar{\mathbf{t}}_{\mathbf{r}}$, with
    \begin{equation}\label{eq:tbarx_tbarr_def}
        \bar{\mathbf{t}}_{\mathbf{x}}\triangleq\sqrt{\overline{\mathbf{T}}^{2}-N},\text{~and~}\bar{\mathbf{t}}_{\mathbf{r}}\triangleq\left(\mathbf{\bar{\mathbf{t}}_{\mathbf{x}}}^{2}+3\right)\bar{\mathbf{z}},
    \end{equation}
    where
    \begin{equation}\label{eq:TBar_def}
        \TBar\triangleq\sup\left\{ \left\Vert \mathcal{X}\right\Vert _{2}\mid\mathcal{X}\in\mathcal{K}\right\}
    \end{equation}
    and
    \begin{equation}\label{eq:zbar_def}
        \bar{\mathbf{z}}\triangleq\underset{\left(i,j\right)\in\mathcal{E}}{\max}\left\{ \left\Vert \zij\right\Vert _{2}\right\}.
    \end{equation}
\end{lemma}
\emph{Proof:} Because $\K$ is compact, it is valid to define $\TBar$ as in~\eqref{eq:TBar_def}. It then follows that for all $\X\in\K$, we have
\begin{equation}\label{eq:X_T_bound}
    \left\Vert \mathcal{X}\right\Vert _{2}	=\sqrt{\sum_{i\in\mathcal{V}}\left\Vert \mathbf{x}_{i}\right\Vert _{2}^{2}}\leq\overline{\mathbf{T}}.
\end{equation}
Since $\mathbf{x}_{i}=[\mathbf{x}_{i,r}^{\top},\mathbf{x}_{i,d}^{\top}]^{\top}$, we have $\left\Vert \mathbf{x}_{i}\right\Vert _{2}^{2}=1+\left\Vert \mathbf{x}_{i,d}\right\Vert _{2}$, and therefore
\begin{equation}
    \sum_{i\in\mathcal{V}}\left\Vert \mathbf{x}_{i}\right\Vert _{2}^{2}	=N+\sum_{i\in\mathcal{V}}\left\Vert \mathbf{x}_{i,d}\right\Vert _{2}^{2}.
\end{equation}
Applying this to equation~\eqref{eq:X_T_bound} and simplifying yields
\begin{equation}\label{eq:tbarx_bound}
    \left\Vert \mathbf{x}_{i,d}\right\Vert _{2}	\leq\sqrt{\overline{\mathbf{T}}^{2}-N}\triangleq\bar{\mathbf{t}}_{\mathbf{x}},
\end{equation}
which gives the left side of~\eqref{eq:tbarx_tbarr_def}. We now address the translational component of $\mathbf{r}_{ij}=[\mathbf{r}_{ij,r}^{\top},\mathbf{r}_{ij,d}^{\top}]^{\top}$. Applying~\eqref{eq:Q_L_R} and~\eqref{eq:Q_R_MM} to the definition of $\rij$ given in~\eqref{eq:r_ij_def} yields
\begin{equation}
    \rij=\zij^{-1}\boxplus\mathbf{x}_{i}^{-1}\boxplus\mathbf{x}_{j}=Q_{R}\left(\mathbf{x}_{j}\right)Q_{R}^{--}\left(\mathbf{x}_{i}\right)\zij.
\end{equation}
It then holds that
\begin{equation}\label{eq:rij_normbound_1}
    \left\Vert \mathbf{r}_{ij}\right\Vert _{2}	\leq\left\Vert Q_{R}\left(\mathbf{x}_{j}\right)\right\Vert _{2}\left\Vert Q_{R}^{--}\left(\mathbf{x}_{i}\right)\right\Vert _{2}\left\Vert \zij\right\Vert _{2},
\end{equation}
where $\Vert Q_{R}(\cdot)\Vert_{2}$ and $\Vert Q_{R}^{--}(\cdot)\Vert_{2}$ denote the matrix 2-norm given by~\eqref{eq:matrix_two_norm}. To simplify~\eqref{eq:rij_normbound_1}, we first derive a bound on $\norm{Q_{R}(\cdot)}$. For any $\mathbf{x}=\left[\mathbf{x}_{r}^{\transpose},\mathbf{x}_{d}^{\transpose}\right]^{\transpose}\in\mathcal{M}$, applying~\eqref{eq:matrix_two_norm} to~\eqref{eq:Q_L_R} and simplifying gives
\begin{equation}
    \left\Vert Q_{R}\left(\mathbf{x}\right)\right\Vert _{2}=\sqrt{1+\frac{1}{2}\left({\left\Vert \mathbf{x}_{d}\right\Vert}_{2}^{2}+\sqrt{\left\Vert \mathbf{x}_{d}\right\Vert_{2}^{2}\left(\left\Vert \mathbf{x}_{d}\right\Vert_{2}^{2}+4\right)}\right)}.
\end{equation}
Therefore, letting $\mathbf{x}=[x_0,x_1,x_2,x_3]^{\top}$, we have
\begin{align}
    \left\Vert Q_{R}\left(\mathbf{x}\right)\right\Vert _{2}&\leq\sqrt{1+\frac{1}{2}\left(\left\Vert \mathbf{x}_{d}\right\Vert_{2}^{2}+\sqrt{\left(\left\Vert \mathbf{x}_{d}\right\Vert_{2}^{2}+4\right)^{2}}\right)}\\
    &=\sqrt{x_{2}^{2}+x_{3}^{2}+3}.\label{eq:Q_R_bound_1}
\end{align}
Since $\left\Vert \mathbf{x}\right\Vert _{2}^{2}=x_{2}^{2}+x_{3}^{2}+1$, equation~\eqref{eq:Q_R_bound_1} implies that
\begin{equation}\label{eq:Q_R_bound_2}
    \left\Vert Q_{R}\left(\mathbf{x}\right)\right\Vert _{2}	\leq\sqrt{\left\Vert \mathbf{x}\right\Vert _{2}^{2}+2}\leq\left\Vert \mathbf{x}\right\Vert _{2}+\sqrt{2},
\end{equation}
which holds for all $\mathbf{x}\in\M$. Noting that $\left\Vert Q_{L}\left(\mathbf{x}\right)\right\Vert _{2}=\left\Vert Q_{R}\left(\mathbf{x}\right)\right\Vert _{2}=\Vert Q_{R}^{--}\left(\mathbf{x}\right)\Vert_{2}$, we can apply~\eqref{eq:Q_R_bound_2} to~\eqref{eq:rij_normbound_1} to write
\begin{equation}\label{eq:rij_normbound_2}
    \left\Vert \mathbf{r}_{ij}\right\Vert _{2}\leq\left(\left\Vert \mathbf{x}_{i}\right\Vert _{2}+\sqrt{2}\right)\left(\left\Vert \mathbf{x}_{j}\right\Vert _{2}+\sqrt{2}\right)\left\Vert \zij\right\Vert _{2}.
\end{equation}
Now, we apply~\eqref{eq:tbarx_bound} and the fact that $\left\Vert\rij\right\Vert_{2}^{2}=\left\Vert \mathbf{r}_{ij,d}\right\Vert _{2}^{2}+1$ to~\eqref{eq:rij_normbound_2} to obtain
\begin{equation}\label{eq:rij_normbound_3}
    \sqrt{\left\Vert\mathbf{r}_{ij,d}\right\Vert _{2}^{2}+1}\leq\left(\mathbf{\bar{\mathbf{t}}_{\mathbf{x}}}^{2}+3\right)\left\Vert\zij\right\Vert_{2},
\end{equation}
and applying~\eqref{eq:zbar_def} to~\eqref{eq:rij_normbound_3} yields, for all $(i,j)\in\E$,
\begin{equation}\label{eq:tbarr_bound}
    \left\Vert \mathbf{r}_{ij,d}\right\Vert _{2}\leq\left(\mathbf{\bar{\mathbf{t}}_{\mathbf{x}}}^{2}+3\right)\bar{\mathbf{z}}=\bar{\mathbf{t}}_{\mathbf{r}},
\end{equation}
with $\tbarr$ from~\eqref{eq:tbarx_tbarr_def}, completing the proof.~\hfill $\blacksquare$

Using Lemma~\ref{lem:compact_trans_bounds}, we now derive a set of preliminary bounds that will aid in the forthcoming analysis. Given a pose graph $\mathcal{G}=\left(\mathcal{V},\mathcal{E}\right)$ as defined in Lemma~\ref{lem:compact_trans_bounds}, we denote $\mathbf{x}_{i}=[x_{i,0},x_{i,1},x_{i,2},x_{i,3}]^{\top}$ and $\mathbf{x}_{j}=[x_{j,0},x_{j,1},x_{j,2},x_{j,3}]^{\top}$, with $i,j\in\V$, to be the poses corresponding to relative measurement $\zij=[z_{0},z_{1},z_{2},z_{3}]^{\top}$, with $(i,j)\in\E$, and let $\mathbf{r}_{ij}=[r_{0},r_{1},r_{2},r_{3}]^{\top}$ be the manifold residual computed via~\eqref{eq:r_ij_def}. Noting that $\mathbf{x}_{i},\mathbf{x}_{j},\zij,\mathbf{r}_{ij}\in\mathcal{M}$, we denote $\phi_{i}$ and $\phi_{j}$ to be the rotation half-angles associated with $\mathbf{x}_{i}$ and $\mathbf{x}_{j}$ such that
\begin{equation}\label{eq:xi_xj_0_1}
    x_{i,0}=\cos\left(\phi_{i}\right),\ x_{i,1}=\sin\left(\phi_{i}\right),\ x_{j,0}=\cos\left(\phi_{j}\right),\ x_{j,1}=\sin\left(\phi_{j}\right),
\end{equation}
and we denote $\phi_{z}$ and $\phi_{r}$ to be the rotation half-angles associated with $\zij$ and $\rij$ such that
\begin{equation}\label{eq:zij_rij_0_1}
    z_{0}=\cos\left(\phi_{z}\right),\ z_{1}=\sin\left(\phi_{z}\right),\ r_{0}=\cos\left(\phi_{r}\right),\ r_{1}=\sin\left(\phi_{r}\right). 
\end{equation}
From~\eqref{eq:xi_xj_0_1}-\eqref{eq:zij_rij_0_1}, we can immediately write
\begin{equation}\label{eq:sin_cos_bounds}
    \left|x_{i,0}\right|,\left|x_{i,1}\right|,\left|x_{j,0}\right|,\left|x_{j,1}\right|,\left|z_{0}\right|,\left|z_{1}\right|,\left|r_{0}\right|,\left|r_{1}\right|	\leq1.
\end{equation}
We now define constants $\bar{\mathbf{z}}_{2}$, $\bar{\mathbf{z}}_{3}$, and $\bar{\mathbf{z}}_{23}$ such that 
\begin{equation}\label{eq:z_2_3_23_def}
    \bar{\mathbf{z}}_{2}\triangleq\underset{\left(i,j\right)\in\mathcal{E}}{\max}\left|z_{2}\right|,~\mathbf{z}_{3}\triangleq\underset{\left(i,j\right)\in\mathcal{E}}{\max}\left|z_{3}\right|,~\bar{\mathbf{z}}_{23}\triangleq\bar{\mathbf{z}}_{2}+\bar{\mathbf{z}}_{3}.
\end{equation}
It then follows from~\eqref{eq:z_2_3_23_def} that 
\begin{equation}\label{eq:z_2_3_23_bound}
    \left|z_{2}\right|\leq\bar{\mathbf{z}}_{2},~\left|z_{3}\right|\leq\bar{\mathbf{z}}_{3},~\left|z_{2}\right|+\left|z_{3}\right|\leq\bar{\mathbf{z}}_{23}
\end{equation}
for all $\left(i,j\right)\in\mathcal{E}$. Furthermore, the function $\operatorname{sinc}(\phi)$ is maximized at $\phi=0$, so $\gamma(\phi(\mathbf{x}))$, as defined in~\eqref{eq:gamma_app}, is bounded by
\begin{equation}\label{eq:gamma_bound}
    \left|\gamma\left(\phi\left(\mathbf{x}\right)\right)\right|\leq\gamma\left(0\right)=1
\end{equation}
for all $\mathbf{x}\in\M$. Additionally, the reciprocal $(\mathrm{sinc}(\phi))^{-1}$ is maximized at $\phi=\nicefrac{\pi}{2}$ over the domain $\phi\in\left(-\frac{\pi}{2},\frac{\pi}{2}\right]$. Applying this fact to~\eqref{eq:gamma_app} and~\eqref{eq:wrap_app} yields
\begin{equation}\label{eq:gamma_recip_bound}
    \left|\frac{1}{\gamma\left(\phi\left(\mathbf{x}\right)\right)}\right|	\leq\left|\frac{1}{\gamma\left(\nicefrac{\pi}{2}\right)}\right|=\frac{\pi}{2}
\end{equation}
for all $\mathbf{x}\in\M$. Because the function $f_{1}\left(\phi\right)$ from~\eqref{eq:f_1_def} takes on values within the range $\left(-1,1\right]$ over the domain $\phi\in\left(-\frac{\pi}{2},\frac{\pi}{2}\right]$, it holds that
\begin{equation}\label{eq:f1_bound}
    \left|f_{1}\left(\phi\left(\mathbf{x}\right)\right)\right|\leq\left|f_{1}\left(\nicefrac{\pi}{2}\right)\right|=1
\end{equation}
for all $\mathbf{x}\in\M$. Since $\left\Vert \cdot\right\Vert _{1}\leq\sqrt{2}\left\Vert \cdot\right\Vert _{2}$, it holds from~\eqref{eq:tbarx_bound} that
\begin{equation}\label{eq:xi_23_bound}
    \left|x_{i,2}\right|,\left|x_{i,3}\right|\leq\left|x_{i,2}\right|+\left|x_{i,3}\right|=\left\Vert \mathbf{x}_{i,d}\right\Vert _{1}\leq\sqrt{2}\left\Vert \mathbf{x}_{i,d}\right\Vert _{2}\leq\tbarx\sqrt{2}
\end{equation}
and
\begin{equation}\label{eq:xj_23_bound}
    \left|x_{j,2}\right|,\left|x_{j,3}\right|\leq\left|x_{j,2}\right|+\left|x_{j,3}\right|=\left\Vert \mathbf{x}_{j,d}\right\Vert _{1}\leq\sqrt{2}\left\Vert \mathbf{x}_{j,d}\right\Vert _{2}\leq\tbarx\sqrt{2}.
\end{equation}
From~\eqref{eq:tbarr_bound}, we have
\begin{equation}\label{eq:rij_23_bound}
    \left|r_{2}\right|,\left|r_{3}\right|\leq\left|r_{2}\right|+\left|r_{3}\right|=\left\Vert\mathbf{r}_{ij,d}\right\Vert_{1}\leq\sqrt{2}\left\Vert\mathbf{r}_{ij,d}\right\Vert _{2}\leq\tbarr\sqrt{2}.
\end{equation}
We now bound entries of the matrix $Q_{i}$ from~\eqref{eq:Q_i_j}, which correspond to~\eqref{eq:mu_i}-\eqref{eq:beta_3}. Substituting~\eqref{eq:sin_cos_bounds} into~\eqref{eq:mu_i}-\eqref{eq:kappa_i},~\eqref{eq:xi_1}-\eqref{eq:zeta_1} and applying angle sum and difference identities yields
\begin{align}
    \mu_{i}&=z_{0}x_{j,0}+z_{1}x_{j,1}=\cos\left(\phi_{z}\right)\cos\left(\phi_{j}\right)+\sin\left(\phi_{z}\right)\sin\left(\phi_{j}\right)=\cos\left(\phi_{j}-\phi_{z}\right),\label{eq:mui_trig_form} \\
    \omega_{i}&=-z_{1}x_{j,0}+z_{0}x_{j,1}=-\sin\left(\phi_{z}\right)\cos\left(\phi_{j}\right)+\cos\left(\phi_{z}\right)\sin\left(\phi_{j}\right)=\sin\left(\phi_{j}-\phi_{z}\right),\label{eq:omegai_trig_form} \\
    \eta_{i}&=-x_{j,0}z_{1}+x_{j,1}z_{0}=-\cos\left(\phi_{j}\right)\sin\left(\phi_{z}\right)+\sin\left(\phi_{j}\right)\cos\left(\phi_{z}\right)=\sin\left(\phi_{j}-\phi_{z}\right),\label{eq:etai_trig_form} \\
    \kappa_{i}&=-x_{j,0}z_{0}-x_{j,1}z_{1}=-\cos\left(\phi_{j}\right)\cos\left(\phi_{z}\right)-\sin\left(\phi_{j}\right)\sin\left(\phi_{z}\right)=-\cos\left(\phi_{j}-\phi_{z}\right),\label{eq:kappai_trig_form} \\
    \xi_{1}&=-x_{j,0}z_{0}+x_{j,1}z_{1}=-\cos\left(\phi_{j}\right)\cos\left(\phi_{z}\right)+\sin\left(\phi_{j}\right)\sin\left(\phi_{z}\right)=-\cos\left(\phi_{j}+\phi_{z}\right), \label{eq:xi1_trig_form} \\
    \zeta_{1}&=-x_{j,0}z_{1}-x_{j,1}z_{0}=-\cos\left(\phi_{j}\right)\sin\left(\phi_{z}\right)-\sin\left(\phi_{j}\right)\cos\left(\phi_{z}\right)=-\sin\left(\phi_{j}+\phi_{z}\right).\label{eq:zeta1_trig_form} 
\end{align}
\noeqref{eq:omegai_trig_form}\noeqref{eq:etai_trig_form}\noeqref{eq:xi1_trig_form}It then follows from~\eqref{eq:mui_trig_form}-\eqref{eq:zeta1_trig_form} that
\begin{equation}\label{eq:mui-zeta1_bounds}
    \left|\mu_{i}\right|,\left|\omega_{i}\right|,\left|\eta_{i}\right|,\left|\kappa_{i}\right|,\left|\xi_{1}\right|,\left|\zeta_{1}\right|	\leq1.
\end{equation}
Next, we apply the triangle inequality and~\eqref{eq:sin_cos_bounds} to the absolute value of~\eqref{eq:alpha_1} yields
\begin{equation}\label{eq:alpha1_tri_eq}
    \left|\alpha_{1}\right|	=\left|-x_{j,0}z{}_{2}-x_{j,1}z_{3}+x_{j,2}z_{0}+x_{j,3}z_{1}\right|\leq\left|z_{2}\right|+\left|z_{3}\right|+\left|x_{j,2}\right|+\left|x_{j,3}\right|,
\end{equation}
and further applying~\eqref{eq:z_2_3_23_bound} and~\eqref{eq:rij_23_bound} to~\eqref{eq:alpha1_tri_eq},~\eqref{eq:beta_1}, and~\eqref{eq:alpha_3}-\eqref{eq:beta_3} yields
\begin{equation}\label{eq:alpha1_bound}
    \left|\alpha_{1}\right|,\left|\beta_{1}\right|,\left|\alpha_{2}\right|,\left|\beta_{2}\right|\leq\bar{\mathbf{z}}_{23}+\tbarx\sqrt{2}.
\end{equation}
We now bound entries of the matrix $Q_{j}$ from~\eqref{eq:Q_i_j}, which correspond to~\eqref{eq:mu_j}-\eqref{eq:beta_2}. Substituting~\eqref{eq:sin_cos_bounds} into~\eqref{eq:mu_j}-\eqref{eq:kappa_j} and applying angle sum and difference identities yields
\begin{align}
    \mu_{j}&=z_{0}x_{i,0}-z_{1}x_{i,1}=\cos\left(\phi_{z}\right)\cos\left(\phi_{i}\right)-\sin\left(\phi_{z}\right)\sin\left(\phi_{i}\right)=\cos\left(\phi_{i}+\phi_{z}\right), \label{eq:muj_trig_form} \\
    \omega_{j}&=z_{1}x_{i,0}+z_{0}x_{i,1}=\sin\left(\phi_{z}\right)\cos\left(\phi_{i}\right)+\cos\left(\phi_{z}\right)\sin\left(\phi_{i}\right)=\sin\left(\phi_{i}+\phi_{z}\right), \label{eq:omegaj_trig_form} \\
    \eta_{j}&=-x_{i,0}z_{1}-x_{i,1}z_{0}=-\cos\left(\phi_{i}\right)\sin\left(\phi_{z}\right)-\sin\left(\phi_{i}\right)\cos\left(\phi_{z}\right)=-\sin\left(\phi_{i}+\phi_{z}\right), \label{eq:etaj_trig_form} \\
    \kappa_{j}&=x_{i,0}z_{0}-x_{i,1}z_{1}=\cos\left(\phi_{i}\right)\cos\left(\phi_{z}\right)-\sin\left(\phi_{i}\right)\sin\left(\phi_{z}\right)=\cos\left(\phi_{i}+\phi_{z}\right). \label{eq:kappaj_trig_form}
\end{align}
\noeqref{eq:omegaj_trig_form}\noeqref{eq:etaj_trig_form}It then follows from~\eqref{eq:muj_trig_form}-\eqref{eq:kappaj_trig_form} that
\begin{equation}\label{eq:muj-kappaj_bounds}
    \left|\mu_{j}\right|,\left|\omega_{j}\right|,\left|\eta_{j}\right|,\left|\kappa_{j}\right|	\leq1.
\end{equation}
Furthermore, applying the derivation of~\eqref{eq:alpha1_bound} to~\eqref{eq:alpha_2}-\eqref{eq:beta_2} yields
\begin{equation}\label{eq:alpha2_beta2_bounds}
    \left|\alpha_{3}\right|,\left|\beta_{3}\right|	\leq\bar{\mathbf{z}}_{23}+\tbarx\sqrt{2}.
\end{equation}
We can also write derivatives of $\phi_{r}$ in trigonometric form by substituting~\eqref{eq:zij_rij_0_1},~\eqref{eq:mui_trig_form}-\eqref{eq:kappai_trig_form}, and~\eqref{eq:muj_trig_form}-\eqref{eq:kappaj_trig_form} and applying angle sum and difference identities, which yields
\begin{align}
    \frac{\partial\phi_{r}}{\partial x_{i,0}}&=\eta_{i}r_{0}-\mu_{i}r_{1}=\sin\left(\phi_{j}-\phi_{z}\right)\cos\left(\phi_{r}\right)-\cos\left(\phi_{j}-\phi_{z}\right)\sin\left(\phi_{r}\right)=\sin\left(\phi_{j}-\phi_{z}-\phi_{ij}\right),\label{eq:dphir_dxi0_trig}\\
    \frac{\partial\phi_{r}}{\partial x_{i,1}}&=\kappa_{i}r_{0}-\omega_{i}r_{1}=-\cos\left(\phi_{j}-\phi_{z}\right)\cos\left(\phi_{r}\right)-\sin\left(\phi_{j}-\phi_{z}\right)\sin\left(\phi_{r}\right)=-\cos\left(\phi_{j}-\phi_{z}-\phi_{r}\right),\label{eq:dphir_dxi1_trig}\\
    \frac{\partial\phi_{r}}{\partial x_{j,0}}&=\eta_{j}r_{0}-\mu_{j}r_{1}=-\sin\left(\phi_{i}+\phi_{z}\right)\cos\left(\phi_{r}\right)-\cos\left(\phi_{i}+\phi_{z}\right)\sin\left(\phi_{r}\right)=-\sin\left(\phi_{i}+\phi_{z}+\phi_{r}\right),\label{eq:dphir_dxi2_trig}\\
    \frac{\partial\phi_{r}}{\partial x_{j,1}}&=\kappa_{j}r_{0}-\omega_{j}r_{1}=\cos\left(\phi_{i}+\phi_{z}\right)\cos\left(\phi_{r}\right)-\sin\left(\phi_{z}+\phi_{i}\right)\sin\left(\phi_{r}\right)=\cos\left(\phi_{i}+\phi_{z}+\phi_{r}\right).\label{eq:dphir_dxi3_trig}
\end{align}
\noeqref{eq:dphir_dxi1_trig}\noeqref{eq:dphir_dxi2_trig}It then follows from~\eqref{eq:dphir_dxi0_trig}-\eqref{eq:dphir_dxi3_trig} that
\begin{equation}\label{eq:dphi_dx_bounds}
    \left|\eta_{i}r_{0}-\mu_{i}r_{1}\right|,\left|\kappa_{i}r_{0}-\omega_{i}r_{1}\right|,\left|\eta_{j}r_{0}-\mu_{j}r_{1}\right|,\left|\kappa_{j}r_{0}-\omega_{j}r_{1}\right|\leq1,
\end{equation}
concluding our preliminary derivations for computing Euclidean gradient bounds.
\subsection{Residual Bounds}\label{app:e_ij_bounds}
We now compute a bound on $\Vert\mathbf{e}_{ij}\Vert_{2}$ for all $(i,j)\in\E$. Applying the definition of $\Vert\cdot\Vert_{2}$ to~\eqref{eq:e_ij_def} yields
\begin{equation}\label{eq:eij_bound_1}
    \left\Vert \mathbf{e}_{ij}\right\Vert _{2}=\left\Vert \mathrm{Log}_{\fatone}\left(\mathbf{r}_{ij}\right)\right\Vert _{2}=\left\Vert \frac{1}{\gamma\left(\phi\left(\mathbf{r}_{ij}\right)\right)}\left[r_{1},r_{2},r_{3}\right]^{\top}\right\Vert _{2}=\left|\frac{1}{\gamma\left(\phi\left(\mathbf{r}_{ij}\right)\right)}\right|\sqrt{r_{1}^{2}+r_{2}^{2}+r_{3}^{2}}.
\end{equation}
Applying the bounds from~\eqref{eq:gamma_recip_bound} and the fact that $r_{1}^{2}=\sin(\phi_{r})^{2}\leq 1$ to~\eqref{eq:eij_bound_1} yields
\begin{equation}\label{eq:eij_bound_2}
    \left\Vert\mathbf{e}_{ij}\right\Vert_{2}\leq\frac{\pi}{2}\sqrt{1+r_{2}^{2}+r_{3}^{2}}=\frac{\pi}{2}\sqrt{1+\left\Vert \mathbf{r}_{ij,d}\right\Vert_{2}^{2}}.
\end{equation}
Finally, applying~\eqref{eq:tbarr_bound} to~\eqref{eq:eij_bound_2} gives, for all $(i,j)\in\E$,
\begin{equation}\label{eq:eij_bound}
    \left\Vert \mathbf{e}_{ij}\right\Vert _{2}\leq\frac{\pi}{2}\sqrt{\bar{\mathbf{t}}_{\mathbf{r}}^{2}+1}\triangleq\bar{\mathbf{e}},
\end{equation}
where we have defined the constant $\bar{\mathbf{e}}$ such that $\left\Vert \mathbf{e}_{ij}\right\Vert _{2}\leq\bar{\mathbf{e}}$.
\subsection{Jacobian Bounds}\label{app:Aij_Bij_bounds}
We now compute a bound on the Frobenius norm of $\mathcal{A}_{ij}$, whose elements are included in equations~\eqref{eq:A_11}-\eqref{eq:A_14},~\eqref{eq:A_21}-\eqref{eq:A_24}, and~\eqref{eq:A_31}-\eqref{eq:A_34}. First, applying the triangle inequality to $\left|\mathcal{A}_{11}\right|$ yields
\begin{equation}\label{eq:A11_bound_stmt}
    \left|\mathcal{A}_{11}\right|	=\left|\frac{\eta_{i}}{\gamma}+r_{1}\left(\eta_{i}r_{0}-\mu_{i}r_{1}\right)f_{1}\right|\leq\left|\frac{\eta_{i}}{\gamma}\right|+\left|r_{1}\right|\left|\eta_{i}r_{0}-\mu_{i}r_{1}\right|\left|f_{1}\right|.
\end{equation}
Applying~\eqref{eq:mui-zeta1_bounds},~\eqref{eq:gamma_recip_bound},~\eqref{eq:sin_cos_bounds},~\eqref{eq:dphi_dx_bounds}, and~\eqref{eq:f1_bound} to \eqref{eq:A11_bound_stmt} yields
\begin{equation}\label{eq:A11_bound}
    \left|\mathcal{A}_{11}\right|\leq\frac{\pi}{2}+1.
\end{equation}
Because $\mathcal{A}_{12}$ has similar structure, applying the same procedure yields
\begin{equation}\label{eq:A12_bound}
    \left|\mathcal{A}_{12}\right|=\left|\frac{\kappa_{i}}{\gamma}+r_{1}\left(\kappa_{i}r_{0}-\omega_{i}r_{1}\right)f_{1}\right|\leq\frac{\pi}{2}+1.
\end{equation}
\noeqref{eq:A12_bound}From~\eqref{eq:A_12}-\eqref{eq:A_13}, we have
\begin{equation}\label{eq:A12_13_bound}
    \left|\mathcal{A}_{12}\right|=\left|\mathcal{A}_{13}\right|=0.
\end{equation}
For $\left|\mathcal{A}_{21}\right|$, we can apply the triangle inequality to write
\begin{equation}\label{eq:A21_bound_stmt}
    \left|\mathcal{A}_{21}\right|	=\left|\frac{\alpha_{1}}{\gamma}+r_{2}\ensuremath{\left(\eta_{i}r_{0}-\mu_{i}r_{1}\right)f_{1}}\right|\leq\frac{\left|\alpha_{1}\right|}{\left|\gamma\right|}+\left|r_{2}\right|\ensuremath{\left|\eta_{i}r_{0}-\mu_{i}r_{1}\right|\left|f_{1}\right|}
\end{equation}
We now define
\begin{equation}\label{eq:rho_def}
    \rho\triangleq\frac{\pi}{2}\left(\bar{\mathbf{z}}_{23}+\tbarx\sqrt{2}\right)+\tbarr\sqrt{2}.
\end{equation}
Then, applying equations~\eqref{eq:alpha1_bound},~\eqref{eq:gamma_recip_bound},~\eqref{eq:rij_23_bound},~\eqref{eq:dphi_dx_bounds}, and~\eqref{eq:f1_bound} to~\eqref{eq:A21_bound_stmt} yields
\begin{equation}\label{eq:A21_bound}
    \left|\mathcal{A}_{21}\right|\leq\frac{\pi}{2}\left(\bar{\mathbf{z}}_{23}+\tbarx\sqrt{2}\right)+\tbarr\sqrt{2}=\rho.
\end{equation}
Applying the same process to $\mathcal{A}_{22}$ from~\eqref{eq:A_22} yields
\begin{equation}\label{eq:A22_bound}
    \left|\mathcal{A}_{22}\right|=\left|\frac{\beta_{1}}{\gamma}+r_{2}\ensuremath{\left(\kappa_{i}r_{0}-\omega_{i}r_{1}\right)}f_{1}\right|\leq\rho.
\end{equation}
The remaining terms have similar structure to $\mathcal{A}_{11}$-$\mathcal{A}_{22}$, so applying the derivations for~\eqref{eq:A11_bound}-\eqref{eq:A12_13_bound},~\eqref{eq:A21_bound}-\eqref{eq:A22_bound} to~\eqref{eq:A_23}, \eqref{eq:A_24}, \eqref{eq:A_31},~\eqref{eq:A_32},~\eqref{eq:A_33}, and~\eqref{eq:A_34} yields
\begin{equation}\label{eq:A23-A34_bound} 
    \left|\mathcal{A}_{23}\right|,\left|\mathcal{A}_{24}\right|,\left|\mathcal{A}_{33}\right|,\left|\mathcal{A}_{34}\right|\leq\frac{\pi}{2}
\end{equation}
and
\begin{equation}\label{eq:A31-A32_bound}
    \left|\mathcal{A}_{31}\right|,\left|\mathcal{A}_{32}\right|\leq\rho.
\end{equation}
Now, we define
\begin{equation}\label{eq:JBar_def}
    \JBar\triangleq\sqrt{2\left(\frac{\pi}{2}+1\right)^{2}+4\rho^{2}+4\left(\frac{\pi}{2}\right)^{2}}.
\end{equation}
\noeqref{eq:A22_bound,eq:A23-A34_bound}Then, substituting~\eqref{eq:A11_bound}-\eqref{eq:A12_13_bound} and~\eqref{eq:A21_bound}-\eqref{eq:A31-A32_bound} into the definition of the Frobenius norm from~\eqref{eq:frob_norm_def} yields
\begin{equation}\label{eq:Aij_bound}
    \left\Vert \mathcal{A}_{ij}\right\Vert _{F}	\leq\sqrt{2\left(\frac{\pi}{2}+1\right)^{2}+4\rho^{2}+4\left(\frac{\pi}{2}\right)^{2}}=\JBar,
\end{equation}
which holds for all $(i,j)\in\E$.\\

We now derive a bound on $\Vert\mathcal{B}_{ij}\Vert_{F}$. Because $\mathcal{A}_{11}$-$\mathcal{A}_{34}$ and $\mathcal{B}_{11}$-$\mathcal{B}_{34}$ share identical structure, we apply~\eqref{eq:sin_cos_bounds},~\eqref{eq:gamma_recip_bound},~\eqref{eq:f1_bound}, \eqref{eq:rij_23_bound}, \eqref{eq:muj-kappaj_bounds}-\eqref{eq:alpha2_beta2_bounds}, and~\eqref{eq:dphi_dx_bounds} to the definitions of $\mathcal{B}_{ij}$ entries in~\eqref{eq:B_11}-\eqref{eq:B_14},~\eqref{eq:B_21}-\eqref{eq:B_24}, and~\eqref{eq:B_31}-\eqref{eq:B_34} to write
\begin{equation}\label{eq:Bij_bound_1}
    \left|\mathcal{B}_{11}\right|,\left|\mathcal{B}_{12}\right|\leq\frac{\pi}{2}+1,
\end{equation}
\begin{equation}\label{eq:Bij_bound_2}
    \left|\mathcal{B}_{13}\right|=\left|\mathcal{B}_{14}\right|=0,
\end{equation}
\begin{equation}\label{eq:Bij_bound_3}
    \left|\mathcal{B}_{21}\right|,\left|\mathcal{B}_{22}\right|,\left|\mathcal{B}_{31}\right|,\left|\mathcal{B}_{32}\right|\leq\rho,
\end{equation}
and
\begin{equation}\label{eq:Bij_bound_4}
    \left|\mathcal{B}_{23}\right|,\left|\mathcal{B}_{24}\right|,\left|\mathcal{B}_{33}\right|,\left|\mathcal{B}_{34}\right|\leq\frac{\pi}{2},
\end{equation}
\noeqref{eq:Bij_bound_2,eq:Bij_bound_3}with $\rho$ given by~\eqref{eq:rho_def}. Therefore, we can substitute~\eqref{eq:Bij_bound_1}-\eqref{eq:Bij_bound_4} into the definition of the Frobenius norm from~\eqref{eq:frob_norm_def} to obtain
\begin{equation}\label{eq:Bij_bound}
    \left\Vert\mathcal{B}_{ij}\right\Vert _{F}\leq\JBar,
\end{equation}
with $\JBar$ given by~\eqref{eq:JBar_def}, which holds for all $(i,j)\in\E$.
\subsection{Euclidean Gradient Bounds}\label{app:gij_bounds}
The proof of Lipschitz continuity of the Riemannian gradient in Appendix~\ref{app:lipschitz} depends on the boundedness of the first two entries of $g_{ij,k}$ from~\eqref{eq:g_ijl} for all $(i,j)\in\E$ and for all $k\in\V$, which we now show. To accomodate the subsequent analysis, we write $\mathcal{A}_{ij}^{\top}\Omega_{ij}\mathbf{e}_{ij}=[g_{i, 0}, g_{i, 1}, g_{i, 2}, g_{i, 3}]^{\top}$ and $\mathcal{B}_{ij}^{\top}\Omega_{ij}\mathbf{e}_{ij}=[g_{j, 0}, g_{j, 1}, g_{j, 2}, g_{j, 3}]^{\top}$ in entry-wise vector form. It then suffices to show that $\left|g_{i,0}\right|$, $\left|g_{i,1}\right|$, $\left|g_{j,0}\right|$, and $\left|g_{j,1}\right|$ are bounded. First, we have
\begin{equation}\label{eq:gi_def}
    \mathcal{A}_{ij}^{\top}\Omega_{ij}\mathbf{e}_{ij}=\mathcal{A}_{ij}^{\top}\left[\begin{array}{c}
        \left\langle \left[\Omega_{ij}\right]_{1},\mathbf{e}_{ij}\right\rangle \\
        \left\langle \left[\Omega_{ij}\right]_{2},\mathbf{e}_{ij}\right\rangle \\
        \left\langle \left[\Omega_{ij}\right]_{3},\mathbf{e}_{ij}\right\rangle 
        \end{array}\right],
\end{equation}
where $\left\langle \cdot,\cdot\right\rangle$ denotes the Euclidean inner product and $\left[\Omega_{ij}\right]_{l}$ denotes the $l$th row of $\Omega_{ij}$. We then have
\begin{equation}\label{eq:gi_expansion}
    \left[\begin{array}{c}
        g_{i,0}\\
        g_{i,1}\\
        g_{i,2}\\ 
        g_{i,3}
        \end{array}\right]	=\left[\begin{array}{c}
        \mathcal{A}_{11}\left[\Omega_{ij}\right]_{1}^{\top}\mathbf{e}_{ij}+\mathcal{A}_{21}\left[\Omega_{ij}\right]_{2}^{\top}\mathbf{e}_{ij}+\mathcal{A}_{31}\left[\Omega_{ij}\right]_{3}^{\top}\mathbf{e}_{ij}\\
        \mathcal{A}_{12}\left[\Omega_{ij}\right]_{1}^{\top}\mathbf{e}_{ij}+\mathcal{A}_{22}\left[\Omega_{ij}\right]_{2}^{\top}\mathbf{e}_{ij}+\mathcal{A}_{32}\left[\Omega_{ij}\right]_{3}^{\top}\mathbf{e}_{ij}\\
        \mathcal{A}_{13}\left[\Omega_{ij}\right]_{1}^{\top}\mathbf{e}_{ij}+\mathcal{A}_{23}\left[\Omega_{ij}\right]_{2}^{\top}\mathbf{e}_{ij}+\mathcal{A}_{32}\left[\Omega_{ij}\right]_{3}^{\top}\mathbf{e}_{ij}\\
        \mathcal{A}_{14}\left[\Omega_{ij}\right]_{1}^{\top}\mathbf{e}_{ij}+\mathcal{A}_{24}\left[\Omega_{ij}\right]_{2}^{\top}\mathbf{e}_{ij}+\mathcal{A}_{34}\left[\Omega_{ij}\right]_{3}^{\top}\mathbf{e}_{ij}
        \end{array}\right].
\end{equation}
Extracting the first two terms from~\eqref{eq:gi_expansion} and taking absolute values yields
\begin{align}
    \left|g_{i,0}\right|&=\left|\mathcal{A}_{11}\left[\Omega_{ij}\right]_{1}^{\top}\mathbf{e}_{ij}+\mathcal{A}_{21}\left[\Omega_{ij}\right]_{2}^{\top}\mathbf{e}_{ij}+\mathcal{A}_{31}\left[\Omega_{ij}\right]_{3}^{\top}\mathbf{e}_{ij}\right|, \label{eq:gi_0_expansion}\\
    \left|g_{i,1}\right|&=\left|\mathcal{A}_{12}\left[\Omega_{ij}\right]_{1}^{\top}\mathbf{e}_{ij}+\mathcal{A}_{22}\left[\Omega_{ij}\right]_{2}^{\top}\mathbf{e}_{ij}+\mathcal{A}_{32}\left[\Omega_{ij}\right]_{3}^{\top}\mathbf{e}_{ij}\right|. \label{eq:gi_1_expansion}
\end{align}
Letting $\mathbf{e}_{ij}=[e_{0},e_{1},e_{2}]^{\top}$, then applying the triangle inequality to~\eqref{eq:gi_0_expansion} and simplifying yields
\begin{equation}\label{eq:gi_0_pre_bound}
    \left|g_{i,0}\right|=\left(\sum_{l=1}^{3}\left|\mathcal{A}_{l1}\right|\left|\Omega_{l1}\right|\right)\left|e_{0}\right|+\left(\sum_{l=1}^{3}\left|\mathcal{A}_{l1}\right|\left|\Omega_{l2}\right|\right)\left|e_{1}\right|+\left(\sum_{l=1}^{3}\left|\mathcal{A}_{l1}\right|\left|\Omega_{l3}\right|\right)\left|e_{2}\right|.
\end{equation}
Now, we define 
\begin{equation}\label{eq:gbar_def}
    \gbar\triangleq\sqrt{2}\left(\left(\frac{\pi}{2}+1\right)\left(\sum_{l=1}^{3}\left|\Omega_{1l}\right|\right)+\rho\left(\sum_{l=1}^{3}\left(\left|\Omega_{2l}\right|+\left|\Omega_{3l}\right|\right)\right)\right)\ebar,
\end{equation}
with $\rho$ defined in~\eqref{eq:rho_def} and $\ebar$ defined in~\eqref{eq:eij_bound}. From~\eqref{eq:eij_bound}, it follows that
\begin{equation}\label{eq:eij_entry_bound}
    \left|e_{0}\right|,\left|e_{1}\right|,\left|e_{2}\right|\leq\Vert\mathbf{e}_{ij}\Vert_{1}\leq\sqrt{2}\Vert\mathbf{e}_{ij}\Vert_{2}\leq\ebar\sqrt{2}.
\end{equation}
Applying~\eqref{eq:eij_entry_bound} and the bounds from~\eqref{eq:A11_bound},~\eqref{eq:A21_bound},~\eqref{eq:A31-A32_bound} into~\eqref{eq:gi_0_pre_bound} yields $|g_{i,0}|\leq\gbar$, and applying the same procedure for $|g_{i,1}|$ yields $|g_{i,1}|\leq\gbar$. Furthermore, repeating the derivation for $|g_{j,0}|$, and $|g_{j,1}|$ using the bounds from~\eqref{eq:Bij_bound_1} and~\eqref{eq:Bij_bound_2} yields $|g_{j,0}|,|g_{j,1}|\leq\gbar$. Summarizing, we have
\begin{equation}\label{eq:gij_bounds}
    \left|g_{i,0}\right|,\left|g_{i,1}\right|,\left|g_{j,0}\right|,\left|g_{j,1}\right|	\leq\gbar,
\end{equation}
with $\gbar$ given by~\eqref{eq:gbar_def}, which holds for all $(i,j)\in\E$. Moreover, we observe from $\mathbf{g}_{ij,k}$ in~\eqref{eq:g_ijl} that~\eqref{eq:gij_bounds} holds for all $k$.
\section{Derivation of Euclidean Hessian Bounds}\label{app:ehess_bounds}
In this appendix, we derive bounds for the Euclidean Hessian tensors derived in Appendix~\ref{app:ehess_tensors} that are necessary for the proof of Lipschitz continuity of the Riemannian gradient in Appendix~\ref{app:lipschitz}. Specifically, we will show that
\begin{equation}
    \left\Vert \frac{\partial\mathcal{A}_{ij}}{\partial x_{i,l}}\right\Vert _{F}^{2}, \left\Vert\frac{\partial\mathcal{A}_{ij}}{\partial x_{j,l}}\right\Vert _{F}^{2}, \left\Vert \frac{\partial\mathcal{B}_{ij}}{\partial x_{j,l}}\right\Vert _{F}^{2},\text{~and~}\left\Vert \frac{\partial\mathcal{B}_{ij}}{\partial x_{i,l}}\right\Vert _{F}^{2}
\end{equation}
are bounded for $k=0\ldots3$ and for all $(i,j)\in\E$, given that $\X\in\K$, where $\K$ is a compact subset of $\MN$.. We begin by providing preliminary derivations that will serve as a reference for the subsequent analysis in this appendix.
\subsection{Preliminaries}\label{app:ehess_prelim_bounds}
The function $f_{2}\left(\phi\right)$ given by equation~\eqref{eq:f_2_def} takes on values within the range $\left(-\frac{\pi}{2},\frac{\pi}{2}\right]$ over the domain $\phi\in\left(-\frac{\pi}{2},\frac{\pi}{2}\right]$, so it holds that
\begin{equation}\label{eq:f2_bound}
    \left|f_{2}\left(\phi\left(\mathbf{x}\right)\right)\right|\leq\left|f_{2}\left(\nicefrac{\pi}{2}\right)\right|=\frac{\pi}{2}
\end{equation}
for all $\mathbf{x}\in\M$. Using the techniques from Appendix~\ref{app:egrad_bounds_prelim}, we now compute the following quantities in trigonometric form.
\begin{align}
    \mu_{i}r_{0}+\eta_{i}r_{1}&=\cos\left(\phi_{j}-\phi_{z}\right)\cos\left(\phi_{r}\right)+\sin\left(\phi_{j}-\phi_{z}\right)\sin\left(\phi_{r}\right)=\cos\left(\phi_{j}-\phi_{z}-\phi_{r}\right), \label{eq:mui_r0_plus_etai_r1} \\
    \omega_{i}r_{0}+\kappa_{i}r_{1}&=\sin\left(\phi_{j}-\phi_{z}\right)\cos\left(\phi_{r}\right)-\cos\left(\phi_{j}-\phi_{z}\right)\sin\left(\phi_{r}\right)=\sin\left(\phi_{j}-\phi_{z}-\phi_{r}\right), \label{eq:omegai_r0_plus_kappai_r1}\\
    \mu_{j}r_{0}+\eta_{j}r_{1}&=\cos\left(\phi_{i}+\phi_{z}\right)\cos\left(\phi_{r}\right)-\sin\left(\phi_{i}+\phi_{z}\right)\sin\left(\phi_{r}\right)=\cos\left(\phi_{i}+\phi_{z}+\phi_{r}\right), \label{eq:muj_r0_plus_etaj_r1}\\
    \omega_{j}r_{0}+\kappa_{j}r_{1}&=\sin\left(\phi_{i}+\phi_{z}\right)\cos\left(\phi_{r}\right)+\cos\left(\phi_{i}+\phi_{z}\right)\sin\left(\phi_{r}\right)=\sin\left(\phi_{i}+\phi_{z}+\phi_{r}\right). \label{eq:omegaj_r0_plus_kappaj_r1}
\end{align}
\noeqref{eq:omegai_r0_plus_kappai_r1,eq:muj_r0_plus_etaj_r1}It follows from~\eqref{eq:mui_r0_plus_etai_r1}-\eqref{eq:omegaj_r0_plus_kappaj_r1} that
\begin{equation}\label{eq:mu_eta_omega_kappa_r01_bounds}
    \left|\mu_{i}r_{0}+\eta_{i}r_{1}\right|,\left|\omega_{i}r_{0}+\kappa_{i}r_{1}\right|,\left|\mu_{j}r_{0}+\eta_{j}r_{1}\right|,\left|\omega_{j}r_{0}+\kappa_{j}r_{1}\right|	\leq1.
\end{equation}
From equations~\eqref{eq:zij_rij_0_1},~\eqref{eq:etai_trig_form}, and~\eqref{eq:mui_r0_plus_etai_r1}-\eqref{eq:omegaj_r0_plus_kappaj_r1} we apply angle sum and difference identities to compute
\begin{align}
    \eta_{i}-r_{1}\left(\mu_{i}r_{0}+\eta_{i}r_{1}\right)&=\cos\left(\phi_{j}-\phi_{z}\right)\sin\left(\phi_{j}-\phi_{z}-\phi_{r}\right),\\
    \eta_{i}-2r_{1}\left(\mu_{i}r_{0}+\eta_{i}r_{1}\right)&=\sin\left(\phi_{j}-\phi_{z}-2\phi_{r}\right), \\
    \kappa_{i}-2r_{1}\left(\omega_{i}r_{0}+\kappa_{i}r_{1}\right)&=\cos\left(\phi_{j}-\phi_{z}-2\phi_{r}\right), \\
    \eta_{j}-r_{1}\left(\mu_{j}r_{0}+\eta_{j}r_{1}\right)&=\cos\left(\phi_{r}\right)\sin\left(\phi_{i}+\phi_{z}+\phi_{r}\right), \\
    \eta_{j}-2r_{1}\left(\mu_{j}r_{0}+\eta_{j}r_{1}\right)&=-\sin\left(\phi_{i}+\phi_{z}+2\phi_{r}\right), \\
    \kappa_{j}-2r_{1}\left(\omega_{j}r_{0}+\kappa_{j}r_{1}\right)&=\cos\left(\phi_{i}+\phi_{z}+2\phi_{r}\right), 
\end{align}
from which it follows that
\begin{equation}\label{eq:hess_trig_bound_1}
    \left|\eta_{i}-r_{1}\left(\mu_{i}r_{0}+\eta_{i}r_{1}\right)\right|,\left|\eta_{i}-2r_{1}\left(\mu_{i}r_{0}+\eta_{i}r_{1}\right)\right|,\left|\kappa_{i}-2r_{1}\left(\omega_{i}r_{0}+\kappa_{i}r_{1}\right)\right|	\leq1
\end{equation}
and
\begin{equation}\label{eq:hess_trig_bound_2}
    \left|\eta_{j}-r_{1}\left(\mu_{j}r_{0}+\eta_{j}r_{1}\right)\right|,\left|\eta_{j}-2r_{1}\left(\mu_{j}r_{0}+\eta_{j}r_{1}\right)\right|,\left|\kappa_{j}-2r_{1}\left(\omega_{j}r_{0}+\kappa_{j}r_{1}\right)\right|	\leq1.
\end{equation}
We also compute
\begin{align}
    \mu_{j}\eta_{i}-\eta_{j}\mu_{i}&=\cos\left(\phi_{i}+\phi_{z}\right)\sin\left(\phi_{j}-\phi_{z}\right)+\sin\left(\phi_{i}+\phi_{z}\right)\cos\left(\phi_{j}-\phi_{z}\right)=\sin\left(\phi_{i}+\phi_{j}\right), \\
    \mu_{j}\kappa_{i}-\eta_{j}\omega_{i}&=-\cos\left(\phi_{i}+\phi_{z}\right)\cos\left(\phi_{j}-\phi_{z}\right)+\sin\left(\phi_{i}+\phi_{z}\right)\sin\left(\phi_{j}-\phi_{z}\right)=-\cos\left(\phi_{i}+\phi_{j}\right),
\end{align}
from which it follows that
\begin{equation}\label{eq:mu_eta_bounds}
    \left|\mu_{j}\eta_{i}-\eta_{j}\mu_{i}\right|,\left|\mu_{j}\kappa_{i}-\eta_{j}\omega_{i}\right|	\leq1.
\end{equation}
Furthermore, it holds that
\begin{align}
    -z_{0}r_{0}+z_{1}r_{1}&=-\cos\left(\phi_{z}\right)\cos\left(\phi_{r}\right)+\sin\left(\phi_{z}\right)\sin\left(\phi_{r}\right)=-\cos\left(\phi_{z}+\phi_{r}\right)\leq1, \\
    z_{1}r_{0}+z_{0}r_{1}&=\cos\left(\phi_{z}\right)\cos\left(\phi_{r}\right)+\sin\left(\phi_{z}\right)\sin\left(\phi_{r}\right)=\cos\left(\phi_{z}-\phi_{r}\right)\leq1,
\end{align}
and, therefore,
\begin{equation}\label{eq:hess_zr_bounds}
    \left|-z_{0}r_{0}+z_{1}r_{1}\right|,\left|z_{1}r_{0}+z_{0}r_{1}\right|\leq1.
\end{equation}
Finally, we have the trigonometric bounds
\begin{equation}\label{eq:hess_trig_bound_3}
    \left|\kappa_{i}\left(\eta_{i}r_{0}-\mu_{i}r_{1}\right)+\left(\eta_{i}-2r_{1}\left(\mu_{i}r_{0}+\eta_{i}r_{1}\right)\right)\left(\kappa_{i}r_{0}-\omega_{i}r_{1}\right)-r_{1}\right|=\left|\cos\left(\phi_{ij}\right)\sin\left(2\left(\phi_{j}-\phi_{z}-\phi_{ij}\right)\right)\right|\leq1
\end{equation}
and
\begin{equation}\label{eq:hess_trig_bound_4}
    \left|\eta_{i}\left(\kappa_{i}r_{0}-\omega_{i}r_{1}\right)+\left(\kappa_{i}-2r_{1}\left(\omega_{i}r_{0}+\kappa_{i}r_{1}\right)\right)\left(\eta_{i}r_{0}-\mu_{i}r_{1}\right)+r_{1}\right|	=\left|-\cos\left(\phi_{ij}\right)\sin\left(2\left(\phi_{j}-\phi_{z}-\phi_{ij}\right)\right)\right|\leq1,
\end{equation}
which concludes our derivation of preliminary bounds for the Euclidean Hessian Tensors.
\subsection{$\mathcal{A}_{ij}$ Tensor Bounds}\label{app:Aij_tensor_bounds}
We first derive bounds for $\left\Vert \frac{\partial\mathcal{A}_{ij}}{\partial x_{i,l}}\right\Vert _{F}^{2}$ for $k=0\ldots3$, starting with $\left\Vert \frac{\partial\mathcal{A}_{ij}}{\partial x_{i,0}}\right\Vert _{F}$. Applying the triangle inequality to~\eqref{eq:dA11_dxi0} yields
\begin{align}
    \left|\frac{\partial\mathcal{A}_{11}}{\partial x_{i,0}}\right|&=\left|2\left(\eta_{i}-r_{1}\left(\mu_{i}r_{0}+\eta_{i}r_{1}\right)\right)\left(\eta_{i}r_{0}-\mu_{i}r_{1}\right)f_{1}+r_{1}\left(\eta_{i}r_{0}-\mu_{i}r_{1}\right)^{2}f_{2}\right| \\
    &\leq 2\left|\left(\eta_{i}-r_{1}\left(\mu_{i}r_{0}+\eta_{i}r_{1}\right)\right)\right|\left|\eta_{i}r_{0}-\mu_{i}r_{1}\right|\left|f_{1}\right|+\left|r_{1}\right|\left|\eta_{i}r_{0}-\mu_{i}r_{1}\right|^{2}\left|f_{2}\right|\label{eq:dA11_bound_1}
\end{align}
Applying~\eqref{eq:hess_trig_bound_1},~\eqref{eq:dphi_dx_bounds},~\eqref{eq:f1_bound},~\eqref{eq:sin_cos_bounds}, and~\eqref{eq:f2_bound} to~\eqref{eq:dA11_bound_1} and simplifying yields
\begin{equation}\label{eq:dA11_dxi0_bound}
    \left|\frac{\partial\mathcal{A}_{11}}{\partial x_{i,0}}\right|\leq\frac{\pi}{2}.
\end{equation}
Next, applying the triangle inequality to~\eqref{eq:dA12_dxi0} yields 
\begin{align}
    \left|\frac{\partial\mathcal{A}_{12}}{\partial x_{i,0}}\right|=&\left|\left(\kappa_{i}\left(\eta_{i}r_{0}-\mu_{i}r_{1}\right)+\left(\eta_{i}-2r_{1}\left(\mu_{i}r_{0}+\eta_{i}r_{1}\right)\right)\left(\kappa_{i}r_{0}-\omega_{i}r_{1}\right)-r_{1}\right)f_{1}+r_{1}\left(\kappa_{i}r_{0}-\omega_{i}r_{1}\right)\left(\eta_{i}r_{0}-\mu_{i}r_{1}\right)f_{2}\right|\\
     \leq&\left|\left(\kappa_{i}\left(\eta_{i}r_{0}-\mu_{i}r_{1}\right)+\left(\eta_{i}-2r_{1}\left(\mu_{i}r_{0}+\eta_{i}r_{1}\right)\right)\left(\kappa_{i}r_{0}-\omega_{i}r_{1}\right)-r_{1}\right)f_{1}\right|\\
     &+\left|r_{1}\right|\left|\kappa_{i}r_{0}-\omega_{i}r_{1}\right|\left|\eta_{i}r_{0}-\mu_{i}r_{1}\right|\left|f_{2}\right|\label{eq:dA12_bound_1}
\end{align}
To simplify~\eqref{eq:dA12_bound_1}, we observe that
\begin{equation}\label{eq:dA12_trig_bound}
    \left|\kappa_{i}\left(\eta_{i}r_{0}-\mu_{i}r_{1}\right)+\left(\eta_{i}-2r_{1}\left(\mu_{i}r_{0}+\eta_{i}r_{1}\right)\right)\left(\kappa_{i}r_{0}-\omega_{i}r_{1}\right)-r_{1}\right|=\left|\cos\left(\phi_{ij}\right)\sin\left(2\left(\phi_{j}-\phi_{z}-\phi_{ij}\right)\right)\right| \leq1.
\end{equation}
Applying~\eqref{eq:dA12_trig_bound},~\eqref{eq:sin_cos_bounds},~\eqref{eq:dphi_dx_bounds}, and~\eqref{eq:f2_bound} to~\eqref{eq:dA12_bound_1} yields
\begin{equation}\label{eq:dA12_dxi0_bound}
    \left|\frac{\partial\mathcal{A}_{12}}{\partial x_{i,0}}\right|\leq\frac{\pi}{2}+1.
\end{equation}
From~\eqref{eq:dA13_dx_all} and~\eqref{eq:dA14_dx_all}, we have
\begin{equation}\label{eq:dA13_dA14_dxi0_bound}
    \left|\frac{\partial\mathcal{A}_{13}}{\partial x_{i,0}}\right|=\left|\frac{\partial\mathcal{A}_{14}}{\partial x_{i,0}}\right|=0.
\end{equation}
From~\eqref{eq:dA21_dxi0}, we have
\begin{align}
    \left|\frac{\partial\mathcal{A}_{21}}{\partial x_{i,0}}\right|&=\left|2\left(\alpha_{1}-r_{2}\left(\mu_{i}r_{0}+\eta_{i}r_{1}\right)\right)\left(\eta_{i}r_{0}-\mu_{i}r_{1}\right)f_{1}+r_{2}\left(\eta_{i}r_{0}-\mu_{i}r_{1}\right)^{2}f_{2}\right| \\
    &\leq 2\left(\left|\alpha_{1}\right|+\left|r_{2}\right|\left|\mu_{i}r_{0}+\eta_{i}r_{1}\right|\right)\left|\eta_{i}r_{0}-\mu_{i}r_{1}\right|\left|f_{1}\right|+\left|r_{2}\right|\left|\eta_{i}r_{0}-\mu_{i}r_{1}\right|^{2}\left|f_{2}\right|,
\end{align}
and further simplifying with~\eqref{eq:dphi_dx_bounds} and~\eqref{eq:f1_bound} yields
\begin{equation}\label{eq:dA21_bound_1}
    \left|\frac{\partial\mathcal{A}_{21}}{\partial x_{i,0}}\right|\leq2\left|\alpha_{1}\right|+\left(\left|f_{2}\right|+2\right)\left|r_{2}\right|.
\end{equation}
We now define
\begin{equation}\label{eq:tau_1_def}
    \bar{\tau}_{1}\triangleq2\left(\bar{\mathbf{z}}_{23}+\sqrt{2}\bar{\mathbf{t}}_{\mathbf{x}}\right)+\sqrt{2}\left(\frac{\pi}{2}+2\right)\bar{\mathbf{t}}_{\mathbf{r}}.
\end{equation}
Applying~\eqref{eq:alpha1_bound},~\eqref{eq:f2_bound} and~\eqref{eq:rij_23_bound} to~\eqref{eq:dA21_bound_1} and simplifying yields
\begin{equation}\label{eq:dA21_dxi0_bound}
    \left|\frac{\partial\mathcal{A}_{21}}{\partial x_{i,0}}\right|\leq2\left(\bar{\mathbf{z}}_{23}+\sqrt{2}\bar{\mathbf{t}}_{\mathbf{x}}\right)+\sqrt{2}\left(\frac{\pi}{2}+2\right)\bar{\mathbf{t}}_{\mathbf{r}}=\bar{\tau}_{1}.
\end{equation}
Applying the triangle inequality to~\eqref{eq:dA22_dxi0} yields
\begin{align}
    \left|\frac{\partial\mathcal{A}_{22}}{\partial x_{i,0}}\right|&=\left|\left(\beta_{1}\left(\eta_{i}r_{0}-\mu_{i}r_{1}\right)+\left(\alpha_{1}-2r_{2}\left(\mu_{i}r_{0}+\eta_{i}r_{1}\right)\right)\left(\kappa_{i}r_{0}-\omega_{i}r_{1}\right)-r_{2}\right)f_{1}+r_{2}\left(\kappa_{i}r_{0}-\omega_{i}r_{1}\right)\left(\eta_{i}r_{0}-\mu_{i}r_{1}\right)f_{2}\right|\\
    &\leq \left(\left|\beta_{1}\right|\left|\eta_{i}r_{0}-\mu_{i}r_{1}\right|+\left(\left|\alpha_{1}\right|+2\left|r_{2}\right|\left|\mu_{i}r_{0}+\eta_{i}r_{1}\right|\right)\left|\kappa_{i}r_{0}-\omega_{i}r_{1}\right|+\left|r_{2}\right|\right)\left|f_{1}\right|+\left|r_{2}\right|\left|\kappa_{i}r_{0}-\omega_{i}r_{1}\right|\left|\eta_{i}r_{0}-\mu_{i}r_{1}\right|\left|f_{2}\right|,
\end{align}
and simplifying with~\eqref{eq:dphi_dx_bounds} and~\eqref{eq:f1_bound} yields
\begin{equation}\label{eq:dA22_bound_1}
    \left|\frac{\partial\mathcal{A}_{22}}{\partial x_{i,0}}\right|\leq\left|\beta_{1}\right|+\left|\alpha_{1}\right|+\left(\left|f_{2}\right|+3\right)\left|r_{2}\right|
\end{equation}
Next, we define
\begin{equation}\label{eq:tau_2_def}
    \bar{\tau}_{2}	\triangleq2\left(\bar{\mathbf{z}}_{23}+\sqrt{2}\bar{\mathbf{t}}_{\mathbf{x}}\right)+\sqrt{2}\left(\frac{\pi}{2}+3\right)\bar{\mathbf{t}}_{\mathbf{r}}.
\end{equation}
Applying~\eqref{eq:alpha1_bound},~\eqref{eq:f2_bound} and~\eqref{eq:rij_23_bound} to~\eqref{eq:dA22_bound_1} and simplifying yields
\begin{equation}\label{eq:dA22_dxi0_bound}
    \left|\frac{\partial\mathcal{A}_{22}}{\partial x_{i,0}}\right|\leq2\left(\bar{\mathbf{z}}_{23}+\sqrt{2}\bar{\mathbf{t}}_{\mathbf{x}}\right)+\sqrt{2}\left(\frac{\pi}{2}+3\right)\bar{\mathbf{t}}_{\mathbf{r}}=\bar{\tau}_{2}.
\end{equation}
Applying the triangle inequality and~\eqref{eq:mui-zeta1_bounds},~\eqref{eq:dphi_dx_bounds},and~\eqref{eq:f1_bound} to~\eqref{eq:dA23_dxi0} and~\eqref{eq:dA24_dxi0} yields
\begin{equation}\label{eq:dA23_dA24_dxi0_bound}
    \left|\frac{\partial\mathcal{A}_{23}}{\partial x_{i,0}}\right|,\left|\frac{\partial\mathcal{A}_{24}}{\partial x_{i,0}}\right|\leq1
\end{equation}
Since $\mathcal{A}_{21}$ and $\mathcal{A}_{31}$ have similar structure, applying the derivation for~\eqref{eq:dA21_dxi0_bound} to~\eqref{eq:dA31_dxi0} yields
\begin{equation}\label{eq:dA31_dxi0_bound}
    \left|\frac{\partial\mathcal{A}_{31}}{\partial x_{i,0}}\right|\leq\bar{\tau}_{1}.
\end{equation}
Similarly, $\mathcal{A}_{22}$ and $\mathcal{A}_{32}$ have similar structure, so applying the derivation for~\eqref{eq:dA22_dxi0_bound} to~\eqref{eq:dA32_dxi0} yields
\begin{equation}\label{eq:dA32_dxi0_bound}
    \left|\frac{\partial\mathcal{A}_{32}}{\partial x_{i,0}}\right|\leq\bar{\tau}_{2}.
\end{equation}
From~\eqref{eq:dA33_dxi0},~\eqref{eq:dA34_dxi0}, and~\eqref{eq:dA23_dA24_dxi0_bound}, it holds that
\begin{align}
    \left|\frac{\partial\mathcal{A}_{33}}{\partial x_{i,0}}\right| & =\left|-\frac{\partial\mathcal{A}_{24}}{\partial x_{i,0}}\right|\leq1\label{eq:dA33_dxi0_bound}\\
    \left|\frac{\partial\mathcal{A}_{34}}{\partial x_{i,0}}\right| & =\left|\frac{\partial\mathcal{A}_{23}}{\partial x_{i,0}}\right|\leq1.\label{eq:dA34_dxi0_bound}
\end{align}
\noeqref{eq:dA31_dxi0_bound,eq:dA32_dxi0_bound,eq:dA33_dxi0_bound}Finally, substituting~\eqref{eq:dA11_dxi0_bound},~\eqref{eq:dA12_dxi0_bound}-\eqref{eq:dA13_dA14_dxi0_bound},~\eqref{eq:dA21_dxi0_bound},~\eqref{eq:dA22_dxi0_bound}, and~\eqref{eq:dA23_dA24_dxi0_bound}-\eqref{eq:dA34_dxi0_bound} into the Frobenius norm definition from~\eqref{eq:frob_norm_def} yields
\begin{equation}\label{eq:dAij_dxi0_bound}
    \left\Vert \frac{\partial\mathcal{A}_{ij}}{\partial x_{i,0}}\right\Vert _{F}^{2}\leq2\left(\bar{\tau}_{1}^{2}+\bar{\tau}_{2}^{2}\right)+\left(\frac{\pi}{2}\right)^{2}+\left(\frac{\pi}{2}+1\right)^{2}+4,
\end{equation}
which holds for all $(i,j)\in\E$. We now address $\left\Vert \frac{\partial\mathcal{A}_{ij}}{\partial x_{i,1}}\right\Vert _{F}$. Applying the triangle inequality to~\eqref{eq:dA11_dxi1} gives
\begin{align}
    \left|\frac{\partial\mathcal{A}_{11}}{\partial x_{i,1}}\right|  \leq&\left|\left(\eta_{i}\left(\kappa_{i}r_{0}-\omega_{i}r_{1}\right)+\left(\kappa_{i}-2r_{1}\left(\omega_{i}r_{0}+\kappa_{i}r_{1}\right)\right)\left(\eta_{i}r_{0}-\mu_{i}r_{1}\right)+r_{1}\right)f_{1}+r_{1}\left(\eta_{i}r_{0}-\mu_{i}r_{1}\right)\left(\kappa_{i}r_{0}-\omega_{i}r_{1}\right)f_{2}\right| \\
    \leq&\left|\eta_{i}\left(\kappa_{i}r_{0}-\omega_{i}r_{1}\right)+\left(\kappa_{i}-2r_{1}\left(\omega_{i}r_{0}+\kappa_{i}r_{1}\right)\right)\left(\eta_{i}r_{0}-\mu_{i}r_{1}\right)+r_{1}\right|\left|f_{1}\right|\\
    &+\left|r_{1}\right|\left|\eta_{i}r_{0}-\mu_{i}r_{1}\right|\left|\kappa_{i}r_{0}-\omega_{i}r_{1}\right|\left|f_{2}\right|.\label{eq:dA11_dxi1_bound_1}
\end{align}
Noting that
\begin{equation}
    \left|\eta_{i}\left(\kappa_{i}r_{0}-\omega_{i}r_{1}\right)+\left(\kappa_{i}-2r_{1}\left(\omega_{i}r_{0}+\kappa_{i}r_{1}\right)\right)\left(\eta_{i}r_{0}-\mu_{i}r_{1}\right)+r_{1}\right|=\left|-\cos\left(\phi_{ij}\right)\sin\left(2\left(\phi_{j}-\phi_{z}-\phi_{ij}\right)\right)\right|\leq1,
\end{equation}
we see that applying~\eqref{eq:f1_bound},~\eqref{eq:sin_cos_bounds},~\eqref{eq:dphi_dx_bounds} and~\eqref{eq:f2_bound} to~\eqref{eq:dA11_dxi1_bound_1} yields
\begin{equation}\label{eq:dA11_dxi1_bound}
    \left|\frac{\partial\mathcal{A}_{11}}{\partial x_{i,1}}\right|\leq\frac{\pi}{2}+1.
\end{equation}
Moreover, applying the triangle inequality,~\eqref{eq:hess_trig_bound_1},~\eqref{eq:dphi_dx_bounds},~\eqref{eq:f1_bound},~\eqref{eq:sin_cos_bounds}, and~\eqref{eq:f2_bound} to~\eqref{eq:dA12_dxi1} yields
\begin{equation}\label{eq:dA12_dxi1_bound}
    \left|\frac{\partial\mathcal{A}_{12}}{\partial x_{i,1}}\right|=\left|2\left(\kappa_{i}-r_{1}\left(\omega_{i}r_{0}+\kappa_{i}r_{1}\right)\right)\left(\kappa_{i}r_{0}-\omega_{i}r_{1}\right)f_{1}+r_{1}\left(\kappa_{i}r_{0}-\omega_{i}r_{1}\right)^{2}f_{2}\right|\leq\frac{\pi}{2}+2.
\end{equation}
From~\eqref{eq:dA13_dx_all} and~\eqref{eq:dA14_dx_all}, we have
\begin{equation}\label{eq:dA13_dxi1_bound}
    \left|\frac{\partial\mathcal{A}_{13}}{\partial x_{i,1}}\right|=\left|\frac{\partial\mathcal{A}_{14}}{\partial x_{i,1}}\right|=0.
\end{equation}
Applying the derivation for~\eqref{eq:dA22_dxi0_bound} to~\eqref{eq:dA21_dxi1} yields
\begin{equation}\label{eq:dA21_dxi1_bound}
    \left|\frac{\partial\mathcal{A}_{21}}{\partial x_{i,1}}\right|\leq\left|\alpha_{1}\right|+\left|\beta_{1}\right|+\left(\left|f_{2}\right|+3\right)\left|r_{2}\right|\leq\bar{\tau}_{2},
\end{equation}
with $\bar{\tau}_{2}$ given by~\eqref{eq:tau_2_def}, and applying the derivation for~\eqref{eq:dA21_dxi0_bound} to~\eqref{eq:dA21_dxi1} gives
\begin{equation}\label{eq:dA22_dxi1_bound}
    \left|\frac{\partial\mathcal{A}_{22}}{\partial x_{i,1}}\right|\leq2\left|\beta_{1}\right|+\left(\left|f_{2}\right|+2\right)\left|r_{2}\right|\leq\bar{\tau}_{1},
\end{equation}
with $\bar{\tau}_{1}$ given by~\eqref{eq:tau_1_def}. Applying the triangle inequality and~\eqref{eq:mui-zeta1_bounds},~\eqref{eq:dphi_dx_bounds},and~\eqref{eq:f1_bound} to~\eqref{eq:dA23_dxi1} and~\eqref{eq:dA24_dxi1} yields
\begin{equation}\label{eq:dA23_dA24_dxi1_bound}
    \left|\frac{\partial\mathcal{A}_{23}}{\partial x_{i,1}}\right|,\left|\frac{\partial\mathcal{A}_{24}}{\partial x_{i,1}}\right|\leq1,
\end{equation}
Now, applying the derivation for~\eqref{eq:dA22_dxi0_bound} to~\eqref{eq:dA31_dxi1} yields
\begin{equation}\label{eq:dA31_dxi1_bound}
    \left|\frac{\partial\mathcal{A}_{31}}{\partial x_{i,1}}\right|\leq\left|\alpha_{2}\right|+\left|\beta_{2}\right|+\left(\left|f_{2}\right|+3\right)\left|r_{3}\right|\leq\bar{\tau}_{2},
\end{equation}
with $\bar{\tau}_{2}$ given by~\eqref{eq:tau_2_def}, and applying the derivation for~\eqref{eq:dA21_dxi0_bound} to~\eqref{eq:dA31_dxi1} gives
\begin{equation}\label{eq:dA32_dxi1_bound}
    \left|\frac{\partial\mathcal{A}_{32}}{\partial x_{i,1}}\right|\leq2\left|\beta_{2}\right|+\left(\left|f_{2}\right|+2\right)\left|r_{3}\right|\leq\bar{\tau}_{1}
\end{equation}
with $\bar{\tau}_{1}$ given by~\eqref{eq:tau_1_def}. Finally, from~\eqref{eq:dA33_dxi1},~\eqref{eq:dA34_dxi1}, and~\eqref{eq:dA23_dA24_dxi1_bound}, we have
\begin{align}
    \left|\frac{\partial\mathcal{A}_{33}}{\partial x_{i,1}}\right| & =\left|-\frac{\partial\mathcal{A}_{24}}{\partial x_{i,1}}\right|\leq1\label{eq:dA33_dxi1_bound}\\
    \left|\frac{\partial\mathcal{A}_{34}}{\partial x_{i,1}}\right| & =\left|\frac{\partial\mathcal{A}_{23}}{\partial x_{i,1}}\right|\leq1\label{eq:dA34_dxi1_bound}
\end{align}
Finally, substituting~\eqref{eq:dA11_dxi1_bound}-\eqref{eq:dA34_dxi1_bound} into~\eqref{eq:frob_norm_def} yields
\begin{equation}\label{eq:dAij_dxi1_bound}
    \left\Vert \frac{\partial\mathcal{A}_{ij}}{\partial x_{i,1}}\right\Vert_{F}^{2}\leq2\left(\bar{\tau}_{1}^{2}+\bar{\tau}_{2}^{2}\right)+\left(\frac{\pi}{2}+1\right)^{2}+\left(\frac{\pi}{2}+2\right)^{2}+4,
\end{equation}
which holds for all $(i,j)\in\E$. We now address $\left\Vert \frac{\partial\mathcal{A}_{ij}}{\partial x_{i,2}}\right\Vert _{F}$. From~\eqref{eq:dA11_dxi2_xi3},~\eqref{eq:dA12_dxi2_xi3},~\eqref{eq:dA13_dx_all},~\eqref{eq:dA14_dx_all},~\eqref{eq:dA23_dxi2_xi3},~\eqref{eq:dA24_dxi2_xi3},~\eqref{eq:dA33_dxi2_xi3} and~\eqref{eq:dA34_dxi2_xi3}, we have
\begin{equation}\label{eq:dAij_dxi2_bound_1}
    \left|\frac{\partial\mathcal{A}_{11}}{\partial x_{i,2}}\right|=\left|\frac{\partial\mathcal{A}_{12}}{\partial x_{i,2}}\right|=\left|\frac{\partial\mathcal{A}_{13}}{\partial x_{i,2}}\right|=\left|\frac{\partial\mathcal{A}_{14}}{\partial x_{i,2}}\right|=\left|\frac{\partial\mathcal{A}_{23}}{\partial x_{i2}}\right|=\left|\frac{\partial\mathcal{A}_{24}}{\partial x_{i2}}\right|=\left|\frac{\partial\mathcal{A}_{33}}{\partial x_{i,2}}\right|=\left|\frac{\partial\mathcal{A}_{34}}{\partial x_{i,2}}\right|=0
\end{equation}
Moreover, applying the triangle inequality and~\eqref{eq:mui-zeta1_bounds},~\eqref{eq:dphi_dx_bounds}, and~\eqref{eq:f1_bound} to~\eqref{eq:dA21_dxi2},~\eqref{eq:dA22_dxi2},~\eqref{eq:dA31_dxi2}, and~\eqref{eq:dA32_dxi2} yields
\begin{equation}\label{eq:dAij_dxi2_bound_2}
    \left|\frac{\partial\mathcal{A}_{21}}{\partial x_{i,2}}\right|,\left|\frac{\partial\mathcal{A}_{22}}{\partial x_{i,2}}\right|,\left|\frac{\partial\mathcal{A}_{31}}{\partial x_{i,2}}\right|,\left|\frac{\partial\mathcal{A}_{32}}{\partial x_{i,2}}\right|\leq1,
\end{equation}
Substituting~\eqref{eq:dAij_dxi2_bound_1}-\eqref{eq:dAij_dxi2_bound_2} into~\eqref{eq:frob_norm_def} yields
\begin{equation}\label{eq:dAij_dxi2_bound}
    \left\Vert \frac{\partial\mathcal{A}_{ij}}{\partial x_{i,2}}\right\Vert _{F}^{2}\leq4,
\end{equation}
which holds for all $(i,j)\in\E$. To address $\left\Vert \frac{\partial\mathcal{A}_{ij}}{\partial x_{i,3}}\right\Vert _{F}$, we first observe from~\eqref{eq:dA11_dxi2_xi3},~\eqref{eq:dA12_dxi2_xi3},~\eqref{eq:dA13_dx_all},~\eqref{eq:dA14_dx_all},~\eqref{eq:dA23_dxi2_xi3},~\eqref{eq:dA24_dxi2_xi3},~\eqref{eq:dA33_dxi2_xi3}, and~\eqref{eq:dA34_dxi2_xi3} that
\begin{equation}\label{eq:dAij_dxi3_bound_1}
    \left|\frac{\partial\mathcal{A}_{11}}{\partial x_{i,3}}\right|=\left|\frac{\partial\mathcal{A}_{12}}{\partial x_{i,3}}\right|=\left|\frac{\partial\mathcal{A}_{13}}{\partial x_{i,3}}\right|=\left|\frac{\partial\mathcal{A}_{14}}{\partial x_{i,3}}\right|=\left|\frac{\partial\mathcal{A}_{23}}{\partial x_{i,3}}\right|=\left|\frac{\partial\mathcal{A}_{24}}{\partial x_{i,3}}\right|=\left|\frac{\partial\mathcal{A}_{33}}{\partial x_{i,3}}\right|=\left|\frac{\partial\mathcal{A}_{34}}{\partial x_{i,3}}\right|=0,
\end{equation}
Furthermore, applying the triangle inequality and~\eqref{eq:mui-zeta1_bounds},~\eqref{eq:dphi_dx_bounds}, and~\eqref{eq:f1_bound} to~\eqref{eq:dA21_dxi3},~\eqref{eq:dA22_dxi3},~\eqref{eq:dA31_dxi3}, and~\eqref{eq:dA32_dxi3} gives
\begin{equation}\label{eq:dAij_dxi3_bound_2}
    \left|\frac{\partial\mathcal{A}_{21}}{\partial x_{i,3}}\right|,\left|\frac{\partial\mathcal{A}_{22}}{\partial x_{i,3}}\right|,\left|\frac{\partial\mathcal{A}_{31}}{\partial x_{i,3}}\right|,\left|\frac{\partial\mathcal{A}_{32}}{\partial x_{i,3}}\right| \leq1.
\end{equation}
Finally, substituting~\eqref{eq:dAij_dxi3_bound_1}-\eqref{eq:dAij_dxi3_bound_2} into~\eqref{eq:frob_norm_def} yields
\begin{equation}\label{eq:dAij_dxi3_bound}
    \left\Vert\frac{\partial\mathcal{A}_{ij}}{\partial x_{i,3}}\right\Vert _{F}^{2}\leq4,
\end{equation}
which holds for all $(i,j)\in\E$. We now derive bounds for $\left\Vert \frac{\partial\mathcal{A}_{ij}}{\partial x_{j,l}}\right\Vert _{F}^{2}$ for $k=0\ldots3$, starting with $\left\Vert \frac{\partial\mathcal{A}_{ij}}{\partial x_{j,0}}\right\Vert _{F}$. First, applying the triangle inequality to~\eqref{eq:dA11_dxj0} yields
\begin{align}
    \left|\frac{\partial\mathcal{A}_{11}}{\partial x_{j,0}}\right|\leq&\left|-\frac{z_{1}}{\gamma}\right|+\left(\left|\eta_{i}\right|\left|\eta_{j}r_{0}-\mu_{j}r_{1}\right|+\left|\eta_{j}-2r_{1}\left(\mu_{j}r_{0}+\eta_{j}r_{1}\right)\right|\left|\eta_{i}r_{0}-\mu_{i}r_{1}\right|+\left|r_{1}\right|\left(\left|\mu_{j}\eta_{i}-\eta_{j}\mu_{i}\right|+\left|z_{1}r_{0}-z_{0}r_{1}\right|\right)\right)\left|f_{1}\right| \\
    &+\left|r_{1}\right|\left|\eta_{i}r_{0}-\mu_{i}r_{1}\right|\left|\eta_{j}r_{0}-\mu_{j}r_{1}\right|\left|f_{2}\right|.\label{eq:dA11_dxj0_bound_1}
\end{align}
To simplify~\eqref{eq:dA11_dxj0_bound_1}, we apply~\eqref{eq:sin_cos_bounds},~\eqref{eq:gamma_recip_bound},~\eqref{eq:mui-zeta1_bounds},~\eqref{eq:dphi_dx_bounds},~\eqref{eq:hess_trig_bound_2},~\eqref{eq:mu_eta_bounds},~\eqref{eq:hess_zr_bounds}, and~\eqref{eq:f1_bound} to obtain
\begin{equation}\label{eq:dA11_dxj0_bound}
    \left|\frac{\partial\mathcal{A}_{11}}{\partial x_{j,0}}\right|\leq\pi+4.
\end{equation}
Similarly, applying the same process to~\eqref{eq:dA12_dxj0} yields
\begin{equation}\label{eq:dA12_dxj0_bound}
    \left|\frac{\partial\mathcal{A}_{12}}{\partial x_{j,0}}\right|\leq\pi+4.
\end{equation}
From~\eqref{eq:dA13_dx_all} and~\eqref{eq:dA14_dx_all}, we have
\begin{equation}\label{eq:dA13-dA14_dxj0_bound}
    \left|\frac{\partial\mathcal{A}_{13}}{\partial x_{j,0}}\right|=\left|\frac{\partial\mathcal{A}_{14}}{\partial x_{j,0}}\right|=0
\end{equation}
Next, applying the triangle inequality to~\eqref{eq:dA21_dxj0} yields
\begin{align}
    \left|\frac{\partial\mathcal{A}_{21}}{\partial x_{j,0}}\right|\leq&\left|\frac{z_{2}}{\gamma}\right|+\left(\left|\alpha_{1}\right|\left|\eta_{j}r_{0}-\mu_{j}r_{1}\right|+(\left|\alpha_{2}\right|+2\left|r_{2}\right|\left|\mu_{j}r_{0}+\eta_{j}r_{1}\right|)\left|\eta_{i}r_{0}-\mu_{i}r_{1}\right|+\left|r_{2}\right|\left(\left|\mu_{j}\eta_{i}-\eta_{j}\mu_{i}\right|+\left|z_{1}r_{0}+z_{0}r_{1}\right|\right)\right)\left|f_{1}\right|\\
    &+\left|r_{2}\right|\left|\eta_{i}r_{0}-\mu_{i}r_{1}\right|\left|\eta_{j}r_{0}-\mu_{j}r_{1}\right|\left|f_{2}\right|\label{eq:dA21_dxj0_bound_1}
\end{align}
Now, we define
\begin{equation}\label{eq:tau_3_def}
    \bar{\tau}_{3}\triangleq\frac{\pi}{2}\bar{\mathbf{z}}_{2}+2\left(\bar{\mathbf{z}}_{23}+\sqrt{2}\bar{\mathbf{t}}_{\mathbf{x}}\right)+\sqrt{2}\left(\frac{\pi}{2}+4\right)\bar{\mathbf{t}}_{\mathbf{r}}.
\end{equation}
Applying~\eqref{eq:z_2_3_23_bound},~\eqref{eq:alpha1_bound},~\eqref{eq:mu_eta_omega_kappa_r01_bounds},~\eqref{eq:rij_23_bound},~\eqref{eq:dphi_dx_bounds},~\eqref{eq:mu_eta_bounds},~\eqref{eq:hess_zr_bounds},~\eqref{eq:f1_bound}, and~\eqref{eq:f2_bound} to~\eqref{eq:dA21_dxj0_bound_1} yields
\begin{equation}\label{eq:dA21_dxj0_bound}
    \left|\frac{\partial\mathcal{A}_{21}}{\partial x_{j,0}}\right|\leq\frac{\pi}{2}\bar{\mathbf{z}}_{2}+2\left(\bar{\mathbf{z}}_{23}+\sqrt{2}\bar{\mathbf{t}}_{\mathbf{x}}\right)+\sqrt{2}\left(\frac{\pi}{2}+4\right)\bar{\mathbf{t}}_{\mathbf{r}}=\bar{\tau}_{3}.
\end{equation}
Applying the triangle inequality to~\eqref{eq:dA22_dxj0} gives
\begin{align}
    \left|\frac{\partial\mathcal{A}_{22}}{\partial x_{j,0}}\right|\leq&\left|\frac{z_{3}}{\gamma}\right|+\left(\left|\beta_{1}\right|\left|\eta_{j}r_{0}-\mu_{j}r_{1}\right|+\left(\left|\alpha_{2}\right|+2\left|r_{2}\right|\left|\mu_{j}r_{0}+\eta_{j}r_{1}\right|\right)\left|\kappa_{i}r_{0}-\omega_{i}r_{1}\right|+\left|r_{2}\right|\left(\left|\kappa_{i}\mu_{j}-\omega_{i}\eta_{j}\right|+\left|-z_{0}r_{0}+z_{1}r_{1}\right|\right)\right)\left|f_{1}\right| \\
    &+\left|r_{2}\right|\left|\kappa_{i}r_{0}-\omega_{i}r_{1}\right|\left|\eta_{j}r_{0}-\mu_{j}r_{1}\right|\left|f_{2}\right|.\label{eq:dA22_dxj0_bound_1}
\end{align}
By letting
\begin{equation}\label{eq:tau_4_def}
    \bar{\tau}_{4}=\frac{\pi}{2}\bar{\mathbf{z}}_{3}+2\left(\bar{\mathbf{z}}_{23}+\sqrt{2}\bar{\mathbf{t}}_{\mathbf{x}}\right)+\sqrt{2}\left(\frac{\pi}{2}+4\right)\bar{\mathbf{t}}_{\mathbf{r}},
\end{equation}
we see that applying~\eqref{eq:z_2_3_23_bound},~\eqref{eq:alpha1_bound},~\eqref{eq:mu_eta_omega_kappa_r01_bounds},~\eqref{eq:rij_23_bound},~\eqref{eq:dphi_dx_bounds},~\eqref{eq:mu_eta_bounds},~\eqref{eq:hess_zr_bounds},~\eqref{eq:f1_bound}, and~\eqref{eq:f2_bound} to~\eqref{eq:dA22_dxj0_bound_1} yields
\begin{equation}\label{eq:dA22_dxj0_bound}
    \left|\frac{\partial\mathcal{A}_{22}}{\partial x_{j,0}}\right|\leq\frac{\pi}{2}\bar{\mathbf{z}}_{3}+2\left(\bar{\mathbf{z}}_{23}+\sqrt{2}\bar{\mathbf{t}}_{\mathbf{x}}\right)+\sqrt{2}\left(\frac{\pi}{2}+4\right)\bar{\mathbf{t}}_{\mathbf{r}}=\bar{\tau}_{4}.
\end{equation}
Now, applying the triangle inequality and~\eqref{eq:sin_cos_bounds},~\eqref{eq:mui-zeta1_bounds},~\eqref{eq:dphi_dx_bounds}, and~\eqref{eq:f1_bound} to~\eqref{eq:dA23_dxj0} and~\eqref{eq:dA24_dxj0} gives
\begin{equation}\label{eq:dA23_dA24_dxj0_bound}
    \left|\frac{\partial\mathcal{A}_{23}}{\partial x_{j,0}}\right|,\left|\frac{\partial\mathcal{A}_{24}}{\partial x_{j,0}}\right|\leq\frac{\pi}{2}+1
\end{equation}
Applying the derivations for~\eqref{eq:dA22_dxj0_bound} and~\eqref{eq:dA21_dxj0_bound} to~\eqref{eq:dA31_dxj0} and~\eqref{eq:dA32_dxj0} yields
\begin{align}
    \left|\frac{\partial\mathcal{A}_{31}}{\partial x_{j,0}}\right|&\leq\bar{\tau}_{4},\label{eq:dA31_dxj0_bound} \\
    \left|\frac{\partial\mathcal{A}_{32}}{\partial x_{j,0}}\right|&\leq\bar{\tau}_{3}.\label{eq:dA32_dxj0_bound}
\end{align}
Finally, from~\eqref{eq:dA33_dxj0},~\eqref{eq:dA34_dxj0},~\eqref{eq:dA23_dA24_dxj0_bound}, and~\eqref{eq:dA23_dxj0}, we have
\begin{align}
    \left|\frac{\partial\mathcal{A}_{33}}{\partial x_{j,0}}\right| & =\left|-\frac{\partial\mathcal{A}_{24}}{\partial x_{j,0}}\right|\leq\frac{\pi}{2}+1,\label{eq:dA33_dxj0_bound}\\
    \left|\frac{\partial\mathcal{A}_{34}}{\partial x_{j,0}}\right| & =\left|\frac{\partial\mathcal{A}_{23}}{\partial x_{j,0}}\right|\leq\frac{\pi}{2}+1\label{eq:dA34_dxj0_bound}
\end{align}
\noeqref{eq:dA12_dxj0_bound,eq:dA23_dA24_dxj0_bound,eq:dA31_dxj0_bound,eq:dA32_dxj0_bound,eq:dA33_dxj0_bound}Finally, substituting~\eqref{eq:dA11_dxj0_bound}-\eqref{eq:dA13-dA14_dxj0_bound},~\eqref{eq:dA21_dxj0_bound}, and~\eqref{eq:dA22_dxj0_bound}-\eqref{eq:dA34_dxj0_bound} into~\eqref{eq:frob_norm_def} yields
\begin{equation}\label{eq:dAij_dxj0_bound}
    \left\Vert \frac{\partial\mathcal{A}_{ij}}{\partial x_{j,0}}\right\Vert _{F}^{2}\leq2\left(\bar{\tau}_{3}^{2}+\bar{\tau}_{4}^{2}\right)+4\left(\frac{\pi}{2}+1\right)^{2}+2\left(\pi+4\right)^{2},
\end{equation}
which holds for all $(i,j)\in\E$. To address $\left\Vert \frac{\partial\mathcal{A}_{ij}}{\partial x_{j,1}}\right\Vert _{F}$, we apply the derivations for~\eqref{eq:dA11_dxj0_bound}-\eqref{eq:dA13-dA14_dxj0_bound},~\eqref{eq:dA21_dxj0_bound}, and~\eqref{eq:dA22_dxj0_bound}-\eqref{eq:dA34_dxj0_bound} to~\eqref{eq:dA11_dxj1},~\eqref{eq:dA12_dxj1},~\eqref{eq:dA13_dx_all},~\eqref{eq:dA14_dx_all},~\eqref{eq:dA21_dxj1},~\eqref{eq:dA22_dxj1},~\eqref{eq:dA23_dxj1},~\eqref{eq:dA24_dxj1},~\eqref{eq:dA31_dxj1},~\eqref{eq:dA32_dxj1},~\eqref{eq:dA33_dxj1}, and~\eqref{eq:dA34_dxj1} to compute
\begin{equation}\label{eq:dAij_dxj1_1}
    \left|\frac{\partial\mathcal{A}_{11}}{\partial x_{j,1}}\right|,\left|\frac{\partial\mathcal{A}_{12}}{\partial x_{j,1}}\right|\leq\pi+4,
\end{equation}
\begin{equation}\label{eq:dAij_dxj1_2}
    \left|\frac{\partial\mathcal{A}_{13}}{\partial x_{j,1}}\right|=\left|\frac{\partial\mathcal{A}_{14}}{\partial x_{j,1}}\right|=0,
\end{equation}
\begin{equation}\label{eq:dAij_dxj1_3}
    \left|\frac{\partial\mathcal{A}_{21}}{\partial x_{j,1}}\right|,\left|\frac{\partial\mathcal{A}_{32}}{\partial x_{j,1}}\right|\leq\bar{\tau}_{4},
\end{equation}
\begin{equation}\label{eq:dAij_dxj1_4}
    \left|\frac{\partial\mathcal{A}_{22}}{\partial x_{j,1}}\right|,\left|\frac{\partial\mathcal{A}_{31}}{\partial x_{j,1}}\right|\leq\bar{\tau}_{3},
\end{equation}
and
\begin{equation}\label{eq:dAij_dxj1_5}
    \left|\frac{\partial\mathcal{A}_{23}}{\partial x_{j,1}}\right|,\left|\frac{\partial\mathcal{A}_{24}}{\partial x_{j,1}}\right|,\left|\frac{\partial\mathcal{A}_{33}}{\partial x_{j,1}}\right|,\left|\frac{\partial\mathcal{A}_{34}}{\partial x_{j,1}}\right|\leq\frac{\pi}{2}+1.
\end{equation}
\noeqref{eq:dAij_dxj1_2,eq:dAij_dxj1_3,eq:dAij_dxj1_4}Substituting~\eqref{eq:dAij_dxj1_1}-\eqref{eq:dAij_dxj1_5} into~\eqref{eq:frob_norm_def} yields
\begin{equation}\label{eq:dAij_dxj1_bound}
    \left\Vert \frac{\partial\mathcal{A}_{ij}}{\partial x_{j,1}}\right\Vert _{F}^{2}\leq2\left(\bar{\tau}_{3}^{2}+\bar{\tau}_{4}^{2}\right)+4\left(\frac{\pi}{2}+1\right)^{2}+2\left(\pi+4\right)^{2},
\end{equation}
which holds for all $(i,j)\in\E$. We now address $\left\Vert \frac{\partial\mathcal{A}_{ij}}{\partial x_{j,2}}\right\Vert _{F}$. From~\eqref{eq:dA11_dxj2_xj3},~\eqref{eq:dA12_dxj2_xj3},~\eqref{eq:dA13_dx_all},~\eqref{eq:dA14_dx_all},~\eqref{eq:dA23_dxj2_xj3},~\eqref{eq:dA24_dxj2_xj3},~\eqref{eq:dA33_dxj2_xj3}, and~\eqref{eq:dA34_dxj2_xj3}, we observe that
\begin{equation}\label{eq:dA11-dA34_dxj2_bound}
    \left|\frac{\partial\mathcal{A}_{11}}{\partial x_{j,2}}\right|=\left|\frac{\partial\mathcal{A}_{12}}{\partial x_{j,2}}\right|=\left|\frac{\partial\mathcal{A}_{13}}{\partial x_{j,2}}\right|=\left|\frac{\partial\mathcal{A}_{14}}{\partial x_{j,2}}\right|=\left|\frac{\partial\mathcal{A}_{23}}{\partial x_{j,2}}\right|=\left|\frac{\partial\mathcal{A}_{24}}{\partial x_{j,2}}\right|=\left|\frac{\partial\mathcal{A}_{33}}{\partial x_{j,2}}\right|=\left|\frac{\partial\mathcal{A}_{34}}{\partial x_{j,2}}\right|=0.
\end{equation}
Furthermore, applying the triangle inequality and~\eqref{eq:sin_cos_bounds},~\eqref{eq:gamma_recip_bound},~\eqref{eq:muj-kappaj_bounds},~\eqref{eq:dphi_dx_bounds}, and~\eqref{eq:f1_bound} to~\eqref{eq:dA12_dxj2_xj3},~\eqref{eq:dA22_dxj2},~\eqref{eq:dA31_dxj2}, and~\eqref{eq:dA32_dxj2} gives
\begin{equation}\label{eq:dA21-dA32_dxj2_bound}
    \left|\frac{\partial\mathcal{A}_{21}}{\partial x_{j,2}}\right|,\left|\frac{\partial\mathcal{A}_{22}}{\partial x_{j,2}}\right|,\left|\frac{\partial\mathcal{A}_{31}}{\partial x_{j,2}}\right|,\left|\frac{\partial\mathcal{A}_{32}}{\partial x_{j,2}}\right|\leq\frac{\pi}{2}+1.
\end{equation}
Substituting~\eqref{eq:dA11-dA34_dxj2_bound}-\eqref{eq:dA21-dA32_dxj2_bound} into~\eqref{eq:frob_norm_def} yields
\begin{equation}\label{eq:dAij_dxj2_bound}
    \left\Vert \frac{\partial\mathcal{A}_{ij}}{\partial x_{j,2}}\right\Vert _{F}^{2}\leq4\left(\frac{\pi}{2}+1\right)^{2}.
\end{equation}
which holds for all $(i,j)\in\E$. We now address $\left\Vert \frac{\partial\mathcal{A}_{ij}}{\partial x_{j,3}}\right\Vert _{F}$. From~\eqref{eq:dA11_dxj2_xj3},~\eqref{eq:dA12_dxj2_xj3},~\eqref{eq:dA13_dx_all}, ~\eqref{eq:dA14_dx_all},~\eqref{eq:dA23_dxj2_xj3},~\eqref{eq:dA24_dxj2_xj3},~\eqref{eq:dA33_dxj2_xj3}, and~\eqref{eq:dA34_dxj2_xj3}, it holds that
\begin{equation}\label{eq:dA13-dA34_dxj3_bound}
    \left|\frac{\partial\mathcal{A}_{13}}{\partial x_{j,3}}\right|=\left|\frac{\partial\mathcal{A}_{11}}{\partial x_{j,3}}\right|=\left|\frac{\partial\mathcal{A}_{12}}{\partial x_{j,3}}\right|=\left|\frac{\partial\mathcal{A}_{14}}{\partial x_{j,3}}\right|=\left|\frac{\partial\mathcal{A}_{23}}{\partial x_{j,3}}\right|=\left|\frac{\partial\mathcal{A}_{24}}{\partial x_{j,3}}\right|=\left|\frac{\partial\mathcal{A}_{33}}{\partial x_{j,2}}\right|=\left|\frac{\partial\mathcal{A}_{34}}{\partial x_{j,2}}\right|=0.
\end{equation}
Moreover, applying the triangle inequality and~\eqref{eq:sin_cos_bounds},~\eqref{eq:gamma_recip_bound},~\eqref{eq:muj-kappaj_bounds},~\eqref{eq:dphi_dx_bounds}, and~\eqref{eq:f1_bound} to~\eqref{eq:dA12_dxj2_xj3},~\eqref{eq:dA22_dxj3},~\eqref{eq:dA31_dxj3}, and~\eqref{eq:dA32_dxj3} yields
\begin{equation}\label{eq:dA21-dA32_dxj3_bound}
    \left|\frac{\partial\mathcal{A}_{21}}{\partial x_{j,3}}\right|,\left|\frac{\partial\mathcal{A}_{22}}{\partial x_{j,3}}\right|,\left|\frac{\partial\mathcal{A}_{31}}{\partial x_{j,3}}\right|,\left|\frac{\partial\mathcal{A}_{32}}{\partial x_{j,3}}\right| \leq\frac{\pi}{2}+1,
\end{equation}
and substituting~\eqref{eq:dA13-dA34_dxj3_bound}-\eqref{eq:dA21-dA32_dxj3_bound} into~\eqref{eq:frob_norm_def} yields
\begin{equation}\label{eq:dAij_dxj3_bound}
    \left\Vert \frac{\partial\mathcal{A}_{ij}}{\partial x_{j,3}}\right\Vert _{F}^{2} \leq4\left(\frac{\pi}{2}+1\right)^{2},
\end{equation}
which concludes our derivation of tensor bounds involving $\mathcal{A}_{ij}$.
\subsection{$\mathcal{B}_{ij}$ Tensor Bounds}\label{app:Bij_tensor_bounds}
We now derive bounds for $\left\Vert \frac{\partial\mathcal{B}_{ij}}{\partial x_{i,l}}\right\Vert _{F}^{2}$ and $\left\Vert \frac{\partial\mathcal{B}_{ij}}{\partial x_{j,l}}\right\Vert _{F}^{2}$ for $k=0\ldots3$. Due to symmetries between $\mathcal{A}_{ij}$ and $\mathcal{B}_{ij}$, the derivations for these bounds are identical to those in Appendix~\ref{app:Aij_tensor_bounds}, so we omit them here and summarize our findings. First, following from the derivations of~\eqref{eq:dAij_dxi0_bound} and~\eqref{eq:dAij_dxi1_bound}, we have
\begin{align}
    \left\Vert \frac{\partial\mathcal{B}_{ij}}{\partial x_{j,0}}\right\Vert _{F}^{2}&\leq2\left(\bar{\tau}_{1}^{2}+\bar{\tau}_{2}^{2}\right)+\left(\frac{\pi}{2}\right)^{2}+\left(\frac{\pi}{2}+1\right)^{2}+4,\label{eq:dBij_dxj0_bound}\\
    \left\Vert \frac{\partial\mathcal{B}_{ij}}{\partial x_{j,1}}\right\Vert _{F}^{2}&\leq2\left(\bar{\tau}_{1}^{2}+\bar{\tau}_{2}^{2}\right)+\left(\frac{\pi}{2}+1\right)^{2}+\left(\frac{\pi}{2}+2\right)^{2}+4.\label{eq:dBij_dxj1_bound}
\end{align}
Next, following from the derivations of~\eqref{eq:dAij_dxi2_bound} and~\eqref{eq:dAij_dxi3_bound}, we have
\begin{equation}\label{eq:dBij_dxj2_dxj3}
    \left\Vert \frac{\partial\mathcal{B}_{ij}}{\partial x_{j,2}}\right\Vert _{F}^{2},\left\Vert \frac{\partial\mathcal{B}_{ij}}{\partial x_{j,3}}\right\Vert _{F}^{2}\leq4.
\end{equation}
Furthermore, following from the derivations of~\eqref{eq:dAij_dxj0_bound} and~\eqref{eq:dAij_dxj1_bound}, we have
\begin{equation}\label{eq:dBij_dxi0_dxi1}
    \left\Vert \frac{\partial\mathcal{B}_{ij}}{\partial x_{i,0}}\right\Vert _{F}^{2}, \left\Vert \frac{\partial\mathcal{B}_{ij}}{\partial x_{i,1}}\right\Vert _{F}^{2}\leq2\left(\bar{\tau}_{3}^{2}+\bar{\tau}_{4}^{2}\right)+4\left(\frac{\pi}{2}+1\right)^{2}+2\left(\pi+4\right)^{2}.
\end{equation}
Finally, following from the derivations of~\eqref{eq:dAij_dxj2_bound} and~\eqref{eq:dAij_dxj3_bound}, we have
\begin{equation}\label{eq:dBij_dxi2_dxi3}
    \left\Vert \frac{\partial\mathcal{B}_{ij}}{\partial x_{i,2}}\right\Vert _{F}^{2},\left\Vert \frac{\partial\mathcal{B}_{ij}}{\partial x_{i,3}}\right\Vert _{F}^{2}\leq4\left(\frac{\pi}{2}+1\right)^{2},
\end{equation}
\noeqref{eq:dBij_dxj1_bound}which concludes our derivation of bounds for the Frobenius norms of Euclidean Hessian tensors.
\subsection{Euclidean Hessian Bounds}\label{app:ehess_combined_bounds}
We now utilize the bounds derived in Appendices~\ref{app:e_ij_bounds},~\ref{app:Aij_Bij_bounds},~\ref{app:gij_bounds},~\ref{app:Aij_tensor_bounds}, and~\ref{app:Bij_tensor_bounds} to derive bounds for $\fnorm{\hii}$, $\fnorm{\hij}$, $\fnorm{\hji}$, and $\fnorm{\hjj}$, which appear in the Euclidean Hessian definition in~\eqref{eq:ehess_F}. First, letting $\mathbf{x}_{i}=[x_{i,0},x_{i,1},x_{i,2},x_{i,3}]^{\top}$ we apply the definitions of $\Cii$-$\Cjj$ given in Appendix~\ref{app:ehess} to compute
\begin{align}
    \mathcal{C}_{ii}=\left(\frac{\partial\mathcal{A}_{ij}}{\partial\mathbf{x}_{i}}\right)^{\top}\Omega_{ij}\mathbf{e}_{ij}=\left[\begin{array}{c|c|c|c}
        \left(\frac{\partial\mathcal{A}_{ij}}{\partial x_{i,0}}\right)^{\top}\Omega_{ij}\mathbf{e}_{ij} & \left(\frac{\partial\mathcal{A}_{ij}}{\partial x_{i,1}}\right)^{\top}\Omega_{ij}\mathbf{e}_{ij} & \left(\frac{\partial\mathcal{A}_{ij}}{\partial x_{i,2}}\right)^{\top}\Omega_{ij}\mathbf{e}_{ij} & \left(\frac{\partial\mathcal{A}_{ij}}{\partial x_{i,3}}\right)^{\top}\Omega_{ij}\mathbf{e}_{ij}\end{array}\right],\label{eq:Cii_exp}\\
        \mathcal{C}_{ij}	=\left(\frac{\partial\mathcal{A}_{ij}}{\partial\mathbf{x}_{j}}\right)^{\top}\Omega_{ij}\mathbf{e}_{ij}=\left[\begin{array}{c|c|c|c}
        \left(\frac{\partial\mathcal{A}_{ij}}{\partial x_{j,0}}\right)^{\top}\Omega_{ij}\mathbf{e}_{ij} & \left(\frac{\partial\mathcal{A}_{ij}}{\partial x_{j,1}}\right)^{\top}\Omega_{ij}\mathbf{e}_{ij} & \left(\frac{\partial\mathcal{A}_{ij}}{\partial x_{j,2}}\right)^{\top}\Omega_{ij}\mathbf{e}_{ij} & \left(\frac{\partial\mathcal{A}_{ij}}{\partial x_{j,3}}\right)^{\top}\Omega_{ij}\mathbf{e}_{ij}\end{array}\right],\label{eq:Cij_exp}\\
        \mathcal{C}_{ji}	=\left(\frac{\partial\mathcal{B}_{ij}}{\partial\mathbf{x}_{i}}\right)^{\top}\Omega_{ij}\mathbf{e}_{ij}=\left[\begin{array}{c|c|c|c}
        \left(\frac{\partial\mathcal{B}_{ij}}{\partial x_{i,0}}\right)^{\top}\Omega_{ij}\mathbf{e}_{ij} & \left(\frac{\partial\mathcal{B}_{ij}}{\partial x_{i,1}}\right)^{\top}\Omega_{ij}\mathbf{e}_{ij} & \left(\frac{\partial\mathcal{B}_{ij}}{\partial x_{i,2}}\right)^{\top}\Omega_{ij}\mathbf{e}_{ij} & \left(\frac{\partial\mathcal{B}_{ij}}{\partial x_{i,3}}\right)^{\top}\Omega_{ij}\mathbf{e}_{ij}\end{array}\right],\label{eq:Cji_exp}\\
        \mathcal{C}_{jj}	=\left(\frac{\partial\mathcal{B}_{ij}}{\partial\mathbf{x}_{j}}\right)^{\top}\Omega_{ij}\mathbf{e}_{ij}=\left[\begin{array}{c|c|c|c}
        \left(\frac{\partial\mathcal{B}_{ij}}{\partial x_{j,0}}\right)^{\top}\Omega_{ij}\mathbf{e}_{ij} & \left(\frac{\partial\mathcal{B}_{ij}}{\partial x_{j,1}}\right)^{\top}\Omega_{ij}\mathbf{e}_{ij} & \left(\frac{\partial\mathcal{B}_{ij}}{\partial x_{j,2}}\right)^{\top}\Omega_{ij}\mathbf{e}_{ij} & \left(\frac{\partial\mathcal{B}_{ij}}{\partial x_{j,3}}\right)^{\top}\Omega_{ij}\mathbf{e}_{ij}\end{array}\right].\label{eq:Cjj_exp}
\end{align}
\noeqref{eq:Cji_exp}Next, taking the Frobenius norm of $\mathbf{h}_{ii}$ from~\eqref{eq:h_ii-h_jj}, applying the triangle inequality, then simplifying, yields
\begin{equation}\label{eq:h_ii_frob}
    \left\Vert \mathbf{h}_{ii}\right\Vert _{F}	=\left\Vert \mathcal{C}_{ii}+\mathcal{A}_{ij}^{\top}\Omega_{ij}\mathcal{A}_{ij}\right\Vert _{F}\leq\left\Vert \mathcal{C}_{ii}\right\Vert _{F}+\left\Vert \mathcal{A}_{ij}^{\top}\Omega_{ij}\mathcal{A}_{ij}\right\Vert _{F}\leq\left\Vert \mathcal{C}_{ii}\right\Vert _{F}+\left\Vert \Omega_{ij}\right\Vert _{F}\left\Vert \mathcal{A}_{ij}\right\Vert _{F}^{2}.
\end{equation}
We now take the Frobeinus norm of~\eqref{eq:Cii_exp} and apply the triangle and Cauchy-Schwarz inequalities to obtain
\begin{equation}\label{eq:Cii_frob_bound}
    \left\Vert \mathcal{C}_{ii}\right\Vert _{F}\leq\left(\sum_{l=0}^{3}\left\Vert \left(\frac{\partial\mathcal{A}_{ij}}{\partial x_{i,l}}\right)^{\top}\Omega_{ij}\mathbf{e}_{ij}\right\Vert _{F}^{2}\right)^{\frac{1}{2}}\leq\left(\sum_{l=0}^{3}\left\Vert \frac{\partial\mathcal{A}_{ij}}{\partial x_{i,l}}\right\Vert _{F}^{2}\right)^{\frac{1}{2}}\left\Vert \Omega_{ij}\right\Vert _{F}\left\Vert \mathbf{e}_{ij}\right\Vert _{2}
\end{equation}
Substituting~\eqref{eq:Cii_frob_bound} into~\eqref{eq:h_ii_frob} and simplifying yields
\begin{equation}\label{eq:h_ii_frob_bound}
    \left\Vert \mathbf{h}_{ii}\right\Vert _{F}\leq\left(\sum_{l=0}^{3}\left\Vert \frac{\partial\mathcal{A}_{ij}}{\partial x_{i,l}}\right\Vert _{F}^{2}\right)^{\frac{1}{2}}\left(\left\Vert \mathbf{e}_{ij}\right\Vert _{2}+\left\Vert \mathcal{A}_{ij}\right\Vert _{F}^{2}\right)\left\Vert \Omega_{ij}\right\Vert _{F}.
\end{equation}
Furthermore, applying the derivation of~\eqref{eq:h_ii_frob_bound} to $\left\Vert \mathbf{h}_{ij}\right\Vert _{F}$, $\left\Vert \mathbf{h}_{ji}\right\Vert _{F}$, and $\left\Vert \mathbf{h}_{jj}\right\Vert _{F}$ using~\eqref{eq:Cij_exp}-\eqref{eq:Cjj_exp} yields
\begin{align}
    \left\Vert \mathbf{h}_{ij}\right\Vert _{F}&\leq\left(\sum_{l=0}^{3}\left\Vert \frac{\partial\mathcal{A}_{ij}}{\partial x_{j,l}}\right\Vert _{F}^{2}\right)^{\frac{1}{2}}\left(\left\Vert \mathbf{e}_{ij}\right\Vert _{2}+\left\Vert \mathcal{A}_{ij}\right\Vert _{F}\left\Vert \mathcal{B}_{ij}\right\Vert _{F}\right)\left\Vert \Omega_{ij}\right\Vert _{F},\label{eq:h_ij_frob_bound}\\
    \left\Vert \mathbf{h}_{ji}\right\Vert _{F}&\leq\left(\sum_{l=0}^{3}\left\Vert \frac{\partial\mathcal{B}_{ij}}{\partial x_{i,l}}\right\Vert _{F}^{2}\right)^{\frac{1}{2}}\left(\left\Vert \mathbf{e}_{ij}\right\Vert _{2}+\left\Vert \mathcal{A}_{ij}\right\Vert _{F}\left\Vert \mathcal{B}_{ij}\right\Vert _{F}\right)\left\Vert \Omega_{ij}\right\Vert _{F},\label{eq:h_ji_frob_bound}\\
    \left\Vert \mathbf{h}_{jj}\right\Vert _{F}&\leq\left(\sum_{l=0}^{3}\left\Vert \frac{\partial\mathcal{B}_{ij}}{\partial x_{j,l}}\right\Vert _{F}^{2}\right)^{\frac{1}{2}}\left(\left\Vert \mathbf{e}_{ij}\right\Vert _{2}+\left\Vert \mathcal{B}_{ij}\right\Vert _{F}^{2}\right)\left\Vert \Omega_{ij}\right\Vert _{F}.\label{eq:h_jj_frob_bound}
\end{align}
Now, from~\eqref{eq:dAij_dxi0_bound},~\eqref{eq:dAij_dxi1_bound},~\eqref{eq:dAij_dxi2_bound},~\eqref{eq:dAij_dxi3_bound}, and~\eqref{eq:dBij_dxj0_bound}-\eqref{eq:dBij_dxj2_dxj3}, we have
\begin{equation}\label{eq:hii_hjj_sum_bounds}
    \sum_{l=0}^{3}\left\Vert \frac{\partial\mathcal{A}_{ij}}{\partial x_{i,l}}\right\Vert _{F}^{2},\sum_{l=0}^{3}\left\Vert \frac{\partial\mathcal{B}_{ij}}{\partial x_{j,l}}\right\Vert _{F}^{2}\leq4\left(\bar{\tau}_{1}^{2}+\bar{\tau}_{2}^{2}\right)+\left(\frac{\pi}{2}\right)^{2}+2\left(\frac{\pi}{2}+1\right)^{2}+\left(\frac{\pi}{2}+2\right)^{2}+16,
\end{equation}
with $\bar{\tau}_{1}$ from~\eqref{eq:tau_1_def} and $\bar{\tau}_{2}$ from~\eqref{eq:tau_2_def}. Similarly,~\eqref{eq:dAij_dxj0_bound},~\eqref{eq:dAij_dxj1_bound},~\eqref{eq:dAij_dxj2_bound},~\eqref{eq:dBij_dxi0_dxi1} and~\eqref{eq:dBij_dxi2_dxi3} give
\begin{equation}\label{eq:hij_hji_sum_bounds}
    \sum_{l=0}^{3}\left\Vert \frac{\partial\mathcal{A}_{ij}}{\partial x_{j,l}}\right\Vert _{F}^{2},\sum_{l=0}^{3}\left\Vert \frac{\partial\mathcal{B}_{ij}}{\partial x_{i,l}}\right\Vert _{F}^{2}\leq4\left(\bar{\tau}_{3}^{2}+\bar{\tau}_{4}^{2}\right)+16\left(\frac{\pi}{2}+1\right)^{2}+4\left(\pi+4\right)^{2},
\end{equation}
with $\bar{\tau}_{3}$ from~\eqref{eq:tau_3_def} and $\bar{\tau}_{4}$ from~\eqref{eq:tau_4_def}. To aid in formulating bounds for~\eqref{eq:h_ii_frob_bound}-\eqref{eq:h_jj_frob_bound}, we define
\begin{equation}\label{eq:hiibar_def}
    \overline{\mathbf{h}}_{ii}\triangleq\left(4\left(\bar{\tau}_{1}^{2}+\bar{\tau}_{2}^{2}\right)+\left(\frac{\pi}{2}\right)^{2}+2\left(\frac{\pi}{2}+1\right)^{2}+\left(\frac{\pi}{2}+2\right)^{2}+16\right)^{\frac{1}{2}}\left(\ebar+\JBar^{2}\right)
\end{equation}
and
\begin{equation}\label{eq:hijbar_def}
    \overline{\mathbf{h}}_{ij}\triangleq\left(4\left(\bar{\tau}_{3}^{2}+\bar{\tau}_{4}^{2}\right)+16\left(\frac{\pi}{2}+1\right)^{2}+4\left(\pi+4\right)^{2}\right)^{\frac{1}{2}}\left(\ebar+\JBar^{2}\right),
\end{equation}
with $\ebar$ from~\eqref{eq:eij_bound} and $\JBar$ from~\eqref{eq:JBar_def}. Finally, applying~\eqref{eq:eij_bound},~\eqref{eq:Aij_bound},~\eqref{eq:Bij_bound},~\eqref{eq:hii_hjj_sum_bounds} and~\eqref{eq:hij_hji_sum_bounds} to~\eqref{eq:h_ii_frob_bound}-\eqref{eq:h_jj_frob_bound} yields
\begin{align}
    \fnorm{\hii},\fnorm{\hjj}&\leq\overline{\mathbf{h}}_{ii}\fnorm{\Omegaij}\label{eq:hii_hjj_bound}\\
    \fnorm{\hij},\fnorm{\hji}&\leq\overline{\mathbf{h}}_{ij}\fnorm{\Omegaij}\label{eq:hij_hji_bound},
\end{align}
with $\overline{\mathbf{h}}_{ii}$ from~\eqref{eq:hiibar_def} and $\overline{\mathbf{h}}_{ij}$ from~\eqref{eq:hijbar_def}, which hold for all $\X\in\K$, where $\K\subset\MN$ is compact. This concludes the derivation of bounds for the Euclidean Hessian.

\end{appendices}

\bibliographystyle{unsrt}
\bibliography{references}

\end{document}